\documentclass[11pt,letterpaper]{article}
\usepackage[margin=1in]{geometry}
\usepackage[T1]{fontenc}
\usepackage[utf8]{inputenc}
\usepackage{lmodern}
\usepackage[american]{babel}
\usepackage{amsmath,amssymb,amsthm}
\usepackage{graphicx,xcolor,booktabs,array}
\usepackage[shortlabels]{enumitem}
\setlist{topsep=5pt plus 2pt minus 2pt,itemsep=2pt plus 1pt,
  parsep=0pt,partopsep=0pt}
\usepackage[titles]{tocloft}
\usepackage[hidelinks]{hyperref}
\usepackage{mathtools,mathrsfs,bm}
\usepackage{microtype,longtable}
\usepackage{tikz}
\makeatletter
\newcommand{\bibalias}[1]{%
  \hypertarget{cite.#1}{}%
  \protected@write\@auxout{}{\string\bibcite{#1}{\the\value{\@listctr}}}}
\makeatother
\hypersetup{
  pdftitle={Nonlinear Inviscid Damping of Couette Flow with Vorticity Perturbations Supported up to the Boundary},
  pdfauthor={Siming He, Sameer Iyer, and Weiren Zhao},
  pdfsubject={Nonlinear inviscid damping in a finite channel},
  pdflang={en-US},
  pdfkeywords={2D Euler, Couette flow, inviscid damping, Gevrey regularity, boundary effects}
}

\usepackage{hyperref}
\hypersetup{
    colorlinks=true,
    linkcolor=blue,
    filecolor=magenta,      
    urlcolor=blue,
    citecolor=blue,
}

\numberwithin{equation}{section}
\allowdisplaybreaks

\makeatletter
\g@addto@macro\normalsize{%
  \setlength{\abovedisplayskip}{5.5pt plus 2pt minus 2pt}%
  \setlength{\belowdisplayskip}{5.5pt plus 2pt minus 2pt}%
  \setlength{\abovedisplayshortskip}{0pt plus 2pt}%
  \setlength{\belowdisplayshortskip}{3pt plus 2pt minus 1pt}%
}
\def\thm@space@setup{%
  \thm@preskip=6pt plus 2pt minus 2pt
  \thm@postskip=6pt plus 2pt minus 2pt
}
\makeatother
\theoremstyle{plain}
\newtheorem{theorem}{Theorem}[section]
\newtheorem{proposition}[theorem]{Proposition}
\newtheorem{lemma}[theorem]{Lemma}
\newtheorem{corollary}[theorem]{Corollary}
\theoremstyle{definition}

\theoremstyle{definition}
\newtheorem{remark}[theorem]{Remark}

\newcommand{\T}{\mathbb T}
\newcommand{\R}{\mathbb R}
\newcommand{\N}{\mathbb N}
\newcommand{\eps}{\varepsilon}
\newcommand{\dd}{\,\mathrm d}
\newcommand{\ii}{\mathrm i}
\newcommand{\avg}[1]{\langle #1\rangle}
\newcommand{\norm}[1]{\left\lVert #1\right\rVert}
\newcommand{\abs}[1]{\left\lvert #1\right\rvert}
\newcommand{\supp}{\operatorname{supp}}
\newcommand{\dist}{\operatorname{dist}}
\newcommand{\Pzero}{\mathbb P_0}
\newcommand{\Pneq}{\mathbb P_{\neq}}
\newcommand{\Gs}[3]{\mathcal G^{#1,#2;#3}}
\newcommand{\CK}{\mathcal{CK}}
\newcommand{\Ext}{\mathrm{Ext}}
\newcommand{\Int}{\mathrm{Int}}
\newcommand{\NR}{\mathrm{NR}}
\newcommand{\Res}{\mathrm{R}}
\newcommand{\pa}{\partial}

\newcommand{\A}{\mathcal{A}}
\newcommand{\J}{\mathcal{J}}

\newcommand{\E}{\mathcal{E}}

\newcommand{\be}{\begin{eqnarray*}}
\newcommand{\bel}{\begin{eqnarray}}

\newcommand{\ee}{\end{eqnarray*}}
\newcommand{\eel}{\end{eqnarray}}
\newcommand{\ba}{\begin{aligned}}
\newcommand{\ea}{\end{aligned}}
\newcommand{\al}{\alpha}
\newcommand{\na}{\nabla}

\newcommand{\pw}[2]{#1^{(#2)}}

\newcommand{\pvt}{\partial_v^t}

\newcommand{\pex}[1]{\Phi_{\rm Ext}^{(#1)}}
\newcommand{\pin}[2]{\Phi_{\rm Int}^{(#2;#1)}}

\newcommand{\bgs}{\begin{split}}
\newcommand{\eds}{\end{split}}

\usepackage{needspace}

\begin{document}

\title{Nonlinear Inviscid Damping of Couette Flow
with Vorticity Perturbations Supported up to the Boundary}
\author{%
Siming He\thanks{Department of Mathematics, University of South Carolina,
Columbia, SC 29208. Email: \texttt{siming@mailbox.sc.edu}.}
\and
Sameer Iyer\thanks{Department of Mathematics,
University of California, Davis, CA 95616, USA.
Email: \texttt{samiyer@ucdavis.edu}.}
\and
Weiren Zhao\thanks{Department of Mathematics, New York University Abu Dhabi,
Abu Dhabi, UAE. Email: \texttt{wz19@nyu.edu}.}}
\date{\today}

\maketitle

\begin{abstract}
We establish nonlinear inviscid damping of the near-Couette flows in a finite channel, with vorticity perturbations vanishing to infinite order at the boundary. 
To study the boundary effect, we decompose the vorticity into dyadic layers whose supports expand toward the boundary, and use an inductive Newton scheme to track nonlinear interactions across layers. Finally, precise elliptic estimates are established to quantify the velocity generated by different layers. 
\end{abstract}



\setcounter{tocdepth}{1}
\hypersetup{bookmarksdepth=3}
\tableofcontents

\addcontentsline{toc}{part}{Main text}

\section{Introduction}
\label{sec:setup}

We consider the two-dimensional Euler equation that governs the dynamics of the perturbation $U=(U^x,U^y)$ of the Couette flow $(y,0)$: 
\begin{align}\label{eq:euler}
\begin{cases}    &\pa_t \omega +y\pa_x \omega+U\cdot \na \omega=0,\\
    &U=\na^\perp\Delta^{-1}\omega,\\
    &U^y(t,x,y=\pm 1)\equiv 0, \\
    &\omega(t=0)=\omega_{\rm in}.
    \end{cases}
\end{align}
Here, the domain is the periodic channel
$\Omega=\mathbb{T}_{2\pi}\times[-1,1]$, and the vorticity
$\omega:=\pa_xU^y-\pa_yU^x$ determines the fluid velocity $U$ through the
Biot--Savart law, in which we use $\Delta^{-1}$ to denote the inverse of the
Dirichlet Laplacian and $\na^\perp$ to denote $(-\pa_y,\pa_x)$.  No
penetration makes the stream function constant on each wall, but does not
force the two constants to coincide.  Throughout the paper we fix the
zero-flux (equivalently, zero harmonic-flow) normalization and subtract the
remaining additive constant, so that $\psi|_{y = \pm 1}=0$ and the
Dirichlet inverse is the correct Biot--Savart operator.

Classical global well-posedness results for the two-dimensional Euler equation, following Wolibner \cite{Wolibner33}, Yudovich \cite{Yudovich63,Yudovich95}, and subsequent work, contrast with the limited understanding of its long-time dynamics in many settings.

A natural route towards a better understanding of the long-time behavior of the Euler flows is to study the dynamics of perturbations to given Euler steady states. Parallel shear flows have served as
a basic testing ground since the nineteenth century. Rayleigh's normal-mode
analysis introduced the Rayleigh equation and the inflection-point necessary
condition for spectral instability \cite{Rayleigh80}. The Couette flow is one of the simplest nontrivial Euler steady states and has been a central model in hydrodynamic stability. The work of Kelvin \cite{Kelvin87}, Orr \cite{Orr07}, Rayleigh \cite{Rayleigh80}, and later contributions \cite{faddeev71,arnold98} laid the foundations for this line of research.  For Couette flow the linear mechanism is particularly transparent.  For every
nonzero horizontal mode \(k\), the linearized vorticity satisfies
\[
  \omega_k(t,y)=e^{-ikty}\omega_{{\rm in},k}(y).
\]
Thus the vorticity is transported to progressively finer cross-stream scales
without losing its \(L^2\) size.  Elliptic inversion in the Biot--Savart law
converts this filamentation into decay of the velocity.  For the linearized
problem, the anisotropic rates, with their datum norms, are
\[
 \begin{aligned}
  \|P_{\neq}U^x(t)\|_{L^2}
    &\lesssim \langle t\rangle^{-1}
                     \|P_{\neq}\omega_{\rm in}\|_{L^2_xH^1_y},\\
  \|U^y(t)\|_{L^2}
    &\lesssim \langle t\rangle^{-2}
                     \|P_{\neq}\omega_{\rm in}\|_{L^2_xH^2_y}.
 \end{aligned}
\] 
This damping phenomenon, which does not involve any fluid viscosity, is called \emph{inviscid damping}. 
At a high level, this damping effect arises from the small-scale creation/phase mixing of the fluid vorticity, which in turn dampens the velocity through the inverse Laplacian in the Biot-Savart law.  

The modern interpretation of inviscid damping is closely related to Landau
damping for the Vlasov--Poisson equation.  Landau predicted collisionless
decay of the electric field in 1946 \cite{Landau46}; van Kampen and Case
subsequently explained the linear effect through phase mixing in the continuous spectrum \cite{VKampen55, Case1959Plasma}. The analogy is structural: kinetic free transport mixes a distribution in velocity, while Couette advection mixes vorticity in the cross-stream variable; the electric field and the fluid velocity are recovered by smoothing operators and therefore decay even though the underlying density need not converge strongly. A landmark nonlinear stability theorem was proved by Mouhot and Villani \cite{MouhotVillani11}.  Bedrossian, Masmoudi, and Mouhot extended the admissible Gevrey range and developed a paraproduct-based proof \cite{BedrossianMasmoudiMouhot16}, also see \cite{grenier2022landau} for a short proof and \cite{ionescu2024nonlinear} for the sharp regularity result.  Plasma echoes show that such regularity assumptions are not merely technical: Bedrossian \cite{Bedrossian21} constructed long echo chains that invalidate a uniform Sobolev analogue of nonlinear Landau damping in the relevant settings.  The same resonance geometry underlies the time-frequency weights used in nonlinear inviscid damping. 

Away from the affine Couette profile, the term
\(u''\partial_x\Delta^{-1}\omega\) makes the linearized Euler equation
genuinely nonlocal and spectral.  Resolvent singularities at critical layers
\(U(y)=c\), discrete modes, and embedded eigenvalues must all be controlled. For monotone shears, Zillinger established damping and scattering by assuming the smallness of $u''$ \cite{Zillinger2014}. Wei, Zhang, and Zhao removed the smallness assumption and proved damping and scattering for a class of monotone shear flows, which have no embedded eigenvalues or eigenvalues, and obtained sharp decay estimates \cite{WZZ18}. Jia \cite{Jia2020siam, jia2020linear} developed complementary pointwise and Gevrey estimates. Non-monotone shear flows exhibit a second mechanism.  At a stationary streamline \(y_*\), where \(u'(y_*)=0\), direct stationary phase suggests weaker mixing. Bouchet and Morita \cite{BouchetMorita10} predicted that the self-consistent dynamics instead depletes the asymptotic vorticity at \(y_*\), cancelling the leading stationary-phase contribution. Wei, Zhang, and Zhao proved linear vorticity depletion for a general class of shear flows with nondegenerate stationary points and obtained sharp damping for important symmetric examples including Poiseuille-type flows \cite{WeiZhangZhao2019}. They also treated the Kolmogorov flow and its enhanced dissipation counterpart
\cite{WeiZhangZhao20}.  Ionescu, Iyer, and Jia subsequently removed the symmetry assumption and proved quantitative depletion and optimal velocity damping rates for channel shears with one nondegenerate critical point and no discrete or generalized embedded eigenvalues \cite{ionescu2024stability}.  The nonlinear inviscid damping and vorticity depletion for non-monotone shears remain open. 

The nonlinear problem contains a new obstruction: resonant interactions can partially reverse mixing and create an echo cascade.  Lin and Zeng constructed nonparallel steady structures arbitrarily close to Couette flow in vorticity \(H^s\) for \(s<3/2\), ruling out a uniform damping theorem at that regularity \cite{LinZeng11}.  Bedrossian and Masmoudi \cite{BM13} proved a nonlinear asymptotic-stability theorem for Couette flow on \(\mathbb T_{2\pi}\times\mathbb R\) for small perturbations in Gevrey classes strictly smoother than Gevrey class \(2\). The main challenges in developing the nonlinear stability result come from two aspects. First, the proof uses sharp pointwise-in-time inviscid damping estimates, whose derivation requires regularity. Second, at the nonlinear level, another weakly nonlinear effect, so-called ``echoes'', can cause a loss of regularity at specific time--frequency configurations. If one cannot control nonlinear echoes, establishing stability in high-regularity spaces will be difficult. In the work \cite{BM13}, the authors introduce a delicate imbalanced \emph{Fourier multiplier} $\mathcal{J}$ to capture all potentially dangerous echo-cascade growth in the system. By carefully tuning the Gevrey regularity space over time, the authors show that the solution is sufficiently regular in a suitable nonlinear coordinate system, thereby yielding optimal inviscid damping. For the inhomogeneous incompressible Euler equation, \cite{chen2025nonlinear} treats perturbations of Couette flow with constant reference density on $\mathbb T\times\mathbb R$. On the instability side, Deng and Masmoudi constructed nonlinear echo cascades producing long-time instability at regularities weaker than Gevrey class \(2\) \cite{deng2023long}. 

The establishment of nonlinear stability of Couette flow on the unbounded
domain motivates the corresponding stability problem in a bounded region.
Let us focus on the channel $\mathbb{T}_{2\pi}\times[-1,1]$, which also
supports Couette flow. As far as positive stability results go, Ionescu and Jia extended the nonlinear Couette damping proven by Bedrossian-Masmoudi to a finite channel in Gevrey regularity, assuming that the perturbation is supported a fixed distance from the walls \cite{IonescuJia20inviscid}.   Nonlinear inviscid damping near a class of spectrally stable monotone shears was proved by Ionescu-Jia \cite{IJ2020} and  Masmoudi-Zhao \cite{MasmoudiZhao24} independently, under the assumption that both $u''$ and the perturbation are supported a fixed distance from the walls. {With the similar compact support assumptions on the monotone shear flows, the background density, and the initial perturbation, \cite{zhao2025inviscid} treats variable-density monotone flow in a finite channel for the inhomogeneous incompressible Euler equation. The traveling and steady structures in \cite{zhao2026traveling} obstruct damping in velocity $H^{5/2-\tau}$ and density $H^{3/2-\tau}$ neighborhoods for small $\tau>0$; this lower-regularity obstruction does not contradict the Gevrey results. }

In all these positive results, one of the primary ``structural" or qualitative assumptions is that the perturbation is compactly supported away from the boundaries (which is then propagated for all time due to the integrability 
$\|U^y(t)\|_{L^\infty}\leq C\eps\langle t\rangle^{-2}$). Indeed, the interplay between the boundary and inviscid damping is delicate, as discovered by Zillinger in a series of \textit{boundary driven instability} results, \cite{Zillinger2016Boundary}. For smooth shears with uniformly positive derivative
on the closed channel, Zillinger identified a boundary obstruction to
uniform Sobolev control of the
vorticity in coordinates moving with the background shear.  More precisely,
if $u''(b)\omega_{{\rm in},k}(b)\neq0$ for a nonzero horizontal mode $k$
at a wall $b\in\{-1,1\}$, a uniform-in-time $H^2$ bound for that
transported vorticity mode would force its normal derivative at the wall
to grow logarithmically, contradicting the bound
\cite[Theorem~B.1]{Zillinger2014}.  This linear mechanism is absent for
the exact Couette profile: $u''=0$, and the transported linearized
vorticity is constant in time. Nevertheless, the nonlinear problem near the Couette flow requires treating the evolving shear $y + U^x_0(t, y)$ as the background shear, and therefore one expects the nonlinear problem to generically fall into the Zillinger instability regime.  Zillinger also proved a related nonlinear
obstruction under additional nondegeneracy assumptions on the evolving
mean shear and a boundary vorticity mode, together with a decay assumption
on the differentiated nonlinear boundary term
\cite[Theorem~5]{Zillinger2016Boundary}.    

We now observe that there remains a gap between Zillinger's obstruction and the compact-support results \cite{IonescuJia20inviscid, IJ23, MasmoudiZhao24}. The remaining question is whether one can relax the compact support assumption stated above to an infinite-order vanishing assumption at the boundary while retaining nonlinear inviscid damping, which is the purpose of this work. 

 \subsection{Statement and discussion}

\subsubsection{Regularity and boundary vanishing}
\label{sec:fourier-conventions}

Denote the distance to the channel walls by
\[
  d(y)=1-|y|.
\]
We call a smooth function $q$ on $\overline\Omega$ \emph{boundary-flat} if
\[
  \partial_y^j q(x,\pm1)=0
  \qquad\text{for every }j\geq0.
\]
For such a function, let the same symbol denote its zero extension to
$\T_{2\pi}\times\R$.  We use the mixed Fourier transform
\[
 \widehat q_k(\eta)
 =\frac1{2\pi}\int_0^{2\pi}\int_{\R}
       q(x,y)e^{-\ii(kx+\eta y)}\dd y\dd x,
\]
and use the same convention with $(x,y)$ replaced by $(z,v)$ for profile
variables.  This transform is unitary from $L^2(\mathbb T_{2\pi}\times\mathbb R,
\dd z\dd v)$ to $\ell_k^2L_\eta^2$.  All unsubscripted $L^2$ norms and
inner products in the profile multiplier energies use this flat cylindrical
measure.  For scalar fields, viewed as independent of $z$, the energy
pairing is
\[
 \langle u_1,u_2\rangle
 =\int_{\mathbb R}\widehat u_1(\eta)
                \overline{\widehat u_2(\eta)}\,\dd\eta
 =2\pi\int_{\mathbb R}u_1(v)\overline{u_2(v)}\,\dd v.
\]
In particular, scalar energy norms have the same Parseval normalization
as their displayed Cauchy--Kovalevskaya ($\CK$) energy terms below.  Explicit one-dimensional $L^p(\mathbb R)$
and $H^r(\mathbb R)$ bounds use ordinary Lebesgue measure; the fixed factors
relating them to the energy norms are retained in equalities and absorbed
only in inequality constants.  The vertical product convention is
\[
 \widehat{u_1u_2}(\eta)
 =\frac1{2\pi}\int_{\mathbb R}
    \widehat u_1(\eta-\xi)\widehat u_2(\xi)\,\dd\xi.
\]
We write
\[
 q_0=\Pzero q=\frac1{2\pi}\int_0^{2\pi}q(x,\cdot)\dd x,
 \qquad q_{\neq}=\Pneq q=q-q_0.
\]
Thus the subscripts $0$ and $\neq$ denote the zero and nonzero horizontal
Fourier modes.  Following \cite[Section~1.2]{BM15}, frequency lengths
are taken in the sum norm:
\[
 |k,\eta|=|(k,\eta)|=|k|+|\eta|,\qquad
 \avg{k,\eta}=\bigl(1+(|k|+|\eta|)^2\bigr)^{1/2}.
\]
For a scalar frequency the same convention gives
$\avg\eta=(1+|\eta|^2)^{1/2}$.  This sum norm is used in the
Gevrey weights and full-frequency Littlewood--Paley projections.
The symbols $k^2+(\eta-kt)^2$ of $-\Delta_L$ remain quadratic.
For $s\in(0,1)$, $\lambda>0$, and
$\sigma\ge0$, define
\begin{equation}
  \norm{q}_{\Gs{s}{\lambda}{\sigma}}^2
  =\sum_{k\in\mathbb Z}\int_{\R}
    e^{2\lambda\abs{(k,\eta)}^s}
    \avg{k,\eta}^{2\sigma}
    \abs{\widehat q_k(\eta)}^2\dd\eta.
  \label{eq:gevrey-norm}
\end{equation}
We fix the basic Gevrey regularity index $s$ and the associated vanishing-order index $\alpha_*$ satisfying
\begin{equation}
  \frac12<s<1,
  \qquad
  \alpha_*>\frac{s}{1-s},
  \label{eq:parameter-range}
\end{equation}
and, for fixed $\lambda_0>0$ and sufficiently large $\sigma$,
assume
\begin{equation}
  \norm{\omega_{\mathrm{in}}}_{\Gs{s}{3\lambda_0}{2\sigma+8}}
  +\norm{e^{c_0d(y)^{-\alpha_*}}\omega_{\mathrm{in}}}_{L^2(\Omega)}
  \le \eps
  \label{eq:data-assumption}
\end{equation}
for some $c_0>0$.  In particular, the vorticity vanishes to infinite order at
both walls.  
The second term in
\eqref{eq:data-assumption} is stronger than a
qualitative flatness condition: on $d(y)\leq2^{-n}$ it gives the explicit bound
\[
  \norm{\omega_{\mathrm{in}}}_{L^2(\left\{(x,y)\in\Omega:\ d(y)\leq2^{-n}\right\})}
  \leq\eps e^{-c_0 2^{n\alpha_*}},\qquad n\geq0,
\]
which is a crucial, scale-specific gain used to close our nonlinear estimates. 

Throughout, $C>0$ denotes a generic constant that may change from
line to line. The positive constant $c$ may differ between computations,
but is held fixed within each computation, subject to the stated local
restrictions. Any reduction of an exponential decay coefficient is
displayed explicitly, for example as $c/2$. 

\subsubsection{Main result}

Write \(\avg{U^x}=\Pzero U^x\) for the zero mode.

\begin{theorem}[Main theorem]
\label{thm:main}
Assume \eqref{eq:parameter-range}, and let $\omega_{\rm in}$ be real-valued
and boundary-flat.
There exists
$\sigma_0=\sigma_0(s,\alpha_*)$ such that the following holds for every
$\sigma\geq\sigma_0$.  Fix $\lambda_0,c_0>0$.  There are constants
\begin{equation}
  \eps_0>0,
  \qquad
  \lambda_\infty\in(0,\lambda_0/2),
  \qquad
  c_\infty\in(0,c_0),
  \label{eq:main-parameter-choices}
\end{equation}
depending only on the displayed parameters, such that if
$0<\eps\leq\eps_0$ and \eqref{eq:data-assumption} holds, then the solution of
\eqref{eq:euler} exists globally and satisfies the following conclusions.
The generic constant $C$ below may depend on
$s,\alpha_*,\sigma,\lambda_0,c_0$, but is independent of $\eps$ and $t$.
\begin{enumerate}[(i)]
  \item The boundary vanishing is propagated uniformly:
  \begin{equation}
    \sup_{t\geq0}
    \norm{e^{c_\infty d(y)^{-\alpha_*}}\omega(t)}_{L^2(\Omega)}
    \leq C\eps.
    \label{eq:main-boundary}
  \end{equation}

  \item The vorticity remains uniformly Gevrey and scatters in a lower norm after
  translation by the shear. Set $\Phi(t,y)=\int_0^t\avg{U^x}(\tau,y)\dd\tau$.
  There exists $\omega_\infty\in\Gs{s}{\lambda_\infty}{\sigma-6}$ such that
  \begin{equation}
    \begin{aligned}
    \sup_{t\geq1}\norm{\omega(t,x+ty+\Phi(t,y),y)}_{\Gs{s}{\lambda_\infty}{\sigma}}
      &\leq C\eps,\\
    \lim_{t\to\infty}
      \norm{\omega(t,x+ty+\Phi(t,y),y)-\omega_\infty(x,y)}_{\Gs{s}{\lambda_\infty}{\sigma-6}}&=0.
    \end{aligned}
    \label{eq:main-scattering}
  \end{equation}

  \item The velocity perturbation obeys the anisotropic damping estimates
  \begin{equation}
    \norm{\Pneq U^x(t)}_{L^2(\Omega)}
      \leq\frac{C\eps}{\avg{t}},
    \qquad
    \norm{U^y(t)}_{L^2(\Omega)}
      \leq\frac{C\eps}{\avg{t}^2}.
    \label{eq:main-damping}
  \end{equation}
  Consequently, $\Pneq\omega(t)\rightharpoonup0$ weakly in $L^2(\Omega)$.
\end{enumerate}
\end{theorem}
Remarks are listed in order: 
\begin{remark}
    We expect an analogue of Theorem~\ref{thm:main} for stable monotone shears $u(y)$ satisfying the hypotheses of \cite{MZ24}, with perturbations supported up to the boundary. The required wave-operator estimates and nonlinear closure are not proved here. The assumption that $u''$ has compact support in \cite{IJ23, MZ24} could be weakened only after specifying quantitative boundary hypotheses and establishing the corresponding operator estimates.
\end{remark}
\begin{remark}
It is an interesting question to study the optimal boundary-flatness.  
The boundary-flatness assumption $\alpha_{*}>\frac{s}{1-s}$ might be improved to $\alpha_{*}\geq \frac{s}{1-s}$ by a careful computation and adjusting $c_0$ in \eqref{eq:parameter-range}. 
    One can also improve Theorem \ref{thm:main} to $s=1/2$, by applying the Fourier multipliers introduced in \cite{IJ20}, correspondingly, the assumptions \eqref{eq:parameter-range} might be changed to $\alpha_*\geq 1$. 

    On the negative side, whether the finite vanishing order of perturbation will cause instability remains open. Indeed, at the linear level, suppose that the initial vorticity \(\omega_{\mathrm{in}}(x,y)=(1-y^2)^m\cos x\), then \(\omega(t,x,y)=(1-y^2)^m\cos(x-ty).\) Let $z=x-ty$ and \(\omega(t, x, y)=f(t,z,y)=(1-y^2)^m\cos z\). Let \(\phi(t,z,y)=(\Delta^{-1}\omega)(t,z+ty,y)=e^{iz}\hat{\phi}(t,y),\) then \[
\begin{aligned}
\hat{\phi}(t,y)
=-\frac{e^{ity}}{\sinh2}\Bigg[
&\sinh(1-y)\int_{-1}^{y}
 \sinh(\eta+1)e^{-it\eta}(1-\eta^2)^m\,d\eta\\
+&\sinh(y+1)\int_y^1
 \sinh(1-\eta)e^{-it\eta}(1-\eta^2)^m\,d\eta
\Bigg].
\end{aligned}
\]
Let \(p_m(t,y)
=\frac12\sum_{j=0}^{2m}
\left[
\frac{1}{(it-1)^{j+1}}
-\frac{1}{(it+1)^{j+1}}
\right]q_m^{(j)}(y),\) with \(q_m^{(j)}(y)
=
\sum_{\ell=\lceil j/2\rceil}^{m}
(-1)^\ell\binom{m}{\ell}
\frac{(2\ell)!}{(2\ell-j)!}\,
y^{2\ell-j}.\)
Then 
\[
\begin{aligned}
\hat{\phi}(t,y)=p_m(t,y)
-\frac{e^{ity}}{\sinh2}\Big[
&e^{it}p_m(t,-1)\sinh(1-y)
+e^{-it}p_m(t,1)\sinh(1+y)
\Big].
\end{aligned}
\]
Thus we have the linear instability $\|\phi\|_{H^{m+3}}\gtrsim  t$. 
\end{remark}

\section{Outline}
\label{sec:outline}

\subsection{Abstract argument and key observation}
\label{sec:abstract-argument}

On $\Omega=\mathbb T_{2\pi}\times[-1,1]$, with $\mathbb T_{2\pi}=\mathbb R/(2\pi\mathbb Z)$, consider
\[
 \partial_tF+y\partial_xF=N(F,F),\qquad
 N(F,F)=\mathcal E(F)\cdot\nabla F,\qquad F|_{t=0}=F_{\mathrm{in}}.
\]
Here the linear velocity-recovery operator $\mathcal E$ gains at least one spatial derivative, and $\mathcal E(F)$ is divergence-free and tangent to the walls. This model illustrates the exact decomposition and support preservation.

Decompose the initial datum into smooth pieces supported away from the walls:
\[
 F_{\mathrm{in}}=\sum_{n\geq1}F_{\mathrm{in}}^{(n)},\qquad
 \supp F_{\mathrm{in}}^{(n)}
 \subset\left\{(x,y)\in\Omega:1-|y|\geq2^{-n}\right\}.
\]
Set $F^{[n]}=\sum_{m=1}^nF^{(m)}$, $F^{[0]}=0$, and solve
\begin{equation}
 \bigl(\partial_t+y\partial_x-\mathcal E(F^{[n]})\cdot\nabla\bigr)F^{(n)}
   =\mathcal E(F^{(n)})\cdot\nabla F^{[n-1]},\qquad
 F^{(n)}(0)=F_{\mathrm{in}}^{(n)}.
 \label{eq:abstract-increment}
\end{equation}

Linearity of $\mathcal E$ gives
\[
 \mathcal E(F^{[n]})\cdot\nabla F^{(n)}
   +\mathcal E(F^{(n)})\cdot\nabla F^{[n-1]}
   =\mathcal E(F^{[n]})\cdot\nabla F^{[n]}
   -\mathcal E(F^{[n-1]})\cdot\nabla F^{[n-1]}.
\]
Thus each partial sum solves the toy equation with datum $\sum_{m=1}^nF_{\mathrm{in}}^{(m)}$. The increment equation remains nonlinear through its transport field. Its source differentiates an earlier partial sum, which we estimate at a slightly larger Gevrey radius. The derivative gain of $\mathcal E$ helps estimate the velocity; nonlinear energy estimates must still give summable increment bounds before we pass to the infinite sum.

To preserve separation from the walls, suppose inviscid damping gives, for small data size $\eps>0$ and $C$ independent of $n,t$,
\begin{equation}
 \bigl|\bigl(\mathcal E(F^{[n]})\bigr)^y(t,x,y)\bigr|
 \leq\frac{C\eps(1-|y|)}{(1+t)^2},\qquad t\geq0.
 \label{eq:key-ob}
\end{equation}
Here the superscript $y$ denotes the vertical component. This decay is an additional hypothesis, not a consequence of smoothing. Along characteristics of $(y,0)-\mathcal E(F^{[n]})$, integration of \eqref{eq:key-ob} gives
\[
 1-|y(t)|\geq (1-|y(\tau)|)
 \exp\left[-C\eps\left(\frac{t}{1+t}-\frac{\tau}{1+\tau}\right)\right],
 \qquad 0\leq\tau\leq t.
\]
Indeed, $|\partial_t\log(1-|y(t)|)|\leq C\eps(1+t)^{-2}$ almost everywhere.

Although $\mathcal E(F^{(n)})$ need not vanish near the walls, its product with $\nabla F^{[n-1]}$ is supported on the earlier increments. Induction gives
\begin{equation}
 \supp F^{(n)}(t)\subset
 \left\{(x,y)\in\Omega:
 1-|y|\geq2^{-n}\exp\left(-C\eps\frac{t}{1+t}\right)\right\},
 \qquad t\geq0.
 \label{eq:abstract-profile-drift}%
\end{equation}
For $n=1$ the source vanishes. For $m<n$, $2^{-m}\geq2^{-n}$, and propagation from any source time $\tau$ gives the same bound \eqref{eq:abstract-profile-drift}. Thus the separation loses at most the factor $e^{-C\eps}$, uniformly in $n$, even as the initial distances tend to zero.

For Euler, we must prove velocity decay and Gevrey bounds for elliptic recovery with explicit dependence on each layer's distance from the boundary. In nonlinear coordinates the elliptic coefficients depend on the solution, so we also take exact differences of the elliptic equations.

\begin{remark}
    A boundary-flat analogue of Theorem \ref{thm:main} for the Stokes-transport equation is deferred to forthcoming work; it is not proved here. For the interior-supported setting, see \cite{sinambela2026asymptotic}.
\end{remark}

The boundary-vanishing condition in \eqref{eq:parameter-range}
compensates for the localization loss in the elliptic estimates. The reason for this choice is a competition
between two Fourier exponentials. First, fix a regularity paramter $\rho$ that is smoother than our working Gevrey-$s$ space, namely $s<\rho<1$. In the nonlinear estimates, a key ingredient for controlling the stream function is the Fourier decay of the localized Green kernels at scale $2^{-m}$, which is defined in terms of a cutoff function of regularity $G^{\rho}$, \eqref{eq:exterior-green-operator}, \eqref{eq:interior-green-operator}. At cutoff scale $2^{-m}$, this kernel provides the factor
\[
 e^{-\lambda_m^\Delta |\eta-\xi|^\rho},\qquad
  \lambda_m^\Delta=c_\Delta2^{-m\rho};
\]
for the full kernel estimate, see \eqref{eq:localized-green-kernel}.  Here $\eta$ and $\xi$ are the
vertical output and source frequencies at a fixed nonzero horizontal mode.
Transferring the Gevrey weight from
$\eta$ to $\xi$, as in a weighted elliptic estimate or a commutator,
costs at most $e^{C|\eta-\xi|^s}$ at the level of the Gevrey exponential.
The inequality $\rho>s$ makes this cost finite at each fixed scale, but
the cutoff radius $\lambda_m^\Delta$ shrinks as the layer approaches the
wall.  Optimizing the competing exponents gives
\[
 e^{C|\eta-\xi|^s-\lambda_m^\Delta |\eta-\xi|^\rho}
 \le e^{C_1 2^{m\alpha_\dagger}}
       e^{-\frac12\lambda_m^\Delta |\eta-\xi|^\rho},
 \qquad \alpha_\dagger=\frac{\rho s}{\rho-s}.
\]
The remaining negative exponential is kept when integrating the kernel;
its polynomial moments are absorbed by the same localization loss.
The complete estimate, including the elliptic denominators and multiplier
corrections, is given in
\eqref{eq:common-mlambda-transfer} and
\eqref{eq:kernel-absorption} below.

On the other hand, the boundary-vanishing assumption \eqref{eq:data-assumption} supplies the gain
$e^{-c2^{n\alpha_*}}$ on the $n$th boundary layer.
Section~\ref{sec:decomposition} converts it into compatible initial
increments of size $\eps e^{-c2^{n\alpha_*}}$; see \eqref{eq:shell-data}.
For $1\le m\le n$, the competing scale factors satisfy
\[
 e^{C 2^{m\alpha_\dagger}}e^{-c2^{n\alpha_*}}
 \lesssim e^{-\frac c2 2^{n\alpha_*}},
\]
which decays double exponentially in $n$ for $\alpha_*>\alpha_\dagger$.

A convenient way to keep track of $\rho$ and $\alpha_{\dagger}$ is the following. Choose $K>0$ and set $\alpha_\dagger$ and $\rho$ as follows: 
\begin{equation}
  \alpha_\dagger
    =\frac{K+s}{K}\frac{s}{1-s},\qquad \rho=\frac{K+s}{K+1}.
  \label{eq:alpha-dagger}
\end{equation}
Since we assume in \eqref{eq:parameter-range} that $\alpha_*>s/(1-s)$, one can choose $K$ so that
$\alpha_\dagger<\alpha_*$. 

\Needspace{14\baselineskip}
\subsection{Nonlinear coordinates}
\label{sec:boundary-newton}%
\label{sec:euler-roadmap}%
For times $0 \le t \le 1$, we appeal to Appendix \ref{sec:short-time} for the construction of the Euler solution. For the long-time
argument, we will use the $(z, v)$ coordinate system introduced in \cite{BM13}. First, we will need a preliminary lemma. 
\label{sec:coordinates}
Let $\avg{U^x}(t,y)=(2\pi)^{-1}\int_0^{2\pi}U^x(t,x,y)\dd x$ denote the horizontal mean of the perturbation velocity.

\begin{lemma}[Boundary shear]
\label{lem:boundary-shear}
Let $\omega$ be a smooth solution of the Euler system \eqref{eq:euler} with
$\omega(t,x,\pm1)=0$. Then
\begin{equation}
  \avg{U^x}(t,\pm1)=c_\pm:=\avg{U^x}(0,\pm1).
  \label{eq:boundary-shear}
\end{equation}
\end{lemma}

\begin{proof}
At $y=\pm1$, $U^y=0$, so averaging the horizontal momentum
equation over $x$ gives $\partial_t\avg{U^x}=-\avg{U^y\,\partial_yU^x}=0$. 
\end{proof}

We now recall the $(z, v)$ coordinates: 
\begin{equation}
  z=x-tv(t,y),
  \qquad
  v(t,y)=y+\frac1t\int_0^t\avg{U^x}(\tau,y)\dd\tau.
  \label{eq:coordinates}
\end{equation}
Lemma~\ref{lem:boundary-shear} implies that the range of $v$ is the fixed
interval
\begin{equation}
  I=[L,R]=[-1+c_-,1+c_+].
  \label{eq:fixed-interval}
\end{equation}
For $\psi=\Delta^{-1}\omega$ in \eqref{eq:euler}, define the profiles
\begin{equation}
\begin{split}
  f(t,z,v)&=\omega(t,x,y),\qquad \phi(t,z,v)=\psi(t,x,y),\\
  g(t,v)&=\partial_t v(t,y),\qquad h(t,v)=\partial_yv(t,y)-1.
\end{split}
\label{eq:profile-f-phi}
\end{equation}
With $\partial_v^t=\partial_v-t\partial_z$, these quantities satisfy
\begin{equation}
\left\{
\begin{aligned}
  \partial_tf
    &=(1+h)\bigl(\partial_v\phi_{\neq}\,\partial_zf
      -\partial_z\phi_{\neq}\,\partial_vf\bigr)-g\partial_vf,\\
  \partial_th+g\partial_vh&=(1+h)\partial_vg,\\
  \partial_tg+\frac2t g+g\partial_vg
    &=\frac1t\Pzero\bigl(f_{\neq}\partial_z\phi_{\neq}\bigr),\\
  \Delta_t\phi&=f,
\end{aligned}
\right.
\label{eq:profile-system}
\end{equation}
where
\begin{equation}
  \Delta_t
  =\partial_z^2+(1+h)^2(\partial_v^t)^2
    +(1+h)(\partial_vh)\partial_v^t.
  \label{eq:variable-laplacian}
\end{equation}

\subsection{Gevrey cutoffs}\label{sec:decomposition}

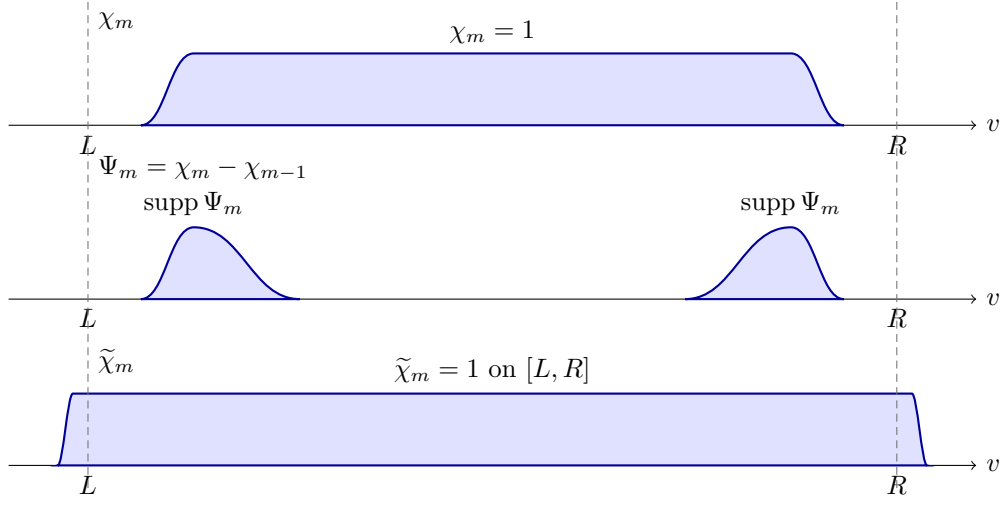
\begin{figure}[tbp]
\centering
\begin{tikzpicture}[x=1cm,y=1cm,font=\small,
  cutoff/.style={draw=blue!60!black,fill=blue!12,thick},
  wall/.style={draw=black!55,densely dashed},
  axis/.style={draw=black,thin,->}]
 \begin{scope}[yshift=4.5cm]
  \node[anchor=west] at (1.3,1.38) {$\chi_m$};
  \draw[axis] (.25,0)--(13.05,0) node[right] {$v$};
  \node[below] at (1.3,0) {$L$};
  \node[below] at (12.0,0) {$R$};
  \filldraw[cutoff]
   (2.0,0) .. controls (2.32,0) and (2.38,.95) .. (2.7,.95)
   --(10.6,.95)
   .. controls (10.92,.95) and (10.98,0) .. (11.3,0)--cycle;
  \node[above] at (6.65,.95) {$\chi_m=1$};
 \end{scope}
 \begin{scope}[yshift=2.2cm]
  \node[anchor=west] at (1.3,1.73) {$\Psi_m=\chi_m-\chi_{m-1}$};
  \draw[axis] (.25,0)--(13.05,0) node[right] {$v$};
  \node[below] at (1.3,0) {$L$};
  \node[below] at (12.0,0) {$R$};
  \filldraw[cutoff]
   (2.0,0) .. controls (2.32,0) and (2.38,.95) .. (2.7,.95)
   .. controls (3.34,.95) and (3.46,0) .. (4.1,0)--cycle;
  \filldraw[cutoff]
   (9.2,0) .. controls (9.84,0) and (9.96,.95) .. (10.6,.95)
   .. controls (10.92,.95) and (10.98,0) .. (11.3,0)--cycle;
  \node[above] at (2.7,.98) {$\supp\Psi_m$};
  \node[above] at (10.6,.98) {$\supp\Psi_m$};

 \end{scope}
 \begin{scope}
  \node[anchor=west] at (1.3,1.38) {$\widetilde\chi_m$};
  \draw[axis] (.25,0)--(13.05,0) node[right] {$v$};
  \node[below] at (1.3,0) {$L$};
  \node[below] at (12.0,0) {$R$};
  \filldraw[cutoff]
   (.9,0) .. controls (.97,0) and (1.03,.95) .. (1.1,.95)
   --(12.2,.95)
   .. controls (12.27,.95) and (12.33,0) .. (12.4,0)--cycle;
  \node[above] at (6.65,.95) {$\widetilde\chi_m=1$ on $[L,R]$};
 \end{scope}
 \foreach \x in {1.3,12.0}{
   \draw[wall] (\x,-.3)--(\x,6.15);
 }
\end{tikzpicture}
\caption{Cutoffs on the channel,
schematically for large $m$ and a nondecreasing choice of $\gamma$.
The dashed lines mark the walls; boundary strips are enlarged and not
to scale. The blue profiles show $\chi_m$, $\Psi_m$, and
$\widetilde\chi_m$ from \eqref{eq:exact-shell-partition} and
\eqref{eq:exterior-cap-geometry}.}
\label{fig:cutoff-profiles}
\label{fig:compatible-data}
\end{figure}

Recall $s$ and the boundary exponent $\alpha_*$ from
\eqref{eq:parameter-range}, and the localization exponent
$\alpha_\dagger$ and cutoff exponent $\rho$ from \eqref{eq:alpha-dagger}.
\stepcounter{equation}
We will define a Gevrey-$\rho$ partition of unity
$\left\{\Psi_m\right\}_{m\ge0}$ on $(L,R)$ as follows.  Put
\begin{equation}
  \mathbf d(v)=\dist(v,\left\{L,R\right\}),\qquad d_m=2^{-m}.
  \label{eq:wall-distance-scales}
\end{equation}
Let
$\gamma\in G^\rho(\mathbb R)$ satisfy $\gamma(r)=0$ for $r\leq1$ and
$\gamma(r)=1$ for $r\geq2$.  
Here $G^\rho$ means that $\gamma$ is bounded and
$\|\partial_r^j\gamma\|_\infty\leq C^{j+1}(j!)^{1/\rho}$ for $j\geq1$.
We then define the interior and shell cutoffs as follows: 
\begin{equation}
  \chi_m(v)=
  \begin{cases}
    \gamma(2^{m+4}\mathbf d(v)),&v\in(L,R),\\
    0,&v\notin(L,R),
  \end{cases}
  \qquad
  \Psi_0=\chi_0,\qquad \Psi_m=\chi_m-\chi_{m-1}\quad(m\geq1).
  \label{eq:exact-shell-partition}
\end{equation}
We also define exterior extensions via 
\begin{align}
\widetilde\chi_m(v)=\begin{cases}
    1,&v\in(L,R),\\
    \gamma(3 - 256 d_m^{-1} \mathbf d(v)),&v\notin(L,R).
  \end{cases}
\end{align}
Figure~\ref{fig:cutoff-profiles} shows these cutoffs on the channel.

Now, we will collect some useful properties of these cutoff functions. Clearly, the shell cutoffs form a partition of unity $\sum_{m\geq0}\Psi_m=1$ pointwise on $(L,R)$. Regarding the supports, we have  \begin{subequations} 
\begin{align}
  &\supp\Psi_m\subset\left\{v\in(L,R)\mid c2^{-m}\leq \mathbf d\leq C2^{-m}\right\},
  && \chi_{m+2}=1
  \quad\hbox{on }\supp\Psi_m,
  \label{eq:psi-support} \\
   &\supp\chi_{m+2}\subset[L+d_m/64,R-d_m/64],
 &&\chi_{m+2}=1\quad\hbox{on }[L+d_m/32,R-d_m/32], \label{eq:interior-cutoff-geometry} \\
 &\supp\widetilde\chi_m\subset[L-d_m/128,R+d_m/128]  &&\widetilde\chi_m=1\quad\hbox{if }\dist(v,[L,R])\le d_m/256.
 \label{eq:exterior-cap-geometry}
\end{align}
\end{subequations}

Now, we turn to quantitative regularity bounds on these cutoff functions. First, 
For $m\geq0$,
\begin{equation}
  \abs{\widehat{\Psi_m}(\eta)}
  \le Cd_m\exp\left(-c(d_m\abs{\eta})^{\rho}\right).
  \label{eq:cutoff-fourier}
\end{equation}
The two cutoff families satisfy
\begin{equation}
 \|\partial_v^j\chi_{m}\|_\infty
 +\|\partial_v^j\widetilde\chi_m\|_\infty
 \le C_0^{j+1}d_m^{-j}(j!)^{1/\rho}.
 \label{eq:extension-cutoff-derivatives}
\end{equation}

\begin{lemma}[Boundary localization]
\label{lem:localization}
Fix $\lambda>0$, $\sigma\geq0$, $\alpha>\alpha_\dagger$, and a 
coefficient $c_0>0$, as in \eqref{eq:data-assumption}.
Assume that $q\in\Gs{s}{2\lambda}{2\sigma}$ is boundary-flat and
$\|e^{c_0\mathbf d(v)^{-\alpha}}q\|_{L^2(\mathbb T_{2\pi}\times(L,R))}\leq M_0.$ 
Then, 
\begin{equation}
  \norm{q\Psi_n}_{\Gs{s}{\lambda}{\sigma}}
  \leq C\bigl(\norm{q}_{\Gs{s}{2\lambda}{2\sigma}}+M_0\bigr)e^{-c2^{n\alpha}},\quad \text{for }n\geq0.
  \label{eq:localization-estimate}
\end{equation}
Here $C,c>0$ depend only on the fixed parameters and cutoffs, not on $q$, $M_0$, or $n$.
\end{lemma}

\begin{proof}
Cauchy--Schwarz and Young's inequality give, respectively,
\begin{align*}
 \|q\Psi_n\|_{\Gs{s}{\lambda}{\sigma}}
 \leq& C\|q\|_{L^2(\mathbb T_{2\pi}\times\supp\Psi_n)}^{\frac{1}{2}}
          \|q\Psi_n\|_{\Gs{s}{2\lambda}{2\sigma}}^{\frac{1}{2}},\\
 \leq &C\|q\|_{L^2(\mathbb T_{2\pi}\times\supp\Psi_n)}^{\frac{1}{2}} \bigl\|e^{2\lambda|\xi|^s}\avg\xi^{2\sigma}
              \widehat\Psi_n(\xi)\bigr\|_{L^1_\xi}^{\frac12}\|q\|_{\Gs{s}{2\lambda}{2\sigma}}^{\frac12}\\
  \leq{}& Ce^{C2^{\alpha_\dagger n}}e^{-\frac{c_0}{2} 2^{\al n}}M_0^{1/2}\|q\|_{\Gs{s}{2\lambda}{2\sigma}}^\frac12.
\end{align*}
The conclusion follows from $\alpha > \alpha_{\dagger}$ and an application of Young's inequality for products. 
\end{proof}

\subsection{Initial data for the Newton iteration}

We will begin the Newton evolution at $t=1$  with the profiles \eqref{eq:profile-f-phi}:
\begin{align} \label{fingindef}
  (f_{\mathrm{in}},h_{\mathrm{in}},g_{\mathrm{in}}):=(f,h,g)|_{t=1}.
\end{align}
Fix a terminal index $\mathcal{N} \in \mathbb{N}$.  We will assume 
\begin{equation}
 \chi_{\mathcal{N}}(f_{\mathrm{in}},h_{\mathrm{in}},g_{\mathrm{in}})
       =(f_{\mathrm{in}},h_{\mathrm{in}},g_{\mathrm{in}}).
 \label{eq:terminal-cutoff-identity}
\end{equation}
Even though no such $\mathcal{N}$ exists for the true datum, we can perform an approximation procedure on the initial data and pass to the limit as $\mathcal{N} \rightarrow \infty$, as all of our estimates will be uniform in $\mathcal{N}$. This limiting procedure will be detailed in Section \ref{sec:convergence}. Our Newton iteration will be conducted for $1 \le n \le \mathcal{N}$, and ultimately be uniform in $\mathcal{N}$.

For each $1 \le n \le \mathcal{N}$, our aim is to  decompose these profiles into initial increments $(f_{\mathrm{in}}^{(n)},h_{\mathrm{in}}^{(n)},g_{\mathrm{in}}^{(n)})$,  whose support is approximately $2^{-n}$ away from the walls, to initialize each step of the iteration.  
We localize all three profiles directly by defining
\begin{align}
 (f_{\mathrm{in}}^{(1)},h_{\mathrm{in}}^{(1)},g_{\mathrm{in}}^{(1)})
   :=\chi_1(f_{\mathrm{in}},h_{\mathrm{in}},g_{\mathrm{in}});
   \quad (f_{\mathrm{in}}^{(n)},h_{\mathrm{in}}^{(n)},g_{\mathrm{in}}^{(n)})
   &:=\Psi_n(f_{\mathrm{in}},h_{\mathrm{in}},g_{\mathrm{in}}),
       \qquad 2\leq n\leq {\mathcal{N}}.
\end{align}

Every increment vanishes near both walls, and
\begin{equation*}
 \supp(f_{\mathrm{in}}^{(n)},h_{\mathrm{in}}^{(n)},g_{\mathrm{in}}^{(n)})
 \subset\left\{v\in[L,R]:\mathbf d(v)\geq2^{-n-4}\right\},
       \qquad 1\leq n\leq \mathcal{N}.
\end{equation*}
\begin{equation}
 \norm{f_{\mathrm{in}}^{(n)}}_{\Gs{s}{\lambda_0}{\sigma+2}}
 +\norm{h_{\mathrm{in}}^{(n)}}_{\Gs{s}{\lambda_0}{\sigma+2}}
 +\norm{g_{\mathrm{in}}^{(n)}}_{\Gs{s}{\lambda_0}{\sigma+2}}
 \leq C\eps e^{-c_{\rm in}2^{n\alpha_*}},\qquad 1\leq n\leq \mathcal{N}.
 \label{eq:shell-data}
\end{equation}
Here $c_{\rm in}>0$ is chosen to equal $c_2$ in Proposition \ref{prop:short-time-weighted-propagation}. 

\subsection{Newton iteration and stream function decomposition}
\label{sec:newton}
Throughout the finite-system construction, all solution
layers satisfy $1\leq n\leq \mathcal{N}$; the estimates and parameter choices are
independent of $\mathcal{N}$.  Infinite sums of fixed shell weights may be used
as numerical majorants.  The local dyadic indices $N,N'$ in Fourier
decompositions are unrelated to this terminal layer.
We derive the increment equations from consecutive finite profile systems and construct full-line extensions of their nonzero-mode stream functions.

For a finite sequence \(q^{(n)}\), write
\begin{equation}
  q^{[n]}=\sum_{m=1}^n q^{(m)},
  \qquad q^{[0]}=0,\qquad 1\leq n\leq \mathcal{N}.
  \label{eq:partial-sum-notation}
\end{equation}
The terminal partial solution is
\begin{equation}
\begin{gathered}
  f^{[{ \mathcal N}]}=\sum_{n=1}^{ \mathcal N}f^{(n)},
  \qquad
  h^{[{ \mathcal N}]}=\sum_{n=1}^{ \mathcal N}h^{(n)},
  \qquad
  g^{[{ \mathcal N}]}=\sum_{n=1}^{ \mathcal N}g^{(n)},
\end{gathered}
  \label{eq:newton-series}
\end{equation}
with initial data for the \(n\)th increment given by
\begin{equation}
  \bigl(f^{(n)},h^{(n)},g^{(n)}\bigr)(1)
  =\bigl(f_{\mathrm{in}}^{(n)},h_{\mathrm{in}}^{(n)},g_{\mathrm{in}}^{(n)}\bigr).
  \label{eq:newton-data}
\end{equation}
Replacing $h$ by $h^{[n]}$ in \eqref{eq:variable-laplacian} gives
\begin{equation}
  \Delta_t^{[n]}
  =\partial_z^2+(1+h^{[n]})^2(\pvt)^2
   +(1+h^{[n]})(\partial_vh^{[n]})\pvt.
  \label{eq:partial-elliptic-operator}
\end{equation}
We define \(\Phi^{[n]}\) by the Dirichlet problem
\begin{equation}
  \Delta_t^{[n]}\Phi^{[n]}=f^{[n]},
  \qquad
  \Phi^{[n]}|_{v=L,R}=0,
  \qquad
  \Phi^{[0]}=0,
  \label{eq:partial-stream-function}
\end{equation}
and set
\begin{equation}
  \phi^{(n)}=\Phi^{[n]}-\Phi^{[n-1]}.
  \label{eq:stream-increment}
\end{equation}

To express the difference between consecutive elliptic operators, define the coefficient increments
\begin{align}
  Z^{(n)}
    &=-h^{(n)}\bigl(2+h^{(n)}+2h^{[n-1]}\bigr),
    \label{eq:Z-increment}\\
  H^{(n)}
    &=(1+h^{[n]})\partial_vh^{(n)}
      +h^{(n)}\partial_vh^{[n-1]}.
    \label{eq:H-increment}
\end{align}
Thus the consecutive elliptic operators satisfy
\begin{equation}
  \Delta_t^{[n]}-\Delta_t^{[n-1]}
  =-Z^{(n)}(\pvt)^2+H^{(n)}\pvt.
  \label{eq:elliptic-telescoping}
\end{equation}
Subtracting the two problems in \eqref{eq:partial-stream-function} therefore
gives the exact increment equation
\begin{equation}
  \Delta_t^{[n]}\phi^{(n)}
    =f^{(n)}+Z^{(n)}(\pvt)^2\Phi^{[n-1]}
      -H^{(n)}\pvt\Phi^{[n-1]},\qquad
  \phi^{(n)}|_{v=L,R}=0.
  \label{eq:elliptic-increment}
\end{equation}

\subsubsection{The \texorpdfstring{\(n\)-system}{n-system}}

Writing $\Phi_{\neq}^{[n]}=\Pneq\Phi^{[n]}$, define the partial velocity 
and its increment
\begin{align}
  \mathbf V^{[n]}
  =&\left(
      -(1+h^{[n]})\partial_v\Phi_{\neq}^{[n]},
      (1+h^{[n]})\partial_z\Phi_{\neq}^{[n]}+g^{[n]}
    \right)
  , \quad 
  \mathbf P^{(n)}=\mathbf V^{[n]}-\mathbf V^{[n-1]}.
  \label{eq:partial-transport-field}
\end{align}
Then \eqref{eq:stream-increment} gives
\[
 \mathbf P^{(n)}
 =(1+h^{[n]})\nabla_{z,v}^{\perp}\phi_{\neq}^{(n)}
   +h^{(n)}\nabla_{z,v}^{\perp}\Phi_{\neq}^{[n-1]}
   +\binom{0}{g^{(n)}},
 \qquad \nabla_{z,v}^{\perp}=(-\partial_v,\partial_z).
\]
The increments satisfy the following equations
\begin{subequations}\label{eq:n-system}
\begin{align}
  \partial_tf^{(n)}
    +\mathbf V^{[n]}\cdot\nabla_{z,v}f^{(n)}
    &=-\mathbf P^{(n)}\cdot\nabla_{z,v}f^{[n-1]},
    \label{eq:n-system-f}\\
  \partial_th^{(n)}
    +g^{[n]}\partial_vh^{(n)}
    -(\partial_vg^{[n]})h^{(n)}
    &=(1+h^{[n-1]})\partial_vg^{(n)}
      -g^{(n)}\partial_vh^{[n-1]},
    \label{eq:n-system-h}\\
  \partial_tg^{(n)}+\frac2t g^{(n)}
    +g^{[n]}\partial_vg^{(n)}
    &=-g^{(n)}\partial_vg^{[n-1]}
      +\frac1t\mathcal L_g^{(n)},
    \label{eq:n-system-g}
\end{align}
\end{subequations}
where 
\begin{align}
  \mathcal L_g^{(n)}
&=\Pzero\bigl(f_{\neq}^{(n)}\partial_z\Phi_{\neq}^{[n]}\bigr)
    +\Pzero\bigl(f_{\neq}^{[n-1]}\partial_z\phi_{\neq}^{(n)}\bigr).
    \label{eq:Lg-expanded}
\end{align}
Together with \eqref{eq:newton-data} and
\eqref{eq:elliptic-increment}, these equations form the \(n\)-system, and we have: 

\begin{lemma}[Equivalence]
\label{prop:newton-equivalence}
Let ${\mathcal N}\geq1$. If the increments $(f^{(n)},h^{(n)},g^{(n)})$ at levels $1\leq n\leq {\mathcal N}$ satisfy
\eqref{eq:newton-data}, \eqref{eq:elliptic-increment}, and \eqref{eq:n-system},
then $(f^{[{\mathcal N}]},\allowbreak h^{[{\mathcal N}]},\allowbreak g^{[{\mathcal N}]},\allowbreak\Phi^{[{\mathcal N}]})$ solves \eqref{eq:profile-system}
with initial datum $\sum_{n=1}^{\mathcal N}(f_{\mathrm{in}}^{(n)},h_{\mathrm{in}}^{(n)},g_{\mathrm{in}}^{(n)})$.
\end{lemma}

\subsubsection{Interior and exterior Green operators}
We will need full-line extensions of the stream functions to perform Fourier analysis. We define, for $k\ne0$:
\begin{equation}
 G_k^{}(v,w)=-\frac{\sinh\!\bigl(|k|(\min\left\{v,w\right\}-L)\bigr)
 \sinh\!\bigl(|k|(R-\max\left\{v,w\right\})\bigr)}
 {|k|\sinh(|k|(R-L))}.
 \label{eq:full-line-green-continuation}
\end{equation}
We incorporate
the shear by setting
\begin{equation}
  G_{k,t}^{}(v,w)=e^{\mathrm i kt(v-w)}G_k^{}(v,w).
  \label{eq:sheared-green-kernel}
\end{equation}
For a function $q\in C_c^\infty((L,R))$, we define the exterior and interior Green operators by\begin{subequations}
\begin{align}
  \mathscr G_{m,k}^{\Ext}q(v)
    &=\int_{\mathbb R}\widetilde\chi_m(v)G_{k,t}^{}(v,w)
       \chi_{m+2}(w)q(w)\dd w,
    \label{eq:exterior-green-operator}\\
  \mathscr G_{m,k}^{\Int}q(v)
    &=\int_{\mathbb R}\chi_{m+2}(v)G_{k,t}^{}(v,w)
       \widetilde\chi_m(w)q(w)\dd w.
    \label{eq:interior-green-operator}
\end{align}
\end{subequations}
These operators are only defined for $k \neq 0$, and when $k = 0$, we use the convention that these operators are identically zero. The extensions therefore contain only nonzero horizontal modes.

To state the channel equation compactly, we recall
$\pvt=\partial_v-t\partial_z$ from \eqref{eq:profile-f-phi} (or $\partial_v-\mathrm i kt$ on the
$k$th horizontal Fourier mode) and set
$\Delta_L=\partial_z^2+(\pvt)^2$. 
We define
\begin{equation}
 S^{(n)}=f^{(n)}+\mathcal C^{(n)}[\phi^{[n-1]}],\qquad
 \mathcal C^{(j)}[Q]=Z^{(j)}(\pvt)^2Q-H^{(j)}\pvt Q.
 \label{eq:coefficient-source}
\end{equation}
Here $Z^{(j)}$ and $H^{(j)}$ are defined in
\eqref{eq:Z-increment} and \eqref{eq:H-increment}, respectively.
We write $\mathcal C_k^{(j)}[Q_k]=(\mathcal C^{(j)}[Q])_k$ for the
$k$th horizontal Fourier coefficient.
The elliptic equation \eqref{eq:elliptic-increment} then becomes
\begin{equation}
 \Delta_L\phi^{(n)}=S^{(n)}+\sum_{j=1}^n\mathcal C^{(j)}[\phi^{(n)}].
 \label{eq:flat-increment-equation}
\end{equation}

We denote the inversion operator by
\begin{equation}
 \mathscr R_{k,t}q(v)=\int_L^R G_{k,t}^{}(v,w)q(w)\dd w.
 \label{eq:flat-dirichlet-resolvent}
\end{equation}
Using these notations, the relation \eqref{eq:flat-increment-equation} can be written as
\begin{equation}
 \phi_k^{(n)}=\mathscr R_{k,t}\Bigl(S_k^{(n)}+\sum_{j=1}^n\mathcal C_{k}^{(j)}[\phi_k^{(n)}]\Bigr).
 \label{eq:increment-fixed-point}
\end{equation}
We observe the following support properties:
\begin{equation}
 \supp S^{(n)}\Subset\left\{v\in I:\mathbf d(v)>d_n/32\right\},
 \qquad
 \supp Z^{(j)}\cup\supp H^{(j)}
 \Subset\left\{v\in I:\mathbf d(v)>d_j/32\right\}.
 \label{eq:extension-support-plateaux}
\end{equation}
Here $\Subset$ denotes compact inclusion. 
These conditions are guaranteed by the bootstrap support assumption
\eqref{eq:bootstrap-support-assumption}.

Define the exterior extension by
\begin{equation}
 \Phi_{\Ext,k}^{(n)}
 =\mathscr G_{n,k}^{\Ext}S_k^{(n)}
  +\sum_{j=1}^n\mathscr G_{j,k}^{\Ext}
       \mathcal C_{k}^{(j)}[\phi_k^{(n)}].
 \label{eq:exterior-extension-definition}
\end{equation}
At the endpoint, define
\begin{equation}
 \Phi_{\Int}^{(n;n)}:=\Phi_{\Ext}^{(n)}.
 \label{eq:exact-extension-endpoint}
\end{equation}
For $1\leq m<n$, define the mixed interior extension by
\begin{equation}
 \Phi_{\Int,k}^{(n;m)} ={}\mathscr G_{m,k}^{\Int}S_k^{(n)} +\sum_{j\leq m}\mathscr G_{j,k}^{\Ext} \mathcal C_{k}^{(j)}[\phi_k^{(n)}] +\sum_{m<j\leq n}\mathscr G_{j,k}^{\Int} \mathcal C_{k}^{(j)}[\phi_k^{(n)}].
 \label{eq:interior-extension-definition}
\end{equation}
Set $\Phi_{\Ext,0}^{(n)}=\Phi_{\Int,0}^{(n;m)}=0$.
The terms involving $\phi^{[n-1]}$ are already included in $S^{(n)}$.

\begin{lemma}[Extension consistency]
\label{lem:extension-consistency}
Assume \eqref{eq:extension-support-plateaux} and that
$\|h^{[n]}\|_\infty+\|\partial_vh^{[n]}\|_\infty$ is sufficiently small.
For $S_k^{(n)}\in L^2(I)$, the channel problem \eqref{eq:increment-fixed-point}
has a unique solution in the norm
\begin{equation}
 \|u\|_{H^2_{t,k}}=\|(\pvt)^2u\|_2+|k|\|\pvt u\|_2+k^2\|u\|_2,
 \label{eq:sheared-graph-norm}
\end{equation}
with norms over $I$ and bounds uniform in $t$ and $k\ne0$. Moreover,
\begin{equation}
\begin{aligned}  \Phi_{\Ext}^{(n)}&=\Pneq\phi^{(n)}&&\quad\hbox{on }[L,R],\\
\Phi_{\Int}^{(n;m)}&=\Pneq\phi^{(n)}&&\quad\hbox{on }\left\{v\in I:\mathbf d(v)>d_m/32\right\}\quad(m<n).
\end{aligned}
 \label{eq:extension-compatibility}
\end{equation} 
For $k\ne0$, the extensions satisfy the following identities, with $1\leq m<n$ in the second:
\begin{align}
 \Phi_{\Ext,k}^{(n)}
 &=\mathscr G_{n,k}^{\Ext}S_k^{(n)}
   +\sum_{j=1}^n\mathscr G_{j,k}^{\Ext}\mathcal C_{k}^{(j)}[\Phi_{\Ext,k}^{(n)}],
 \label{eq:exterior-derived-recursion}\\
 \Phi_{\Int,k}^{(n;m)}
 &=\mathscr G_{m,k}^{\Int}S_k^{(n)}
   +\sum_{j\leq m}\mathscr G_{j,k}^{\Ext}\mathcal C_{k}^{(j)}[\Phi_{\Int,k}^{(n;m)}]
   +\sum_{m<j\leq n}\mathscr G_{j,k}^{\Int}\mathcal C_{k}^{(j)}[\Phi_{\Int,k}^{(n;j)}].
 \label{eq:interior-derived-recursion}
\end{align}
\end{lemma}

\begin{proof}The argument is standard, and we omit the details.
\end{proof}

\begin{figure}[htbp]
\centering
\includegraphics[scale=0.5]{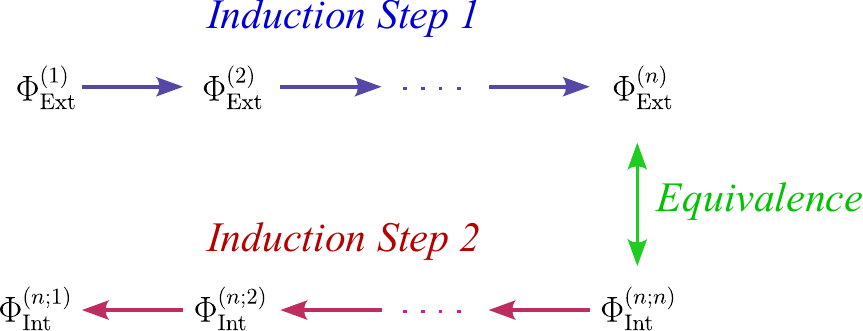}
\caption{Inductively defined elliptic operators $\pex{n}$ and $\pin{m}{n}$.}
\label{fig:ext-int-induction}
\end{figure}

Figure~\ref{fig:ext-int-induction} summarizes the induction process to define the elliptic extension operators.

\subsection{Parameters, weights, and bootstrap assumptions}
\label{sec:weights}
\label{sec:energy-setup}

We now choose Gevrey radii, Fourier multipliers, and energy functionals used in the long-time
argument.  The constants that appear below are chosen in the following order.  First fix \(s\),
\(\alpha_*\), \(K\), and hence \(\alpha_\dagger<\alpha_*\), as in
\eqref{eq:parameter-range} and \eqref{eq:alpha-dagger}.  Next choose parameters $a$ and $\mu$ satisfying 
\begin{equation}
  0<a\ll s-\frac12,
  \qquad a\leq1-s,
  \qquad \mu=54.
  \label{eq:weight-parameters}
\end{equation}
Take the Sobolev correction \(\sigma\) sufficiently large in terms of
these parameters.  Next fix \(c_*\in(0,c_{\rm in})\), where $c_{\rm in}$ is the constant from
\eqref{eq:shell-data}. We then choose the bootstrap constant \(M\gg1\) and, finally,
\(\eps>0\) sufficiently small. The parameter $\lambda_0>0$ is the Gevrey radius in the initial-data hypothesis.

\subsubsection{\texorpdfstring{Fourier multipliers}{Essential Fourier multipliers}}

For each $n\geq1$ choose the radius
\begin{equation}
  \lambda_n(t)=\lambda_{\mathrm{in},n}
    \left(\frac12+2^{a-1}(1+t)^{-a}\right),
  \qquad
  \lambda_{\mathrm{in},n}
    =\lambda_0\left(1-\sum_{j=1}^n2^{-j-2}\right),
  \label{eq:lambda-n}
\end{equation}
If $m<n$, then
\begin{equation}
  \lambda_m(t)-\lambda_n(t)
  \geq\frac{\lambda_0}{16}2^{-m}.
  \label{eq:radius-gap}
\end{equation}
Moreover, we note 
\begin{equation}
  -\dot\lambda_n(t)
  =a2^{a-1}\lambda_{\mathrm{in},n}(1+t)^{-1-a}>0.
  \label{eq:lambda-derivative}
\end{equation}

In this work, we will use the weights $w_k,w_{\NR},w_{\Res}$
introduced by Bedrossian--Masmoudi \cite[Section~3]{BM15}. Their explicit construction is given in Appendix~\ref{sec:orr-weights}.
Define 
\begin{align}
  {\J}_k(t,\eta)
    &=\frac{e^{\mu|\eta|^{1/2}}}{w_k(t,\eta)}
      +e^{\mu|k|^{1/2}},
  &\widetilde {\J}_k(t,\eta)
    &=\frac{e^{\mu|\eta|^{1/2}}}{w_k(t,\eta)},
  \label{eq:J-mode}\\
  {\J}_{\NR}(t,\eta)
    &=\frac{e^{\mu|\eta|^{1/2}}}{w_{\NR}(t,\eta)},
  &{\J}_{\Res}(t,\eta)
    &=\frac{e^{\mu|\eta|^{1/2}}}{w_{\Res}(t,\eta)},
  \label{eq:J-zero}
\end{align}
and
\begin{equation}
\begin{alignedat}{2}
  {\A}_k^{(n)}(t,\eta)
    &=e^{\lambda_n(t)|(k,\eta)|^s}
      \avg{k,\eta}^{\sigma}{\J}_k(t,\eta),
  &\quad\widetilde {\A}_k^{(n)}(t,\eta)
    &=e^{\lambda_n(t)|(k,\eta)|^s}
      \avg{k,\eta}^{\sigma}\widetilde {\J}_k(t,\eta),\\
  {\A}_{\NR}^{(n)}(t,\eta)
    &=e^{\lambda_n(t)|\eta|^s}\avg{\eta}^{\sigma}
      {\J}_{\NR}(t,\eta),
  &\quad{\A}_{\Res}^{(n)}(t,\eta)
    &=e^{\lambda_n(t)|\eta|^s}\avg{\eta}^{\sigma}
      {\J}_{\Res}(t,\eta).
\end{alignedat}
  \label{eq:main-multiplier}
\end{equation}

\subsubsection{Energy and Cauchy--Kovalevskaya functionals}
Define the energies on the $n$'th level: 
\begin{subequations}\label{eq:component-energies}
\begin{align}
  \E_f^{(n)}(t)
    &=\norm{{\A}^{(n)}f^{(n)}(t)}_2^2,
  \label{eq:energy-f}\\
  \E_{\partial g}^{(n)}(t)
    &=\avg{t}^{2+2s}
      \norm{{\A}_{\NR}^{(n)}\avg{\partial_v}^{-s}
        \partial_vg^{(n)}(t)}_2^2,
  \label{eq:energy-dg}\\
  \E_{g,\mathrm{low}}^{(n)}(t)
    &=\avg{t}^{\beta_g}
      \norm{g^{(n)}(t)}_{\Gs{s}{\lambda_n(t)}{\sigma-6}}^2,
  \label{eq:energy-g-low}\\
  \E_h^{(n)}(t)
    &=\norm{{\A}_{\Res}^{(n)}h^{(n)}(t)}_2^2.
  \label{eq:energy-h}
\end{align}
\end{subequations}
where the parameter $\beta_g=4-\frac a2$. We define the total energy at the $n$'th level: 
\begin{equation}
  \E^{(n)}
  =\E_f^{(n)}+\E_{\partial g}^{(n)}+\E_{g,\mathrm{low}}^{(n)}
    +M^{-2}\E_h^{(n)}.
  \label{eq:combined-energy}
\end{equation}

Define the $\CK$ terms associated to these energies via  \begin{subequations}
\begin{align}
  \CK_{f;\lambda}^{(n)}(t)
  &=\sum_{k\in\mathbb Z}\int_{\mathbb R}
    (-\dot\lambda_n)|(k,\eta)|^s
    \left|{\A}_k^{(n)}\widehat f_k^{(n)}(\eta)\right|^2\dd\eta,
  \label{eq:CK-f-lambda}\\
  \CK_{f;w}^{(n)}(t)
  &=\sum_{k\in\mathbb Z}\int_{\mathbb R}
    \frac{\partial_tw_k}{w_k}
    \left|\sqrt{\widetilde {\A}_k^{(n)}{\A}_k^{(n)}}\,
      \widehat f_k^{(n)}(\eta)\right|^2\dd\eta,
  \label{eq:CK-f-w}\\
  \CK_{f,+}^{(n)}
  &=\left[
      \CK_{f;\lambda}^{(n)}+\CK_{f;w}^{(n)}
      +(-\dot\lambda_n)\E_f^{(n)}\right],
  \label{eq:CK-f} \\
  \CK_{\partial g}^{(n)}(t)
  &= \avg{t}^{2+2s}\int_{\mathbb R}
 \left[(-\dot\lambda_n)|\eta|^s
       +\frac{\partial_tw_{\NR}}{w_{\NR}}\right]
 \left|\A_{\NR}^{(n)}\avg{\eta}^{-s}(\ii\eta)
       \widehat g^{(n)}(\eta)\right|^2\dd\eta.
 \label{eq:CK-dg} \\
 \CK_h^{(n)}(t) &=\int_{\mathbb R}
 \left[(-\dot\lambda_n)|\eta|^s
       +\frac{\partial_tw_{\Res}}{w_{\Res}}\right]
 \left|\A_{\Res}^{(n)}\widehat h^{(n)}(\eta)\right|^2\dd\eta.
 \label{eq:CK-h} \\
 \CK_{\mathrm{aux}}^{(n)}(t)
  &= \avg{t}^{2+2s}\int_{\mathbb R}
    \left[(-\dot\lambda_n)|\eta|^s(\A_{\NR}^{(n)})^2
      +\frac{\partial_tw_{\NR}}{w_{\NR}}(\widetilde\A_0^{(n)})^2\right]
    \frac{|\widehat g^{(n)}(\eta)|^2}{\avg\eta^{2s}}\dd\eta.
  \label{eq:CK-aux}
\end{align}
\end{subequations}
We note that since $\avg{k,\eta}^s\,\approx1+|(k,\eta)|^s$, we have the coercivity bound 
\begin{equation}
 (-\dot\lambda_n)\sum_{k\in\mathbb Z}\int_{\mathbb R}
 \avg{k,\eta}^{s}\left|{\A}_k^{(n)}(t,\eta)
             \widehat f_k^{(n)}(t,\eta)\right|^2\dd\eta
 \lesssim \CK_{f,+}^{(n)}.
 \label{eq:augmented-radius-coercivity}
\end{equation}

 We combine the $\CK$ terms defined above into the weighted $\CK$-term
\begin{equation}
  \CK_{f,\partial g,h}^{(n)}
  =\CK_{f,+}^{(n)}+\CK_{\partial g}^{(n)}
    +M^{-2}\CK_h^{(n)}.
  \label{eq:combined-CK}
\end{equation}

\subsubsection{The bootstrap proposition}

Define 
\begin{equation}
  \delta_n=e^{-c_*2^{n\alpha_*}}
  \label{eq:shell-envelope}
\end{equation}
For \(T\geq1\), define the total energy 
\begin{equation}
\begin{aligned}
  \mathcal B_n(T)
  :={}&\sup_{1\leq t\leq T}
    \left({\E_f}^{(n)}+\E_{g,\mathrm{low}}^{(n)}
      +\E_{\partial g}^{(n)}\right)(t)+\int_1^T\left(\CK_{f,+}^{(n)}
      +\CK_{\partial g}^{(n)}
      +t^{-1}\E_{\partial g}^{(n)}
      +t^{-1}\E_{g,\mathrm{low}}^{(n)}\right)\dd t\\
  &+M^{-2}\left[
      \sup_{1\leq t\leq T}\E_h^{(n)}(t)
      +\int_1^T\left(\CK_h^{(n)}
        +\CK_{\mathrm{aux}}^{(n)}\right)(t)\dd t
    \right].
\end{aligned}
  \label{eq:bootstrap-envelope}
\end{equation}

There is a constant $C_{\mathrm{in}}\geq1$,
depending only on the fixed parameters and independent of
$n,\mathcal N,M\geq1$, and $\eps$, such that
\[
 \mathcal B_n(1)\leq C_{\mathrm{in}}\eps^2\delta_n^2.
\]
Fix once and for all
\begin{equation}
 C_{\rm B}:=\frac{32}{a}e^{2\lambda_0}.
 \label{eq:bootstrap-constant}
\end{equation}
 Choose $M\geq1$ sufficiently large and such that 
$M^2\geq100C_{\rm B}C_{\mathrm{in}}$. Then
\begin{align}
 \mathcal B_n(1)
 \leq &\frac{M^2\eps^2}{100C_{\rm B}}\delta_n^2,
 \qquad 1\leq  n\leq\mathcal N, 
 \label{eq:initial-bootstrap-margin} \\
 \supp X_{\mathrm{in}}^{(n)}\subset &\left\{\mathbf d\geq d_n/16\right\},
 \qquad X\in\left\{f,g,h\right\}.
 \label{eq:support-scale}
\end{align}
With $C_{\rm B}$ and $M$ fixed in this order, choose
$\eps>0$ sufficiently small depending on $M$.
All these choices are independent of $\mathcal N$.

For solutions of the $n$-systems on $[1,T]$, assume for every $1\leq n\leq {\mathcal{N}}$:
\begin{align}
  \mathcal B_n(T)&\leq4M^2\eps^2\delta_n^2,
  \label{eq:bootstrap-energy-assumption}\\
  \supp X^{(n)}(t)&\subset
    \left\{v:\ \mathbf d(v)\geq\frac{5d_n}{96}\right\},
  \quad X\in\left\{f,g,h\right\},\quad 1\leq t\leq T.
  \label{eq:bootstrap-support-assumption}
\end{align}  

The following are consequences of \eqref{eq:bootstrap-energy-assumption}.

\begin{lemma}[Bootstrap consequences]
\label{lem:energy-bootstrap-consequences}
Under \eqref{eq:bootstrap-energy-assumption}, for every \(1\leq j\leq {\mathcal{N}}\)
and \(1\leq t\leq T\), for every fixed $c\in(0,1)$,
\begin{align}
  \norm{f^{(j)}}_{\Gs{s}{c\lambda_j}{\sigma-4}}
  +\avg t^{1+s}\norm{\avg{\partial_v}^{1-s}g^{(j)}}_{
       \Gs{s}{c\lambda_j}{\sigma-4}}
  +\avg t^{\beta_g/2}\norm{g^{(j)}}_{
       \Gs{s}{\lambda_j}{\sigma-6}}
  &\lesssim M\eps\delta_j,
  \label{eq:bootstrap-consequences-fg}\\
  \norm{h^{(j)}}_{\Gs{s}{c\lambda_j}{\sigma-4}}
  +\norm{Z^{(j)}}_{\Gs{s}{c\lambda_j}{\sigma-5}}
  +\norm{H^{(j)}}_{\Gs{s}{c\lambda_j}{\sigma-6}}
  &\lesssim M^C\eps\delta_j.
  \label{eq:bootstrap-consequences-coefficients}
\end{align}
Moreover, the corresponding partial sums with $1\leq n\leq {\mathcal{N}}$ satisfy
\begin{align}
  \norm{f^{[n]}}_{\Gs{s}{c\lambda_n}{\sigma-4}}
  +\norm{h^{[n]}}_{\Gs{s}{c\lambda_n}{\sigma-4}}
  &\lesssim M^C\eps,
  \label{eq:bootstrap-consequences-partial-fh}\\
  \norm{\partial_vg^{[n]}}_{L^\infty}
  &\lesssim \frac{M^C\eps}{\avg t^{1+s}},
  \qquad
  \norm{g^{[n]}}_{\Gs{s}{c\lambda_n}{\sigma-6}}
  \lesssim\frac{M^C\eps}{\avg t^{\beta_g/2}}.
  \label{eq:bootstrap-consequences-partial-g}
\end{align}
All constants are independent of \(j,n,{\mathcal{N}}\) and \(T\).
The weaker-radius constants may depend on the chosen fixed $c$.
\end{lemma}

\begin{proof}
The proof follows immediately from the definitions. 
\end{proof}

\begin{proposition}[Bootstrap improvement]
\label{prop:bootstrap-improvement}
There exist \(M_0\gg1\) and, for every \(M\geq M_0\), an
\(\eps_0=\eps_0(M)>0\) with the following property.  Let \(T\geq1\), and
suppose that the \(n\)-systems have solutions on \([1,T]\) satisfying \eqref{eq:bootstrap-energy-assumption}--\eqref{eq:bootstrap-support-assumption} for every $1\leq n\leq {\mathcal{N}}$.
If \(0<\eps\leq\eps_0\), then the strict improvements
\begin{align}
  \mathcal B_n(T)&\leq2M^2\eps^2\delta_n^2,
  \label{eq:bootstrap-energy-improvement}\\
  \supp X^{(n)}(t)&\subset
    \left\{v:\ \mathbf d(v)\geq\frac{5d_n}{88}\right\}
  \label{eq:bootstrap-support-improvement}
\end{align}
hold for every \(1\leq n\leq {\mathcal{N}}\), every \(X\in\left\{f,g,h\right\}\), and every
\(1\leq t\leq T\). The constants $M_0$ and $\eps_0(M)$ are independent of ${\mathcal{N}},T$.
\end{proposition}

The improvement from \(4\) to \(2\), together with
\eqref{eq:initial-bootstrap-margin}, implies by a continuity argument that
\eqref{eq:bootstrap-energy-improvement} holds for every \(T\in[1,\infty)\). Standard arguments in Section~\ref{sec:convergence} then yield
Theorem~\ref{thm:main}.

\section{Elliptic estimates}
\label{sec:elliptic}
In this section, we develop the necessary estimates for the exterior
and interior stream functions $\Phi_{\Ext}^{(n)}$ and $\Phi_{\Int}^{(n;m)}$, \eqref{eq:exterior-extension-definition} - \eqref{eq:interior-extension-definition}.
We work with the finite family $1\leq n\leq {\mathcal{N}}$. All constants below are independent of ${\mathcal{N}}$.

Recall that $s\in(1/2,1)$ is the Gevrey exponent and that
$a>0$ is chosen with $5a/4\leq2s-1$ and $a\leq1-s$;
the low-frequency time exponent is $\beta_g=4-a/2$.
Fix, throughout this section and Section~\ref{sec:energy},
\begin{equation}
  0<\delta<\min\left\{a,\frac{2s-1}{4},
    \frac{4-\beta_g}{8}\right\}.
  \label{eq:energy-decay-delta}
\end{equation}

Fix once and for all a sufficiently large constant $C_{\mathrm E}$
to dominate the exponential localization factors in
Lemmas~\ref{lem:lossy-green-calculus},
\ref{lem:precision-green-calculus}, and
\ref{lem:newton-velocity-output} below.  With
$\delta_j=e^{-c_*2^{j\alpha_*}}$ from \eqref{eq:shell-envelope} and
$\alpha_\dagger$ from \eqref{eq:alpha-dagger}, put
\begin{equation}
\varsigma_j=\exp\bigl(C_{\mathrm E}2^{j\alpha_\dagger}\bigr)
  .
  \label{eq:elliptic-shell-weights}
\end{equation}
The following will be used repeatedly, 
\begin{equation}
  \sum_{j\geq1}\varsigma_j\delta_j
  <\infty.
  \label{eq:elliptic-shell-weight-sums}
\end{equation}

We use the following low and high-regularity norms 
for the coefficient increments $Z^{(j)},H^{(j)}$ of
\eqref{eq:Z-increment}--\eqref{eq:H-increment},
\begin{equation}
 c_{Z,H;\mathrm{low}}^{(j)}
 =\|Z^{(j)}\|_{\Gs{s}{\lambda_n}{2}}
   +\|H^{(j)}\|_{\Gs{s}{\lambda_n}{2}},\qquad
 c_{Z,H;\mathrm{high}}^{(j)}
 =\|Z^{(j)}\|_{\Gs{s}{\lambda_n}{\sigma-1}}
   +\|H^{(j)}\|_{\Gs{s}{\lambda_n}{\sigma-1}}.
 \label{eq:elliptic-loss-sequences}
\end{equation}

\subsection{\texorpdfstring{Main estimates and proof strategy}{Main estimates and proof strategy}}
\label{sec:elliptic-statements}

The goal of this section is to state the lossy bounds in
Proposition~\ref{prop:lossy-elliptic} and the precision bounds in
Proposition~\ref{prop:precision-elliptic}.  Recall from \eqref{eq:coefficient-source} that the source term has the form $S^{(n)}=f^{(n)}+Z^{(n)}(\partial_v^t)^2\phi^{[n-1]}
-H^{(n)}\partial_v^t\phi^{[n-1]}$, with coefficients
\eqref{eq:Z-increment}--\eqref{eq:H-increment}. To measure this source in low regularity, we use
\begin{equation}
  \,\mathfrak S_n(t)
  :=\norm{\Pneq S^{(n)}(t)}_{\Gs{s}{\lambda_n(t)}{\sigma-1}}.
  \label{eq:elliptic-source-norm}
\end{equation}
Our first main proposition reads as follows.  
\begin{proposition}[Lossy exterior and interior estimates]
\label{prop:lossy-elliptic}
Under the bootstrap assumptions~\eqref{eq:bootstrap-energy-assumption}--\eqref{eq:bootstrap-support-assumption}, for $\eps$ sufficiently small and every
$1\leq m\leq n\leq {\mathcal{N}}$,
\begin{align}
  \norm{\Delta_L\Pneq\Phi_{\Ext}^{(n)}}_{
      \Gs{s}{\lambda_n}{\sigma-1}}
  &\lesssim e^{C2^{n\alpha_\dagger}}\mathfrak S_n,
  \label{eq:lossy-source-ext}\\
  \norm{\Delta_L\Pneq\Phi_{\Int}^{(n;m)}}_{
      \Gs{s}{\lambda_n}{\sigma-1}}
  &\lesssim e^{C2^{m\alpha_\dagger}}\mathfrak S_n.
  \label{eq:lossy-source-int}
\end{align}
Consequently,
\begin{align}
  \norm{\Pneq\Phi_{\Ext}^{(n)}(t)}_{
      \Gs{s}{\lambda_n(t)}{\sigma-3}}
  &\lesssim
    \frac{e^{C2^{n\alpha_\dagger}}}{1+t^2}
    \mathfrak S_n(t),\label{eq:lossy-ext}\\
  \norm{\Pneq\Phi_{\Int}^{(n;m)}(t)}_{
      \Gs{s}{\lambda_n(t)}{\sigma-3}}
  &\lesssim
    \frac{e^{C2^{m\alpha_\dagger}}}{1+t^2}
    \mathfrak S_n(t).\label{eq:lossy-int}
\end{align}
The corresponding exterior and interior velocity bounds are
\begin{equation}
\begin{aligned}
  \norm{\pvt\Pneq\Phi_{\Ext}^{(n)}}_{\Gs{s}{\lambda_n}{\sigma-4}}
  +\avg t\norm{\partial_z\Pneq\Phi_{\Ext}^{(n)}}_{\Gs{s}{\lambda_n}{\sigma-4}}
  &\lesssim \frac{e^{C2^{n\alpha_\dagger}}}{\avg t}\mathfrak S_n,\\
  \norm{\pvt\Pneq\Phi_{\Int}^{(n;m)}}_{\Gs{s}{\lambda_n}{\sigma-4}}
  +\avg t\norm{\partial_z\Pneq\Phi_{\Int}^{(n;m)}}_{\Gs{s}{\lambda_n}{\sigma-4}}
  &\lesssim \frac{e^{C2^{m\alpha_\dagger}}}{\avg t}\mathfrak S_n,
\end{aligned}
\label{eq:lossy-velocity}
\end{equation}
and
\begin{equation}
\begin{aligned}
  \norm{\partial_v\Pneq\Phi_{\Ext}^{(n)}}_{\Gs{s}{\lambda_n}{\sigma-5}}
  &\lesssim \frac{e^{C2^{n\alpha_\dagger}}}{1+t^2}\mathfrak S_n, \qquad \norm{\partial_v\Pneq\Phi_{\Int}^{(n;m)}}_{\Gs{s}{\lambda_n}{\sigma-5}}
  \lesssim \frac{e^{C2^{m\alpha_\dagger}}}{1+t^2}\mathfrak S_n.
\end{aligned}
\label{eq:lossy-profile-velocity}
\end{equation}
\end{proposition}

By the Gevrey product estimate \cite[Lemma~A.3]{BM15},
\begin{align}
  \mathfrak S_n
  &\lesssim \norm{\Pneq f^{(n)}}_{\Gs{s}{\lambda_n}{\sigma-1}}
   +\norm{Z^{(n)}}_{\Gs{s}{\lambda_n}{\sigma-1}}
      \norm{(\pvt)^2\Phi_{\Ext,\neq}^{[n-1]}}_{\Gs{s}{\lambda_n}{\sigma-1}}
  \notag\\
  &\qquad
   +\norm{H^{(n)}}_{\Gs{s}{\lambda_n}{\sigma-1}}
      \norm{\pvt\Phi_{\Ext,\neq}^{[n-1]}}_{\Gs{s}{\lambda_n}{\sigma-1}}.
  \label{eq:elliptic-source-expanded}
\end{align}
Along with the proposition,
we will prove the source bound
\begin{equation}
  \mathfrak S_n(t)\lesssim M^C\eps\delta_n.
  \label{eq:elliptic-source-bootstrap}
\end{equation}

For the precision estimate, recall ${\A}_k^{(n)},\widetilde {\A}_k^{(n)}$
from \eqref{eq:main-multiplier}, $w_k$ from
\eqref{eq:mode-dependent-weight}, and $-\dot\lambda_n$ from
\eqref{eq:lambda-derivative}.  Define two multipliers on nonzero modes by\begin{subequations}
\begin{align}
  \mathcal M_{\lambda,k}^{(n)}(t,\eta)
  &=\sqrt{-\dot\lambda_n(t)}\,|(k,\eta)|^{s/2}
    \left\langle\frac{\eta}{kt}\right\rangle^{-1}
    {\A}_k^{(n)}(t,\eta),
  \label{eq:elliptic-M-lambda}\\
  \mathcal M_{w,k}^{(n)}(t,\eta)
  &=\sqrt{\frac{\partial_tw_k(t,\eta)}{w_k(t,\eta)}}
    \left\langle\frac{\eta}{kt}\right\rangle^{-1}
    \widetilde {\A}_k^{(n)}(t,\eta).
  \label{eq:elliptic-M-w}
\end{align}
\end{subequations}

\begin{proposition}[Precision elliptic control]
\label{prop:precision-elliptic}
Under the hypotheses of Proposition~\ref{prop:lossy-elliptic}, the following bounds hold:
\begin{equation}
  \sum_{\iota=\lambda,w}\norm{\mathcal M_\iota^{(n)}\Delta_L \Phi_{\Ext}^{[n]}}_2^2
  \lesssim
    \sum_{j\leq n}\delta_j^{-1}\varsigma_j\CK_{f,\partial g,h}^{(j)}
    +\frac{M^C\eps^2}{\avg t^{1+\delta}}.
  \label{eq:precision-elliptic-partial-sum}
\end{equation} \space For $1\leq m\leq n$,
\begin{equation}
  \sum_{\iota=\lambda,w}\norm{\mathcal M_\iota^{(n)}\Delta_L \Phi_{\Int}^{(n;m)}}_2^2
  \lesssim
    \varsigma_m\CK_{f,\partial g,h}^{(n)}
    +\varsigma_m\delta_n^2\sum_{j<n}
       \delta_j^{-1}\varsigma_j\CK_{f,\partial g,h}^{(j)}
    +\frac{M^C\eps^2\varsigma_m\delta_n^2}
       {\avg t^{1+\delta}}.
  \label{eq:precision-elliptic-increment}
\end{equation} \space The sources satisfy
\begin{align}
  \sum_{\iota=\lambda,w}\left\|\mathcal M_\iota^{(n)}\Pneq
                  \sum_{j\leq n}S^{(j)}\right\|_2^2
  &\lesssim
    \sum_{j\leq n}\delta_j^{-1}\varsigma_j\CK_{f,\partial g,h}^{(j)}
    +\frac{M^C\eps^2}{\avg t^{1+\delta}},
  \label{eq:precision-source-partial-sum}\\
  \sum_{\iota=\lambda,w}\|\mathcal M_\iota^{(n)}\Pneq S^{(n)}\|_2^2
  &\lesssim
    \varsigma_n\CK_{f,\partial g,h}^{(n)}
    +\varsigma_n\delta_n^2\sum_{j<n}
       \delta_j^{-1}\varsigma_j\CK_{f,\partial g,h}^{(j)}
    +\frac{M^C\eps^2\varsigma_n\delta_n^2}{\avg t^{1+\delta}}.
  \label{eq:precision-source-increment}
\end{align} 
\end{proposition}

The section is organized as follows.
Subsection~\ref{sec:elliptic-kernel} collects the localized Dirichlet kernel bounds.
Subsection~\ref{sec:elliptic-green-bounds} proves the lossy Green estimates.
Subsection~\ref{sec:elliptic-precision-green-proof} proves the precision Green estimates.
Subsection~\ref{sec:elliptic-lossy-proof} proves \eqref{eq:elliptic-source-bootstrap} and Proposition~\ref{prop:lossy-elliptic}.
Subsection~\ref{sec:elliptic-precision-proof} proves Proposition~\ref{prop:precision-elliptic}.

\subsection{\texorpdfstring{The Dirichlet kernel and its extension}{The Dirichlet kernel and its extension}}
\label{sec:elliptic-kernel}

We establish the localized kernel and weighted operator bounds needed
for Lemmas~\ref{lem:lossy-green-calculus} and~\ref{lem:precision-green-calculus}.
These estimates are independent of the Newton recursion.  

 Let $\widehat{\mathscr G}_{m,k}^{*}(t,\eta,\xi)$ denote the Fourier kernel
in $v,w$ of either of the operators in
\eqref{eq:exterior-green-operator}--\eqref{eq:interior-green-operator},
where $*\in\left\{\Ext,\Int\right\}$:
\[
 \widehat{\mathscr G_{m,k}^{*}q_k}(\eta)
 =\int_{\R}\widehat{\mathscr G}_{m,k}^{*}(t,\eta,\xi)
                  \widehat q_k(\xi)\dd\xi.
\]
Write
\begin{equation}
  \lambda_m^\Delta=c_\Delta2^{-m\rho},
  \qquad 0<c_\Delta\ll1,
  \qquad \rho=\frac{K+s}{K+1}.
  \label{eq:kernel-localization-radius}
\end{equation}

\begin{lemma}[Localized Green-kernel estimate]
\label{lem:localized-green-kernel}
For $k\neq0$ and $m\geq1$,
\begin{align}
  |\widehat{\mathscr G}_{m,k}^{*}(t,\eta,\xi)|
  &\lesssim 2^{2m}
  \min\left\{
    \frac1{k^2+|\xi-kt|^2},
    \frac1{k^2+|\eta-kt|^2}
  \right\}
  e^{-\lambda_m^\Delta|\eta-\xi|^\rho},
  \label{eq:localized-green-kernel}\\
 (k^2+|\eta-kt|^2) |\widehat{\mathscr G}_{m,k}^{*}(t,\eta,\xi)| &\lesssim 2^{2m} \left(1+\frac{|\eta-\xi|^2} {k^2+|\xi-kt|^2}\right) e^{-\lambda_m^\Delta|\eta-\xi|^\rho}.
 \label{eq:localized-green-kernel-delta}
\end{align}
\end{lemma}

\begin{proof} With an explicit expression of the Green's function \eqref{eq:full-line-green-continuation}, the proof is similar to the localized Fourier-kernel argument of
Sinambela--Zhao--Zi \cite{sinambela2026asymptotic}, with the dependence on $d_m$ kept explicit, see also \cite{IonescuJia20inviscid}. We omit the details for the sake of brevity.
\end{proof}

\Needspace{12\baselineskip}
\begin{lemma}
\label{lem:weighted-common-scale-green}
Recall $\mathscr G_{m,k}^{*}$ from
\eqref{eq:exterior-green-operator}--\eqref{eq:interior-green-operator}
and $\Delta_{L,k}=(\partial_v-ikt)^2-k^2$.
We have the following identity: 
\begin{align}
 \Delta_{L,k}\mathscr G_{m,k}^{\Ext}q_{k}
 &=\widetilde\chi_m\chi_{m+2}q_{k}
   +2\widetilde\chi_m'(\partial_v-ikt)\!
     \int_{\mathbb R}G_{k,t}^{}(\,\cdot,w)\chi_{m+2}(w)q_{k}(w)\dd w
 \notag\\*
 &\quad
   +\widetilde\chi_m''
     \int_{\mathbb R}G_{k,t}^{}(\,\cdot,w)\chi_{m+2}(w)q_{k}(w)\dd w,
 \label{eq:exact-exterior-cap-residual}\\
 \Delta_{L,k}\mathscr G_{m,k}^{\Int}q_{k}
 &=\chi_{m+2}\widetilde\chi_mq_{k}
   +2\chi_{m+2}'(\partial_v-ikt)\!
     \int_{\mathbb R}G_{k,t}^{}(\,\cdot,w)\widetilde\chi_m(w)q_{k}(w)\dd w
 \notag\\*
 &\quad
   +\chi_{m+2}''
     \int_{\mathbb R}G_{k,t}^{}(\,\cdot,w)\widetilde\chi_m(w)q_{k}(w)\dd w.
 \label{eq:exact-interior-cap-residual}
\end{align}
The following bounds are satisfied (recall the definitions \eqref{eq:elliptic-M-lambda}, \eqref{eq:elliptic-M-w})\begin{subequations}
\begin{align}
 \norm{\mathcal M_\lambda^{(n)}\Delta_L\mathscr G_m^*q}_2
 &\leq Ce^{C2^{m\alpha_\dagger}}
       \norm{\mathcal M_\lambda^{(n)}q}_2,
 \label{eq:weighted-common-lambda}\\
 \norm{\mathcal M_w^{(n)}\Delta_L\mathscr G_m^*q}_2
 &\leq Ce^{C2^{m\alpha_\dagger}}
   \left(\norm{\mathcal M_w^{(n)}q}_2
        +\norm{\mathcal M_\lambda^{(n)}q}_2\right).
 \label{eq:weighted-common-w}
\end{align}
\end{subequations}
\end{lemma}

\begin{proof}
Identities \eqref{eq:exact-exterior-cap-residual} and
\eqref{eq:exact-interior-cap-residual} follow from Leibniz' rule. Lemma~\ref{lem:localized-green-kernel} gives,
\begin{equation}
 |\widehat{\Delta_{L,k}\mathscr G_{m,k}^{*}}(t,\eta,\xi)|
 \leq C2^{2m}\avg{\eta-\xi}^{2}
       e^{-\lambda_m^\Delta|\eta-\xi|^\rho}.
 \label{eq:common-delta-kernel-envelope}
\end{equation}
By a standard commutator estimate (see \cite{BM13}) with a
polynomial $P$ whose degree and coefficients are independent of $t,m,n,k$,
\begin{equation}
 \mathcal M_\lambda^{(n)}(k,\eta)
 \leq P(\avg{\eta-\xi})e^{\lambda_n|\eta-\xi|^s+C\mu|\eta-\xi|^{1/2}}
       \mathcal M_\lambda^{(n)}(k,\xi).
 \label{eq:common-mlambda-transfer}
\end{equation}

Next, we move on to the bound of the  $\mathcal{M}_w^{(n)}$-multiplier.  First suppose that $|\eta|\approx|\xi|$.  We apply \eqref{eq:raw-w-endpoint-transfer} and combine it with
Lemma~\ref{lem:reduced-J}, which gives
\begin{align}
 \mathcal M_w^{(n)}(k,\eta)
 &\leq P(\avg{\eta-\xi})e^{\lambda_n|\eta-\xi|^s+C\mu|\eta-\xi|^{1/2}}
 \left[\mathcal M_w^{(n)}(k,\xi)
   +\avg{\xi/(kt)}^{-1}\frac{|(k,\xi)|^{s/2}}{\avg t^s}
        {\A}_k^{(n)}(t,\xi)\right]
 \notag\\
 &\leq P(\avg{\eta-\xi})e^{\lambda_n|\eta-\xi|^s+C\mu|\eta-\xi|^{1/2}}
       \bigl(\mathcal M_w^{(n)}(k,\xi)
             +\mathcal M_\lambda^{(n)}(k,\xi)\bigr),
 \label{eq:common-mw-transfer}
\end{align}
where the second line uses \eqref{eq:w-endpoint-paid-by-lambda}.

If the $|\eta|$ and $|\xi|$ are not comparable, then
$|\eta-\xi|\geq c\max\left\{|\eta|,|\xi|\right\}$.  The assertion is immediate when
$\frac{\partial_tw_k(t,\eta)}{w_k(t,\eta)}=0$.  Otherwise, $\frac{\partial_tw_k(t,\eta)}{w_k(t,\eta)}\lesssim1$, and the active time range gives
$t\leq2|\eta|\lesssim\avg{\eta-\xi}$.  Lemma~\ref{lem:reduced-J} transfers
$\widetilde {\A}_k(t,\eta)$ to ${\A}_k(t,\xi)$ with the same exponential.
Moreover,
\[
 \avg{\xi/(kt)}\leq C\avg{\eta-\xi},
 \qquad
 (-\dot\lambda_n)^{-1/2}
 \lesssim\avg t^{(1+a)/2}\lesssim\avg{\eta-\xi}^{(1+a)/2}.
\]
Since $|(k,\xi)|\geq1$, these polynomial factors allow us to control the output
by $\mathcal M_\lambda^{(n)}(k,\xi)$ and establish
\eqref{eq:common-mw-transfer} in all cases.

Combining \eqref{eq:common-delta-kernel-envelope} and
\eqref{eq:common-mlambda-transfer} or
\eqref{eq:common-mw-transfer} gives for the $w$ multiplier,
\begin{align*}
 &\left|\mathcal M_w^{(n)}(k,\eta)
       \widehat{\Delta_L\mathscr G_m^*q}_k(\eta)\right|\\
 &\quad\leq C2^{2m}\int_{\R}\avg{\eta-\xi}^N
 \exp\!\left[-\lambda_m^\Delta|\eta-\xi|^\rho
       +\lambda_n|\eta-\xi|^s+C\mu|\eta-\xi|^{1/2}\right]\\
 &\qquad\qquad\times
       \bigl[\mathcal M_w^{(n)}(k,\xi)+\mathcal M_\lambda^{(n)}(k,\xi)\bigr]
       |\widehat q_k(\xi)|\dd\xi.
\end{align*}
For the $\lambda$ multiplier, the bracket contains only
$\mathcal M_\lambda^{(n)}(k,\xi)$.
The optimization
\[
 \lambda_n r^s+C \mu r^{1/2}-\lambda_m^\Delta{}r^\rho+N\log\avg r
 \leq C2^{m\alpha_\dagger}
      -\frac14\lambda_m^\Delta{}r^\rho
\]
uses $s>1/2$ and
$\alpha_\dagger=\rho s/(\rho-s)$.  The remaining kernel has
$L^1(\mathbb R)$ norm $O(2^m)$, absorbed in $e^{C2^{m\alpha_\dagger}}$.
Young's inequality at each fixed $k$, followed by the $\ell_k^2$ norm, gives
\begin{align*}
 \left(\sum_{k\ne0}
   \|\mathcal M_w^{(n)}(k,\cdot)
      \widehat{\Delta_L\mathscr G_m^*q}_k\|_{L^2_\eta}^2\right)^{1/2}
 &\lesssim e^{C2^{m\alpha_\dagger}}
   \left(\sum_{k\ne0}
   \|[\mathcal M_w^{(n)}(k,\cdot)+\mathcal M_\lambda^{(n)}(k,\cdot)]
      |\widehat q_k|\|_{L^2_\xi}^2\right)^{1/2}\\
 &\lesssim e^{C2^{m\alpha_\dagger}}
       \bigl(\|\mathcal M_w^{(n)}q\|_2+\|\mathcal M_\lambda^{(n)}q\|_2\bigr).
\end{align*}
The same argument with only $\mathcal M_\lambda^{(n)}q$ proves the radius bound.
These are
\eqref{eq:weighted-common-lambda}--
\eqref{eq:weighted-common-w}.
\end{proof}

\begin{remark} The above proof used the following observation, which we now restate in a form that will be used in the future. The choice $\rho>s$ allows us to convert kernel regularity into the following
explicit shell loss.  For $r\geq0$,
\begin{equation}
 \lambda_n r^s-\frac12\lambda_m^\Delta r^\rho \leq \sup_{x\geq0} \left(\lambda_0x^s-\frac12c_\Delta2^{-m\rho}x^\rho\right) \leq C \left(2^{m\rho}\right)^{s/(\rho-s)} =C2^{m\rho s/(\rho-s)} =C2^{m\alpha_\dagger}.
 \label{eq:kernel-to-shell-loss}
\end{equation}
Here
\(\rho s/(\rho-s)=s(K+s)/(K(1-s))=\alpha_\dagger\).
Consequently,
\begin{equation}
  e^{\lambda_n r^s+C\mu r^{1/2}}
    e^{-\lambda_m^\Delta r^\rho}
  \lesssim
  e^{C2^{m\alpha_\dagger}}
    e^{-\frac14\lambda_m^\Delta r^\rho},
  \label{eq:kernel-absorption}
\end{equation}
where $s>1/2$ also absorbs the square-root exponential in the multiplier
comparison.
\end{remark}

\subsection{\texorpdfstring{Lossy and precision Green estimates}{Lossy and precision Green estimates}}
\label{sec:elliptic-green-bounds}

We reduce the operator estimates for
Propositions~\ref{prop:lossy-elliptic} and~\ref{prop:precision-elliptic}
to three model terms, whose lossy and precision bounds are stated in
Lemmas~\ref{lem:lossy-green-calculus} and~\ref{lem:precision-green-calculus}.
For Proposition~\ref{prop:precision-elliptic}, we also need the coefficient
$\CK$-estimate stated below in Lemma \ref{lem:coefficient-CK-product}.

Recall the operators
\eqref{eq:exterior-green-operator}--\eqref{eq:interior-green-operator},
the coefficients \eqref{eq:Z-increment}--\eqref{eq:H-increment}, and
$\partial_v^t=\partial_v-t\partial_z$.
We introduce the three model terms
\begin{equation}
  \mathcal T_{m;0}^{*}[q]=\mathscr G_m^{*}q,\qquad
  \mathcal T_{m;1}^{(j),*}[Z,Q]=\mathscr G_m^{*}\bigl(Z^{(j)}(\pvt)^2Q\bigr),\qquad
  \mathcal T_{m;2}^{(j),*}[H,Q]=\mathscr G_m^{*}\bigl(H^{(j)}\pvt Q\bigr).
  \label{eq:model-terms}
\end{equation}
The subscript $m$ denotes the spatial cutoff scale, the superscript $(j)$
denotes the coefficient layer, and $*$ specifies the exterior or interior
operator. 

The coefficient $\CK$ terms are
\begin{align}
 \mathcal{CCK}_Z^{(j;n)}
 &=\left[(-\dot\lambda_n)\norm{|\partial_v|^{s/2}{\A}_{\Res}^{(n)}Z^{(j)}}_2^2
 +\norm{\sqrt{\frac{\partial_tw_{\Res}}{w_{\Res}}}{\A}_{\Res}^{(n)}Z^{(j)}}_2^2\right],
 \label{eq:CCK-Z}\\
 \mathcal{CCK}_H^{(j;n)}
 &=\left[(-\dot\lambda_n)\norm{|\partial_v|^{s/2}{\A}_{\Res}^{(n)}\avg{\partial_v}^{-1}H^{(j)}}_2^2
 +\norm{\sqrt{\frac{\partial_tw_{\Res}}{w_{\Res}}}{\A}_{\Res}^{(n)}\avg{\partial_v}^{-1}H^{(j)}}_2^2\right],
 \label{eq:CCK-H}\\
 \mathcal{CCK}_{Z,H}^{(j;n)}
 &=\mathcal{CCK}_Z^{(j;n)}+\mathcal{CCK}_H^{(j;n)}.
 \label{eq:CCK-ZH}
\end{align}

\begin{lemma}[Coefficient $\CK$ product estimate]
\label{lem:coefficient-CK-product}
Under \eqref{eq:bootstrap-energy-assumption}, for $1\leq j\leq n$,
\begin{equation}
 \mathcal{CCK}_{Z,H}^{(j;n)}
 \lesssim M^C\left(
 \CK_{f,\partial g,h}^{(j)}+\eps^2\delta_j^2
 \sum_{r<j}\varsigma_r\CK_{f,\partial g,h}^{(r)}\right).
 \label{eq:coefficient-CK-product-bound}
\end{equation}
\end{lemma}
\begin{proof} This is a straightforward adaptation from \cite{BM13}, and we omit the proof for brevity.
\end{proof}

\begin{lemma}[Lossy Green bounds]
\label{lem:lossy-green-calculus}
For arbitrary $m\geq1$,
\begin{gather}
  \norm{\Delta_L\mathcal T_{m;0}^{*}[q]}_{\Gs{s}{\lambda_n}{\sigma-1}}
  \lesssim e^{C2^{m\alpha_\dagger}}
    \norm{q}_{\Gs{s}{\lambda_n}{\sigma-1}},
  \label{eq:T0-lossy}\\
  \begin{aligned}
  &\norm{\Delta_L\mathcal T_{m;1}^{(j),*}[Z,Q]}_{\Gs{s}{\lambda_n}{\sigma-1}}
   +\norm{\Delta_L\mathcal T_{m;2}^{(j),*}[H,Q]}_{\Gs{s}{\lambda_n}{\sigma-1}}
  \\
  &\quad\lesssim e^{C2^{m\alpha_\dagger}}
    \bigl(\norm{Z^{(j)}}_{\Gs{s}{\lambda_n}{\sigma-1}}
          +\norm{H^{(j)}}_{\Gs{s}{\lambda_n}{\sigma-1}}\bigr)
    \norm{\Delta_LQ}_{\Gs{s}{\lambda_n}{\sigma-1}}.
  \end{aligned}
  \label{eq:T12-lossy}
\end{gather}
\end{lemma}

\begin{lemma}[Precision Green bounds]
\label{lem:precision-green-calculus}
Each sum below contains exactly two terms, indexed by
${\iota}\in\left\{\lambda,w\right\}$, with the precision multipliers
$\mathcal M_\lambda^{(n)}$ and $\mathcal M_w^{(n)}$ defined
in \eqref{eq:elliptic-M-lambda}--\eqref{eq:elliptic-M-w}.  We have
\begin{align}
 &\sum_{{\iota}\in\left\{\lambda,w\right\}} \norm{\mathcal M_{\iota}^{(n)}\Delta_L\mathcal T_{m;0}^{*}[q]}_2^2 \lesssim e^{C2^{m\alpha_\dagger}} \sum_{{\iota}\in\left\{\lambda,w\right\}}\norm{\mathcal M_{\iota}^{(n)} q}_2^2,
 \label{eq:T0-precision}\\
  &\sum_{{\iota}\in\left\{\lambda,w\right\}}
  \norm{\mathcal M_{\iota}^{(n)}\Delta_L\mathcal T_{m;1}^{(j),*}[Z,Q]}_2^2
  \notag\\
  &\quad\lesssim e^{C2^{m\alpha_\dagger}}
  \left[
    \norm{Z^{(j)}}_{\Gs{s}{\lambda_n}{2}}^2
       \sum_{\iota=\lambda,w}\norm{\mathcal M_\iota^{(n)}\Delta_L Q}_2^2
    +\mathcal{CCK}_Z^{(j;n)}\,\avg t^4
       \norm{Q}_{\Gs{s}{\lambda_n}{\sigma-3}}^2
  \right],
  \label{eq:T1-precision}\\
  &\sum_{{\iota}\in\left\{\lambda,w\right\}}
  \norm{\mathcal M_{\iota}^{(n)}\Delta_L\mathcal T_{m;2}^{(j),*}[H,Q]}_2^2
  \notag\\
  &\quad\lesssim e^{C2^{m\alpha_\dagger}}
  \left[
    \norm{H^{(j)}}_{\Gs{s}{\lambda_n}{2}}^2
       \sum_{\iota=\lambda,w}\norm{\mathcal M_\iota^{(n)}\Delta_L Q}_2^2
    +\mathcal{CCK}_H^{(j;n)}\,\avg t^4
       \norm{Q}_{\Gs{s}{\lambda_n}{\sigma-3}}^2
  \right].
  \label{eq:T2-precision}
\end{align}
The observation index and the restriction $j\leq n$ are essential.
We do not assert this bound for a later coefficient measured at an earlier
Gevrey radius.
\end{lemma}

\begin{proof}[Proof of Lemma~\ref{lem:lossy-green-calculus}]
Taking Fourier transform of the second term in \eqref{eq:model-terms} and applying
\eqref{eq:localized-green-kernel-delta} give
\begin{align}
  &e^{\lambda_n|(k,\eta)|^s}\avg{k,\eta}^{\sigma-1}
    |\widehat{\Delta_L\mathcal T_{m;1}^{(j),*}[Z,Q]}(k,\eta)|
  \notag\\*
  &\lesssim e^{Cm}\iint_{\R^2}
   e^{\lambda_n|(k,\eta)|^s}\avg{k,\eta}^{\sigma-1}
   \left(1+\frac{|\eta-\xi|^2}
        {k^2+|\xi-kt|^2}\right)
   e^{-\lambda_m^\Delta|\eta-\xi|^\rho}
    |\widehat{Z^{(j)}}(\xi-\xi')|
    |\xi'-kt|^2|\widehat Q_k(\xi')|
    \dd\xi'\dd\xi.
  \label{eq:lossy-kernel-chain-1}
\end{align}
Subadditivity, the polynomial triangle inequality, and
\eqref{eq:kernel-absorption} imply
\begin{align*}
  &e^{\lambda_n|(k,\eta)|^s}\avg{k,\eta}^{\sigma-1}
    e^{-\lambda_m^\Delta|\eta-\xi|^\rho}\lesssim e^{C2^{m\alpha_\dagger}}
    e^{-\frac14\lambda_m^\Delta|\eta-\xi|^\rho}
    e^{\lambda_n|\xi-\xi'|^s+\lambda_n|(k,\xi')|^s}
    \left(\avg{\xi-\xi'}^{\sigma-1}+\avg{k,\xi'}^{\sigma-1}\right).
\end{align*}
Substitution in \eqref{eq:lossy-kernel-chain-1}, with
$1+|\eta-\xi|^2/(k^2+|\xi-kt|^2)\leq\avg{\eta-\xi}^2$, yields
\begin{align*}
 &e^{\lambda_n|(k,\eta)|^s}\avg{k,\eta}^{\sigma-1}
        |\widehat{\Delta_L\mathcal T_{m;1}^{(j),*}[Z,Q]}(k,\eta)|\\
 &\quad\lesssim e^{C2^{m\alpha_\dagger}}
 \iint_{\R^2}\avg{\eta-\xi}^2e^{-\frac14\lambda_m^\Delta|\eta-\xi|^\rho}
 e^{\lambda_n|\xi-\xi'|^s+\lambda_n|(k,\xi')|^s}\\
 &\qquad\times
 \left(\avg{\xi-\xi'}^{\sigma-1}+\avg{k,\xi'}^{\sigma-1}\right)
 |\widehat{Z^{(j)}}(\xi-\xi')|\,|\xi'-kt|^2
 |\widehat Q_k(\xi')|\dd\xi'\dd\xi.
\end{align*}
We observe that the kernel satisfies the following bound 
\[
 \int_{\R}\avg\xi^2e^{-\frac14\lambda_m^\Delta|\xi|^\rho}\dd\xi
 \lesssim2^{3m}\lesssim e^{C2^{m\alpha_\dagger}}.
\]
Therefore, we have 
\begin{align}
  \norm{\Delta_L\mathcal T_{m;1}^{(j),*}[Z,Q]}_{\Gs{s}{\lambda_n}{\sigma-1}}
  &\lesssim e^{C2^{m\alpha_\dagger}}
    \norm{Z^{(j)}}_{\Gs{s}{\lambda_n}{\sigma-1}}
    \norm{\Delta_LQ}_{\Gs{s}{\lambda_n}{\sigma-1}}.
  \label{eq:lossy-kernel-chain-3}
\end{align}
The $\mathcal T_{m;0}^*[q]$ and $\mathcal T_{m;2}^{(j),*}[H,Q]$ estimates follow in a nearly identical manner. This proves
Lemma~\ref{lem:lossy-green-calculus}.
\end{proof}

\subsection{\texorpdfstring{Proof of the precision Green estimate}{Proof of the precision Green estimate}}
\label{sec:elliptic-precision-green-proof}

First of all, we present the multiplier estimates used to prove
Lemma~\ref{lem:precision-green-calculus}, which provides the precision
Green bounds needed for Proposition~\ref{prop:precision-elliptic}.

\begin{lemma}
\label{lem:stream-frequency-multiplier-bounds}
Recall the multipliers \eqref{eq:elliptic-M-lambda}--\eqref{eq:elliptic-M-w}
and the radii \eqref{eq:lambda-n}, with $1/2<s<1$, $\sigma\geq0$, and
$0<a\leq2s-1$. Let $n\geq1$, $k\in\mathbb Z\setminus\{0\}$, and
$\xi,\xi'\in\mathbb R$. The following bounds hold for almost every $t\geq1$.

\emph{(i)}
Fix $b>0$ and choose $c'=(1+b)^s-b^s<c<1$.
If $|(k,\xi')|\geq b|\xi-\xi'|$, then
\begin{align}
 \mathcal M_\lambda^{(n)}(k,\xi)
 &\lesssim e^{c\lambda_n|\xi-\xi'|^s}
             \mathcal M_\lambda^{(n)}(k,\xi'),
 \label{eq:precision-shared-stream-transfer-lambda}\\*
 \mathcal M_w^{(n)}(k,\xi)
 &\lesssim e^{c\lambda_n|\xi-\xi'|^s}
       \bigl[\mathcal M_w^{(n)}(k,\xi')
                   +\mathcal M_\lambda^{(n)}(k,\xi')\bigr].
 \label{eq:precision-shared-stream-transfer}
\end{align}

\emph{(ii)}
Fix $s3^{s-1}<c<1$. If $|(k,\xi')|\leq\tfrac14|\xi-\xi'|$, then
\begin{align}
 \mathcal M_\lambda^{(n)}(k,\xi)
 &\lesssim \avg{\xi/(kt)}^{-1}e^{c\lambda_n|(k,\xi')|^s}
       \sqrt{-\dot\lambda_n}|\xi-\xi'|^{s/2}{\A}_{\Res}^{(n)}(\xi-\xi'),
 \label{eq:precision-shared-coefficient-transfer-lambda}\\*
 \mathcal M_w^{(n)}(k,\xi)
 &\lesssim \avg{\xi/(kt)}^{-1}e^{c\lambda_n|(k,\xi')|^s}
       \left[\sqrt{\frac{\partial_tw_{\Res}(t,\xi-\xi')}
                              {w_{\Res}(t,\xi-\xi')}}
                +\sqrt{-\dot\lambda_n}|\xi-\xi'|^{s/2}\right]
       {\A}_{\Res}^{(n)}(\xi-\xi').
 \label{eq:precision-shared-coefficient-transfer}
\end{align}
In each part, the constants may depend on the chosen $c$ and the fixed
weight and radius parameters, and also on $b$ in part~(i),
but are independent of $n,k,t,\xi,\xi'$.
\end{lemma}

\begin{proof}
The proof is similar to the corresponding estimates in \cite{BM13},
hence we omit it for the sake of brevity.
\end{proof}

\begin{proof}[Proof of Lemma~\ref{lem:precision-green-calculus}]
The bound \eqref{eq:T0-precision} follows from the
two bounds \eqref{eq:weighted-common-lambda}, \eqref{eq:weighted-common-w}. We now turn to the bounds \eqref{eq:T1-precision}--\eqref{eq:T2-precision}. We use the projections fixed in
\eqref{eq:explicit-LP-symbols}--\eqref{eq:explicit-LP-low-high-support}.
Here $\pi_N$ denotes the vertical projection symbol:
$\widehat{u_N}(\xi)=\pi_N(\xi)\widehat u(\xi)$, and
$\pi_{<N/8}(\xi)=\sum_{N'<N/8}\pi_{N'}(\xi)$.

Let $(r,F)\in\{(1,H^{(j)}),(2,Z^{(j)})\}$.
Throughout the proof $k\ne0$ and $t\geq1$, and we suppress $t$ in the
Fourier transforms and multipliers.

On the Fourier side, we have for $\iota\in\left\{\lambda,w\right\}$
\begin{align*}
 &\left|\mathcal M_\iota^{(n)}(k,\eta)
       \widehat{\Delta_L\mathscr G_m^*(F(\pvt)^rQ)}_k(\eta)\right|\\
 &=\frac1{2\pi}\left|\iint_{\R^2}
       \mathcal M_\iota^{(n)}(k,\eta)(k^2+|\eta-kt|^2)
       \widehat{\mathscr G}_{m,k}^*(t,\eta,\xi)\right.\left.\widehat F(\xi-\xi')
       [\ii(\xi'-kt)]^r\widehat Q_k(\xi')\dd\xi'\dd\xi\right|.
\end{align*}
By \eqref{eq:localized-green-kernel-delta} and
\eqref{eq:common-mlambda-transfer}--\eqref{eq:common-mw-transfer}, we have that
\begin{align}
 &\sum_{\iota\in\{\lambda,w\}}
   \left|\mathcal M_\iota^{(n)}(k,\eta)
       \widehat{\Delta_L\mathscr G_m^*(F(\pvt)^rQ)}_k(\eta)\right|
 \notag\\
 &\quad\lesssim e^{C2^{m\alpha_\dagger}}
       \iint_{\R^2}e^{-\frac14\lambda_m^\Delta|\eta-\xi|^\rho}
       |\widehat F(\xi-\xi')|\,|\xi'-kt|^r
       \sum_{\iota\in\{\lambda,w\}}
       \left|\mathcal M_\iota^{(n)}(k,\xi)\widehat Q_k(\xi')\right|
       \dd\xi'\dd\xi.
 \label{eq:precision-kernel-integral-transfer}
\end{align}

In \eqref{eq:precision-kernel-integral-transfer}, at fixed $k$, we have
two convolutions: the one involving the Fourier variables
$(\eta-\xi,\xi)$, integrated in $\xi$, is called the outer convolution;
the one involving $(\xi-\xi',\xi')$, integrated in $\xi'$, is called
the inner convolution.

For the precision bounds, all Littlewood--Paley projections in this
subsection act only on the vertical Fourier variable. Decompose
\begin{align}
 F (\pvt)^rQ &=\sum_{N\geq8}F_N(\pvt)^rQ_{<N/8} +\sum_{N\geq8}F_{<N/8}(\pvt)^rQ_N + \sum_N\sum_{N/8\leq N'\leq8N}F_{N'}(\pvt)^rQ_N
 \notag\\
  &=:(F(\pvt)^rQ)_{HL}+(F(\pvt)^rQ)_{LH}+(F(\pvt)^rQ)_R.
  \label{eq:elliptic-paraproduct}
\end{align}

To apply \eqref{eq:precision-kernel-integral-transfer}, we insert
\eqref{eq:elliptic-paraproduct} into the Fourier representation and write
\begin{equation}
 \widehat{\Delta_L\mathscr G_m^*(F(\pvt)^rQ)}_k(\eta)
 =T_{HL}(k,\eta)+T_{LH}(k,\eta)+T_R(k,\eta),
 \label{eq:precision-green-T-decomposition}
\end{equation}
where
{\small
\begin{equation}
 \begin{aligned}
 \begin{pmatrix}T_{HL}(k,\eta)\\T_{LH}(k,\eta)\end{pmatrix}
 &:=-\frac1{2\pi}\sum_{N\geq8}\iint_{\R^2}
 (k^2+|\eta-kt|^2)\widehat{\mathscr G}_{m,k}^*(t,\eta,\xi)\\
 &\quad\times[\ii(\xi'-kt)]^r
 \begin{pmatrix}
 \widehat F_N(\xi-\xi')\widehat Q_{<N/8,k}(\xi')\\
 \widehat F_{<N/8}(\xi-\xi')\widehat Q_{N,k}(\xi')
 \end{pmatrix}\dd\xi'\dd\xi,
 \end{aligned}
 \label{eq:precision-green-T-HL-LH}
\end{equation}
}
{\small
\begin{equation}
 \begin{aligned}
 T_R(k,\eta)&:=-\frac1{2\pi}
 \sum_N\sum_{N/8\leq N'\leq8N}\iint_{\R^2}
 (k^2+|\eta-kt|^2)\widehat{\mathscr G}_{m,k}^*(t,\eta,\xi)\\
 &\quad\times\widehat F_{N'}(\xi-\xi')[\ii(\xi'-kt)]^r
 \widehat Q_{N,k}(\xi')\dd\xi'\dd\xi.
 \end{aligned}
 \label{eq:precision-green-T-R}
\end{equation}
}
We estimate $T_{HL}$ first, then $T_{LH}$ and $T_R$.

\medskip

\medskip
We begin with $T_{HL}$ in \eqref{eq:precision-green-T-HL-LH}.
   We split according to the relative sizes of $k$ and $\xi$, and write $T_{HL}=T_{HL,z}+T_{HL,v}$, where
{\small
\begin{equation}
 \begin{aligned}
 \begin{pmatrix}T_{HL,z}(k,\eta)\\T_{HL,v}(k,\eta)\end{pmatrix}
 &:=-\frac1{2\pi}\sum_{N\geq8}\iint_{\R^2}
 (k^2+|\eta-kt|^2)\widehat{\mathscr G}_{m,k}^*(t,\eta,\xi)\\
 &\quad\times\widehat F_N(\xi-\xi')[\ii(\xi'-kt)]^r
 \widehat Q_{<N/8,k}(\xi')
 \begin{pmatrix}\mathbf1_{|\xi|\leq32|k|}\\
                         \mathbf1_{|\xi|>32|k|}\end{pmatrix}\dd\xi'\dd\xi.
 \end{aligned}
 \label{eq:precision-green-HL-z-v}
\end{equation}
}
We call $32|k|\geq|\xi|$ the $\partial_z$-dominant region and
$32|k|<|\xi|$ the $\partial_v$-dominant region.  

For $T_{HL,z}$ in \eqref{eq:precision-green-HL-z-v}, the full expression
to estimate is
\begin{equation}
 \begin{aligned}
 &\sum_{\iota=\lambda,w}
   \left|\mathcal M_\iota^{(n)}(k,\eta)T_{HL,z}(k,\eta)\right|
 =\frac1{2\pi}\sum_{\iota=\lambda,w}
   \Biggl|\sum_{N\geq8}\iint_{\R^2}
   \mathcal M_\iota^{(n)}(k,\eta)(k^2+|\eta-kt|^2)\\
 &\qquad\qquad\times\widehat{\mathscr G}_{m,k}^*(t,\eta,\xi)
   \widehat F_N(\xi-\xi')[\ii(\xi'-kt)]^r
   \widehat Q_{<N/8,k}(\xi')\mathbf1_{|\xi|\leq32|k|}
   \dd\xi'\dd\xi\Biggr|.
 \end{aligned}
 \label{eq:precision-HL-z-weighted-integral}
\end{equation}
We collect three estimates to bound
\eqref{eq:precision-HL-z-weighted-integral}.
First, we estimate the factor involving $Q$.
On the support of the integrand, $|\xi'|\leq3|\xi-\xi'|/16$ and
$13|\xi-\xi'|/16\leq|\xi|\leq19|\xi-\xi'|/16$.
Thus
\[
 |(k,\xi')|\geq|k|\geq|\xi|/32
 \geq13|\xi-\xi'|/512>|\xi-\xi'|/62.
\]
In \eqref{eq:precision-HL-z-weighted-integral}, first apply
\eqref{eq:precision-kernel-integral-transfer} to transfer the multipliers
from $(k,\eta)$ to $(k,\xi)$.
Lemma~\ref{lem:stream-frequency-multiplier-bounds}(i), applied with
$b=1/62$, together with
\eqref{eq:precision-shared-stream-transfer-lambda}--\eqref{eq:precision-shared-stream-transfer},
yields the bound
\begin{equation}
 \sum_{\iota=\lambda,w}
   \left|\mathcal M_\iota^{(n)}(k,\xi)
                    \widehat Q_{<N/8,k}(\xi')\right|
 \lesssim e^{c\lambda_n|\xi-\xi'|^s}
   \sum_{\iota=\lambda,w}
   \left|\mathcal M_\iota^{(n)}(k,\xi')
                    \widehat Q_{<N/8,k}(\xi')\right|.
 \label{eq:precision-HL-z-stream-transfer}
\end{equation}
Secondly, by \eqref{eq:localized-green-kernel-delta}, for $r=1,2$ we have
\begin{equation}
\begin{gathered}
 \begin{aligned}
 \mathscr{F}&:=\bigl(|k|^2+|\eta-kt|^2\bigr)
       |\widehat{\mathscr G}_{m,k}^*(t,\eta,\xi)|
       \frac{|\xi'-kt|^r}{|k|^2+|\xi'-kt|^2}\\
 &\lesssim2^{2m}
       \left(1+\frac{|\eta-\xi|^2}{|k|^2+|\xi-kt|^2}\right)
       e^{-\lambda_m^\Delta|\eta-\xi|^\rho}
       \frac{|\xi'-kt|^r}{|k|^2+|\xi'-kt|^2}\\
 &\lesssim2^{2m}\avg{\eta-\xi}^2e^{-\lambda_m^\Delta|\eta-\xi|^\rho},
 \end{aligned}\\
 \frac{|\xi'-kt|}{|k|^2+|\xi'-kt|^2}\leq\frac1{2|k|}\leq\frac12.
\end{gathered}
 \label{eq:precision-shared-denominator}
\end{equation}

Thirdly, we recall the kernel estimate
\begin{equation}
 \int_{\R}e^{-\frac19\lambda_m^\Delta|\eta-\xi|^\rho}\dd\eta
 \lesssim2^m.
 \label{eq:precision-kernel-L1}
\end{equation}
Finally, combining \eqref{eq:precision-kernel-integral-transfer},
\eqref{eq:precision-HL-z-stream-transfer}, and
\eqref{eq:precision-shared-denominator}, we obtain
\begin{equation}
 \begin{aligned}
 &\sum_{\iota=\lambda,w}\left|
  \mathcal M_\iota^{(n)}(k,\eta)
  T_{HL,z}(k,\eta)\right|
  \lesssim e^{C2^{m\alpha_\dagger}}\sum_{N\geq8}
  \iint_{\R^2}e^{-\frac19\lambda_m^\Delta|\eta-\xi|^\rho}\\
 &\quad\times\mathbf1_{|\xi|\leq32|k|}
   e^{c\lambda_n|\xi-\xi'|^s}|\widehat F_N(\xi-\xi')|
   \sum_{\iota=\lambda,w}
  |\mathcal M_\iota^{(n)}(k,\xi')
       \widehat{\Delta_LQ_{<N/8}}_k(\xi')|\dd\xi'\dd\xi.
 \end{aligned}
 \label{eq:precision-HL-z-pointwise-bound}
\end{equation}
In \eqref{eq:precision-HL-z-pointwise-bound}, sum the nonnegative
projection products using \eqref{eq:explicit-LP-symbols}.
Apply Young's $L^1*L^2\to L^2$ inequality to the inner and outer
convolutions, with $e^{c\lambda_n|\cdot|^s}\widehat F$ in $L^1$,
$\mathcal M_\iota^{(n)}\widehat{\Delta_LQ}$ in $\ell^2L^2$, and the
kernel bound \eqref{eq:precision-kernel-L1}. This gives
\begin{equation}
 \sum_{\iota=\lambda,w}
  \|\mathcal M_\iota^{(n)}T_{HL,z}\|_2^2
 \lesssim e^{C2^{m\alpha_\dagger}}\|F\|_{\Gs{s}{\lambda_n}{2}}^2
   \sum_{\iota=\lambda,w}\|\mathcal M_\iota^{(n)}\Delta_LQ\|_2^2.
 \label{eq:precision-high-low-z-conclusion}
\end{equation}

For $T_{HL,v}$ in \eqref{eq:precision-green-HL-z-v}, the full quantity
we estimate is
\begin{align}
 &\sum_{\iota=\lambda,w}
   \left|\mathcal M_\iota^{(n)}(k,\eta)T_{HL,v}(k,\eta)\right|\notag\\*
 &=\frac1{2\pi}\sum_{\iota=\lambda,w}
   \left|\sum_{N\geq8}\iint_{\R^2}
   \mathcal M_\iota^{(n)}(k,\eta)
   (k^2+|\eta-kt|^2)\widehat{\mathscr G}_{m,k}^*(t,\eta,\xi)\right.\notag\\*
 &\hspace{23mm}\left.\times
   \widehat F_N(\xi-\xi')[\ii(\xi'-kt)]^r
   \widehat Q_{<N/8,k}(\xi')\mathbf1_{|\xi|>32|k|}
   \dd\xi'\dd\xi\right|.
 \label{eq:precision-HL-v-weighted-integral}
\end{align}
We collect three estimates to bound
\eqref{eq:precision-HL-v-weighted-integral}.

First of all, on the support of this integral, $|k|<|\xi|/32$.
By \eqref{eq:explicit-LP-low-high-support},
$|\xi'|\leq3|\xi-\xi'|/16$.  Since $\xi=\xi-\xi'+\xi'$, this gives
\[
 |\xi|\leq\tfrac{19}{16}|\xi-\xi'|,\qquad
 |(k,\xi')|=|k|+|\xi'|
 <\left(\tfrac{19}{512}+\tfrac3{16}\right)|\xi-\xi'|
 =\tfrac{115}{512}|\xi-\xi'|<\tfrac14|\xi-\xi'|.
\]
\Needspace{9\baselineskip}
Thus Lemma~\ref{lem:stream-frequency-multiplier-bounds}(ii) applies.
Fix $s3^{s-1}<c<1$. Multiplying
\eqref{eq:precision-shared-coefficient-transfer-lambda}--\eqref{eq:precision-shared-coefficient-transfer}
by $|\widehat F_N(\xi-\xi')|$ and adding, we obtain
\begin{align}
 &\sum_{\iota=\lambda,w}
    \left|\mathcal M_\iota^{(n)}(k,\xi)\widehat F_N(\xi-\xi')\right|\notag\\*
 &\quad\lesssim
   \avg{\xi/(kt)}^{-1}e^{c\lambda_n|(k,\xi')|^s}
   \left[\sqrt{-\dot\lambda_n}|\xi-\xi'|^{s/2}
        +\sqrt{\frac{\partial_tw_{\Res}(t,\xi-\xi')}
                         {w_{\Res}(t,\xi-\xi')}}\right]\notag\\*
 &\hspace{25mm}\times
   \A_{\Res}^{(n)}(\xi-\xi')|\widehat F_N(\xi-\xi')|.
 \label{eq:precision-HL-v-coefficient-factor}
\end{align}

Secondly, we use the following two straightforward bounds:
\begin{align}
 \avg{\xi/(kt)}^{-1}|\xi'-kt|^2
       \left(1+\frac{(\eta-\xi)^2}{\bigl(|k|^2+|\xi-kt|^2\bigr)}\right)
 &\lesssim\avg t^2\avg{k,\xi'}^2\avg{\eta-\xi}^2,
 \label{eq:T1-coefficient-high-kernel-payment}\\
 \avg{\xi-\xi'}\avg{\xi/(kt)}^{-1}|\xi'-kt|
       \left(1+\frac{(\eta-\xi)^2}{\bigl(|k|^2+|\xi-kt|^2\bigr)}\right)
 &\lesssim\avg t^2\avg{k,\xi'}^2\avg{\eta-\xi}^2.
 \label{eq:T2-coefficient-high-kernel-payment}
\end{align}
For the first, $|\xi'-kt|\leq C\avg t\avg{k,\xi'}$ and $\bigl(|k|^2+|\xi-kt|^2\bigr)\geq1$.
For the second, $\avg{\xi-\xi'}\approx\avg\xi$ and
$\avg{\xi-\xi'}\avg{\xi/(kt)}^{-1}\lesssim|kt|
\lesssim\avg t\avg{k,\xi'}$; here $|\xi-\xi'|\geq4|(k,\xi')|\geq4$.

The decay in \eqref{eq:precision-kernel-integral-transfer} and
$\lambda_m^\Delta=c_\Delta2^{-m\rho}$ give
\begin{equation}
 \avg{\eta-\xi}^2e^{-\frac14\lambda_m^\Delta|\eta-\xi|^\rho}
 \lesssim2^{2m}e^{-\frac19\lambda_m^\Delta|\eta-\xi|^\rho}.
 \label{eq:precision-HL-v-kernel-payment}
\end{equation}
The remaining outer kernel has the $L^1$ bound
\eqref{eq:precision-kernel-L1}.

Thirdly, the definitions \eqref{eq:CCK-Z}--\eqref{eq:CCK-H} give
\begin{equation}
\begin{aligned}
 &\left\|\left(\sqrt{-\dot\lambda_n}|\partial_v|^{s/2}
       +\sqrt{\frac{\partial_tw_{\Res}}{w_{\Res}}}\right)
     \A_{\Res}^{(n)}\avg{\partial_v}^{r-2}F\right\|_2^2
 \leq2
 \begin{cases}
  \mathcal{CCK}_Z^{(j;n)},&r=2,\\
  \mathcal{CCK}_H^{(j;n)},&r=1.
 \end{cases}
\end{aligned}
\label{eq:precision-HL-v-coefficient-norm}
\end{equation}
Here $\avg{\xi-\xi'}^{r-2}$ is one for $F=Z^{(j)}$ when $r=2$, and
$\avg{\xi-\xi'}^{-1}$ for $F=H^{(j)}$ when $r=1$.
For the stream factor, with $\sigma$ chosen as in
\eqref{eq:weight-parameters}, Cauchy--Schwarz gives, uniformly in $k$,
\begin{equation}
 \begin{aligned}
 &\|e^{c\lambda_n|(k,\xi')|^s}\avg{k,\xi'}^2
                       \widehat Q_k(\xi')\|_{L^1_{\xi'}}
 \leq\|\avg{\xi'}^{5-\sigma}e^{-(1-c)\lambda_n|(k,\xi')|^s}\|_{L^2_{\xi'}}
       \|e^{\lambda_n|(k,\xi')|^s}\avg{k,\xi'}^{\sigma-3}\widehat Q_k\|_{L^2_{\xi'}}\\
 &\quad\lesssim\|e^{\lambda_n|(k,\xi')|^s}\avg{k,\xi'}^{\sigma-3}\widehat Q_k\|_{L^2_{\xi'}}.
 \end{aligned}
 \label{eq:precision-HL-v-stream-norm}
\end{equation}

\Needspace{9\baselineskip}
Finally, we combine the three estimates. Apply
\eqref{eq:precision-kernel-integral-transfer} to
\eqref{eq:precision-HL-v-weighted-integral}, then insert
\eqref{eq:precision-HL-v-coefficient-factor},
\eqref{eq:T1-coefficient-high-kernel-payment}--\eqref{eq:T2-coefficient-high-kernel-payment},
and \eqref{eq:precision-HL-v-kernel-payment}.
Sum the nonnegative projection products using
\eqref{eq:explicit-LP-symbols}. This gives the concrete bound
\begin{align}
 &\sum_{\iota=\lambda,w}
   \left|\mathcal M_\iota^{(n)}(k,\eta)T_{HL,v}(k,\eta)\right|\notag\\*
 &\quad\lesssim e^{C2^{m\alpha_\dagger}}\avg t^2
   \iint_{\R^2}e^{-\frac19\lambda_m^\Delta|\eta-\xi|^\rho}
   \left[\sqrt{-\dot\lambda_n}|\xi-\xi'|^{s/2}
        +\sqrt{\frac{\partial_tw_{\Res}(t,\xi-\xi')}
                         {w_{\Res}(t,\xi-\xi')}}\right]\notag\\*
 &\hspace{17mm}\times\A_{\Res}^{(n)}(\xi-\xi')
   \avg{\xi-\xi'}^{r-2}|\widehat F(\xi-\xi')|
   e^{c\lambda_n|(k,\xi')|^s}\avg{k,\xi'}^2|\widehat Q_k(\xi')|
                 \dd\xi'\dd\xi.
 \label{eq:precision-HL-v-pointwise-bound}
\end{align}
Apply Young's $L^2*L^1\to L^2$ inequality to the inner convolution in $\xi'$,
using \eqref{eq:precision-HL-v-coefficient-norm} and
\eqref{eq:precision-HL-v-stream-norm}. Apply its $L^1*L^2\to L^2$
form to the outer convolution in $\xi$, using
\eqref{eq:precision-kernel-L1}.
Squaring and summing over $k\ne0$ yields
\begin{equation}
 \sum_{{\iota}\in\left\{\lambda,w\right\}}
 \|\mathcal M_{\iota}^{(n)}T_{HL,v}\|_2^2
 \lesssim e^{C2^{m\alpha_\dagger}}
 \avg t^4\|Q\|_{\Gs{s}{\lambda_n}{\sigma-3}}^2
 \begin{cases}
 \mathcal{CCK}_Z^{(j;n)},&r=2,\\
 \mathcal{CCK}_H^{(j;n)},&r=1
 \end{cases}.
 \label{eq:precision-shared-coefficient-conclusion}
\end{equation}
Since $T_{HL}=T_{HL,z}+T_{HL,v}$ by
\eqref{eq:precision-green-HL-z-v}, apply
$\|u+v\|_2^2\leq2\|u\|_2^2+2\|v\|_2^2$ and combine
\eqref{eq:precision-high-low-z-conclusion} with
\eqref{eq:precision-shared-coefficient-conclusion} to obtain
\begin{equation}
\begin{aligned}
 \sum_{\iota=\lambda,w}\|\mathcal M_\iota^{(n)}T_{HL}\|_2^2
 &\lesssim e^{C2^{m\alpha_\dagger}}\|F\|_{\Gs{s}{\lambda_n}{2}}^2
     \sum_{\iota=\lambda,w}\|\mathcal M_\iota^{(n)}\Delta_LQ\|_2^2\\
 &\qquad+e^{C2^{m\alpha_\dagger}}\avg t^4
 \|Q\|_{\Gs{s}{\lambda_n}{\sigma-3}}^2
 \begin{cases}
  \mathcal{CCK}_Z^{(j;n)},&r=2,\\
  \mathcal{CCK}_H^{(j;n)},&r=1
 \end{cases}.
\end{aligned}
\label{eq:precision-high-low-conclusion}
\end{equation}
For $(r,F)=(2,Z^{(j)})$ and $(r,F)=(1,H^{(j)})$, this gives the
high--low contribution to \eqref{eq:T1-precision} and
\eqref{eq:T2-precision}, respectively, in
Lemma~\ref{lem:precision-green-calculus}.

\medskip
We now estimate $T_{LH}$ in \eqref{eq:precision-green-T-HL-LH}.
The projection supports give
\[
 |\xi-\xi'|\leq\frac{3N}{32},\qquad
 \frac N2\leq|\xi'|\leq\frac{3N}{2},\qquad
 |\xi-\xi'|\leq\frac3{16}|\xi'|.
\]
On the present low--high support, $|(k,\xi')|\geq(16/3)|\xi-\xi'|$,
so \eqref{eq:precision-shared-stream-transfer} applies.

By the same multiplier-transfer, denominator, and two-convolution
argument used for \eqref{eq:precision-high-low-z-conclusion}, with
\eqref{eq:precision-shared-denominator} and \eqref{eq:precision-kernel-L1},
we obtain
\begin{equation}
 \sum_{\iota=\lambda,w}
 \|\mathcal M_\iota^{(n)}\Delta_L\mathscr G_m^*
       (F(\pvt)^rQ)_{LH}\|_2^2
 \lesssim e^{C2^{m\alpha_\dagger}}\|F\|_{\Gs{s}{\lambda_n}{2}}^2
       \sum_{\iota=\lambda,w}\norm{\mathcal M_\iota^{(n)}\Delta_L Q}_2^2.
 \label{eq:precision-low-high-conclusion}
\end{equation}
Here the nonnegative symbols satisfy
$\sum_N\pi_N(\xi')\pi_{<N/8}(\xi-\xi')\leq1$;
thus the projected sum is taken inside the convolution before applying
Young's inequality.

\medskip
Finally, we estimate $T_R$ in \eqref{eq:precision-green-T-R}.
The comparable vertical frequencies allow us to put the multiplier on
$Q$, including when the output frequency is small.
In the annular comparable remainder,
\begin{equation}
 b|\xi-\xi'|\leq|\xi'|\leq b^{-1}|\xi-\xi'|,
 \qquad b=1/24.
 \label{eq:remainder-vertical-only-geometry}
\end{equation}
A pair with $\min(N,N')=1/2$ has $|\xi-\xi'|,|\xi'|\leq6$;
there $|(k,\xi')|\geq1\geq|\xi-\xi'|/6$.
Thus \eqref{eq:precision-shared-stream-transfer} applies with
$b=1/64$ to every pair, including the lowest-frequency block, for every $k\ne0$.

The nonnegative projection symbols satisfy
\[
 \sum_N\sum_{N/8\leq N'\leq8N}
    \pi_{N'}(\xi-\xi')\pi_N(\xi')\leq1.
\]
By the same argument
used for \eqref{eq:precision-high-low-z-conclusion}, the two Young
inequalities and \eqref{eq:precision-kernel-L1} give
\[
 \sum_{\iota=\lambda,w}
  \|\mathcal M_\iota^{(n)}T_R\|_2^2
 \lesssim e^{C2^{m\alpha_\dagger}}\|F\|_{\Gs{s}{\lambda_n}{2}}^2
  \sum_{\iota=\lambda,w}\|\mathcal M_\iota^{(n)}\Delta_LQ\|_2^2.
\]
This includes arbitrarily large horizontal frequencies in the
bounded vertical block; no orthogonality of the Green outputs is used.

The low--high, $\partial_z$-dominant high--low, and remainder terms
give the first terms on the right of \eqref{eq:T1-precision}--\eqref{eq:T2-precision}.
The $\partial_v$-dominant high--low term gives their coefficient CK terms
by \eqref{eq:precision-shared-coefficient-conclusion}.
Together with the direct-source bound \eqref{eq:T0-precision}, this proves
Lemma~\ref{lem:precision-green-calculus}.
\end{proof}

\subsection{Proof of the lossy estimates}
\label{sec:elliptic-lossy-proof}

We prove Proposition~\ref{prop:lossy-elliptic} and
\eqref{eq:elliptic-source-bootstrap}.

\begin{proof}[Proof of Proposition~\ref{prop:lossy-elliptic} and
\eqref{eq:elliptic-source-bootstrap}]
Assume the source bound \eqref{eq:elliptic-source-bootstrap} and the
lossy bounds in Proposition~\ref{prop:lossy-elliptic} for $r<n$.
We first reduce the decay estimates
\eqref{eq:lossy-ext}--\eqref{eq:lossy-profile-velocity} to
\eqref{eq:lossy-source-ext}--\eqref{eq:lossy-source-int}, then prove the
source bound \eqref{eq:elliptic-source-bootstrap} and the exterior bound
\eqref{eq:lossy-source-ext}, and finally establish the interior bound
\eqref{eq:lossy-source-int} by backward induction on $m$.

It is straightforward to see that, for $k\ne0$, the following bound holds:
\begin{equation}
 \avg{k,\eta}^4\bigl(|k|^2+|\eta-kt|^2\bigr)^2\geq k^8\avg t^4/64.
 \label{eq:elliptic-time-denominator}
\end{equation}
Thus the exact nonzero-mode Fourier identity
$\widehat\Psi_k=-\bigl(|k|^2+|\eta-kt|^2\bigr)^{-1}\widehat{\Delta_L\Psi}_k$ gives
\begin{equation}
\|\Pneq\Psi\|_{\Gs{s}{\lambda_n}{\sigma-3}}^2 \lesssim\avg t^{-4} \||\partial_z|^{-4}\Delta_L\Pneq\Psi\|_{\Gs{s}{\lambda_n}{\sigma-1}}^2 \lesssim\avg t^{-4} \|\Delta_L\Pneq\Psi\|_{\Gs{s}{\lambda_n}{\sigma-1}}^2.
 \label{eq:flat-lossy-inversion}
\end{equation}
Here $|\partial_z|^{-4}$ has symbol $|k|^{-4}$ on $k\ne0$ and zero
on $k=0$.  
The same symbol bounds give
\eqref{eq:lossy-velocity}--\eqref{eq:lossy-profile-velocity} once
\eqref{eq:lossy-source-ext}--\eqref{eq:lossy-source-int} are proved.

We next estimate the source in \eqref{eq:coefficient-source} and then the exterior extension in \eqref{eq:exterior-extension-definition}.
We have $\dot\lambda_n\leq0$ by \eqref{eq:lambda-derivative} and
$\lambda_n(t)\leq\lambda_r(t)$ for $r<n$ by \eqref{eq:radius-gap}.
Using the symbol bound \eqref{eq:precision-shared-denominator},
the outer induction hypotheses \eqref{eq:lossy-source-ext} and
\eqref{eq:elliptic-source-bootstrap} at levels $r<n$, and the
summability bound \eqref{eq:elliptic-shell-weight-sums}, we obtain
\begin{align}
 &\|\Pneq(\pvt)^2\Phi_{\Ext}^{[n-1]}\|_{\Gs{s}{\lambda_n}{\sigma-1}}
 +\|\Pneq\pvt\Phi_{\Ext}^{[n-1]}\|_{\Gs{s}{\lambda_n}{\sigma-1}}
 \lesssim\sum_{r<n}\varsigma_r\mathfrak S_r\lesssim M^C\eps.
 \label{eq:previous-field-lossy-outer-induction}
\end{align}
Substituting this estimate and the coefficient bootstrap bounds into
\eqref{eq:elliptic-source-expanded}, we obtain
$\mathfrak S_n\lesssim M\eps\delta_n+M^C\eps^2\delta_n
\lesssim M^C\eps\delta_n$, proving the source bound \eqref{eq:elliptic-source-bootstrap} at level $n$.
For the coefficient norms \eqref{eq:elliptic-loss-sequences},
\eqref{eq:bootstrap-energy-assumption} and the coefficient formulas give
\begin{equation}
 \mathfrak S_n\leq M^C\eps\delta_n,\qquad
 c_{Z,H;\mathrm{high}}^{(j)}+c_{Z,H;\mathrm{low}}^{(j)}\leq M^C\eps\delta_j.
 \label{eq:common-scale-coefficient-inputs}
\end{equation}
After fixing $M$ and the constants in the estimates, choose
$\varepsilon_{\mathrm{abs}}>0$ sufficiently small for the absorptions below.
By \eqref{eq:common-scale-coefficient-inputs} and
$\alpha_*>\alpha_\dagger$, we can then choose $\eps>0$ sufficiently small,
independently of $n$ and $\mathcal N$, so that
\begin{equation}
 \sum_{j=1}^{n}\left[
 \varsigma_jc_{Z,H;\mathrm{high}}^{(j)}+\varsigma_j^2c_{Z,H;\mathrm{high}}^{(j)}
 +\varsigma_j(c_{Z,H;\mathrm{low}}^{(j)})^2+\varsigma_j^2(c_{Z,H;\mathrm{low}}^{(j)})^2
 \right]\leq\varepsilon_{\mathrm{abs}}.
 \label{eq:common-scale-small-sums}
\end{equation}

Recall the exterior extension \eqref{eq:exterior-extension-definition} and
the operators \eqref{eq:model-terms}. The exterior stream function has the decomposition
\[
 \Phi_{\Ext}^{(n)}
 =\mathcal T_{n;0}^{\Ext}[\Pneq S^{(n)}]
   +\sum_{j\leq n}\bigl(
     \mathcal T_{j;1}^{(j),\Ext}[Z,\Phi_{\Ext}^{(n)}]
       -\mathcal T_{j;2}^{(j),\Ext}[H,\Phi_{\Ext}^{(n)}]\bigr).
\]
Apply \eqref{eq:T0-lossy} to the source term and
\eqref{eq:T12-lossy} to the $Z^{(j)}$ and $H^{(j)}$ terms:
\begin{equation}
\begin{aligned}
 &\|\Delta_L\Pneq\Phi_{\Ext}^{(n)}\|_{\Gs{s}{\lambda_n}{\sigma-1}}\leq
   \|\Delta_L\mathcal T_{n;0}^{\Ext}[\Pneq S^{(n)}]\|_{\Gs{s}{\lambda_n}{\sigma-1}}\\
 &\quad+\sum_{j\leq n}\left(
   \|\Delta_L\mathcal T_{j;1}^{(j),\Ext}[Z,\Phi_{\Ext}^{(n)}]\|_{\Gs{s}{\lambda_n}{\sigma-1}}
  +\|\Delta_L\mathcal T_{j;2}^{(j),\Ext}[H,\Phi_{\Ext}^{(n)}]\|_{\Gs{s}{\lambda_n}{\sigma-1}}
    \right)\\
 &\leq C_1\varsigma_n\mathfrak S_n
             +C_1\sum_{j\leq n}\varsigma_jc_{Z,H;\mathrm{high}}^{(j)}\|\Delta_L\Pneq\Phi_{\Ext}^{(n)}\|_{\Gs{s}{\lambda_n}{\sigma-1}}.
\end{aligned}
 \label{eq:exterior-lossy-absorption-chain}
\end{equation}
Choose $C_1\varepsilon_{\mathrm{abs}}\leq1/2$ and absorb the last term:
\[
 \|\Delta_L\Pneq\Phi_{\Ext}^{(n)}\|_{\Gs{s}{\lambda_n}{\sigma-1}}
 \leq2C_1\varsigma_n\mathfrak S_n.
\]

It remains to prove the interior bound \eqref{eq:lossy-source-int}.
By \eqref{eq:exact-extension-endpoint} and
\eqref{eq:interior-derived-recursion}, for $m<n$ we have
\begin{align*}
 \Phi_{\Int}^{(n;m)}
 ={}&\mathscr G_m^{\Int}S^{(n)}
   +\sum_{j\leq m}\mathscr G_j^{\Ext}
       \bigl[Z^{(j)}(\pvt)^2\Phi_{\Int}^{(n;m)}
                  -H^{(j)}\pvt\Phi_{\Int}^{(n;m)}\bigr]\\*
 &+\sum_{m<j\leq n}\mathscr G_j^{\Int}
       \bigl[Z^{(j)}(\pvt)^2\Phi_{\Int}^{(n;j)}
                  -H^{(j)}\pvt\Phi_{\Int}^{(n;j)}\bigr].
\end{align*}
Applying the Green estimates \eqref{eq:T0-lossy} and
\eqref{eq:T12-lossy} to \eqref{eq:interior-derived-recursion} gives
\begin{equation}
 \begin{aligned}
 \|\Delta_L\Pneq\Phi_{\Int}^{(n;m)}\|_{\Gs{s}{\lambda_n}{\sigma-1}}
 &\leq C_1\varsigma_m\mathfrak S_n
  +C_1\sum_{j\leq m}\varsigma_jc_{Z,H;\mathrm{high}}^{(j)}\|\Delta_L\Pneq\Phi_{\Int}^{(n;m)}\|_{\Gs{s}{\lambda_n}{\sigma-1}}\\
 &\quad+C_1\sum_{m<j\leq n}\varsigma_jc_{Z,H;\mathrm{high}}^{(j)}\|\Delta_L\Pneq\Phi_{\Int}^{(n;j)}\|_{\Gs{s}{\lambda_n}{\sigma-1}}.
 \end{aligned}
 \label{eq:interior-lossy-backward-recurrence}
\end{equation}
After absorption, enlarge $C_2$ to include the exterior bound and take
$2C_2\varepsilon_{\mathrm{abs}}\leq1$.
Assuming $\|\Delta_L\Pneq\Phi_{\Int}^{(n;j)}\|_{\Gs{s}{\lambda_n}{\sigma-1}}\leq2C_2\varsigma_j\mathfrak S_n$
for $m<j\leq n$, \eqref{eq:common-scale-small-sums} and $\varsigma_m\geq1$ give
\[
 \begin{aligned}
 \|\Delta_L\Pneq\Phi_{\Int}^{(n;m)}\|_{\Gs{s}{\lambda_n}{\sigma-1}}
 &\leq C_2\varsigma_m\mathfrak S_n
       +2C_2^2\mathfrak S_n\sum_{m<j\leq n}
                       \varsigma_j^2c_{Z,H;\mathrm{high}}^{(j)}\\
 &\leq\bigl(C_2+2C_2^2\varepsilon_{\mathrm{abs}}\bigr)
           \varsigma_m\mathfrak S_n
 \leq2C_2\varsigma_m\mathfrak S_n.
 \end{aligned}
\]
This closes the backward induction uniformly in $n$.
Equation~\eqref{eq:flat-lossy-inversion} yields
\eqref{eq:lossy-ext}--\eqref{eq:lossy-int}; the corresponding velocity symbol
bounds yield \eqref{eq:lossy-velocity}--\eqref{eq:lossy-profile-velocity}.
This proves Proposition~\ref{prop:lossy-elliptic} and
\eqref{eq:elliptic-source-bootstrap} at level $n$.
\end{proof}

\subsection{Proof of the precision estimates}
\label{sec:elliptic-precision-proof}

We prove Proposition~\ref{prop:precision-elliptic}.

\begin{proof}[Proof of Proposition~\ref{prop:precision-elliptic}]
The proof of this proposition is relatively long, so we now outline the basic strategy. We first prove the source estimate \eqref{eq:precision-source-increment}
and deduce its Green-operator bound.  We then estimate the exterior and
interior stream-function increments simultaneously to obtain
\eqref{eq:precision-elliptic-increment}.  Finally, we sum these estimates
to prove \eqref{eq:precision-elliptic-partial-sum} and
\eqref{eq:precision-source-partial-sum}.
 
Assume the precision estimates are valid for $r<n$. To collect the right-hand side of
\eqref{eq:precision-elliptic-increment}, set
\begin{equation}
 \mathcal Q_{\rm obs}^{(n)}
 :=\CK_{f,\partial g,h}^{(n)}
   +\delta_n^2\sum_{r<n}\delta_r^{-1}\varsigma_r\CK_{f,\partial g,h}^{(r)}
   +\frac{M^C\eps^2\delta_n^2}{\avg t^{1+\delta}}.
 \label{eq:interior-precision-source}
\end{equation}
The subscript $\mathrm{obs}$ stands for ``observation".

\medskip
\noindent\emph{1. The source estimate \eqref{eq:precision-source-increment}.}
We first record the inputs from the earlier levels.  By
\eqref{eq:radius-gap}, \eqref{eq:lambda-derivative},
\begin{equation}
 \mathcal M_{\iota,k}^{(n)}(t,\eta)
 \leq C\mathcal M_{\iota,k}^{(r)}(t,\eta),
 \qquad 1\leq r\leq n,\quad k\ne0,\quad
 \iota\in\left\{\lambda,w\right\}.
 \label{eq:precision-radius-comparison}
\end{equation}
For $n\geq2$, apply \eqref{eq:precision-elliptic-partial-sum} at level
$n-1$ and \eqref{eq:precision-radius-comparison} to obtain
\begin{equation}
 \sum_{\iota=\lambda,w}
 \|\mathcal M_\iota^{(n)}\Delta_L\Phi_{\Ext}^{[n-1]}\|_2^2
 \lesssim \left[
 \sum_{r<n}\delta_r^{-1}\varsigma_r\CK_{f,\partial g,h}^{(r)}
 +\frac{M^C\eps^2}{\avg t^{1+\delta}}\right].
 \label{eq:previous-field-precision-input}
\end{equation}
By \eqref{eq:lossy-ext}, \eqref{eq:elliptic-source-bootstrap}, and
\eqref{eq:elliptic-shell-weight-sums},
\begin{equation}
 \avg t^4\|\Phi_{\Ext}^{[n-1]}\|_{\Gs{s}{\lambda_n}{\sigma-3}}^2
 \lesssim\left(\sum_{r<n}\varsigma_r\mathfrak S_r\right)^2
 \lesssim M^C\eps^2\left(\sum_{r<n}\varsigma_r\delta_r\right)^2
 \lesssim M^C\eps^2.
 \label{eq:previous-field-lossy-input}
\end{equation}
Both left-hand sides vanish when $n=1$.
Lemma~\ref{lem:coefficient-CK-product} gives, for $j\leq n$,
\begin{equation}
 \mathcal{CCK}_{Z,H}^{(j;n)}
 \leq CM^C\left(\CK_{f,\partial g,h}^{(j)}
       +\eps^2\delta_j^2\sum_{r<j}\varsigma_r\CK_{f,\partial g,h}^{(r)}\right).
 \label{eq:Gamma-common-scale-bound}
\end{equation}

On the coefficient supports, the extension agreement
\eqref{eq:extension-compatibility} identifies the first two
$\pvt$ derivatives of $\Phi_{\Ext}^{[n-1]}$ with those of
$\Pneq\phi^{[n-1]}$.  Since $Z^{(n)}$ and $H^{(n)}$ are independent of $z$,
\begin{equation}
 \Pneq S^{(n)}
 =\Pneq f^{(n)}+Z^{(n)}(\pvt)^2\Phi_{\Ext}^{[n-1]}
   -H^{(n)}\pvt\Phi_{\Ext}^{[n-1]}
 =:T_1+T_2+T_3.
 \label{eq:complete-source-identity-kernel}
\end{equation}
For $T_1$ in \eqref{eq:complete-source-identity-kernel}, the multiplier
definitions and \eqref{eq:CK-f} give
\[
 \sum_{\iota=\lambda,w}\|\mathcal M_\iota^{(n)}T_1\|_2^2
 \lesssim\CK_{f,+}^{(n)}\leq\CK_{f,\partial g,h}^{(n)}.
\]
For $T_2,T_3$, their Fourier representations are
\[
 \begin{aligned}
 \mathcal M_\iota^{(n)}(k,\eta)\widehat T_{2,k}(\eta)
 &=\frac1{2\pi}\int_{\R}\mathcal M_\iota^{(n)}(k,\eta)
   \widehat Z^{(n)}(\eta-\xi')[\ii(\xi'-kt)]^2
   \widehat{\Phi_{\Ext}^{[n-1]}}_k(\xi')\dd\xi',\\
 \mathcal M_\iota^{(n)}(k,\eta)\widehat T_{3,k}(\eta)
 &=-\frac1{2\pi}\int_{\R}\mathcal M_\iota^{(n)}(k,\eta)
   \widehat H^{(n)}(\eta-\xi')\ii(\xi'-kt)
   \widehat{\Phi_{\Ext}^{[n-1]}}_k(\xi')\dd\xi'.
 \end{aligned}
\]
We now apply the paraproduct decomposition \eqref{eq:elliptic-paraproduct}
and the inner-convolution estimates in the proof of
Lemma~\ref{lem:precision-green-calculus}, with $\xi=\eta$.
The low--high term and the remainder use
\eqref{eq:precision-shared-stream-transfer-lambda}--\eqref{eq:precision-shared-stream-transfer};
the high--low term uses the same $\partial_z$- and $\partial_v$-dominant
split as \eqref{eq:precision-green-HL-z-v} and
\eqref{eq:precision-shared-coefficient-transfer-lambda}--\eqref{eq:precision-shared-coefficient-transfer}. Clearly,
\begin{equation}
 \begin{aligned}
 \sum_{\iota=\lambda,w}
 \bigl(\|\mathcal M_\iota^{(n)}T_2\|_2^2
       +\|\mathcal M_\iota^{(n)}T_3\|_2^2\bigr)
 &\lesssim (c_{Z,H;\mathrm{low}}^{(n)})^2
       \sum_{\iota=\lambda,w}
       \|\mathcal M_\iota^{(n)}\Delta_L\Phi_{\Ext}^{[n-1]}\|_2^2\\
 &\quad+\mathcal{CCK}_{Z,H}^{(n;n)}\avg t^4
       \|\Phi_{\Ext}^{[n-1]}\|_{\Gs{s}{\lambda_n}{\sigma-3}}^2.
 \end{aligned}
 \label{eq:direct-source-precision-ledger}
\end{equation}
Combining these bounds with
\eqref{eq:common-scale-coefficient-inputs},
\eqref{eq:previous-field-precision-input}--\eqref{eq:Gamma-common-scale-bound},
and choosing $M$ large and $\eps$ small, we obtain
\begin{equation}
 \begin{aligned}
 \sum_{\iota=\lambda,w}\|\mathcal M_\iota^{(n)}\Pneq S^{(n)}\|_2^2
 &\lesssim (1+M^C\eps^2)\CK_{f,\partial g,h}^{(n)}
       +M^C\eps^2\delta_n^2
          \left[\sum_{r<n}\delta_r^{-1}\varsigma_r\CK_{f,\partial g,h}^{(r)}
       +\frac{M^C\eps^2}{\avg t^{1+\delta}}\right]\lesssim\mathcal Q_{\rm obs}^{(n)}.
 \end{aligned}
 \label{eq:source-precision-direct-full-line}
\end{equation}
Since $\varsigma_n\geq1$, this proves
\eqref{eq:precision-source-increment}. For $(q,*)=(m,\Int)$ or $(n,\Ext)$, apply
\eqref{eq:T0-precision} and \eqref{eq:source-precision-direct-full-line}:
\begin{equation}
 \sum_{\iota=\lambda,w}
 \|\mathcal M_\iota^{(n)}\Delta_L\mathscr G_q^*\Pneq S^{(n)}\|_2^2
 \lesssim\varsigma_q\sum_{\iota=\lambda,w}
       \|\mathcal M_\iota^{(n)}\Pneq S^{(n)}\|_2^2
 \lesssim\varsigma_q\mathcal Q_{\rm obs}^{(n)}.
 \label{eq:direct-source-precision-bound}
\end{equation}

\medskip
\noindent\emph{2. The stream-function increment estimate
\eqref{eq:precision-elliptic-increment}.}
On the support of $Z^{(j)},H^{(j)}$, the agreement
\eqref{eq:extension-compatibility} gives
\begin{equation}
 \Pneq\mathcal C^{(j)}[\phi^{(n)}]
 =\mathcal C^{(j)}[\Pneq\phi^{(n)}]
 =\mathcal C^{(j)}[\Phi_{\Int}^{(n;j)}].
 \label{eq:whole-product-j-view-identity}
\end{equation}
By \eqref{eq:exterior-extension-definition}--\eqref{eq:interior-extension-definition}
and \eqref{eq:whole-product-j-view-identity}, for $1\leq m\leq n$,
\begin{equation}
 \begin{aligned}
 \Phi_{\Int}^{(n;m)}
 &=\begin{cases}
  \mathscr G_m^{\Int}\Pneq S^{(n)},&m<n,\\
  \mathscr G_n^{\Ext}\Pneq S^{(n)},&m=n,
 \end{cases}\\
 &\quad+\sum_{j\leq m}\mathscr G_j^{\Ext}
       \bigl[Z^{(j)}(\pvt)^2\Phi_{\Int}^{(n;j)}
                      -H^{(j)}\pvt\Phi_{\Int}^{(n;j)}\bigr]\\
 &\quad+\sum_{m<j\leq n}\mathscr G_j^{\Int}
       \bigl[Z^{(j)}(\pvt)^2\Phi_{\Int}^{(n;j)}
                      -H^{(j)}\pvt\Phi_{\Int}^{(n;j)}\bigr]\\
 &=:T_4+T_5+T_6.
 \end{aligned}
 \label{eq:precision-stream-term-definitions}
\end{equation}
Here $\Phi_{\Int}^{(n;n)}=\Phi_{\Ext}^{(n)}$ by
\eqref{eq:exact-extension-endpoint}, and $T_6=0$ when $m=n$.
The $T_4$ bound is \eqref{eq:direct-source-precision-bound}.

For $T_5,T_6$ in \eqref{eq:precision-stream-term-definitions}, first
apply the triangle inequality in the two-component $L^2$ norm, then
weighted Cauchy--Schwarz:
\begin{equation}
 \begin{aligned}
 &\sum_{\iota=\lambda,w}
   \bigl(\|\mathcal M_\iota^{(n)}\Delta_LT_5\|_2^2
        +\|\mathcal M_\iota^{(n)}\Delta_LT_6\|_2^2\bigr)\\
 &\quad\leq\left(\sum_{j=1}^n\varsigma_j^{-1/2}\right)
 \Biggl[\sum_{j\leq m}\varsigma_j^{1/2}
   \sum_{\iota=\lambda,w}
   \|\mathcal M_\iota^{(n)}\Delta_L\mathscr G_j^{\Ext}
                 \mathcal C^{(j)}[\Phi_{\Int}^{(n;j)}]\|_2^2\\
 &\hspace{42mm}
 +\sum_{m<j\leq n}\varsigma_j^{1/2}
   \sum_{\iota=\lambda,w}
   \|\mathcal M_\iota^{(n)}\Delta_L\mathscr G_j^{\Int}
                 \mathcal C^{(j)}[\Phi_{\Int}^{(n;j)}]\|_2^2\Biggr].
 \end{aligned}
 \label{eq:precision-shell-sum-before-square}
\end{equation}
The first factor is bounded independently of $n$.
Apply \eqref{eq:T1-precision}--\eqref{eq:T2-precision} at scale $j$,
with coefficient layer $j$ and stream $\Phi_{\Int}^{(n;j)}$.
The choice of $C_{\mathrm E}$ in $\varsigma_j$ gives
\begin{equation}
 \begin{aligned}
 &\sum_{\iota=\lambda,w}
   \bigl(\|\mathcal M_\iota^{(n)}\Delta_LT_5\|_2^2
        +\|\mathcal M_\iota^{(n)}\Delta_LT_6\|_2^2\bigr)\\
 &\quad\lesssim\sum_{j\leq n}\varsigma_j
 \Bigl[(c_{Z,H;\mathrm{low}}^{(j)})^2
       \sum_{\iota=\lambda,w}
       \|\mathcal M_\iota^{(n)}\Delta_L\Phi_{\Int}^{(n;j)}\|_2^2
       +\mathcal{CCK}_{Z,H}^{(j;n)}\avg t^4
        \|\Phi_{\Int}^{(n;j)}\|_{\Gs{s}{\lambda_n}{\sigma-3}}^2\Bigr].
 \end{aligned}
 \label{eq:whole-product-coefficient-budget}
\end{equation}
We first estimate the second term above in
\eqref{eq:whole-product-coefficient-budget}.
By \eqref{eq:lossy-int}, \eqref{eq:elliptic-source-bootstrap}, and
\eqref{eq:Gamma-common-scale-bound},
\begin{equation}
 \begin{aligned}
 &\sum_{j\leq n}\varsigma_j\mathcal{CCK}_{Z,H}^{(j;n)}\avg t^4
       \|\Phi_{\Int}^{(n;j)}\|_{\Gs{s}{\lambda_n}{\sigma-3}}^2\\
 &\quad\lesssim\mathfrak S_n^2
       \sum_{j\leq n}\varsigma_j^3\mathcal{CCK}_{Z,H}^{(j;n)}\\
 &\quad\lesssim M^C\eps^2\delta_n^2
       \sum_{j\leq n}\varsigma_j^3\CK_{f,\partial g,h}^{(j)}
       +M^C\eps^4\delta_n^2
       \sum_{j\leq n}\varsigma_j^3\delta_j^2
                         \sum_{r<j}\varsigma_r\CK_{f,\partial g,h}^{(r)} \\
& \quad \lesssim \CK_{f,\partial g,h}^{(n)}
       +\delta_n^2\sum_{j<n}\delta_j^{-1}\varsigma_j
                                         \CK_{f,\partial g,h}^{(j)} \\
&\quad \lesssim \mathcal Q_{\rm obs}^{(n)}
 \end{aligned}
\end{equation}

Combining the $T_4$ estimate \eqref{eq:direct-source-precision-bound}
with \eqref{eq:whole-product-coefficient-budget}, and using $\varsigma_m\geq1$, gives for $1\leq m\leq n$,
\begin{equation}
 \begin{aligned}
 \sum_{\iota=\lambda,w}
 \|\mathcal M_\iota^{(n)}\Delta_L\Phi_{\Int}^{(n;m)}\|_2^2
 &\leq C\varsigma_m\mathcal Q_{\rm obs}^{(n)}+C\sum_{j\leq n}\varsigma_j(c_{Z,H;\mathrm{low}}^{(j)})^2
       \sum_{\iota=\lambda,w}
       \|\mathcal M_\iota^{(n)}\Delta_L\Phi_{\Int}^{(n;j)}\|_2^2.
 \end{aligned}
 \label{eq:precision-whole-product-coupled-system}
\end{equation}
From here, we can get:
\[
 \begin{aligned}
 &\max_{1\leq m\leq n}\frac1{\varsigma_m}
       \sum_{\iota=\lambda,w}
       \|\mathcal M_\iota^{(n)}\Delta_L\Phi_{\Int}^{(n;m)}\|_2^2\\
 &\quad\leq C\mathcal Q_{\rm obs}^{(n)}
       +C\left(\sum_{j\leq n}\varsigma_j^2(c_{Z,H;\mathrm{low}}^{(j)})^2\right)
       \max_{1\leq m\leq n}\frac1{\varsigma_m}
       \sum_{\iota=\lambda,w}
       \|\mathcal M_\iota^{(n)}\Delta_L\Phi_{\Int}^{(n;m)}\|_2^2.
 \end{aligned}
\]
The constant $C$ is defined here to be the common constant to these two inequalities.
By \eqref{eq:common-scale-small-sums}, we can choose
$C\varepsilon_{\mathrm{abs}}\leq1/2$ in order to absorb the last term to the left-hand side.
At $m=n$ and $m<n$, respectively, this yields
\begin{align}
 \sum_{\iota=\lambda,w}
 \|\mathcal M_\iota^{(n)}\Delta_L\Phi_{\Ext}^{(n)}\|_2^2
 &\leq2C\varsigma_n\mathcal Q_{\rm obs}^{(n)},
 \label{eq:precision-exterior-common-closed}\\
 \sum_{\iota=\lambda,w}
 \|\mathcal M_\iota^{(n)}\Delta_L\Phi_{\Int}^{(n;m)}\|_2^2
 &\leq2C\varsigma_m\mathcal Q_{\rm obs}^{(n)}.
 \label{eq:interior-precision-final}
\end{align}
This proves \eqref{eq:precision-elliptic-increment}.

\medskip
\noindent\emph{3. The partial-sum estimates
\eqref{eq:precision-elliptic-partial-sum} and
\eqref{eq:precision-source-partial-sum}.}
For the stream function, \eqref{eq:precision-radius-comparison} and
\eqref{eq:precision-exterior-common-closed} at each level $m$ implies:
\[
 \begin{aligned}
 \left(\sum_{\iota=\lambda,w}
       \|\mathcal M_\iota^{(n)}\Delta_L\Phi_{\Ext}^{[n]}\|_2^2\right)^{1/2}
 &\lesssim\sum_{m\leq n}\sqrt{\varsigma_m\CK_{f,\partial g,h}^{(m)}}+\left(\sum_{m\leq n}\delta_m\sqrt{\varsigma_m}\right)
       \left[\left(\sum_{r\leq n}\delta_r^{-1}\varsigma_r
                              \CK_{f,\partial g,h}^{(r)}\right)^{1/2}
             +\frac{M^C\eps}{\avg t^{(1+\delta)/2}}\right].
 \end{aligned}
\]
Weighted Cauchy--Schwarz and \eqref{eq:elliptic-shell-weight-sums} give
\begin{equation}
 \begin{aligned}
 \left(\sum_{m\leq n}\sqrt{\varsigma_m\CK_{f,\partial g,h}^{(m)}}\right)^2
 &\leq\left(\sum_{m\leq n}\delta_m\varsigma_m\right)
       \left(\sum_{m\leq n}\delta_m^{-1}\CK_{f,\partial g,h}^{(m)}\right)\lesssim\sum_{m\leq n}\delta_m^{-1}\varsigma_m
                              \CK_{f,\partial g,h}^{(m)}.
 \end{aligned}
 \label{eq:precision-weighted-CK-sum}
\end{equation}
Also $\sum_{m\leq n}\delta_m\sqrt{\varsigma_m}
\leq\sum_{m\geq1}\delta_m\varsigma_m<\infty$.
Squaring the preceding stream estimate therefore proves
\begin{equation}
 \sum_{\iota=\lambda,w}
 \|\mathcal M_\iota^{(n)}\Delta_L\Phi_{\Ext}^{[n]}\|_2^2
 \lesssim\sum_{m\leq n}\delta_m^{-1}\varsigma_m
                              \CK_{f,\partial g,h}^{(m)}
       +\frac{M^C\eps^2}{\avg t^{1+\delta}},
 \label{eq:precision-partial-sum-chain}
\end{equation}
which is \eqref{eq:precision-elliptic-partial-sum}. By applying \eqref{eq:precision-radius-comparison} and
\eqref{eq:source-precision-direct-full-line} at each level $m$, we get 
\begin{equation}
 \begin{aligned}
 \sum_{\iota=\lambda,w}
 \left\|\mathcal M_\iota^{(n)}\Pneq\sum_{m\leq n}S^{(m)}\right\|_2^2
 &\lesssim\left(\sum_{m\leq n}\sqrt{\mathcal Q_{\rm obs}^{(m)}}\right)^2\\
 &\lesssim\left(\sum_{m\leq n}\sqrt{\CK_{f,\partial g,h}^{(m)}}\right)^2
       +\left(\sum_{m\leq n}\delta_m\right)^2
        \left[\sum_{r\leq n}\delta_r^{-1}\varsigma_r
                              \CK_{f,\partial g,h}^{(r)}
              +\frac{M^C\eps^2}{\avg t^{1+\delta}}\right]\\
 &\lesssim\sum_{m\leq n}\delta_m^{-1}\varsigma_m
                              \CK_{f,\partial g,h}^{(m)}
       +\frac{M^C\eps^2}{\avg t^{1+\delta}}.
 \end{aligned}
 \label{eq:precision-source-partial-sum-chain}
\end{equation}
This proves \eqref{eq:precision-source-partial-sum} and completes
the strong induction for Proposition~\ref{prop:precision-elliptic}.
\end{proof}

\subsection{\texorpdfstring{Consequences for the transport estimates}{Consequences for the transport estimates}}
\label{subsec:transport-consequences}

Propositions~\ref{prop:lossy-elliptic} and~\ref{prop:precision-elliptic}
are now proved.  We record consequences for the transport
estimates.

\begin{corollary}
\label{cor:elliptic-profile-output}
Define the elliptic part of
\eqref{eq:partial-transport-field} via 
\begin{equation}
  \mathbf V_{\phi}^{[n]}
  =\left(-(1+h^{[n]})\partial_v\phi_{\neq}^{[n]},
          (1+h^{[n]})\partial_z\phi_{\neq}^{[n]}\right).
  \label{eq:elliptic-profile-transport-field}
\end{equation}
Then, for every fixed radius factor \(c\in(0,1)\),
\begin{align}
  \norm{\avg t^2\mathbf V_{\phi}^{[n]}}_{
       \Gs{s}{c\lambda_n}{\sigma-6}}
  &\lesssim M^C\eps,
  \label{eq:summed-profile-transport}\\
  \norm{\pvt\Phi_{\Ext,\neq}^{[n]}}_{
       \Gs{s}{c\lambda_n}{\sigma-5}}
  +\avg t\norm{\partial_z\Phi_{\Ext,\neq}^{[n]}}_{
       \Gs{s}{c\lambda_n}{\sigma-5}}
  &\lesssim\frac{M^C\eps}{\avg t}.
  \label{eq:summed-physical-velocity}
\end{align}
\end{corollary}

\begin{proof} Equations \eqref{eq:lossy-velocity} and
\eqref{eq:lossy-profile-velocity}, the radius gap
\eqref{eq:radius-gap}, and \eqref{eq:elliptic-source-bootstrap} give
\[
  \norm{\avg t^2\partial_v\Phi_{\Ext,\neq}^{(j)}}_{
       \Gs{s}{c\lambda_n}{\sigma-6}}
   +\norm{\avg t^2\partial_z\Phi_{\Ext,\neq}^{(j)}}_{
       \Gs{s}{c\lambda_n}{\sigma-6}}
  \lesssim
    M^C\eps\,e^{C2^{j\alpha_\dagger}}\delta_j,
\]
and the same argument gives \eqref{eq:summed-physical-velocity}.  Since
\(\alpha_\dagger<\alpha_*\), the series
\(\sum_j e^{C2^{j\alpha_\dagger}}\delta_j\) converges.  Finally,
the product estimate and the bootstrap bound on \(h^{[n]}\) allow us to
include the factor \(1+h^{[n]}\) in \eqref{eq:elliptic-profile-transport-field}.
\end{proof}

\begin{lemma}
\label{lem:newton-velocity-output}
For $1\leq m\leq n$, for every fixed $c\in(0,1)$,
\begin{align}
  \norm{\chi_{m+2}\mathbf P^{(n)}}_{
       \Gs{s}{c\lambda_n}{\sigma-7}}
  &\lesssim M^C\eps\delta_n\varsigma_m
    \left(\frac1{\avg t^2}
          +\frac1{\avg t^{\beta_g/2}}\right),
  \label{eq:newton-velocity-interior}\\
  \norm{\mathbf P^{(n)}}_{
       \Gs{s}{c\lambda_n}{\sigma-7}}
  &\lesssim M^C\eps\delta_n
    \left(\frac{\varsigma_n}{\avg t^2}
          +\frac1{\avg t^{\beta_g/2}}\right).
  \label{eq:newton-velocity-global}
\end{align}
\end{lemma}

\begin{proof}
We first prove \eqref{eq:newton-velocity-interior}. Fix $0<c<c'<1$.
By \eqref{eq:exact-shell-partition} and
\eqref{eq:interior-cutoff-geometry},
$\{\chi_{m+2}\ne0\}\subset\{\mathbf d>d_m/64=d_{m+1}/32\}$.
Thus \eqref{eq:exact-extension-endpoint} and
\eqref{eq:extension-compatibility} identify the derivatives of
$\Phi_{\Int}^{(n;m+1)}$ with those of $\phi_{\neq}^{(n)}$ on this set
when $m<n$; for $m=n$ we use $\Phi_{\Ext}^{(n)}$.
Expanding \eqref{eq:partial-transport-field}, we obtain three components
\begin{equation}
\begin{aligned}
 T_1&:=h^{(n)}
   \left(-\partial_v\Phi_{\Ext}^{[n-1]},
                 \partial_z\Phi_{\Ext}^{[n-1]}\right),\\
 T_2&:=
 \begin{cases}
 (1+h^{[n]})\left(-\partial_v\Phi_{\Int}^{(n;m+1)},
                         \partial_z\Phi_{\Int}^{(n;m+1)}\right),&m<n,\\
 (1+h^{[n]})\left(-\partial_v\Phi_{\Ext}^{(n)},
                         \partial_z\Phi_{\Ext}^{(n)}\right),&m=n,
 \end{cases}\\
 T_3&:=(0,g^{(n)}),\qquad
 \chi_{m+2}\mathbf P^{(n)}=\chi_{m+2}(T_1+T_2+T_3).
\end{aligned}
\label{eq:newton-velocity-term-decomposition}
\end{equation}
Here, recall that $\Phi_{\Ext}^{[n-1]}=\sum_{r<n}\Phi_{\Ext}^{(r)}$. We discuss these three terms separately. 

For the $T_1$-term in \eqref{eq:newton-velocity-term-decomposition}, we combine the Gevrey
product inequality \cite[Lemma~A.3]{BM15}, the $\lambda_n$-comparison
\eqref{eq:radius-gap}, the $(\pw hn,\pw Zn,\pw Hn)$-bound
\eqref{eq:bootstrap-consequences-coefficients}, and the lossy elliptic estimates 
\eqref{eq:lossy-velocity}--\eqref{eq:lossy-profile-velocity} (supplemented by 
\eqref{eq:elliptic-shell-weight-sums} and \eqref{eq:elliptic-source-bootstrap}) to obtain that for a fixed constant $C_v>0$\footnote{
During the proof, we will introduce two local constants $C_v,C_\chi$, which will not be used elsewhere in the paper. },
\begin{equation}
\begin{aligned}
 \norm{T_1}_{\Gs{s}{c'\lambda_n}{\sigma-7}}
 &\lesssim\norm{h^{(n)}}_{\Gs{s}{c'\lambda_n}{\sigma-7}}
       \sum_{r<n}\norm{\nabla\Phi_{\Ext}^{(r)}}_{
                                     \Gs{s}{c'\lambda_r}{\sigma-7}}\\
 &\lesssim M^C\eps^2\delta_n\avg t^{-2}
                    \sum_{r<n}e^{C_v2^{r\alpha_\dagger}}\delta_r
 \lesssim M^C\eps^2\delta_n\avg t^{-2}.
\end{aligned}
\label{eq:newton-velocity-T1-bound}
\end{equation}

For the $T_2$-term in \eqref{eq:newton-velocity-term-decomposition}, we distinguish between the $m<n$ case and the $m=n$ case. For the first case, applications of the
Gevrey product inequality, the $(f^{[n]}, h^{[n]})$-bound \eqref{eq:bootstrap-consequences-partial-fh}
and the lossy elliptic estimates \eqref{eq:lossy-velocity}--\eqref{eq:lossy-profile-velocity} (supplemented with 
\eqref{eq:elliptic-source-bootstrap}) yields
\begin{equation}
\begin{aligned}
 \norm{T_2}_{\Gs{s}{c'\lambda_n}{\sigma-7}}
 &\lesssim\left(1+\norm{h^{[n]}}_{\Gs{s}{c'\lambda_n}{\sigma-7}}\right)
       \norm{\nabla\Phi_{\Int}^{(n;m+1)}}_{
                                      \Gs{s}{c'\lambda_n}{\sigma-7}}\\
 &\lesssim M^C\eps\delta_n
            e^{C_v2^{(m+1)\alpha_\dagger}}\avg t^{-2},\qquad m<n.
\end{aligned}
\label{eq:newton-velocity-T2-interior-bound}
\end{equation}
For the $m=n$ case, the lossy elliptic estimate \eqref{eq:lossy-velocity}--\eqref{eq:lossy-profile-velocity} (supplemented with 
\eqref{eq:elliptic-source-bootstrap}) associated with the $\pex n$  gives 
\begin{equation}
 \norm{T_2}_{\Gs{s}{c'\lambda_n}{\sigma-7}}
 \lesssim\left(1+\norm{h^{[n]}}_{\Gs{s}{c'\lambda_n}{\sigma-7}}\right)
       \norm{\nabla\Phi_{\Ext}^{(n)}}_{\Gs{s}{c'\lambda_n}{\sigma-7}}
 \lesssim M^C\eps\delta_n e^{C_v2^{n\alpha_\dagger}}\avg t^{-2}.
\label{eq:newton-velocity-T2-exterior-bound}
\end{equation}

Finally, for the $T_3$-term in \eqref{eq:newton-velocity-term-decomposition}, the $\pw gn$-estimate 
\eqref{eq:bootstrap-consequences-fg} gives
\begin{equation}
 \norm{T_3}_{\Gs{s}{c'\lambda_n}{\sigma-7}}
 =\norm{g^{(n)}}_{\Gs{s}{c'\lambda_n}{\sigma-7}}
 \lesssim M\eps\delta_n\avg t^{-\beta_g/2}.
\label{eq:newton-velocity-T3-bound}
\end{equation} Here, we recall the $\beta_g=4-a/2$ from \eqref{eq:component-energies}.

Before collecting all the $T_1,T_2,T_3$ bounds together, we invoke the bounds
\eqref{eq:alpha-dagger} and \eqref{eq:extension-cutoff-derivatives} to obtain a quantitative description of the cutoff function $\chi_{m+2}$,
\begin{equation}
 \bigl\|e^{c\lambda_n|\eta|^s}\avg\eta^{\sigma-7}
                  \widehat\chi_{m+2}(\eta)\bigr\|_{L^1_\eta}
 \lesssim e^{C_\chi2^{m\alpha_\dagger}}.
\label{eq:newton-velocity-cutoff-bound}
\end{equation}
Take the fixed constant $C_{\mathrm E}$ in
\eqref{eq:elliptic-shell-weights} so that
$C_{\mathrm E}\geq C_\chi+2^{\alpha_\dagger}C_v$. Then
\begin{equation*}
 e^{C_\chi2^{m\alpha_\dagger}}
 e^{C_v2^{(m+1)\alpha_\dagger}}
 =e^{(C_\chi+2^{\alpha_\dagger}C_v)2^{m\alpha_\dagger}}
 \leq\varsigma_m.
\end{equation*}
We then apply Young's convolution inequality, the $\chi$-bound \eqref{eq:newton-velocity-cutoff-bound} and the
$T_i$-controls 
\eqref{eq:newton-velocity-T1-bound}, \eqref{eq:newton-velocity-T2-interior-bound}, \eqref{eq:newton-velocity-T2-exterior-bound}, \eqref{eq:newton-velocity-T3-bound}.
Recalling \eqref{eq:newton-velocity-term-decomposition},
we obtain
\begin{equation*}
 \norm{\chi_{m+2}\mathbf P^{(n)}}_{\Gs{s}{c\lambda_n}{\sigma-7}}
 \lesssim e^{C_\chi2^{m\alpha_\dagger}}
               \sum_{i=1}^3\norm{T_i}_{\Gs{s}{c'\lambda_n}{\sigma-7}}
 \lesssim M^C\eps\delta_n\varsigma_m
                 \left(\avg t^{-2}+\avg t^{-\beta_g/2}\right),
\end{equation*}
which proves \eqref{eq:newton-velocity-interior}.

For \eqref{eq:newton-velocity-global}, recall the decomposition 
\eqref{eq:newton-velocity-term-decomposition} with $m=n$.
The bounds \eqref{eq:newton-velocity-T1-bound},
\eqref{eq:newton-velocity-T2-exterior-bound}, and
\eqref{eq:newton-velocity-T3-bound} give
$
 \norm{\mathbf P^{(n)}}_{\Gs{s}{c\lambda_n}{\sigma-7}}
 \lesssim M^C\eps\delta_n
          \left(e^{C_v2^{n\alpha_\dagger}}\avg t^{-2}
                                  +\avg t^{-\beta_g/2}\right)
 \lesssim M^C\eps\delta_n
          \left(\varsigma_n\avg t^{-2}+\avg t^{-\beta_g/2}\right).
$
\end{proof}

\section{Reaction and transport multiplier estimates}
\label{sec:multipliers}

We prove Propositions~\ref{prop:reaction-remainder-multiplier}
and~\ref{prop:transport-multiplier}, which provide the reaction and transport
bounds used in Section~\ref{sec:energy}. In the Fourier expansion of $\nabla^\perp\phi\cdot\nabla f$, write
\begin{equation}
 \underbrace{(k,\eta)}_{\text{output}}
 =\underbrace{(\ell,\xi)}_{\text{stream function }\phi}
 +\underbrace{(p,\zeta)}_{\text{transported }f}.
\label{eq:nonlinear-frequency-convention}
\end{equation}
Define
\begin{align}
  \mathsf{HL}
    &=\left\{((k,\eta),(\ell,\xi)): |(p,\zeta)|\leq\tfrac18|(\ell,\xi)|\right\},
  \label{eq:HL-region}\\
  \mathsf{LH}
    &=\left\{((k,\eta),(\ell,\xi)): |(\ell,\xi)|\leq\tfrac18|(p,\zeta)|\right\},
  \label{eq:LH-region}\\
  \mathsf{HH}
    &=\left\{((k,\eta),(\ell,\xi)): \tfrac18|(\ell,\xi)|<|(p,\zeta)|<8|(\ell,\xi)|\right\}.
  \label{eq:HH-region}
\end{align}
These are, respectively, the reaction, transport, and
comparable-frequency regions.  

\subsection{Reaction and comparable-frequency estimates}

\begin{proposition}[Reaction and remainder multipliers]
\label{prop:reaction-remainder-multiplier}
There is \(c=c(s)\in(0,1)\) such that the following estimates hold for every
\(n\geq1\) and every fixed \(((k,\eta),(\ell,\xi))\in\mathsf{HL}\cup\mathsf{HH}\),
for almost every \(t\geq1\).
For $\ell\ne0$, the following estimate holds:
\begin{align}
 &\frac{\A_k^{(n)}(t,\eta)}{\A_\ell^{(n)}(t,\xi)}
 \frac{1}{\ell^2(1+|t-\xi/\ell|^2)}
 \begin{cases}\avg{\ell,\xi},&\ell\ne k,\\ |\ell|,&\ell=k,\end{cases}
 \notag\\
 &\quad\lesssim e^{c\lambda_n(t)|(p,\zeta)|^s}
 \Biggl[\sqrt{\frac{\partial_tw_k(t,\eta)}{w_k(t,\eta)}
              \frac{\partial_tw_\ell(t,\xi)}{w_\ell(t,\xi)}
              \frac{\widetilde\A_k^{(n)}(t,\eta)}{\A_k^{(n)}(t,\eta)}}
              \frac{\widetilde\A_\ell^{(n)}(t,\xi)}{\A_\ell^{(n)}(t,\xi)}
       +(-\dot\lambda_n)\avg{k,\eta}^{s/2}\avg{\ell,\xi}^{s/2}\Biggr]
 \left\langle\frac{\xi}{\ell t}\right\rangle^{-1}.
 \label{eq:reaction-multiplier-off-diagonal}
\end{align}
\begingroup
\expandafter\edef\csname @currentlabel\endcsname{\theequation}
\label{eq:reaction-multiplier-diagonal}
\endgroup
If $\ell=0$, the following estimate holds:
\begin{align}
  \frac{{\A}_k^{(n)}(t,\eta)}{{\A}_0^{(n)}(t,\xi)}
  \bm 1_{\left\{\ell=0\right\}}
  &\lesssim
  e^{c\lambda_n(t)|(k,\eta-\xi)|^s}
  \Bigg[
    \avg t^{1+s}\sqrt{\frac{\partial_tw_{k}(\eta)}{w_{k}(\eta)} \frac{\partial_tw_{0}(\xi)}{w_{0}(\xi)}
        \frac{\widetilde\A_k^{(n)}(t,\eta)}{\A_k^{(n)}(t,\eta)}}
        \frac{\widetilde\A_0^{(n)}(t,\xi)}{\A_0^{(n)}(t,\xi)}
    +\frac{\avg{k,\eta}^{s/2}}{\avg{\xi}^{s/2}}
  \Bigg].
  \label{eq:reaction-multiplier-zero-mode}
\end{align}
\end{proposition}

\begin{remark}
\label{rem:reaction-quotients}
We retain the quotients $\widetilde\A_k^{(n)}/\A_k^{(n)}$ in
\eqref{eq:reaction-multiplier-off-diagonal}--\eqref{eq:reaction-multiplier-zero-mode},
in contrast to \cite{BM13}, to simplify
the presentation in Section~\ref{sec:energy}. In
\eqref{eq:f-reaction-stream-majorants}, they combine with the two
high-frequency multipliers to give exactly the $\CK$ and elliptic weights
in $\CK_{f;w}^{(n)}$ and $\mathcal M_w^{(n)}$. The passage to
\eqref{eq:f-reaction-stream-branch-bounds} therefore requires no additional
frequency restrictions for comparing $\A^{(n)}$ with $\widetilde\A^{(n)}$.
\end{remark}

\begin{proof}[Proof of Proposition~\ref{prop:reaction-remainder-multiplier}] 
We first prove \eqref{eq:reaction-multiplier-off-diagonal} for $\ell\ne0$,
treating $\mathsf{HL}$ and $\mathsf{HH}$ separately. We then prove
\eqref{eq:reaction-multiplier-zero-mode} for $\ell=0$, using the preceding
Gevrey and Sobolev bounds. We begin with two elementary estimates.
For $|k|\leq|\eta|$,
the definitions \eqref{eq:J-mode} and \eqref{eq:main-multiplier},
together with $w_k\leq1$, give
\begin{equation}
 \frac12\leq\frac{\widetilde\A_k^{(n)}(t,\eta)}{\A_k^{(n)}(t,\eta)}\leq1
 \qquad(|k|\leq|\eta|).
 \label{eq:theta-vertical-lower-bound}
\end{equation}
If $((k,\eta),(\ell,\xi))\in\mathsf{HL}$ and $|\ell|>|\xi|$, then
\begin{equation}
 \frac{\avg{\ell,\xi}}
 {\ell^2+|\xi-\ell t|^2}
 \lesssim\frac1{\avg t^2\avg{\ell,\xi}}
 \lesssim(-\dot\lambda_n(t))
   \avg{k,\eta}^{s/2}\avg{\ell,\xi}^{s/2},
 \qquad t\geq1,
 \label{eq:horizontal-elliptic-paid-by-radius}
\end{equation}

\noindent{\bf Step 1: Proof of
\eqref{eq:reaction-multiplier-off-diagonal} for $\ell\ne0$.}

We first record the Gevrey and Sobolev bounds used in all four
$\mathsf{HL}$ cases. The triangle inequality and the defining restriction
$|(p,\zeta)|\leq|(\ell,\xi)|/8$ give
\begin{equation}
 \frac{\avg{k,\eta}^{\sigma}}{\avg{\ell,\xi}^{\sigma}}
 =\left(\frac{1+|(k,\eta)|^2}{1+|(\ell,\xi)|^2}\right)^{\sigma/2}
 \leq\left(\frac98\right)^\sigma.
 \label{eq:reaction-HL-Sobolev-ratio}
\end{equation}
We also observe the following consequence of the strict concavity of
$r\mapsto r^s$, where $1/2<s<1$ and $c'=s8^{s-1}<1$:
\begin{equation}
 |(k,\eta)|^s-|(\ell,\xi)|^s
 \leq s|(\ell,\xi)|^{s-1}|(p,\zeta)|
 \leq c'|(p,\zeta)|^s.
 \label{eq:reaction-HL-concavity}
\end{equation}
For any fixed $c\in(c',1)$ and $B\geq0$, the lower bound
$\lambda_n(t)\geq3\lambda_0/8$ from \eqref{eq:lambda-n} gives
\begin{equation}
\begin{aligned}
 e^{c'\lambda_n(t)|(p,\zeta)|^s+C\mu|(p,\zeta)|^{1/2}}
   \avg{p,\zeta}^{B}
 &=e^{c\lambda_n(t)|(p,\zeta)|^s}
   e^{C\mu|(p,\zeta)|^{1/2}-(c-c')\lambda_n(t)|(p,\zeta)|^s}
   \avg{p,\zeta}^{B}\lesssim_B e^{c\lambda_n(t)|(p,\zeta)|^s}.
\end{aligned}
 \label{eq:reaction-HL-exponential-absorption}
\end{equation}

Now we distinguish between several cases to prove the result.

\noindent {\bf Case 1a)} $\mathsf{HL}$ as defined in
\eqref{eq:HL-region}, $k\ne\ell$, and
$t\in\mathbf I_{k,\eta}\cap\mathbf I_{k,\xi}$.
The intervals $\mathbf I_{k,\eta}$ and $\mathbf I_{k,\xi}$ have the same
horizontal index $k$. By the definitions \eqref{eq:critical-endpoints}
and \eqref{eq:restricted-critical-intervals}, if their intersection is
nonempty, we necessarily have $|\eta|>4k^2, |\xi|>4k^2$. In particular, $k\ne0$ and $|k|<|\xi|/4$.
Using the $\mathsf{HL}$ restriction \eqref{eq:HL-region}, we further obtain
$7|\ell|\leq8|k|+|\xi|<3|\xi|$ and
$|(p,\zeta)|\leq(|\ell|+|\xi|)/8<5|\xi|/28$.
Consequently,
\begin{equation}
 \frac{23}{28}|\xi|<|\eta|<\frac{33}{28}|\xi|,
 \qquad |\ell|<|\xi|,\qquad |k|\leq|\eta|.
 \label{eq:reaction-1a-frequency-bounds}
\end{equation}
These inequalities verify the frequency comparability required by the
estimate \eqref{eq:raw-J-two-CK-ledger}.
Application of \eqref{eq:raw-J-two-CK-ledger} yields
\[
 \frac{\J_k(t,\eta)}{\J_\ell(t,\xi)}
 \lesssim \frac{|\eta|}{k^2}
 \sqrt{\frac{\partial_tw_k(t,\eta)}{w_k(t,\eta)}
       \frac{\partial_tw_\ell(t,\xi)}{w_\ell(t,\xi)}}
 e^{C\mu|(p,\zeta)|^{1/2}}.
\]
By the multiplier definition \eqref{eq:main-multiplier}, the Sobolev
ratio \eqref{eq:reaction-HL-Sobolev-ratio}, and the Gevrey bounds
\eqref{eq:reaction-HL-concavity}--\eqref{eq:reaction-HL-exponential-absorption},
we therefore have
\begin{equation}
 \frac{\A_k^{(n)}(t,\eta)}{\A_\ell^{(n)}(t,\xi)}
 \lesssim e^{c\lambda_n(t)|(p,\zeta)|^s}
 \frac{|\eta|}{k^2}
 \sqrt{\frac{\partial_tw_k(t,\eta)}{w_k(t,\eta)}
       \frac{\partial_tw_\ell(t,\xi)}{w_\ell(t,\xi)}}.
 \label{eq:reaction-1a-multiplier-bound}
\end{equation}
We next estimate
$\frac{|\eta|}{k^2}\frac{\avg{\ell,\xi}}{\ell^2+|\xi-\ell t|^2}$.
Since $t\in I_{k,\xi}$ and $\ell\ne k$, the definition
\eqref{eq:critical-endpoints} gives
\[
 \left|t-\frac\xi\ell\right|
 \geq\frac12\frac{|\xi||k-\ell|}{|k\ell|}.
\]
To see this, first suppose that $\ell$ and $\xi$ have the same sign.
For $|k|\geq2$, the two endpoints of $I_{k,\xi}$ are the midpoints
between $\xi/k$ and the adjacent critical times. Thus each critical
time $\xi/\ell$, $\ell\ne k$, is at least half its distance from
$\xi/k$ away from this interval. For $|k|=1$, the lower endpoint
$3|\xi|/4$ gives the same conclusion for every $|\ell|\geq2$.
If $\ell$ and $\xi$ have opposite signs, then
$|t-\xi/\ell|=t+|\xi|/|\ell|$ and
$t\geq3|\xi|/(4|k|)$ give the same bound.
Using this separation, $|\eta|\approx|\xi|$, and
$t\geq3|\xi|/(4|k|)$, we obtain
\begin{equation}
\begin{aligned}
 \frac{|\eta|}{k^2}
 \frac{\avg{\ell,\xi}}{\ell^2+|\xi-\ell t|^2}
 &\lesssim\frac{|\xi|^2}{k^2\ell^2+|\xi|^2(k-\ell)^2},\\
 \frac{|\xi|^2\avg{\xi/(\ell t)}}{k^2\ell^2+|\xi|^2(k-\ell)^2}
 &\lesssim\frac{1+|k|/|\ell|}{|k-\ell|^2}
 \leq\frac{2+|k-\ell|}{|k-\ell|^2}\lesssim1.
\end{aligned}
 \label{eq:reaction-1a-elliptic-bounds}
\end{equation}
Here $|\ell|\geq1$ and $|k-\ell|\geq1$ are used in the last two steps.

Finally, the bounds $|k|\leq|\eta|$ and $|\ell|\leq|\xi|$ in
\eqref{eq:reaction-1a-frequency-bounds} and
the bound \eqref{eq:theta-vertical-lower-bound} give
$\frac{\widetilde\A_k^{(n)}(t,\eta)}{\A_k^{(n)}(t,\eta)}\geq\frac12$
and
$\frac{\widetilde\A_\ell^{(n)}(t,\xi)}{\A_\ell^{(n)}(t,\xi)}\geq\frac12$.
Hence the quotients can be inserted with the precise placement in
\eqref{eq:reaction-multiplier-off-diagonal}, since
\begin{equation}
\begin{aligned}
 &\sqrt{\frac{\partial_tw_k(t,\eta)}{w_k(t,\eta)}
        \frac{\partial_tw_\ell(t,\xi)}{w_\ell(t,\xi)}}
 \leq2\sqrt2
 \sqrt{\frac{\partial_tw_k(t,\eta)}{w_k(t,\eta)}
        \frac{\partial_tw_\ell(t,\xi)}{w_\ell(t,\xi)}
        \frac{\widetilde\A_k^{(n)}(t,\eta)}{\A_k^{(n)}(t,\eta)}}
 \frac{\widetilde\A_\ell^{(n)}(t,\xi)}{\A_\ell^{(n)}(t,\xi)}.
\end{aligned}
 \label{eq:reaction-1a-quotient-insertion}
\end{equation}
Combining \eqref{eq:reaction-1a-multiplier-bound},
\eqref{eq:reaction-1a-elliptic-bounds}, and
\eqref{eq:reaction-1a-quotient-insertion} proves
\eqref{eq:reaction-multiplier-off-diagonal} in Case 1a.

In Cases 1b, 1c, and 1d, we always retain the restriction
\begin{equation}
 k=\ell\qquad\text{or}\qquad
 t\notin\mathbf I_{k,\eta}\cap\mathbf I_{k,\xi}.
 \label{eq:reaction-HL-remaining-cases}
\end{equation}
Indeed, \eqref{eq:reaction-1a-multiplier-bound}--\eqref{eq:reaction-1a-quotient-insertion}
cover the alternative in Case 1a.

\noindent {\bf Case 1b)} $\mathsf{HL}$ as defined in
\eqref{eq:HL-region}, and $|\ell|>|\xi|$.
The $\mathsf{HL}$ restriction gives
$|k-\ell|<|\ell|/4$, $|k|>3|\ell|/4$, and
$|\eta|<5|\ell|/4<4k^2$, where the last inequality uses $|\ell|\geq1$.
Since a nonempty restricted
interval $\mathbf I_{k,\eta}$ requires $|\eta|>4k^2$ by the definitions
\eqref{eq:critical-endpoints}--\eqref{eq:restricted-critical-intervals},
we have $\mathbf I_{k,\eta}=\varnothing$.
Thus the parenthetical factor in \eqref{eq:J-ratio} is one, so
$\frac{\J_k(t,\eta)}{\J_\ell(t,\xi)}
\lesssim e^{C\mu|(p,\zeta)|^{1/2}}$.
Together with the bound \eqref{eq:reaction-HL-Sobolev-ratio} and
the bounds \eqref{eq:reaction-HL-concavity}--\eqref{eq:reaction-HL-exponential-absorption},
this gives
$\frac{\A_k^{(n)}(t,\eta)}{\A_\ell^{(n)}(t,\xi)}
\lesssim e^{c\lambda_n(t)|(p,\zeta)|^s}$.
Also, $|\xi|<|\ell|$ and $t\geq1$ give
$\avg{\xi/(\ell t)}\leq\sqrt2$.
Since $|\ell|t\leq|\xi-\ell t|+|\xi|$ and $|\xi|<|\ell|$, we have
$\ell^2\avg t^2\leq3(\ell^2+|\xi-\ell t|^2)$ and
$\avg{\ell,\xi}^2\leq5\ell^2$.
Hence \eqref{eq:horizontal-elliptic-paid-by-radius} gives
\[
 \frac{\avg{\ell,\xi}\avg{\xi/(\ell t)}}{\ell^2+|\xi-\ell t|^2}
 \lesssim\frac1{\avg t^2\avg{\ell,\xi}}
 \lesssim(-\dot\lambda_n(t))\avg{k,\eta}^{s/2}\avg{\ell,\xi}^{s/2},
\]
where the last inequality uses the radius derivative
\eqref{eq:lambda-derivative} and $a<1$.
Combining the multiplier and elliptic estimates gives
\[
 \frac{\A_k^{(n)}(t,\eta)}{\A_\ell^{(n)}(t,\xi)}
 \frac{\avg{\ell,\xi}\avg{\xi/(\ell t)}}{\ell^2+|\xi-\ell t|^2}
 \lesssim e^{c\lambda_n(t)|(p,\zeta)|^s}
 (-\dot\lambda_n(t))\avg{k,\eta}^{s/2}\avg{\ell,\xi}^{s/2}.
\]
For $k=\ell$, the numerator $|\ell|$ in
\eqref{eq:reaction-multiplier-off-diagonal} is no larger than
$\avg{\ell,\xi}$. Dividing by $\avg{\xi/(\ell t)}$ proves
\eqref{eq:reaction-multiplier-off-diagonal} in Case 1b.

\noindent {\bf Case 1c)} $\mathsf{HL}$, $|\ell|\leq|\xi|$, and
$t\leq8\sqrt{\max\left\{|\eta|,|\xi|,1\right\}}$.
Since $\ell\ne0$ and $|\ell|\leq|\xi|$, we have $|\xi|\geq1$.
The $\mathsf{HL}$ restriction gives
$|\eta-\xi|\leq|(p,\zeta)|\leq(|\ell|+|\xi|)/8\leq|\xi|/4$,
so $3|\xi|/4\leq|\eta|\leq5|\xi|/4$ and
$\max\{|\eta|,|\xi|,1\}\leq5|\xi|/4$.
By the radius derivative formula \eqref{eq:lambda-derivative},
the upper bound on $t$, and $a\leq2s-1$, we have
\[
 (-\dot\lambda_n(t))\avg{k,\eta}^{s/2}\avg{\ell,\xi}^{s/2}
 \gtrsim\max\{|\eta|,|\xi|,1\}^{s-(1+a)/2}\geq1.
\]
If $t\in\mathbf I_{k,\eta}$, then $k\ne0$, and
the definition \eqref{eq:critical-endpoints} gives $t\geq3|\eta|/(4|k|)$ and
\[
 \frac{|\eta|}{k^2}\leq\frac{16t^2}{9|\eta|}
 \lesssim\frac{\max\{|\eta|,|\xi|,1\}}{|\eta|}\lesssim1.
\]
If $t\notin\mathbf I_{k,\eta}$, the indicator in the estimate \eqref{eq:J-ratio}
is zero, so its parenthetical factor is one; this includes $k=0$.
Thus that factor is bounded in both cases, and the estimate \eqref{eq:J-ratio} gives
$\frac{\J_k(t,\eta)}{\J_\ell(t,\xi)}
\lesssim e^{C\mu|(p,\zeta)|^{1/2}}$.
By the definition \eqref{eq:main-multiplier}, the bound \eqref{eq:reaction-HL-Sobolev-ratio}, and
the bounds \eqref{eq:reaction-HL-concavity}--\eqref{eq:reaction-HL-exponential-absorption},
$\frac{\A_k^{(n)}(t,\eta)}{\A_\ell^{(n)}(t,\xi)}
\lesssim e^{c\lambda_n(t)|(p,\zeta)|^s}$.

For the elliptic factor, first suppose
$|\xi-\ell t|\leq|\xi|/2$.
Then $|\ell|t\geq|\xi|/2$, so $\avg{\xi/(\ell t)}\lesssim1$ and
\[
 \frac{\avg{\ell,\xi}\avg{\xi/(\ell t)}}{\ell^2+|\xi-\ell t|^2}
 \lesssim\frac{|\xi|}{\ell^2}
 \lesssim\frac{t^2}{|\xi|}\lesssim1.
\]
If $|\xi-\ell t|>|\xi|/2$, then instead
\[
 \frac{\avg{\ell,\xi}\avg{\xi/(\ell t)}}{\ell^2+|\xi-\ell t|^2}
 \lesssim\frac1{|\xi|}+\frac1{|\ell|t}\lesssim1,
\]
using $|\xi|\geq1$, $|\ell|\geq1$, and $t\geq1$.
Combining the multiplier and elliptic estimates with the lower bound
for the radius term, we obtain
\[
 \frac{\A_k^{(n)}(t,\eta)}{\A_\ell^{(n)}(t,\xi)}
 \frac{\avg{\ell,\xi}\avg{\xi/(\ell t)}}{\ell^2+|\xi-\ell t|^2}
 \lesssim e^{c\lambda_n(t)|(p,\zeta)|^s}
 (-\dot\lambda_n(t))\avg{k,\eta}^{s/2}\avg{\ell,\xi}^{s/2}.
\]
For $k=\ell$, use $|\ell|\leq\avg{\ell,\xi}$.
After division by $\avg{\xi/(\ell t)}$, this proves
\eqref{eq:reaction-multiplier-off-diagonal} in Case 1c.

\noindent {\bf Case 1d)} $\mathsf{HL}$, $|\ell|\leq|\xi|$, and
$t>8\sqrt{\max\left\{|\eta|,|\xi|,1\right\}}$.
We retain the restriction \eqref{eq:reaction-HL-remaining-cases}.
Here $|\eta-\xi|\leq|(p,\zeta)|\leq(|\ell|+|\xi|)/8\leq|\xi|/4$,
so $|\eta|\approx|\xi|$.
By the definition \eqref{eq:main-multiplier}, the left-hand side of
\eqref{eq:reaction-multiplier-off-diagonal}, multiplied by
$\avg{\xi/(\ell t)}$, factors as
\[
 \frac{\A_k^{(n)}(t,\eta)}{\A_\ell^{(n)}(t,\xi)}
 \frac{\avg{\xi/(\ell t)}}{\ell^2+|\xi-\ell t|^2}
 \begin{cases}\avg{\ell,\xi},&k\ne\ell,\\ |\ell|,&k=\ell,\end{cases}
 =:\mathfrak F_1\mathfrak F_2,
\]
where
\begin{align*}
 \mathfrak F_1
 &=e^{\lambda_n(t)(|(k,\eta)|^s-|(\ell,\xi)|^s)}
   \frac{\avg{k,\eta}^{\sigma}}{\avg{\ell,\xi}^{\sigma}},\\*
 \mathfrak F_2
 &=\frac{\J_k(t,\eta)}{\J_\ell(t,\xi)}
   \frac{\avg{\xi/(\ell t)}}{\ell^2+|\xi-\ell t|^2}
   \begin{cases}\avg{\ell,\xi},&k\ne\ell,\\ |\ell|,&k=\ell.\end{cases}
\end{align*}
We first estimate $\mathfrak F_1$, including the multiplier quotients
needed for the $\CK$ term, and then estimate $\mathfrak F_2$ in four subcases.
By the bound \eqref{eq:reaction-HL-Sobolev-ratio} and the estimate \eqref{eq:reaction-HL-concavity},
\begin{equation}
 \mathfrak F_1\leq(9/8)^\sigma e^{c'\lambda_n(t)|(p,\zeta)|^s}.
 \label{eq:reaction-1d-first-factor}
\end{equation}
By the definition \eqref{eq:critical-endpoints}, if $t\in I_{\ell,\xi}$, then 
$|\ell|\leq\sqrt{|\xi|}$.  For $|\xi|\geq4$, the $\mathsf{HL}$ restriction gives
\[
 |(p,\zeta)|\leq\tfrac18(\sqrt{|\xi|}+|\xi|)\leq\tfrac3{16}|\xi|,
 \qquad |k|\leq\tfrac{11}{16}|\xi|<\tfrac{13}{16}|\xi|\leq|\eta|.
\]
Thus the bound \eqref{eq:theta-vertical-lower-bound} gives
$\widetilde\A_k^{(n)}(t,\eta)/\A_k^{(n)}(t,\eta)\geq1/2$ and
$\widetilde\A_\ell^{(n)}(t,\xi)/\A_\ell^{(n)}(t,\xi)\geq1/2$.
If $|\xi|<4$, then $t\in I_{\ell,\xi}$ implies $t\leq2|\xi|<8$,
contradicting $t>8\sqrt{\max\{|\eta|,|\xi|,1\}}$ in Case 1d.

These lower bounds and
\eqref{eq:reaction-1d-first-factor} imply
\begin{equation}
 \mathfrak F_1\leq2\sqrt2(9/8)^\sigma e^{c'\lambda_n(t)|(p,\zeta)|^s}
 \sqrt{\frac{\widetilde\A_k^{(n)}(t,\eta)}{\A_k^{(n)}(t,\eta)}}
 \frac{\widetilde\A_\ell^{(n)}(t,\xi)}{\A_\ell^{(n)}(t,\xi)}
 \qquad(t\in I_{\ell,\xi}).
 \label{eq:reaction-1d-first-factor-CK}
\end{equation}
We next prove
\begin{equation}
 \mathfrak F_2\lesssim e^{C\mu|(p,\zeta)|^{1/2}}\avg{\eta-\xi}^{3}
 \left[\sqrt{\frac{\partial_tw_k(t,\eta)}{w_k(t,\eta)}
             \frac{\partial_tw_\ell(t,\xi)}{w_\ell(t,\xi)}}
       +(-\dot\lambda_n(t))\avg{k,\eta}^{s/2}\avg{\ell,\xi}^{s/2}\right].
 \label{eq:reaction-1d-second-factor}
\end{equation}
For this, we estimate $\J_k(t,\eta)/\J_\ell(t,\xi)$ in four subcases
and insert each bound into $\mathfrak F_2$.

\noindent{\bf Subcase \boldmath$1\mathrm{d}_{1}$:} $t\notin I_{\ell,\xi}$.
For $k=\ell$, $k=0$, or $t\notin I_{k,\eta}$, the parenthetical factor
in \eqref{eq:J-ratio} is one.
Otherwise, $k\ne\ell$ and $t\in I_{k,\eta}$.
By the definition \eqref{eq:critical-endpoints} and
$t>8\sqrt{\max\{|\eta|,|\xi|,1\}}$,
$|k|\leq2|\eta|/t<\sqrt{|\eta|}/4$ and
$t_{|k|,|\eta|}\geq3|\eta|/(4|k|)>2\sqrt{|\eta|}$.
Thus $t\in I_{k,\eta}=\mathbf I_{k,\eta}$ by
the definition \eqref{eq:restricted-critical-intervals}.
Since $t\notin\mathbf I_{k,\eta}\cap\mathbf I_{k,\xi}$,
we have $t\notin\mathbf I_{k,\xi}$.
In either case, the estimate \eqref{eq:J-ratio} gives
\begin{equation}
 \frac{\J_k(t,\eta)}{\J_\ell(t,\xi)}
 \lesssim e^{C\mu|(p,\zeta)|^{1/2}}.
 \label{eq:reaction-1d1-J-ratio}
\end{equation}
In the case $|\xi|\leq2|\ell|t$, we substitute
\eqref{eq:reaction-1d1-J-ratio} into the definition of $\mathfrak F_2$
and apply \eqref{eq:nonresonant-elliptic-symbol}, the bound \eqref{eq:lambda-derivative},
$a\leq2s-1$, and $\avg{k,\eta}\approx\avg{\ell,\xi}$ to obtain
\begin{equation}
\begin{aligned}
 \mathfrak F_2
 &\lesssim e^{C\mu|(p,\zeta)|^{1/2}}
 \frac{\avg{\ell,\xi}\avg{\xi/(\ell t)}}{\ell^2+|\xi-\ell t|^2}
 \lesssim e^{C\mu|(p,\zeta)|^{1/2}}
 \frac{\avg{\ell,\xi}^s}{\avg t^{2s}}\\
 &\lesssim e^{C\mu|(p,\zeta)|^{1/2}}
 (-\dot\lambda_n(t))\avg{k,\eta}^{s/2}\avg{\ell,\xi}^{s/2}.
\end{aligned}
 \label{eq:reaction-1d1-F2-radius}
\end{equation}
If $|\xi|>2|\ell|t$, then
$\avg{\ell,\xi}\avg{\xi/(\ell t)}/(\ell^2+|\xi-\ell t|^2)
\lesssim1/(|\ell|t)\lesssim|\xi|^s/\avg t^{1+a}$.
The last inequality uses $t<|\xi|/(2|\ell|)$ and $a\leq2s-1<s$:
$t^a/(|\ell||\xi|^s)\lesssim|\xi|^{a-s}/|\ell|^{1+a}\leq1$.
Combining this with \eqref{eq:reaction-1d1-J-ratio} gives the final bound
in \eqref{eq:reaction-1d1-F2-radius} also in this range.

\noindent{\bf Subcase \boldmath$1\mathrm{d}_{2}$:} $k\ne\ell$ and
$t\in I_{\ell,\xi}\setminus I_{k,\eta}$.
Equation~\eqref{eq:raw-J-reverse-ledger} gives
\begin{equation}
 \frac{\J_k(t,\eta)}{\J_\ell(t,\xi)}
 \lesssim\frac{\ell^2(1+|t-\xi/\ell|)}{|\xi|}
 e^{C\mu|(p,\zeta)|^{1/2}}.
 \label{eq:reaction-1d2-J-ratio}
\end{equation}
By the definition \eqref{eq:critical-endpoints} and
$t>8\sqrt{\max\{|\eta|,|\xi|,1\}}$,
$|\ell|\leq2|\xi|/t<\sqrt{|\xi|}/4$ and
$t_{|\ell|,|\xi|}\geq3|\xi|/(4|\ell|)>2\sqrt{|\xi|}$.
Hence, the definition \eqref{eq:restricted-critical-intervals} implies
$I_{\ell,\xi}=\mathbf I_{\ell,\xi}$, and the definition
\eqref{eq:critical-endpoints} gives $\avg{\xi/(\ell t)}\lesssim1$.
Substituting \eqref{eq:reaction-1d2-J-ratio} into the definition of
$\mathfrak F_2$ and using $\avg{\ell,\xi}\lesssim|\xi|$, we obtain
\[
 \mathfrak F_2\lesssim e^{C\mu|(p,\zeta)|^{1/2}}
 \frac{1+|t-\xi/\ell|}{1+|t-\xi/\ell|^2}
 \lesssim e^{C\mu|(p,\zeta)|^{1/2}}
 \frac{\partial_tw_\ell(t,\xi)}{w_\ell(t,\xi)},
\]
where the last step uses the estimate \eqref{eq:weight-log-derivative} on the
restricted interval $I_{\ell,\xi}$.

\noindent{\bf Subcase \boldmath$1\mathrm{d}_{3}$:} $k\ne\ell$ and
$t\in(I_{k,\eta}\cap I_{\ell,\xi})\setminus
(\mathbf I_{k,\eta}\cap\mathbf I_{k,\xi})$.
The interval calculation in Subcase $1\mathrm{d}_{1}$ gives
$t\in I_{k,\eta}=\mathbf I_{k,\eta}$.
Since $t\notin\mathbf I_{k,\eta}\cap\mathbf I_{k,\xi}$,
we have $t\notin\mathbf I_{k,\xi}$, and the estimate \eqref{eq:J-ratio} gives
$\J_k(t,\eta)/\J_\ell(t,\xi)\lesssim e^{C\mu|(p,\zeta)|^{1/2}}$.
The interval calculation in Subcase $1\mathrm{d}_{2}$ gives
$I_{\ell,\xi}=\mathbf I_{\ell,\xi}$ and $\avg{\xi/(\ell t)}\lesssim1$.
Since $|\eta|\approx|\xi|$, Lemma~\ref{lem:critical-separation} gives either
$|t-\xi/\ell|\gtrsim|\xi|/\ell^2$ or
$|\eta-\xi|\gtrsim|\xi|/|\ell|$.
In the first and second alternatives, respectively,
\begin{align*}
 \frac{\avg{\ell,\xi}}{\ell^2+|\xi-\ell t|^2}
 &\lesssim\frac{|t-\xi/\ell|}{1+|t-\xi/\ell|^2}
 \lesssim\frac1{1+|t-\xi/\ell|},\\*
 \frac{\avg{\ell,\xi}}{\ell^2+|\xi-\ell t|^2}
 &\lesssim\frac{\avg{\eta-\xi}}{|\ell|(1+|t-\xi/\ell|^2)}
 \lesssim\frac{\avg{\eta-\xi}}{1+|t-\xi/\ell|}.
\end{align*}
Combining either bound with the improved comparison in
\eqref{eq:J-ratio}, $\avg{\xi/(\ell t)}\lesssim1$, and
the estimate \eqref{eq:weight-log-derivative} gives
$\mathfrak F_2\lesssim e^{C\mu|(p,\zeta)|^{1/2}}
\avg{\eta-\xi}\,\partial_tw_\ell(t,\xi)/w_\ell(t,\xi)$.

\noindent{\bf Subcase \boldmath$1\mathrm{d}_{4}$:} $k=\ell\ne0$ and $t\in I_{k,\xi}$.
The same-mode comparison in \eqref{eq:J-ratio} gives
$\J_k(t,\eta)/\J_k(t,\xi)\lesssim e^{C\mu|\eta-\xi|^{1/2}}$.
The interval calculation in Subcase $1\mathrm{d}_{2}$, with $\ell=k$, gives
$I_{k,\xi}=\mathbf I_{k,\xi}$ and $\avg{\xi/(kt)}\lesssim1$.
The numerator in $\mathfrak F_2$ is now $|k|$.
Substitution and the estimate \eqref{eq:weight-log-derivative} give
\begin{align*}
 \mathfrak F_2&\lesssim e^{C\mu|\eta-\xi|^{1/2}}
 \frac{|k|}{k^2+|\xi-kt|^2}
 =e^{C\mu|\eta-\xi|^{1/2}}
 \frac1{|k|(1+|t-\xi/k|^2)}\lesssim e^{C\mu|\eta-\xi|^{1/2}}
 \frac{\partial_tw_k(t,\xi)}{w_k(t,\xi)}.
\end{align*}

Subcases $1\mathrm{d}_{2}$--$1\mathrm{d}_{4}$ give the common bound
$\mathfrak F_2\lesssim e^{C\mu|(p,\zeta)|^{1/2}}
\avg{\eta-\xi}\,\partial_tw_\ell(t,\xi)/w_\ell(t,\xi)$.
We now convert the right-hand side to the form required in
\eqref{eq:reaction-multiplier-off-diagonal}.
Multiply \eqref{eq:raw-w-endpoint-transfer} by
$\sqrt{\partial_tw_\ell(t,\xi)/w_\ell(t,\xi)}$.
Applying Young's inequality to the frequency term and absorbing half of
$\partial_tw_\ell(t,\xi)/w_\ell(t,\xi)$ gives
\[
 \frac{\partial_tw_\ell(t,\xi)}{w_\ell(t,\xi)}
 \lesssim\avg{\eta-\xi}
 \sqrt{\frac{\partial_tw_k(t,\eta)}{w_k(t,\eta)}
       \frac{\partial_tw_\ell(t,\xi)}{w_\ell(t,\xi)}}
 +\avg{\eta-\xi}^2\frac{|\eta|^s}{\avg t^{2s}}.
\]
Adding the bound with $(k,\eta)$ and $(\ell,\xi)$ interchanged and using
the bound \eqref{eq:lambda-derivative}, $|\eta|\approx|\xi|$, and $a\leq2s-1$, we obtain
\[
 \frac{\partial_tw_k(t,\eta)}{w_k(t,\eta)}
 +\frac{\partial_tw_\ell(t,\xi)}{w_\ell(t,\xi)}
 \lesssim\avg{\eta-\xi}^{2}
 \left[\sqrt{\frac{\partial_tw_k(t,\eta)}{w_k(t,\eta)}
             \frac{\partial_tw_\ell(t,\xi)}{w_\ell(t,\xi)}}
       +(-\dot\lambda_n)\avg{k,\eta}^{s/2}\avg{\ell,\xi}^{s/2}\right].
\]
Substitution proves the bound for $\mathfrak F_2$ in
\eqref{eq:reaction-1d-second-factor} in Subcases $1\mathrm{d}_{2}$--$1\mathrm{d}_{4}$.
This bound also holds in Subcase $1\mathrm{d}_{1}$, using only its radius term contributed by  $\dot{\lambda}_n$.

Combine the bound for $\mathfrak F_2$ in \eqref{eq:reaction-1d-second-factor} with
\eqref{eq:reaction-1d-first-factor-CK} for the $\CK$ term in Subcases
$1\mathrm{d}_{2}$--$1\mathrm{d}_{4}$, and with \eqref{eq:reaction-1d-first-factor} for the
term $(-\dot\lambda_n(t))\avg{k,\eta}^{s/2}\avg{\ell,\xi}^{s/2}$. Applying the estimate \eqref{eq:reaction-HL-exponential-absorption} with
$B=3$ and dividing by $\avg{\xi/(\ell t)}$ proves
\eqref{eq:reaction-multiplier-off-diagonal} in Case 1d.

\noindent {\bf Case 1e)} $\mathsf{HH}$.
We prove the estimate under the more general condition
$|(\ell,\xi)|/64\leq|(p,\zeta)|\leq64|(\ell,\xi)|$.
For $1/64\leq r\leq64$,
$\frac{(1+r)^s}{1+r^s}\leq\frac{65^s}{1+64^s}=:\alpha_s<1$:
the quotient is unchanged by $r\mapsto1/r$ and increases for $r\geq1$
since $s<1$. Increase $c<1$, if necessary, so that
$c>\max\{c',\alpha_s\}$.
Applying the preceding inequality with
$r=|(p,\zeta)|/|(\ell,\xi)|$ and using the triangle inequality, we obtain
\begin{equation*}
 |(k,\eta)|^s-|(\ell,\xi)|^s-c|(p,\zeta)|^s
 \leq-(1-\alpha_s)|(\ell,\xi)|^s
      -(c-\alpha_s)|(p,\zeta)|^s.
 \tag{\theequation a}\label{eq:reaction-HH-gevrey-gap}
\end{equation*}
The Sobolev ratio satisfies
$\frac{\avg{k,\eta}^{\sigma}}{\avg{\ell,\xi}^{\sigma}}
=\left(\frac{1+|(k,\eta)|^2}{1+|(\ell,\xi)|^2}\right)^{\sigma/2}
\leq65^\sigma$.
By the bound \eqref{eq:weight-total-growth} and \eqref{eq:J-mode},
$\frac{\J_k(t,\eta)}{\J_\ell(t,\xi)}
\lesssim e^{C\mu(|(\ell,\xi)|^{1/2}+|(p,\zeta)|^{1/2})}$.
Factoring $e^{c\lambda_n(t)|(p,\zeta)|^s}$ out of the Gevrey quotient
in the definition \eqref{eq:main-multiplier}, \eqref{eq:reaction-HH-gevrey-gap} gives
the decay
\[
 e^{-\lambda_n(t)[(1-\alpha_s)|(\ell,\xi)|^s
 +(c-\alpha_s)|(p,\zeta)|^s]}.
\]
Since $s>1/2$ and $\lambda_n(t)\geq3\lambda_0/8$,
\begin{align*}
 &\avg{\ell,\xi}^{4}
 e^{C\mu[|(\ell,\xi)|^{1/2}+|(p,\zeta)|^{1/2}]}\,
 e^{-\lambda_n(t)[(1-\alpha_s)|(\ell,\xi)|^s
                  +(c-\alpha_s)|(p,\zeta)|^s]}
 \lesssim1.
\end{align*}

Since $|\ell|\geq1$ and
$|\ell|t\leq|\xi-\ell t|+|\xi|$, we have
$\avg t^2\lesssim\avg{\ell,\xi}^{2}(\ell^2+|\xi-\ell t|^2)$ and
$\left\langle\frac{\xi}{\ell t}\right\rangle
\lesssim1+\frac{\avg{\ell,\xi}}{\avg t}$.
Thus, using $a<1$ and $t\geq1$ in the last step,
\[
 \frac{\avg{\ell,\xi}}{\ell^2+|\xi-\ell t|^2}
 \left\langle\frac{\xi}{\ell t}\right\rangle
 \lesssim\frac{\avg{\ell,\xi}^{3}}{\avg t^2}
 \left(1+\frac{\avg{\ell,\xi}}{\avg t}\right)
 \lesssim\avg{\ell,\xi}^{4}\avg t^{-1-a}.
\]
Combining this elliptic bound with the preceding exponential estimate
and using the bound \eqref{eq:lambda-derivative} gives
\[
 \frac{\A_k^{(n)}(t,\eta)}{\A_\ell^{(n)}(t,\xi)}
 \frac{\avg{\ell,\xi}}{\ell^2+|\xi-\ell t|^2}
 \left\langle\frac{\xi}{\ell t}\right\rangle
 \lesssim e^{c\lambda_n(t)|(p,\zeta)|^s}(-\dot\lambda_n(t)).
\]
Since $\avg{k,\eta}^{s/2}\avg{\ell,\xi}^{s/2}\geq1$, division by
$\avg{\xi/(\ell t)}$ proves
\eqref{eq:reaction-multiplier-off-diagonal} in Case 1e.
For $k=\ell$, the same estimate applies because
$|\ell|\leq\avg{\ell,\xi}$.
This completes Step 1.

\noindent{\bf Step 2: Proof of
\eqref{eq:reaction-multiplier-zero-mode} for $\ell=0$.}
The argument follows Step 1; we record the adjustments for the zero mode.

\noindent {\bf Case 2a)} $\mathsf{HL}$ and $k\ne0$.
Here $|k|+|\eta-\xi|\leq|\xi|/8$, so
$7|\xi|/8\leq|\eta|\leq9|\xi|/8$, $|k|\leq|\eta|/7$, and
$7|\xi|/8\leq|(k,\eta)|\leq9|\xi|/8$.
Since $w_0=w_{\NR}$, Lemma~\ref{lem:zero-mode-endpoint-dichotomy}
applies with its parameter $c=1/8$.
The Gevrey and Sobolev factors are controlled by
the bounds \eqref{eq:reaction-HL-Sobolev-ratio}--\eqref{eq:reaction-HL-exponential-absorption}
with $\ell=0$.
Also, the bound \eqref{eq:theta-vertical-lower-bound} gives
$\frac{\widetilde\A_k^{(n)}(t,\eta)}{\A_k^{(n)}(t,\eta)}\geq\frac12$
and
$\frac{\widetilde\A_0^{(n)}(t,\xi)}{\A_0^{(n)}(t,\xi)}\geq\frac12$.

If $t\in\mathbf I_{k,\eta}\cap\mathbf I_{k,\xi}$,
the definition \eqref{eq:critical-endpoints} gives $t\geq3|\eta|/(4|k|)$ and
$|\eta|/k^2\leq4t/(3|k|)\lesssim\avg t^{1+s}$.
Use \eqref{eq:zero-mode-J-CK-dichotomy} and insert the quotients as in
\eqref{eq:reaction-1a-quotient-insertion} to obtain
\[
 \frac{\A_k^{(n)}(t,\eta)}{\A_0^{(n)}(t,\xi)}
 \lesssim e^{c\lambda_n(t)|(k,\eta-\xi)|^s}\avg t^{1+s}
 \sqrt{\frac{\partial_tw_k(t,\eta)}{w_k(t,\eta)}
       \frac{\partial_tw_0(t,\xi)}{w_0(t,\xi)}
       \frac{\widetilde\A_k^{(n)}(t,\eta)}{\A_k^{(n)}(t,\eta)}}
 \frac{\widetilde\A_0^{(n)}(t,\xi)}{\A_0^{(n)}(t,\xi)}.
\]
This is the first term in \eqref{eq:reaction-multiplier-zero-mode}.
If $t\notin\mathbf I_{k,\eta}\cap\mathbf I_{k,\xi}$, use
\eqref{eq:zero-mode-J-improved-dichotomy} instead.
Since $\avg{k,\eta}^{s/2}/\avg\xi^{s/2}\gtrsim1$, this gives the second
term in \eqref{eq:reaction-multiplier-zero-mode} and completes Case 2a.

\noindent {\bf Case 2b)} $\mathsf{HL}$ and $k=0$.
Here $|\eta-\xi|\leq|\xi|/8$ and
$7|\xi|/8\leq|\eta|\leq9|\xi|/8$.
For $\xi\ne0$, use the estimate \eqref{eq:J-ratio}, whose parenthetical factor is
one when $k=\ell=0$, and
the bounds \eqref{eq:reaction-HL-Sobolev-ratio}--\eqref{eq:reaction-HL-exponential-absorption}.
Since $\avg\eta^{s/2}/\avg\xi^{s/2}\gtrsim1$, these give the second
term in \eqref{eq:reaction-multiplier-zero-mode}.
If $\xi=0$, the $\mathsf{HL}$ restriction forces $\eta=0$, and the
estimate reduces to $1\lesssim1$.

\noindent {\bf Case 2c)} $\mathsf{HH}$.
Use the Gevrey/Sobolev and weight-growth argument in Case 1e, with
$\ell=0$, the bound \eqref{eq:weight-total-growth}, and
\eqref{eq:reaction-HH-gevrey-gap}.
To obtain the frequency ratio in \eqref{eq:reaction-multiplier-zero-mode},
we also absorb
$\avg\xi^{s/2}/\avg{k,\eta}^{s/2}\leq\avg\xi^{s/2}$
into the decay in \eqref{eq:reaction-HH-gevrey-gap}.
This is allowed just as for the fixed polynomial factors in Case 1e,
since $s>1/2$, $\lambda_n(t)\geq3\lambda_0/8$, and $c>\alpha_s$.
Thus
\[
 \frac{\A_k^{(n)}(t,\eta)}{\A_0^{(n)}(t,\xi)}
 \frac{\avg\xi^{s/2}}{\avg{k,\eta}^{s/2}}
 \lesssim e^{c\lambda_n(t)|(k,\eta-\xi)|^s}.
\]
Rearranging gives the second term in
\eqref{eq:reaction-multiplier-zero-mode}, completing Step 2 and the proof.
\end{proof}

\begin{corollary}
\label{cor:reaction-annular-LP}
The estimates
\eqref{eq:reaction-multiplier-off-diagonal} and
\eqref{eq:reaction-multiplier-zero-mode} also hold, with a possibly larger
fixed $c<1$, on
\begin{equation}
 \mathsf{HH}_{\mathrm{LP}}
 =\left\{((k,\eta),(\ell,\xi)): |(\ell,\xi)|/64\leq|(p,\zeta)|\leq64|(\ell,\xi)|\right\},\qquad (p,\zeta)=(k-\ell,\eta-\xi).
 \label{eq:HH-LP-region}
\end{equation}
\end{corollary}

\begin{proof}
This is a straightforward extension of the proof of
Proposition~\ref{prop:reaction-remainder-multiplier}.
\end{proof}

\subsection{The transport estimate}

\begin{proposition}[Transport multiplier]
\label{prop:transport-multiplier}
There is \(c=c(s)\in(0,1)\) such that, for every \(n\geq1\) and every
fixed \(((k,\eta),(\ell,\xi))\) with
$|(\ell,\xi)|\leq3|(p,\zeta)|/16$, for almost every \(t\geq1\),
\begin{align}
  &\frac{|(p,\zeta)|}{\avg t^2}
  \left|
    \frac{{\A}_k^{(n)}(t,\eta)}
         {{\A}_p^{(n)}(t,\zeta)}
    -1
  \right|
  \notag\\*
  &\qquad\lesssim
  e^{c\lambda_n(t)|(\ell,\xi)|^s}
  \Biggl[\sqrt{\frac{\partial_tw_{k}(\eta)}{w_{k}(\eta)} \frac{\partial_tw_{p}(\zeta)}{w_{p}(\zeta)}
       \frac{\widetilde\A_k^{(n)}(t,\eta)}{\A_k^{(n)}(t,\eta)}
       \frac{\widetilde\A_p^{(n)}(t,\zeta)}{\A_p^{(n)}(t,\zeta)} }
       +(-\dot\lambda_n)\avg{k,\eta}^{s/2}\avg{p,\zeta}^{s/2}\Biggr].
  \label{eq:transport-multiplier}
\end{align}
\end{proposition}

\begin{proof}
We first prove \eqref{eq:transport-multiplier} on $\mathsf{LH}$, and then on
the remaining band $|(p,\zeta)|/8<|(\ell,\xi)|\leq3|(p,\zeta)|/16$.

\noindent{\bf Step 1: Proof of \eqref{eq:transport-multiplier} on $\mathsf{LH}$.}

\noindent {\bf Case 1a)} $|(p,\zeta)|<1$.
The $\mathsf{LH}$ restriction and
$p,\ell,k\in\mathbb Z$ imply $p=\ell=k=0$, with
$|\zeta|<1$, $|\xi|<1/8$, and $|\eta|<9/8$.
The multiplier definition \eqref{eq:main-multiplier},
the radius bound \eqref{eq:lambda-n}, and
the bound \eqref{eq:weight-total-growth} therefore give
\[
 \frac{|\zeta|}{\avg t^2}
 \left|
   \frac{\A_0^{(n)}(t,\eta)}
        {\A_0^{(n)}(t,\zeta)}-1
 \right|
 \lesssim\avg t^{-2}
 \lesssim-\dot\lambda_n(t),
\]
where the last inequality follows from
the bound \eqref{eq:lambda-derivative} and $a<1$.
This proves \eqref{eq:transport-multiplier} in this case.

\noindent {\bf Case 1b)} $|(p,\zeta)|\geq1$.
Following \cite[Section~5]{BM15}, we decompose the difference:

\begin{align}
 \frac{\A_k^{(n)}(t,\eta)}{\A_p^{(n)}(t,\zeta)}-1
 ={}&\left[e^{\lambda_n(t)(|(k,\eta)|^s-|(p,\zeta)|^s)}-1\right]\notag\\
 &+e^{\lambda_n(t)(|(k,\eta)|^s-|(p,\zeta)|^s)}
       \left[\frac{\avg{k,\eta}^{\sigma}}{\avg{p,\zeta}^{\sigma}}-1\right]\notag\\
 &+e^{\lambda_n(t)(|(k,\eta)|^s-|(p,\zeta)|^s)}
       \frac{\avg{k,\eta}^{\sigma}}{\avg{p,\zeta}^{\sigma}}
       \left[\frac{\J_k(t,\eta)}{\J_p(t,\zeta)}-1\right]
 =:T_1+T_2+T_3.
 \label{eq:transport-multiplier-difference-pieces}
\end{align}
For $|(p,\zeta)|\geq1$, $|T_1|+|T_2|$ in
\eqref{eq:transport-multiplier-difference-pieces} is at most
$Ce^{c'\lambda_n|(\ell,\xi)|^s}\allowbreak |(\ell,\xi)|\avg{p,\zeta}^{s-1}$
by the mean value theorem and the $\mathsf{LH}$ restriction, for some
$c'<c<1$. The gap $c-c'>0$ and the uniform lower bound for
$\lambda_n$ from \eqref{eq:lambda-n} absorb $|(\ell,\xi)|$ into the
exponential. Since $\avg{k,\eta}\approx\avg{p,\zeta}$ and
$\avg t^{-2}\lesssim-\dot\lambda_n$ by the bound
\eqref{eq:lambda-derivative}, multiplying by $|(p,\zeta)|/\avg t^2$ gives
\begin{equation}
 \frac{|(p,\zeta)|}{\avg t^2}(|T_1|+|T_2|)
 \lesssim e^{c\lambda_n|(\ell,\xi)|^s}
 (-\dot\lambda_n)\avg{k,\eta}^{s/2}\avg{p,\zeta}^{s/2}.
 \label{eq:transport-first-two-pieces}
\end{equation}

The following subcases estimate only $T_3$ in
\eqref{eq:transport-multiplier-difference-pieces}. We use the common bound for $T_1+T_2$ in
\eqref{eq:transport-first-two-pieces} throughout, together with the identity 
\begin{align*}
 \frac{\J_k(t,\eta)}{\J_p(t,\zeta)}-1
={\frac{\widetilde\A_p^{(n)}(t,\zeta)}{\A_p^{(n)}(t,\zeta)}}
 \left[\frac{w_p(t,\zeta)}{w_k(t,\eta)}
       e^{\mu(\sqrt{|\eta|}-\sqrt{|\zeta|})}-1\right]+\left(1-\frac{\widetilde\A_p^{(n)}(t,\zeta)}{\A_p^{(n)}(t,\zeta)}\right)
       \left[e^{\mu(\sqrt{|k|}-\sqrt{|p|})}-1\right].
\end{align*}
\noindent{\bf Subcase \boldmath$1\mathrm{b}_1$:}
$t\leq\frac12\min\left\{\sqrt{|\eta|},\sqrt{|\zeta|}\right\}$.
The estimate \cite[Lemma~3.7]{BM15} gives
$\left|\frac{\J_k(t,\eta)}{\J_p(t,\zeta)}-1\right|
\lesssim e^{C\mu|(\ell,\xi)|^{1/2}}
\frac{\avg{\ell,\xi}}{\avg{p,\zeta}^{1/2}}$.
After multiplying by $|(p,\zeta)|/\avg t^2$,
we use $s>1/2$ and $a<1$ to obtain
\[
 \frac{|(p,\zeta)|}{\avg t^2}|T_3|
 \lesssim
 e^{c\lambda_n(t)|(\ell,\xi)|^s}
 (-\dot\lambda_n(t))
 \avg{k,\eta}^{s/2}\avg{p,\zeta}^{s/2}.
\]
Together with the common bound for $T_1+T_2$ in \eqref{eq:transport-first-two-pieces}, this proves
\eqref{eq:transport-multiplier} in Subcase $1\mathrm{b}_1$.

\noindent{\bf Subcase \boldmath$1\mathrm{b}_2$:}
$t>\frac12\min\left\{\sqrt{|\eta|},\sqrt{|\zeta|}\right\}$,
and either $|k|>100|\eta|$ or $|p|>100|\zeta|$.
The $\mathsf{LH}$ restriction gives
$|k|\approx|p|\approx|(p,\zeta)|$ and
$\max\{|\eta|,|\zeta|\}\leq\tfrac14|p|$.
By the bound \eqref{eq:weight-total-growth}, $\mu=54$ from
\eqref{eq:weight-parameters}, and the mean value theorem,
\[
\begin{aligned}
 \left|\frac{\J_k(t,\eta)}{\J_p(t,\zeta)}-1\right|
 &\leq e^{60\sqrt{|\eta|}-54\sqrt{|p|}}
       +e^{60\sqrt{|\zeta|}-54\sqrt{|p|}}
       +\left|e^{54(\sqrt{|k|}-\sqrt{|p|})}-1\right|\\
 &\lesssim e^{C\mu|(\ell,\xi)|^{1/2}}
 \left(e^{-c\avg{p,\zeta}^{1/2}}
       +\frac{|(\ell,\xi)|}{\avg{p,\zeta}^{1/2}}\right)
 \lesssim e^{C\mu|(\ell,\xi)|^{1/2}}
          \frac{\avg{\ell,\xi}}{\avg{p,\zeta}^{1/2}}.
\end{aligned}
\]

Multiplying by $|(p,\zeta)|/\avg t^2$ and arguing as in Subcase $1\mathrm{b}_1$
gives the bound for $T_3$.
Combining this bound with the common bound for $T_1+T_2$ in
\eqref{eq:transport-first-two-pieces} proves
\eqref{eq:transport-multiplier} in Subcase $1\mathrm{b}_2$.

\noindent{\bf Subcase \boldmath$1\mathrm{b}_3$:}
$t>\frac12\min\left\{\sqrt{|\eta|},\sqrt{|\zeta|}\right\}$,
$|k|\leq100|\eta|$, and $|p|\leq100|\zeta|$.
Since $7|(p,\zeta)|/8\leq|(k,\eta)|\leq9|(p,\zeta)|/8$, we have
$|(p,\zeta)|/101\leq|\zeta|\leq|(p,\zeta)|$ and
$7|(p,\zeta)|/808\leq|\eta|\leq9|(p,\zeta)|/8$.
Thus $|\eta|\approx|\zeta|\approx|(p,\zeta)|$, so
$t\gtrsim\sqrt{|(p,\zeta)|}$.
In Subcases $1\mathrm{b}_{31}$ and $1\mathrm{b}_{32}$, we assume $k,p\ne0$.
Subcase $1\mathrm{b}_{33}$ treats the zero horizontal modes.

\noindent{\bf Subcase \boldmath$1\mathrm{b}_{31}$:} $k\ne p$ and
$t\in\mathbf I_{k,\eta}\cap\mathbf I_{k,\zeta}$.
By the definition \eqref{eq:critical-endpoints},
$t\geq3|\eta|/(4|k|)$. Together with $|(p,\zeta)|\approx|\eta|$, this gives
$\frac{|(p,\zeta)|}{\avg t^2}\frac{|\eta|}{k^2}\lesssim1$.
The estimate \eqref{eq:weight-log-derivative} also gives 
\[
 \frac{|\eta|}{k^2}
 \sqrt{\frac{\partial_tw_k(t,\eta)}{w_k(t,\eta)}
       \frac{\partial_tw_p(t,\zeta)}{w_p(t,\zeta)}}\gtrsim1,
\]
This lower bound controls the subtracted $1$ when applying
\eqref{eq:raw-J-two-CK-ledger} to
$|\J_k(t,\eta)/\J_p(t,\zeta)-1|$.
Moreover, the definition \eqref{eq:restricted-critical-intervals} and the $\mathsf{LH}$
restriction give $|\eta|>4k^2\geq4|k|$, $|\zeta|>4k^2\geq4|k|$, and
$7|p|\leq8|k|+|\zeta|<3|\zeta|$.
Thus $|k|<|\eta|$ and $|p|<|\zeta|$, so
the bound \eqref{eq:theta-vertical-lower-bound} gives
$\frac{\widetilde\A_k^{(n)}(t,\eta)}{\A_k^{(n)}(t,\eta)}\geq\frac12$
and
$\frac{\widetilde\A_p^{(n)}(t,\zeta)}{\A_p^{(n)}(t,\zeta)}\geq\frac12$.
These lower bounds allow the product
$\frac{\widetilde\A_k^{(n)}(t,\eta)}{\A_k^{(n)}(t,\eta)}
\frac{\widetilde\A_p^{(n)}(t,\zeta)}{\A_p^{(n)}(t,\zeta)}$
to be inserted under the square root in the estimate for
$|\J_k(t,\eta)/\J_p(t,\zeta)-1|$, up to an implicit constant.

Combining these bounds with the estimate \eqref{eq:raw-J-two-CK-ledger} gives
\[
 \frac{|(p,\zeta)|}{\avg t^2}
 \left|\frac{\J_k(t,\eta)}{\J_p(t,\zeta)}-1\right|
 \lesssim e^{C\mu|(\ell,\xi)|^{1/2}}
 \sqrt{\frac{\partial_tw_k(t,\eta)}{w_k(t,\eta)}
       \frac{\partial_tw_p(t,\zeta)}{w_p(t,\zeta)}
       \frac{\widetilde\A_k^{(n)}(t,\eta)}{\A_k^{(n)}(t,\eta)}
       \frac{\widetilde\A_p^{(n)}(t,\zeta)}{\A_p^{(n)}(t,\zeta)}}.
\]
Insert this bound into $T_3$ in
\eqref{eq:transport-multiplier-difference-pieces} and apply the estimate
\eqref{eq:reaction-HL-exponential-absorption}, with $(p,\zeta)$ replaced by
$(\ell,\xi)$, to obtain the $\CK$ term in \eqref{eq:transport-multiplier}.
Together with the common bound for $T_1+T_2$ in \eqref{eq:transport-first-two-pieces}, this proves
\eqref{eq:transport-multiplier} when $k\ne p$ and
$t\in\mathbf I_{k,\eta}\cap\mathbf I_{k,\zeta}$.

\noindent{\bf Subcase \boldmath$1\mathrm{b}_{32}$:} $k=p$ or
$t\notin\mathbf I_{k,\eta}\cap\mathbf I_{k,\zeta}$.
Here
the estimate \eqref{eq:J-ratio} gives
$|\J_k(t,\eta)/\J_p(t,\zeta)-1|
\lesssim e^{C\mu|(\ell,\xi)|^{1/2}}$, and since $a\leq2s-1$, we have 
\[
 \frac{|(p,\zeta)|}{\avg t^2}
 \lesssim(-\dot\lambda_n)\avg{k,\eta}^{s/2}\avg{p,\zeta}^{s/2},
 \qquad
 \frac{\avg{p,\zeta}^{1-s}}{\avg t^{1-a}}
 \lesssim\avg{p,\zeta}^{(1+a-2s)/2}\lesssim1,
\]
Together with the common bound for $T_1+T_2$ in \eqref{eq:transport-first-two-pieces}, this proves
\eqref{eq:transport-multiplier} when $k=p$ or
$t\notin\mathbf I_{k,\eta}\cap\mathbf I_{k,\zeta}$.

\noindent{\bf Subcase \boldmath$1\mathrm{b}_{33}$:} $k=0$ or $p=0$.
For $p=0\ne k$, use $w_0=w_{\NR}$ and apply
Lemma~\ref{lem:zero-mode-endpoint-dichotomy} with $\xi$ replaced by $\zeta$:
use \eqref{eq:zero-mode-J-CK-dichotomy} when
$t\in\mathbf I_{k,\eta}\cap\mathbf I_{k,\zeta}$,
and \eqref{eq:zero-mode-J-improved-dichotomy} otherwise.
For $k=0$, the definition \eqref{eq:restricted-critical-intervals} gives
$\mathbf I_{0,\eta}=\mathbf I_{0,\zeta}=\varnothing$,
so $t\notin\mathbf I_{k,\eta}\cap\mathbf I_{k,\zeta}$. Combining the preceding bounds with the common bound for $T_1+T_2$ in
\eqref{eq:transport-first-two-pieces} proves \eqref{eq:transport-multiplier} in Subcase $1\mathrm{b}_{33}$.
The arguments in Subcases $1\mathrm{b}_{31}$--$1\mathrm{b}_{33}$ establish
\eqref{eq:transport-multiplier} in Subcase $1\mathrm{b}_3$.
Combining Case 1a with Subcases $1\mathrm{b}_1$--$1\mathrm{b}_3$ proves
\eqref{eq:transport-multiplier} on $\mathsf{LH}$, completing Step 1.

\noindent{\bf Step 2: Proof of \eqref{eq:transport-multiplier} on $\mathsf{HH}$.}
Under the transport support restriction, it remains to consider
$|(p,\zeta)|/8<|(\ell,\xi)|\leq3|(p,\zeta)|/16$.
Here $|(p,\zeta)|<8|(\ell,\xi)|$,
$|(k,\eta)|<9|(\ell,\xi)|$, and
$13|(p,\zeta)|/16\leq|(k,\eta)|\leq19|(p,\zeta)|/16$.
Choose fixed $c',c$ such that
$\max_{1/8\leq r\leq3/16}\frac{(1+r)^s-1}{r^s}<c'<c<1$.
The Sobolev ratio is bounded.  The resonant prefactor in
\eqref{eq:J-ratio} is at most $1+|\eta|$ when $k\ne0$ and equals
one when $k=0$.  Therefore
\begin{align}
 \frac{\A_k^{(n)}(t,\eta)}{\A_p^{(n)}(t,\zeta)}
 &\lesssim\avg{\ell,\xi}
 e^{c'\lambda_n|(\ell,\xi)|^s+C\mu|(\ell,\xi)|^{1/2}},\notag\\
 \frac{|(p,\zeta)|}{\avg t^2}
 \left|\frac{\A_k^{(n)}(t,\eta)}{\A_p^{(n)}(t,\zeta)}-1\right|
 &\lesssim\frac{\avg{\ell,\xi}^{2}}{\avg t^2}
 e^{c'\lambda_n|(\ell,\xi)|^s+C\mu|(\ell,\xi)|^{1/2}} \lesssim e^{c\lambda_n|(\ell,\xi)|^s}
 (-\dot\lambda_n)\avg{k,\eta}^{s/2}\avg{p,\zeta}^{s/2}.
 \label{eq:transport-BM-support-band}
\end{align}
Since $s>1/2$ and $\lambda_n$ has a uniform positive lower bound,
the fixed gap $c-c'>0$ gives
\[
 \avg{\ell,\xi}^{2}
 e^{C\mu|(\ell,\xi)|^{1/2}-(c-c')\lambda_n|(\ell,\xi)|^s}
 \lesssim e^{-\frac12(c-c')\lambda_n|(\ell,\xi)|^s}\leq1.
\]

Multiplying this bound by
$e^{c\lambda_n(t)|(\ell,\xi)|^s}\avg t^{-2}$ and using
the bound \eqref{eq:lambda-derivative}, with $\avg t^{-2}\lesssim-\dot\lambda_n$
for $a<1$, gives the last inequality in
\eqref{eq:transport-BM-support-band}, since
$\avg{k,\eta}^{s/2}\avg{p,\zeta}^{s/2}\geq1$.
This proves \eqref{eq:transport-multiplier} on
$|(p,\zeta)|/8<|(\ell,\xi)|\leq3|(p,\zeta)|/16$.
Together with the $\mathsf{LH}$ cases, this completes the proof.

\end{proof}

\section{Energy identities}
\label{sec:energy-identities}

In this section, we derive the four energy identities corresponding to the functionals defined in \eqref{eq:component-energies}. There will be two types of lower-order contributions, for which we now specify notations: 
\begin{align}
  \mathscr R^{(n)}(t)
  &:=\frac{M^C\eps^3\delta_n^2}{\avg t^{1+\delta}},
  \label{eq:integrable-remainder}\\
  \mathscr N^{(n)}(t)
  &:=M^C\eps\delta_n^2
    \sum_{m\leq n}\delta_m^{-1}\varsigma_m\mathscr{Q}^{(m)}(t), \quad\mathscr{Q}^{(m)}(t):=
       \CK_{f,\partial g,h}^{(m)}(t)+t^{-1}\Bigl(\E_{\partial g}^{(m)}(t)+\E_{g,\mathrm{low}}^{(m)}(t)\Bigr).
  \label{eq:newton-CK-coupling}
\end{align}
Here $\CK_{f,\partial g,h}^{(m)}$ is defined in  \eqref{eq:combined-CK}. We observe that the remainder $\mathscr R^{(n)}$ is integrable in time, whereas $\mathscr N^{(n)}$ couples the
Newton levels. 

\subsection{The vorticity energy}
Recall the transport vector-fields $\mathbf V^{[n]}$ and $\mathbf P^{(n)}$ defined in 
\eqref{eq:partial-transport-field}.
Applying \({\A}^{(n)}\) to \eqref{eq:n-system-f}, taking the real
\(L^2\)-pairing with \({\A}^{(n)}f^{(n)}\), and differentiating 
the frequency multiplier gives
\begin{equation}
  \frac12\frac{\dd}{\dd t}\E_f^{(n)}+\CK_{f,+}^{(n)}
  =\mathcal I_f^{(n)}+\mathcal C_f^{(n)}+\mathcal N_f^{(n)}
   +(-\dot\lambda_n)\E_f^{(n)},
  \label{eq:f-energy-identity}
\end{equation}
where\begin{subequations}
\begin{align}
  \mathcal I_f^{(n)}
  &:=\frac12\int_{\T\times\R}
    \operatorname{div}\mathbf V^{[n]}
    \abs{{\A}^{(n)}f^{(n)}}^2\dd z\dd v,
  \label{eq:f-divergence-term}\\
  \mathcal C_f^{(n)}
  &:=-\operatorname{Re}\int_{\T\times\R}
    [\A^{(n)},\mathbf V^{[n]}\cdot\nabla]f^{(n)}
    \overline{{\A}^{(n)}f^{(n)}}\dd z\dd v,
  \label{eq:f-self-commutator}\\
  \mathcal N_f^{(n)}
  &:=-\operatorname{Re}\int_{\T\times\R}
    {\A}^{(n)}\bigl(\mathbf P^{(n)}\cdot\nabla f^{[n-1]}\bigr)
    \overline{{\A}^{(n)}f^{(n)}}\dd z\dd v.
  \label{eq:f-Newton-term}
\end{align}
\end{subequations}

\subsection{The coordinate system energy}
Recall the $\E_h$-energy \eqref{eq:energy-h}. Direct computation yields
\begin{equation}
  \frac12\frac{\dd}{\dd t}\E_h^{(n)}+\CK_h^{(n)}
  =\mathcal T_h^{(n)}+\mathcal S_h^{(n)}
    +\mathcal G_h^{(n)}+\mathcal N_h^{(n)},
  \label{eq:h-energy-identity}
\end{equation}
where \begin{subequations}\label{eq:h-error-terms}
\begin{align}
  \mathcal T_h^{(n)}
  &:=-\operatorname{Re}\left\langle
    {\A}_{\Res}^{(n)}(g^{[n]}\partial_vh^{(n)}),
    {\A}_{\Res}^{(n)}h^{(n)}\right\rangle,
  \label{eq:h-transport-term}\\
  \mathcal S_h^{(n)}
  &:=\operatorname{Re}\left\langle
    {\A}_{\Res}^{(n)}((\partial_vg^{[n]})h^{(n)}),
    {\A}_{\Res}^{(n)}h^{(n)}\right\rangle,
  \label{eq:h-stretching-term}\\
  \mathcal G_h^{(n)}
  &:=\operatorname{Re}\left\langle
    {\A}_{\Res}^{(n)}((1+h^{[n-1]})\partial_vg^{(n)}),
    {\A}_{\Res}^{(n)}h^{(n)}\right\rangle,
  \label{eq:h-g-source}\\
  \mathcal N_h^{(n)}
  &:=-\operatorname{Re}\left\langle
    {\A}_{\Res}^{(n)}(g^{(n)}\partial_vh^{[n-1]}),
    {\A}_{\Res}^{(n)}h^{(n)}\right\rangle.
  \label{eq:h-Newton-source}
\end{align}
\end{subequations}

Next we derive the energy identity for $\pa_v g^{(n)}.$ Recall 
${\A}_{\NR}^{(n)}(t,\eta)$ from \eqref{eq:main-multiplier}.
Set
\begin{equation}
  {\A}_{\partial g}^{(n)}={\A}_{\NR}^{(n)}\avg{\partial_v}^{-s}.
  \label{eq:dg-multiplier}
\end{equation}
With $\mathcal L_g^{(n)}$ defined in \eqref{eq:Lg-expanded},
differentiating \eqref{eq:n-system-g} gives the exact equation
\begin{align}
 &\partial_t\partial_vg^{(n)}+\frac2t\partial_vg^{(n)}
 +g^{[n]}\partial_v^2g^{(n)}
 +(\partial_vg^{[n]})\partial_vg^{(n)}
 \notag\\
 &\quad={}-(\partial_vg^{(n)})\partial_vg^{[n-1]}
      -g^{(n)}\partial_v^2g^{[n-1]}
      +\frac1t\partial_v\mathcal L_g^{(n)}.
 \label{eq:dg-differentiated-equation}
\end{align}
Since
$\partial_t\log\A_{\partial g}^{(n)}
=\dot\lambda_n|\eta|^s-\partial_tw_{\NR}/w_{\NR}$,
applying $\A_{\partial g}^{(n)}$ and testing with
$\avg t^{2+2s}\A_{\partial g}^{(n)}\partial_vg^{(n)}$ gives
\begin{align}
  \frac12\frac{\dd}{\dd t}\E_{\partial g}^{(n)}
  +\CK_{\partial g}^{(n)}
  +\left(\frac2t-(1+s)\frac{t}{\avg t^2}\right)\E_{\partial g}^{(n)}
  =\mathcal T_{\partial g}^{(n)}
    +\mathcal N_{\partial g}^{(n)}
    +\mathcal F_{\partial g}^{(n)},
  \label{eq:dg-energy-identity}
  \end{align} 
  \begin{subequations}
      \begin{align}
  \mathcal T_{\partial g}^{(n)}
 : ={}&-\avg t^{2+2s}\operatorname{Re}
    \left\langle {\A}_{\partial g}^{(n)}(g^{[n]}\partial_v^2g^{(n)}),{\A}_{\partial g}^{(n)}\partial_vg^{(n)}\right\rangle
    \notag\\
  &-\avg t^{2+2s}\operatorname{Re}
    \left\langle {\A}_{\partial g}^{(n)}((\partial_vg^{[n]})\partial_vg^{(n)}),{\A}_{\partial g}^{(n)}\partial_vg^{(n)}\right\rangle,
  \label{eq:dg-transport-term}\\
  \mathcal N_{\partial g}^{(n)}
  :={}&-\avg t^{2+2s}\operatorname{Re}
    \left\langle {\A}_{\partial g}^{(n)}((\partial_vg^{(n)})\partial_vg^{[n-1]}),{\A}_{\partial g}^{(n)}\partial_vg^{(n)}\right\rangle
    \notag\\
  &-\avg t^{2+2s}\operatorname{Re}
    \left\langle {\A}_{\partial g}^{(n)}(g^{(n)}\partial_v^2g^{[n-1]}),{\A}_{\partial g}^{(n)}\partial_vg^{(n)}\right\rangle,
  \label{eq:dg-Newton-term}\\
  \mathcal F_{\partial g}^{(n)}
  :={}&\frac{\avg t^{2+2s}}{t}\operatorname{Re}
    \left\langle {\A}_{\partial g}^{(n)}\partial_v\mathcal L_g^{(n)},{\A}_{\partial g}^{(n)}\partial_vg^{(n)}\right\rangle.
  \label{eq:dg-elliptic-source}
\end{align}\end{subequations}
Note that the damping \(2g^{(n)}/t\) in \eqref{eq:n-system-g} dominates the derivative
of the prefactor \(\avg t^{2+2s}\).  More precisely, 
\begin{equation}
  \frac2t-(1+s)\frac{t}{\avg t^2}
  =\frac{1-s}{t}+\frac{1+s}{t\avg t^2}
  \geq\frac{1-s}{t},\qquad t\geq1.
\end{equation} 

Finally, we consider the low-regularity $g$-energy \eqref{eq:energy-g-low}. Define the multiplier and $\CK$-term
\begin{align}  {\A}_{\mathrm{low}}^{(n)}(t,\eta)
  :=&e^{\lambda_n(t)|\eta|^s}\avg\eta^{\sigma-6},
  \label{eq:g-low-multiplier}\\
\CK_{g,\mathrm{low};\lambda}^{(n)}
  :=&\avg t^{\beta_g}\int_{\mathbb R}
    (-\dot\lambda_n)|\eta|^s
    \left|{\A}_{\mathrm{low}}^{(n)}(t,\eta)
      \widehat g^{(n)}(t,\eta)\right|^2\dd\eta.
  \label{eq:g-low-lambda-CK}
\end{align}
The choice \(\beta_g=4-a/2<4\) retains a small positive damping term after
differentiating the time weight.  More precisely,
\begin{equation}
  d_{g,\mathrm{low}}(t)
  :=\frac2t-\frac{\beta_g}{2}\frac{t}{\avg t^2}
  =\frac{4-\beta_g}{2t}+\frac{\beta_g}{2t\avg t^2}
  \geq\frac{4-\beta_g}{2t},\qquad t\geq1.
  \label{eq:g-low-damping-coefficient}
\end{equation}
A direct computation gives
\begin{align}
  \frac12\frac{\dd}{\dd t}\E_{g,\mathrm{low}}^{(n)}
  +\CK_{g,\mathrm{low};\lambda}^{(n)}
  +d_{g,\mathrm{low}}\E_{g,\mathrm{low}}^{(n)}
  ={}&\mathcal T_{g,\mathrm{low}}^{(n)}
    +\mathcal N_{g,\mathrm{low}}^{(n)}
    +\mathcal F_{g,\mathrm{low}}^{(n)},
  \label{eq:g-low-energy-identity}
\end{align} 
where
\begin{subequations}
\begin{align}\mathcal T_{g,\mathrm{low}}^{(n)}
  :={}&-\avg t^{\beta_g}\operatorname{Re}
    \left\langle {\A}_{\mathrm{low}}^{(n)}
      (g^{[n]}\partial_vg^{(n)}),
      {\A}_{\mathrm{low}}^{(n)}g^{(n)}\right\rangle,
  \label{eq:g-low-transport-term}\\
  \mathcal N_{g,\mathrm{low}}^{(n)}
  :={}&-\avg t^{\beta_g}\operatorname{Re}
    \left\langle {\A}_{\mathrm{low}}^{(n)}
      (g^{(n)}\partial_vg^{[n-1]}),
      {\A}_{\mathrm{low}}^{(n)}g^{(n)}\right\rangle,
  \label{eq:g-low-Newton-term}\\
  \mathcal F_{g,\mathrm{low}}^{(n)}
  :={}&\frac{\avg t^{\beta_g}}{t}\operatorname{Re}
    \left\langle {\A}_{\mathrm{low}}^{(n)}\mathcal L_g^{(n)},
      {\A}_{\mathrm{low}}^{(n)}g^{(n)}\right\rangle.
  \label{eq:g-low-source-term}
\end{align}
\end{subequations}  

\section{Closing the energy estimates}
\label{sec:energy}
In this section, we prove the vorticity bound \eqref{eq:f-energy-closed} and the three
coordinate estimates \eqref{eq:h-energy-bound}, \eqref{eq:dg-energy-bound},
and \eqref{eq:g-low-energy-bound}. We consider the finite regular system \eqref{eq:n-system} 
satisfying the bootstrap assumptions. The constants are uniform in the shell
indices $m,n,\mathcal{N}$; we use the lower-order terms $\mathscr R^{(n)}$ and $\mathscr N^{(n)}$
from \eqref{eq:integrable-remainder} and \eqref{eq:newton-CK-coupling} throughout.

\subsection{Vorticity error terms}
We start from the energy relation \eqref{eq:f-energy-identity} and estimate the divergence term $\mathcal{I}_f^{(n)}$,
the commutator term $\mathcal{C}_f^{(n)}$, and the Newton forcing term $\mathcal{N}_f^{(n)}$ separately.

The commutator term $\mathcal C_f^{(n)}$ \eqref{eq:f-self-commutator} is the most delicate; we apply the Bony paraproduct decomposition to obtain
\begin{subequations}
\begin{align}
  \mathcal C_f^{(n)}
  &=\sum_{N\geq8}\mathcal T_{f;N}^{(n)}
   +\sum_{N\geq8}\mathcal R_{f;N}^{(n)}
   +\mathcal H_f^{(n)},
  \label{eq:f-paraproduct-decomposition}\\
 \mathcal T_{f;N}^{(n)}
 &=-\operatorname{Re}\bigl\langle
 \A^{(n)}(\mathbf V_{<N/8}^{[n]}\cdot\nabla f_N^{(n)})
 -\mathbf V_{<N/8}^{[n]}\cdot\nabla(\A^{(n)}f_N^{(n)}),
 \A^{(n)}f^{(n)}\bigr\rangle,
 \label{eq:f-transport-paraproduct}\\
 \mathcal R_{f;N}^{(n)}
 &=-\operatorname{Re}\bigl\langle
 \A^{(n)}(\mathbf V_N^{[n]}\cdot\nabla f_{<N/8}^{(n)})
 -\mathbf V_N^{[n]}\cdot\nabla(\A^{(n)}f_{<N/8}^{(n)}),
 \A^{(n)}f^{(n)}\bigr\rangle,
 \label{eq:f-reaction-paraproduct}\\
 \mathcal H_f^{(n)}
 &=-\operatorname{Re}\sum_{\substack{N,N'\in\mathbf D\\N/8\leq N'\leq8N}}
 \bigl\langle
 \A^{(n)}(\mathbf V_N^{[n]}\cdot\nabla f_{N'}^{(n)})
 -\mathbf V_N^{[n]}\cdot\nabla(\A^{(n)}f_{N'}^{(n)}),
 \A^{(n)}f^{(n)}\bigr\rangle.
 \label{eq:f-high-high-paraproduct}
\end{align}
\end{subequations}
Here, $\mathcal{T}_f^{(n)}$ is the transport term, $\mathcal{R}_f^{(n)}$ is the reaction term, and $\mathcal{H}_f^{(n)}$ is the remainder. The transport estimate uses Proposition~\ref{prop:transport-multiplier};
the reaction estimate uses Proposition~\ref{prop:reaction-remainder-multiplier}.
These are the two applications of Section~\ref{sec:multipliers} detailed below.

We first record the bound $e^{-c|\eta|^s+C\sqrt{|\eta|}}\langle\eta\rangle^b\lesssim{}e^{-\frac c2|\eta|^s}\leq1$
for $c>0$, $s>1/2$.
Here $c$ and $b$ are fixed throughout this estimate; the implicit
constant may depend on them.
For nonnegative $U_j$ with $1\leq j\leq n\leq {\mathcal{N}}$:
\begin{equation}
 \left(\sum_{j=1}^{n} a_j\sqrt{U_j}\right)^2
 \leq\left(\sum_{j=1}^{n}\frac{\delta_ja_j^2}{\varsigma_j}\right)
       \sum_{j=1}^{n}\delta_j^{-1}\varsigma_jU_j,
 \qquad
 \sum_{j\geq1}2^{B j/s}\varsigma_j\delta_j<\infty\quad(B\geq0\text{ fixed}).
 \label{eq:energy-common-shell-sum}
\end{equation}
Here $\delta_j=e^{-c_*2^{j\alpha_*}}$ \eqref{eq:shell-envelope}, and $\varsigma_j=\exp\bigl(C_{\mathrm E}2^{j\alpha_\dagger}\bigr)
 $ \eqref{eq:elliptic-shell-weights}. 

\begin{lemma}[Vorticity error estimates]
\label{lem:f-energy-errors}
Under the bootstrap assumptions
\eqref{eq:bootstrap-energy-assumption}--\eqref{eq:bootstrap-support-assumption},
\begin{align}
  \abs{\mathcal I_f^{(n)}}
   +\sum_{N\geq8}\abs{\mathcal T_{f;N}^{(n)}}
  &\leq \frac1{32}\CK_{f,+}^{(n)}
    +C M^C\eps\left(\CK_{\partial g}^{(n)}
       +M^{-2}\CK_h^{(n)}\right)+C\mathscr R^{(n)},
  \label{eq:f-transport-bound}\\
  \sum_{N\geq8}\abs{\mathcal R_{f;N}^{(n)}}
   +\abs{\mathcal H_f^{(n)}}
  &\leq \frac1{16}\CK_{f,+}^{(n)}
    +C M^C\eps\left(\CK_{\partial g}^{(n)}
       +M^{-2}\CK_h^{(n)}\right)
    +C\mathscr R^{(n)}+C\mathscr N^{(n)},
  \label{eq:f-reaction-bound}\\
  \abs{\mathcal N_f^{(n)}}
  &\leq \frac1{32}\CK_{f,+}^{(n)}
    +C M^C\eps\left(\CK_{\partial g}^{(n)}
       +M^{-2}\CK_h^{(n)}\right)
    +C\mathscr R^{(n)}+C\mathscr N^{(n)}.
  \label{eq:f-Newton-bound}
\end{align}
Here, the error terms $\mathscr R^{(n)},\ \mathscr N^{(n)}$ are defined in \eqref{eq:integrable-remainder} and \eqref{eq:newton-CK-coupling}.
\end{lemma}

\begin{proof}

We focus on the energy equality \eqref{eq:f-energy-identity}
and divide the proof into steps.

\noindent\emph{Step 1: Estimate the divergence term $\mathcal{I}_{f}^{(n)}$  and the transport term $\mathcal T_{f}^{(n)}$.}

For the divergence term $\mathcal I_f^{(n)}$ \eqref{eq:f-divergence-term},
the identity $\operatorname{div}\mathbf V^{[n]}
  =(\partial_vh^{[n]})\partial_z\phi_{\neq}^{[n]}
    +\partial_vg^{[n]}$ yields
\begin{equation}
|\mathcal I_f^{(n)}|\lesssim
(M^C\eps\langle t\rangle^{-1-s}
+M^C\eps^2\langle t\rangle^{-2})\E_f^{(n)}
\lesssim\mathscr R^{(n)}.
\label{eq:f-divergence-error-bound}
\end{equation}

By Corollary~\ref{cor:elliptic-profile-output}, write
$\mathbf V_\phi^{[n]}=(1+h^{[n]})\nabla_{z,v}^{\perp}\Phi_{\Ext,\neq}^{[n]}$
and $\mathbf V^{[n]}=\mathbf V_\phi^{[n]}+(0,g^{[n]})$.
With $(p,\zeta)=(k-\ell,\eta-\xi)$ as in
\eqref{eq:nonlinear-frequency-convention}, we have
\begin{align}
 \mathcal T_{f;N}^{(n)}={}&-\frac1{2\pi}\operatorname{Re}
 \sum_{k,\ell}\iint_{\R^2}
 \overline{\A_k^{(n)}(\eta)\widehat f_k^{(n)}(\eta)}
 [\A_k^{(n)}(\eta)-\A_p^{(n)}(\zeta)]
 \widehat{\mathbf V}_{\phi,<N/8}^{[n]}(\ell,\xi)
 \cdot i(p,\zeta)\widehat f_N^{(n)}(p,\zeta)\,\dd\xi\dd\eta\notag\\
 &-\frac1{2\pi}\operatorname{Re}\sum_k\iint_{\R^2}
 \overline{\A_k^{(n)}(\eta)\widehat f_k^{(n)}(\eta)}
 [\A_k^{(n)}(\eta)-\A_k^{(n)}(\xi)]
 \widehat{g_{<N/8}^{[n]}}(\eta-\xi)i\xi
 \widehat f_N^{(n)}(k,\xi)\,\dd\xi\dd\eta\notag\\
 =:{}& T_{1;N}+T_{2;N}.
 \label{eq:f-transport-fourier-decomposition}
\end{align}
For the $T_{1;N}$ term in \eqref{eq:f-transport-fourier-decomposition},
we apply the multiplier estimate \eqref{eq:transport-multiplier} in
Proposition~\ref{prop:transport-multiplier} to obtain
\begin{align}
 |T_{1;N}|\lesssim{}&
 \sum_{k,\ell}\iint_{\R^2}
 |\A_k^{(n)}(\eta)\widehat f_k^{(n)}(\eta)|
 \langle t\rangle^2e^{c\lambda_n|(\ell,\xi)|^s}
 |\widehat{\mathbf V}_{\phi,<N/8}^{[n]}(\ell,\xi)|
 |\A_p^{(n)}(\zeta)\widehat f_N^{(n)}(p,\zeta)|\notag\\*
 &\quad\times\left[
 \sqrt{\frac{\partial_tw_k(\eta)}{w_k(\eta)}
              \frac{\widetilde\A_k^{(n)}(\eta)}{\A_k^{(n)}(\eta)}}
 \sqrt{\frac{\partial_tw_p(\zeta)}{w_p(\zeta)}
              \frac{\widetilde\A_p^{(n)}(\zeta)}{\A_p^{(n)}(\zeta)}}
 +(-\dot\lambda_n)\langle k,\eta\rangle^{s/2}
                       \langle p,\zeta\rangle^{s/2}\right]\dd\xi\dd\eta\notag\\
  =:{}&T_{11;N}+T_{12;N}.
 \label{eq:f-transport-elliptic-split}
\end{align}
For the $T_{11;N}$ term, we recall
the $\CK_{f;w}^{(n)}$-term in \eqref{eq:CK-f-w} and Corollary~\ref{cor:elliptic-profile-output} to obtain that \begin{align}
 \sum_{N\geq8}T_{11;N}
 &\lesssim\langle t\rangle^2
 \left\|e^{c\lambda_n|(\ell,\xi)|^s}
       \widehat{\mathbf V}_{\phi}^{[n]}(\ell,\xi)
       \right\|_{\ell^1_\ell L^1_\xi}
\left\|\sqrt{\frac{\partial_tw_k(\eta)}{w_k(\eta)}}
\sqrt{\A_k^{(n)}(\eta)\widetilde{\A}_k^{(n)}(\eta)}
 \widehat f_k^{(n)}(\eta)\right\|_{\ell^2_kL^2_\eta}^{2}
 \notag\\*
 &\approx
 \langle t\rangle^2
 \left\|e^{c\lambda_n|(\ell,\xi)|^s}
       \widehat{\mathbf V}_{\phi}^{[n]}(\ell,\xi)
       \right\|_{\ell^1_\ell L^1_\xi}\CK_{f;w}^{(n)}
 \lesssim
 M^C\eps\CK_{f;w}^{(n)}
 \lesssim M^C\eps\CK_{f,+}^{(n)}.
 \label{eq:f-transport-T11-bound}
\end{align}
For the $T_{12;N}$ term in \eqref{eq:f-transport-elliptic-split}, 
applying the augmented coercivity estimate \eqref{eq:augmented-radius-coercivity} gives
\begin{align}
 \sum_{N\geq8}T_{12;N}
 &\lesssim(-\dot\lambda_n)\langle t\rangle^2
 \left\|e^{c\lambda_n|(\ell,\xi)|^s}
       \widehat{\mathbf V}_{\phi}^{[n]}(\ell,\xi)
       \right\|_{\ell^1_\ell L^1_\xi}\left\|\langle k,\eta\rangle^{s/2}
       \A_k^{(n)}(\eta)\widehat f_k^{(n)}(\eta)
       \right\|_{\ell^2_kL^2_\eta}^{2}
 \notag\\*
 &\lesssim
 \langle t\rangle^2
 \left\|e^{c\lambda_n|(\ell,\xi)|^s}
       \widehat{\mathbf V}_{\phi}^{[n]}(\ell,\xi)
       \right\|_{\ell^1_\ell L^1_\xi}\CK_{f,+}^{(n)}
 \lesssim M^C\eps\CK_{f,+}^{(n)}.
 \label{eq:f-transport-T12-bound}
\end{align}

For the $T_2$ term in \eqref{eq:f-transport-fourier-decomposition}, fix $c<1$ and a nonnegative integer
$p_1$, we have the bounds\begin{equation}
 \begin{aligned}
 \left\|e^{c\lambda_n|\xi|^s}\langle\xi\rangle^{p_1}|\xi|
 \widehat{g^{[n]}}(\xi)\right\|_{L^1_\xi}\lesssim\sum_{j\leq n}\|{\A}_{\partial g}^{(j)}\partial_vg^{(j)}\|_2
 &\lesssim M^C\eps\langle t\rangle^{-1-s},\\
 \left\|e^{c\lambda_n|\xi|^s}\langle\xi\rangle^{p_1}
 \widehat{g^{[n]}}(\xi)\right\|_{L^1_\xi}\lesssim\sum_{j\leq n}\|g^{(j)}\|_{\Gs{s}{\lambda_j}{\sigma-6}}
 &\lesssim M^C\eps\langle t\rangle^{-2+a/4}.
 \end{aligned}
 \label{eq:energy-coordinate-G01}
\end{equation}
Now following the standard argument as in \cite{BM13}, one has \begin{align}
 \sum_{N\geq8}|T_{2;N}|
\lesssim M^C\eps
 \left(\langle t\rangle^{a-s}
       +\langle t\rangle^{1-2s+5a/4}\right)\CK_{f,+}^{(n)}
 +\mathscr R^{(n)}.
 \label{eq:f-transport-coordinate-total}
\end{align}
Both time powers in this bound are nonpositive by
\eqref{eq:weight-parameters}. Combining the divergence bound
\eqref{eq:f-divergence-error-bound}, the elliptic bounds
\eqref{eq:f-transport-T11-bound}--\eqref{eq:f-transport-T12-bound},
and the coordinate bound \eqref{eq:f-transport-coordinate-total} yields 
 \begin{align}
|\mathcal I_f^{(n)}|+\sum_{N\geq8}|\mathcal T_{f;N}^{(n)}|
 \leq\frac1{32}\CK_{f,+}^{(n)}+C\mathscr R^{(n)},
 \end{align}
which is \eqref{eq:f-transport-bound}. 

\noindent\emph{Step 2: Estimate the reaction $\mathcal{ R }_f$ and the  remainder $\mathcal{H}_f$ .}

For $\mathcal R_{f;N}^{(n)}$ in
\eqref{eq:f-reaction-paraproduct}, replace
the  $h^{(j)}\nabla^\perp\phi^{[n]}$ by $h^{(j)}\nabla^\perp\Phi_{\mathrm{mix}}^{[n;j]}$, where
\[
 \Phi_{\mathrm{mix}}^{[n;j]}:=\sum_{b<j}\Phi_{\Ext}^{(b)}+
 \sum_{j\leq b\leq n}\Phi_{\Int}^{(b;j)},\qquad
 \Phi_{\Int}^{(j;j)}=\Phi_{\Ext}^{(j)}.
\] 
We apply the lossy and
precision elliptic estimates to the mixed stream function
$\Phi_{\mathrm{mix}}^{[n;j]}$.

For the reaction term \eqref{eq:f-reaction-paraproduct}, we decompose $\mathcal R_{f;N}^{(n)}$ 
with $(p,\zeta)=(k-\ell,\eta-\xi)$,
\begin{equation}
 \mathcal R_{f;N}^{(n)}=T_{3;N}+T_{4;N}+T_{5;N},
 \label{eq:f-reaction-local-decomposition}
\end{equation} 
where
\begin{align}
 T_{3;N}:={}&-\frac1{2\pi}\operatorname{Re}
 \sum_{\substack{k\in\mathbb Z\\\ell\ne0}}\iint_{\mathbb R^2}
 \overline{\A_k^{(n)}(\eta)\widehat f_k^{(n)}(\eta)}
 \notag\\
 &\quad\times[\A_k^{(n)}(\eta)-\A_p^{(n)}(\zeta)]\pi_N(\ell,\xi)\widehat f_{<N/8}^{(n)}(p,\zeta)
 (p\xi-\ell\zeta)\widehat\Phi_{\Ext}^{[n]}(\ell,\xi)
 \,\dd\xi\dd\eta,\notag\\
 T_{4;N}:={}&-\frac1{2\pi}\operatorname{Re}
 \sum_k\iint_{\mathbb R^2}
 \overline{\A_k^{(n)}(\eta)\widehat f_k^{(n)}(\eta)}
 \notag\\
 &\quad\times[\A_k^{(n)}(\eta)-\A_k^{(n)}(\eta-\xi)]
 i(\eta-\xi)\widehat g_N^{[n]}(\xi)
 \widehat f_{<N/8}^{(n)}(k,\eta-\xi)\,\dd\xi\dd\eta,
 \label{eq:f-reaction-stream-g-terms} \\ \nonumber
 T_{5;N}
 :={}&-\frac1{(2\pi)^2}\operatorname{Re}
 \sum_{j\leq n}\sum_{\substack{k\in\mathbb Z\\\ell\ne0}}
 \iiint_{\mathbb R^3}
 \overline{\A_k^{(n)}(\eta)\widehat f_k^{(n)}(\eta)}
 [\A_k^{(n)}(\eta)-\A_p^{(n)}(\zeta)]\\
 &\quad\times\pi_N(\ell,\xi)\widehat f_{<N/8}^{(n)}(p,\zeta)
 (p\xi'-\ell\zeta)\widehat h^{(j)}(\xi-\xi')
 \widehat\Phi_{\mathrm{mix}}^{[n;j]}(\ell,\xi')
 \dd\xi'\dd\xi\dd\eta. \label{eq:f-reaction-coefficient-terms}
\end{align}

We further decompose $T_{3;N}=T_{31;N}+T_{32;N}$, where
\begin{align}
 \begin{pmatrix}T_{31;N}\\T_{32;N}\end{pmatrix}
 :={}&-\frac1{2\pi}\operatorname{Re}
 \sum_{\substack{k\in\mathbb Z\\\ell\ne0}}\iint_{\mathbb R^2}
 \overline{\A_k^{(n)}(\eta)\widehat f_k^{(n)}(\eta)}
 \notag\\
 &\quad\times\begin{pmatrix}\A_k^{(n)}(\eta)\\-\A_p^{(n)}(\zeta)\end{pmatrix}\widehat f_{<N/8}^{(n)}(p,\zeta)
 (p\xi-\ell\zeta)\pi_N(\ell,\xi)\widehat\Phi_{\Ext}^{[n]}(\ell,\xi)
 \,\dd\xi\dd\eta.\end{align}

We start by estimating $T_{32;N}$. The low-regularity coordinate bounds 
\eqref{eq:bootstrap-consequences-coefficients},
Proposition~\ref{prop:lossy-elliptic}, Corollary~\ref{cor:elliptic-profile-output},
and \eqref{eq:energy-coordinate-G01} give 
 \begin{equation}
\begin{aligned}
 \sum_{N\geq8}|T_{32;N}|
 &\lesssim\sum_{N\geq 8}\sum_{k,\ell}\iint_{\mathbb R^2}
 |\A_k^{(n)}(\eta)\widehat f_k^{(n)}(\eta)|
 |\A_p^{(n)}(\zeta)\widehat f_{<N/8}^{(n)}(p,\zeta)|
 \langle\ell,\xi\rangle
 |(i\ell,i\xi)|\,|\pi_N\widehat\Phi_{\Ext}^{[n]}(\ell,\xi)|
 \dd\xi\dd\eta\\
 &\lesssim\E_f^{(n)}
 \|\langle\ell,\xi\rangle|(i\ell,i\xi)|
           \widehat\Phi_{\Ext}^{[n]}(\ell,\xi)\|_{\ell^1_{\ell\ne0}L^1_\xi}
 \lesssim M^C\eps
 (\langle t\rangle^{-2}+\langle t\rangle^{-2+a/4})\E_f^{(n)}
 \lesssim\mathscr R^{(n)} .
\end{aligned}
 \label{eq:f-reaction-input-bound}
\end{equation}
Next we estimate the $T_{31;N}$ term. By Proposition~\ref{prop:reaction-remainder-multiplier}, 
\begin{equation}
\begin{aligned}
 \sum_{N\geq8}|T_{31;N}|\lesssim {}&\sum_{N\geq8}\sum_{\substack{k\in\mathbb Z\\\ell\ne0}}\iint_{\mathbb R^2}
 |(p,\zeta)|e^{c\lambda_n|(p,\zeta)|^s}|\widehat f_{<N/8}^{(n)}(p,\zeta)|
 |\A_k^{(n)}(\eta)\widehat f_k^{(n)}(\eta)|
 |\A_\ell^{(n)}(\xi)\pi_N\widehat{\Delta_L\Phi_{\Ext}^{[n]}}(\ell,\xi)|\\
 &\quad\times\left\langle\frac{\xi}{\ell t}\right\rangle^{-1}
\left(
 \sqrt{\frac{\partial_tw_k(\eta)}{w_k(\eta)}
       \frac{\partial_tw_\ell(\xi)}{w_\ell(\xi)}
       \frac{\widetilde\A_k^{(n)}(\eta)}{\A_k^{(n)}(\eta)}}
       \frac{\widetilde\A_\ell^{(n)}(\xi)}{\A_\ell^{(n)}(\xi)}
 -\dot\lambda_n\langle k,\eta\rangle^{s/2}
                      \langle\ell,\xi\rangle^{s/2}
 \right)\,\dd\xi\dd\eta .
\end{aligned}
 \label{eq:f-reaction-stream-majorants}
\end{equation}
We use
\eqref{eq:augmented-radius-coercivity},
$\langle\ell,\xi\rangle^{s/2}\lesssim|(\ell,\xi)|^{s/2}$ for $\ell\ne0$, and
Young's convolution inequality with the differentiated vorticity in Fourier
$L^1$ to get

{\begin{equation}
\begin{aligned}
 \sum_{N\geq8}|T_{31;N}|
 &\lesssim M\eps\delta_n\sqrt{\CK_{f;w}^{(n)}}
 \|\mathcal M_w^{(n)}\Delta_L\Phi_{\Ext}^{[n]}\|_2+M\eps\delta_n\sqrt{\CK_{f,+}^{(n)}}
 \|\mathcal M_\lambda^{(n)}\Delta_L\Phi_{\Ext}^{[n]}\|_2.
\end{aligned}
 \label{eq:f-reaction-stream-branch-bounds}
\end{equation}}

By \eqref{eq:precision-elliptic-partial-sum},
\eqref{eq:f-reaction-input-bound}, and
\eqref{eq:f-reaction-stream-branch-bounds}, we obtain
\begin{equation}
\begin{aligned}
 \sum_{N\geq8}|T_{3;N}|
 &\lesssim M\eps\delta_n\sqrt{\CK_{f,+}^{(n)}}
 \left[
 \left(\sum_{b\leq n}\delta_b^{-1}\varsigma_b\mathscr Q^{(b)}\right)^{1/2}
 +M^C\eps\langle t\rangle^{-(1+\delta)/2}\right]+\mathscr R^{(n)}.
\end{aligned}
 \label{eq:f-reaction-T3-bound}
\end{equation}

An adaptation of the arguments in \cite{BM13} gives
\begin{align}
     \sum_{N\geq8}\sum_{j=4}^5|T_{j;N}|
 &\lesssim (M\eps+M^3\eps^2)\delta_n\sqrt{\CK_{f,+}^{(n)}}
 \left(\sum_{b\leq n}\delta_b^{-1}\varsigma_b\mathscr Q^{(b)}\right)^{1/2}
 +M^C\eps^3\delta_n\sqrt{\CK_{f,+}^{(n)}}\langle t\rangle^{-(1+\delta)/2}+\mathscr R^{(n)} . 
 \label{eq:f-reaction-T45-bound}
\end{align}

For $\mathcal{H}_f^{(n)}$, we obtain
\begin{equation}
|\mathcal{H}_f^{(n)}|
 \lesssim\mathscr R^{(n)} .
 \label{eq:f-reaction-T7-bound}
\end{equation} 
Combining \eqref{eq:f-reaction-T3-bound} and \eqref{eq:f-reaction-T45-bound} yields
\begin{equation}
 \sum_{N\geq8}|\mathcal R_{f;N}^{(n)}|+|\mathcal H_f^{(n)}|
 \leq\frac1{16}\CK_{f,+}^{(n)}+C\mathscr N^{(n)}+C\mathscr R^{(n)} .
 \label{eq:f-reaction-final-bound}
\end{equation}
This proves
\eqref{eq:f-reaction-bound}.

\noindent\emph{Step 3:  Estimate the Newton forcing $\mathcal N_f^{(n)}$.} 
For the Newton forcing $\mathcal N_f^{(n)}$ in \eqref{eq:f-Newton-term}, recall the velocity:
\begin{equation}
 \mathbf P^{(n;m)}:=h^{(n)}\nabla^\perp\Phi_{\Ext,\neq}^{[n-1]}
 +(1+h^{[n]})\nabla^\perp\Phi_{\Int,\neq}^{(n;m)}+(0,g^{(n)}).
 \label{eq:f-newton-local-velocity}
\end{equation}
The same dyadic partition as in \eqref{eq:f-paraproduct-decomposition}
now gives an exact decomposition: 
\begin{equation}
\begin{aligned}
 T_{6}:={}&-\operatorname{Re}\sum_{m<n}\sum_{N\geq8}
 \bigl\langle\A^{(n)}(\mathbf P_{<N/8}^{(n;m)}\cdot\nabla f_N^{(m)}),
 \A^{(n)}f^{(n)}\bigr\rangle,\\
 T_{7}:={}&-\operatorname{Re}\sum_{m<n}\sum_{N\geq8}
 \bigl\langle\A^{(n)}(\mathbf P_N^{(n;m)}\cdot\nabla f_{<N/8}^{(m)}),
 \A^{(n)}f^{(n)}\bigr\rangle,\\
 T_{8}:={}&-\operatorname{Re}\sum_{m<n}
 \sum_{\substack{N,N'\in\mathbf D\\N/8\leq N'\leq8N}}
 \bigl\langle\A^{(n)}(\mathbf P_N^{(n;m)}\cdot\nabla f_{N'}^{(m)}),
 \A^{(n)}f^{(n)}\bigr\rangle,\\
 \mathcal N_f^{(n)}={}&T_{6}+T_{7}+T_{8}.
\end{aligned}
 \label{eq:f-newton-local-decomposition}
\end{equation} 

 \noindent\emph{The transport term $T_{6}$ and remainder $T_{8}$.} 

The $T_{6}$ term in
\eqref{eq:f-newton-local-decomposition} has the form\begin{equation} 
 T_{6}={} -\frac1{2\pi}\operatorname{Re}
 \sum_{m<n}\sum_{N\geq8}\sum_{k,\ell}
 \iint_{\mathbb R^2}
 \overline{\A_k^{(n)}(\eta)\widehat f_k^{(n)}(\eta)}
 \A_k^{(n)}(\eta)
 \widehat{\mathbf P_{<N/8}^{(n;m)}}(\ell,\xi)
 \cdot i(p,\zeta)\widehat f_N^{(m)}(p,\zeta)
 \,\dd\xi\dd\eta . 
 \label{eq:f-newton-low-velocity-fourier}
\end{equation} 
On $|(\ell,\xi)|\leq3|(p,\zeta)|/16$, the Gevrey difference has the
strict lower-radius constant $s(16/3)^{s-1}<1$.
The general raw comparison \eqref{eq:J-ratio} contributes at most
$C\langle p,\zeta\rangle e^{C\sqrt{|(\ell,\xi)|}}$.
Thus \eqref{eq:radius-gap} gives, for a fixed $s(16/3)^{s-1}<c<1$,   
\begin{equation} 
 |(p,\zeta)|\A_k^{(n)}(\eta)
  \lesssim\langle p,\zeta\rangle^2
 e^{-(\lambda_m-\lambda_n)|(p,\zeta)|^s}
 e^{c\lambda_n|(\ell,\xi)|^s}\A_p^{(m)}(\zeta) 
 \lesssim 2^{2m/s}
 e^{c\lambda_n|(\ell,\xi)|^s}\A_p^{(m)}(\zeta). 
 \label{eq:f-newton-low-velocity-comparison}
\end{equation}  
Insert it into
\eqref{eq:f-newton-low-velocity-fourier}, sum the nonnegative
projections, and apply Young's inequality:  \begin{align*}
 |T_{6}|
 &\lesssim\sum_{m<n}2^{2m/s}
 \sum_{k,\ell}\iint_{\mathbb R^2}
 |\A_k^{(n)}(\eta)\widehat f_k^{(n)}(\eta)|
 e^{c\lambda_n|(\ell,\xi)|^s}
 |\widehat{\mathbf P^{(n;m)}}_\ell(\xi)|
 |\A_p^{(m)}(\zeta)\widehat f_p^{(m)}(\zeta)|
 \,\dd\xi\dd\eta\\
 &\lesssim\sqrt{\E_f^{(n)}}\sum_{m<n}2^{2m/s}
 \|e^{c\lambda_n|(\ell,\xi)|^s}
        \widehat{\mathbf P^{(n;m)}}(\ell,\xi)\|_{\ell^1_\ell L^1_\xi}
 \sqrt{\E_f^{(m)}} .
\end{align*} The bounds from Proposition~\ref{prop:lossy-elliptic},
Corollary~\ref{cor:elliptic-profile-output}, \eqref{eq:energy-common-shell-sum}, and
\eqref{eq:energy-coordinate-G01} yield 
\begin{equation}
 |T_{6}|\lesssim M^C\eps^3\delta_n^2
 \sum_{m<n}2^{2m/s}\delta_m
 \left(\varsigma_m\langle t\rangle^{-2}
                     +\langle t\rangle^{-2+a/4}\right)
 \lesssim\mathscr R^{(n)}.
 \label{eq:f-newton-T8-bound}
\end{equation}

For the $T_{8}$ term in \eqref{eq:f-newton-local-decomposition}, straightforward estimates  yield \begin{equation}
 |T_{8}|
 \lesssim\mathscr R^{(n)} .
 \label{eq:f-newton-T10-bound}
\end{equation}

\noindent\emph{The reaction term $T_{7}$.} For the $T_7$ term in \eqref{eq:f-newton-local-decomposition}, insert
\eqref{eq:f-newton-local-velocity} and split each coefficient product
at $|\xi-\xi'|=|(\ell,\xi')|/16$, as displayed in
\eqref{eq:f-newton-high-velocity-children}. With
$(p,\zeta)=(k-\ell,\eta-\xi)$, the resulting decomposition is
\[ T_7=T_{71}+T_{72}+T_{73}+T_{74},\qquad
 T_{73}:=T_{731}+T_{732},\qquad T_{74}:=T_{741}+T_{742},
\]
where
\begin{equation}
\begin{aligned}
 T_{71}:={}&-\frac1{2\pi}\operatorname{Re}
 \sum_{m<n}\sum_{N\geq8}\sum_{\substack{k\in\mathbb Z\\\ell\ne0}}
 \iint_{\mathbb R^2}
 \overline{\A_k^{(n)}(\eta)\widehat f_k^{(n)}(\eta)}\\
 &\quad\times
 \A_k^{(n)}(\eta)\pi_N(\ell,\xi)
 \widehat f_{<N/8}^{(m)}(p,\zeta)
 (p\xi-\ell\zeta)\widehat\Phi_{\Int}^{(n;m)}(\ell,\xi)
 \,\dd\xi\dd\eta\\
 T_{72}:={}&-\frac1{2\pi}\operatorname{Re}
 \sum_{m<n}\sum_{N\geq8}\sum_k\iint_{\mathbb R^2}
 \overline{\A_k^{(n)}(\eta)\widehat f_k^{(n)}(\eta)}\\
 &\quad\times
 \A_k^{(n)}(\eta)
 i(\eta-\xi)
 \widehat g_N^{(n)}(\xi)\widehat f_{<N/8}^{(m)}(k,\eta-\xi)
 \,\dd\xi\dd\eta,
\end{aligned}
 \label{eq:f-newton-stream-g-children}
\end{equation}
{\begin{equation}
\begin{aligned}
 &\begin{pmatrix}T_{731}+T_{732}\\T_{741}+T_{742}\end{pmatrix}
 :={}-\frac1{(2\pi)^2}\operatorname{Re}
 \sum_{m<n}\sum_{N\geq8}\sum_{\substack{k\in\mathbb Z\\\ell\ne0}}
 \iiint_{\mathbb R^3}
 \overline{\A_k^{(n)}(\eta)\widehat f_k^{(n)}(\eta)}\\
 &\qquad\times
 \A_k^{(n)}(\eta)\pi_N(\ell,\xi)
 \widehat f_{<N/8}^{(m)}(p,\zeta)(p\xi'-\ell\zeta)\\
 &\qquad\times\begin{pmatrix}
\begin{aligned}
 &\mathbf1_{|\xi-\xi'|\leq|(\ell,\xi')|/16}\widehat h^{[n]}(\xi-\xi')\widehat\Phi_{\Int}^{(n;m)}(\ell,\xi')\\
 &+ \mathbf1_{|\xi-\xi'|\leq|(\ell,\xi')|/16}\widehat h^{(n)}(\xi-\xi')\widehat\Phi_{\Ext}^{[n-1]}(\ell,\xi')
\end{aligned}\\[1ex]
\begin{aligned}
 &\mathbf1_{|\xi-\xi'|>|(\ell,\xi')|/16}\widehat h^{[n]}(\xi-\xi')\widehat\Phi_{\Int}^{(n;m)}(\ell,\xi')\\
 &+ \mathbf1_{|\xi-\xi'|>|(\ell,\xi')|/16}\widehat h^{(n)}(\xi-\xi')\widehat\Phi_{\Ext}^{[n-1]}(\ell,\xi')
\end{aligned}
\end{pmatrix}\,\dd\xi'\dd\xi\dd\eta .
\end{aligned}
 \label{eq:f-newton-high-velocity-children}
\end{equation}}

For the $T_{71}$ term,
and the coefficient-low terms $T_{731},T_{732}$ in
\eqref{eq:f-newton-stream-g-children}--\eqref{eq:f-newton-high-velocity-children},
we apply Proposition~\ref{prop:reaction-remainder-multiplier} directly.
The outer projections give $|(p,\zeta)|\leq3|(\ell,\xi)|/16$.
For the coefficient-low terms,
\[
 |\xi-\xi'|\leq\frac{|(\ell,\xi')|}{16},\qquad
 |(p,\zeta)|\leq\frac{51}{256}|(\ell,\xi')|,\qquad
 |(k-\ell,\eta-\xi')|\leq\frac{67}{256}|(\ell,\xi')|.
\]
Thus the output and stream frequencies belong to
$\mathsf{HL}\cup\mathsf{HH}$ of \eqref{eq:HL-region}--\eqref{eq:HH-region}.
For $k\ne\ell$, use
$|p\xi'-\ell\zeta|\leq\langle\ell,\xi'\rangle|(p,\zeta)|$;
for $k=\ell\ne0$, the derivative is exactly $-\ell\zeta$.
The off-diagonal and diagonal comparisons
\eqref{eq:reaction-multiplier-off-diagonal},
with difference $(p,\zeta+\xi-\xi')$, give
\begin{equation}
\begin{aligned}
 &\frac{\A_k^{(n)}(\eta)}{\A_\ell^{(n)}(\xi')}
 \frac{|p\xi'-\ell\zeta|}{\ell^2+|\xi'-\ell t|^2}
 \lesssim |(p,\zeta)|e^{c\lambda_n(|(p,\zeta)|^s+|\xi-\xi'|^s)}
 \left\langle\frac{\xi'}{\ell t}\right\rangle^{-1}\\
 &\quad\times\Biggl[
 \sqrt{\frac{\partial_tw_k(\eta)}{w_k(\eta)}
       \frac{\partial_tw_\ell(\xi')}{w_\ell(\xi')}
       \frac{\widetilde\A_k^{(n)}(\eta)}{\A_k^{(n)}(\eta)}}
       \frac{\widetilde\A_\ell^{(n)}(\xi')}{\A_\ell^{(n)}(\xi')}
 +(-\dot\lambda_n)\langle k,\eta\rangle^{s/2}
                       \langle\ell,\xi'\rangle^{s/2}\Biggr].
\end{aligned}
\label{eq:f-newton-reaction-comparison}
\end{equation}
For $T_{71}$, take $\xi'=\xi$. Insert
\eqref{eq:f-newton-reaction-comparison} into
\eqref{eq:f-newton-stream-g-children}, using the precision symbols
\eqref{eq:elliptic-M-lambda}--\eqref{eq:elliptic-M-w}.
After summing the nonnegative projection products, whose sum is at most one,
\begin{equation}
\begin{aligned}
 |T_{71}|\lesssim{}&\sum_{m<n}\sum_{\substack{k\in\mathbb Z\\\ell\ne0}}
 \iint_{\substack{(\eta,\xi)\in\mathbb R^2\\|(p,\zeta)|\leq3|(\ell,\xi)|/16}}
 |(p,\zeta)|e^{c\lambda_n|(p,\zeta)|^s}|\widehat f_p^{(m)}(\zeta)|\\
 &\quad\times\Biggl[
 \sqrt{\frac{\partial_tw_k(\eta)}{w_k(\eta)}}
 \left|\sqrt{\A_k^{(n)}(\eta)\widetilde\A_k^{(n)}(\eta)}\widehat f_k^{(n)}(\eta)\right|
 \left|\mathcal M_{w,\ell}^{(n)}(\xi)
 \widehat{\Delta_L\Phi_{\Int}^{(n;m)}}(\ell,\xi)\right|\\
 &\qquad+
 \sqrt{-\dot\lambda_n}\langle k,\eta\rangle^{s/2}
 |\A_k^{(n)}(\eta)\widehat f_k^{(n)}(\eta)|
 \left|\mathcal M_{\lambda,\ell}^{(n)}(\xi)
 \widehat{\Delta_L\Phi_{\Int}^{(n;m)}}(\ell,\xi)\right|
 \Biggr]\,\dd\xi\dd\eta.
\end{aligned}
\label{eq:f-newton-stream-majorant}
\end{equation}
Here $\ell\ne0$ gives
$\langle\ell,\xi\rangle^{s/2}\lesssim|(\ell,\xi)|^{s/2}$,
and \eqref{eq:augmented-radius-coercivity} also covers $k=0$.
Young's inequality puts the output and precision-weighted stream in
$\ell^2L^2$, and the differentiated earlier vorticity in $\ell^1L^1$.
The bootstrap and $\lambda_n\leq\lambda_m$ give

\begin{align*}
 |T_{71}|&\lesssim\sqrt{\CK_{f,+}^{(n)}}\sum_{m<n}
 \|| (p,\zeta)|e^{c\lambda_n|(p,\zeta)|^s}
       \widehat f_p^{(m)}(\zeta)\|_{\ell^1_pL^1_\zeta}
 \left(\sum_{\iota=\lambda,w}
 \|\mathcal M_\iota^{(n)}\Delta_L\Phi_{\Int}^{(n;m)}\|_2^2
 \right)^{1/2}\\
 &\lesssim M\eps\sqrt{\CK_{f,+}^{(n)}}\sum_{m<n}\delta_m
 \left(\sum_{\iota=\lambda,w}
 \|\mathcal M_\iota^{(n)}\Delta_L\Phi_{\Int}^{(n;m)}\|_2^2
 \right)^{1/2}.
\end{align*}
The interior precision elliptic estimate \eqref{eq:interior-precision-final},
with \eqref{eq:interior-precision-source}, gives
\begin{align*}
 &\left(\sum_{\iota=\lambda,w}
 \|\mathcal M_\iota^{(n)}\Delta_L\Phi_{\Int}^{(n;m)}\|_2^2
 \right)^{1/2}\lesssim\sqrt{\varsigma_m}\Biggl[
 \sqrt{\CK_{f,\partial g,h}^{(n)}}
 +\delta_n
 \left(\sum_{b<n}\delta_b^{-1}\varsigma_b\mathscr Q^{(b)}\right)^{1/2}
 +M^C\eps\delta_n\langle t\rangle^{-(1+\delta)/2}\Biggr].
\end{align*}
Since \eqref{eq:energy-common-shell-sum} gives
$\sum_m\delta_m\sqrt{\varsigma_m}<\infty$, $T_{71}$ satisfies
\begin{equation}
\begin{aligned}
 |T_{71}|\lesssim M\eps\sqrt{\CK_{f,+}^{(n)}}\Biggl[
 \sqrt{\CK_{f,\partial g,h}^{(n)}}
 +\delta_n
 \left(\sum_{b<n}\delta_b^{-1}\varsigma_b\mathscr Q^{(b)}\right)^{1/2}
 +M^C\eps\delta_n\langle t\rangle^{-(1+\delta)/2}\Biggr].
\end{aligned}
 \label{eq:f-newton-T91-bound}
\end{equation}

For the $T_{731}$ term in \eqref{eq:f-newton-high-velocity-children},
insert \eqref{eq:f-newton-reaction-comparison} into its first summand.
The second summand in the first row defines $T_{732}$; the two summands
in the second row define $T_{741},T_{742}$ in the same order.
Summing the projections gives the explicit bound
\begin{equation}
\begin{aligned}
 |T_{731}|\lesssim{}&\sum_{m<n}\sum_{\substack{k\in\mathbb Z\\\ell\ne0}}
 \iiint_{\substack{(\eta,\xi,\xi')\in\mathbb R^3\\
 |(p,\zeta)|\leq3|(\ell,\xi)|/16\\|\xi-\xi'|\leq|(\ell,\xi')|/16}}
 |(p,\zeta)|e^{c\lambda_n|(p,\zeta)|^s}|\widehat f_p^{(m)}(\zeta)|
 e^{c\lambda_n|\xi-\xi'|^s}|\widehat h^{[n]}(\xi-\xi')|\\
 &\quad\times\Biggl[
 \sqrt{\frac{\partial_tw_k(\eta)}{w_k(\eta)}}
 \left|\sqrt{\A_k^{(n)}(\eta)\widetilde\A_k^{(n)}(\eta)}\widehat f_k^{(n)}(\eta)\right|
 \left|\mathcal M_{w,\ell}^{(n)}(\xi')
 \widehat{\Delta_L\Phi_{\Int}^{(n;m)}}(\ell,\xi')\right|\\
 &\qquad+
 \sqrt{-\dot\lambda_n}\langle k,\eta\rangle^{s/2}
 |\A_k^{(n)}(\eta)\widehat f_k^{(n)}(\eta)|
 \left|\mathcal M_{\lambda,\ell}^{(n)}(\xi')
 \widehat{\Delta_L\Phi_{\Int}^{(n;m)}}(\ell,\xi')\right|
 \Biggr]\,\dd\xi'\dd\xi\dd\eta.
\end{aligned}
\label{eq:f-newton-coefficient-low-majorant}
\end{equation}
The output and precision-weighted stream take the two $L^2$ norms;
the coefficient and differentiated vorticity take Fourier $L^1$ norms.
The coefficient sum in $h^{[n]}$ is independent of $m$, and the stream
remains localized at scale $m$. Young's inequality and
\eqref{eq:augmented-radius-coercivity} therefore give

\begin{align*}
 |T_{731}|&\lesssim\sqrt{\CK_{f,+}^{(n)}}
 \|e^{c\lambda_n|\xi|^s}\widehat h^{[n]}(\xi)\|_{L^1_\xi}
 \sum_{m<n}
 \|| (p,\zeta)|e^{c\lambda_n|(p,\zeta)|^s}
       \widehat f_p^{(m)}(\zeta)\|_{\ell^1_pL^1_\zeta}
 \left(\sum_{\iota=\lambda,w}
 \|\mathcal M_\iota^{(n)}\Delta_L\Phi_{\Int}^{(n;m)}\|_2^2
 \right)^{1/2}.
\end{align*}

The bootstrap gives
$\|e^{c\lambda_n|\xi|^s}\widehat h^{[n]}\|_1\lesssim M^2\eps$
and
$\|| (p,\zeta)|e^{c\lambda_n|(p,\zeta)|^s}\widehat f_p^{(m)}(\zeta)\|_{\ell^1_pL^1_\zeta}
\lesssim M\eps\delta_m$.
Apply \eqref{eq:interior-precision-source},
\eqref{eq:interior-precision-final}, and
\eqref{eq:energy-common-shell-sum}, retaining
$\sum_{m<n}\delta_m\sqrt{\varsigma_m}<\infty$. Consequently,

\begin{equation}
\begin{aligned}
 |T_{731}|\lesssim M^3\eps^2\sqrt{\CK_{f,+}^{(n)}}\Biggl[
 \sqrt{\CK_{f,\partial g,h}^{(n)}}
 +\delta_n
 \left(\sum_{b<n}\delta_b^{-1}\varsigma_b\mathscr Q^{(b)}\right)^{1/2}
 +M^C\eps\delta_n\langle t\rangle^{-(1+\delta)/2}\Biggr].
\end{aligned}
 \label{eq:f-newton-T931-bound}
\end{equation}

For the $T_{732}$ term in \eqref{eq:f-newton-high-velocity-children},
the same pointwise comparison \eqref{eq:f-newton-reaction-comparison}
gives \eqref{eq:f-newton-coefficient-low-majorant} with
$\widehat h^{[n]}(\xi-\xi')$ replaced by $\widehat h^{(n)}(\xi-\xi')$
and $\widehat{\Delta_L\Phi_{\Int}^{(n;m)}}(\ell,\xi')$ replaced by
$\widehat{\Delta_L\Phi_{\Ext}^{[n-1]}}(\ell,\xi')$.
The stream is independent of $m$, whereas
$\|e^{c\lambda_n|\xi|^s}\widehat h^{(n)}\|_1\lesssim M^2\eps\delta_n$.
For $n\geq2$, \eqref{eq:precision-radius-comparison} bounds its level-$n$ precision norm
by the level-$(n-1)$ norm, to which
\eqref{eq:precision-elliptic-partial-sum} applies. Thus

\begin{equation}
\begin{aligned}
 |T_{732}|&\lesssim\sqrt{\CK_{f,+}^{(n)}}
 \|e^{c\lambda_n|\xi|^s}\widehat h^{(n)}(\xi)\|_{L^1_\xi}
 \left(\sum_{\iota=\lambda,w}
 \|\mathcal M_\iota^{(n)}\Delta_L\Phi_{\Ext}^{[n-1]}\|_2^2
 \right)^{1/2}\\
 &\qquad\times
 \sum_{m<n}
 \|| (p,\zeta)|e^{c\lambda_n|(p,\zeta)|^s}
       \widehat f_p^{(m)}(\zeta)\|_{\ell^1_pL^1_\zeta}\\
 &\lesssim M^3\eps^2\delta_n\sqrt{\CK_{f,+}^{(n)}}
 \left[
 \left(\sum_{b<n}\delta_b^{-1}\varsigma_b\mathscr Q^{(b)}\right)^{1/2}
 +M^C\eps\langle t\rangle^{-(1+\delta)/2}\right].
\end{aligned}
 \label{eq:f-newton-T932-bound}
\end{equation}

For the $T_{741}$ and $T_{742}$ terms in \eqref{eq:f-newton-high-velocity-children},
the support gives
\[
 |(\ell,\xi')|<16|\xi-\xi'|,\qquad
 |(p,\zeta)|<\frac{51}{16}|\xi-\xi'|,
 \qquad |(\ell,\xi')|+|(p,\zeta)|<20|\xi-\xi'|.
\]
We transfer the output multiplier to the coefficient directly.
If $|\xi-\xi'|\geq8$ and
$|(\ell,\xi')|+|(p,\zeta)|\leq|\xi-\xi'|/8$,
Lemma~\ref{lem:elliptic-J-transfer} applies with remaining frequency
$(k,\eta-\xi+\xi')$. The Gevrey increment is at most
$s8^{s-1}(|(\ell,\xi')|+|(p,\zeta)|)^s$, with $s8^{s-1}<1$.
On the complementary range use
\[
 \max_{1/8\leq x\leq20}\frac{(1+x)^s-1}{x^s}<1
\]
and \eqref{eq:weight-total-growth}. Since $s>1/2$ and the radii
have a fixed positive lower bound, the remaining square-root exponential
is absorbed by choosing fixed $c',c$ with
\[
 \max\!\left\{s8^{s-1},
 \max_{1/8\leq x\leq20}\frac{(1+x)^s-1}{x^s}\right\}<c'<c<1.
\]
We have
\[
\begin{aligned}
 &e^{C\mu(|(\ell,\xi')|^{1/2}+|(p,\zeta)|^{1/2})}
 e^{-(c-c')\lambda_n(|(\ell,\xi')|^s+|(p,\zeta)|^s)}
 \lesssim
 e^{-\frac12(c-c')\lambda_n(|(\ell,\xi')|^s+|(p,\zeta)|^s)}
 \lesssim 1.
\end{aligned}
\]
The Sobolev ratio is bounded
because $|(k,\eta)|<21|\xi-\xi'|$. Hence
\begin{equation}
 \A_k^{(n)}(\eta)\lesssim\A_{\Res}^{(n)}(\xi-\xi')
 e^{c\lambda_n(|(\ell,\xi')|^s+|(p,\zeta)|^s)}.
 \label{eq:f-newton-coefficient-high-transfer}
\end{equation}
If $|\xi-\xi'|<8$, the same support gives
$|(\ell,\xi')|<128$ and $|(p,\zeta)|<51/2$; all output multipliers
are uniformly bounded and $\A_{\Res}^{(n)}\geq1$, so
\eqref{eq:f-newton-coefficient-high-transfer} also holds there.
Insert this estimate and
$|p\xi'-\ell\zeta|\leq|(\ell,\xi')||(p,\zeta)|$
into the first summand of the second row of
\eqref{eq:f-newton-high-velocity-children}. Summing the projections yields
\begin{equation}
\begin{aligned}
 |T_{741}|\lesssim{}&\sum_{m<n}\sum_{\substack{k\in\mathbb Z\\\ell\ne0}}
 \iiint_{\mathbb R^3}
 |\A_k^{(n)}(\eta)\widehat f_k^{(n)}(\eta)|
 |\A_{\Res}^{(n)}(\xi-\xi')\widehat h^{[n]}(\xi-\xi')|\\
 &\quad\times |(\ell,\xi')|e^{c\lambda_n|(\ell,\xi')|^s}
 |\widehat\Phi_{\Int}^{(n;m)}(\ell,\xi')|
 |(p,\zeta)|e^{c\lambda_n|(p,\zeta)|^s}
 |\widehat f_p^{(m)}(\zeta)|\,\dd\xi'\dd\xi\dd\eta.
\end{aligned}
\label{eq:f-newton-coefficient-high-majorant}
\end{equation}
For $T_{742}$ replace $h^{[n]}$ by $h^{(n)}$ and
$\Phi_{\Int}^{(n;m)}$ by $\Phi_{\Ext}^{[n-1]}$ in this integral. The bootstrap,
the lossy elliptic bounds, Corollary~\ref{cor:elliptic-profile-output},
 \eqref{eq:energy-common-shell-sum} and Young's convolution inequality give
\begin{equation}
\begin{aligned}
 |T_{741}|&\lesssim\sqrt{\E_f^{(n)}}\|\A_{\Res}^{(n)}h^{[n]}\|_2
 \sum_{m<n}
 \||(\ell,\xi)|e^{c\lambda_n|(\ell,\xi)|^s}
       \widehat\Phi_{\Int,\ell}^{(n;m)}(\xi)\|_{\ell^1_{\ell\ne0}L^1_\xi}
 \|| (p,\zeta)|e^{c\lambda_n|(p,\zeta)|^s}
       \widehat f_p^{(m)}(\zeta)\|_{\ell^1_pL^1_\zeta}\\
 &\qquad\lesssim
 (M\eps\delta_n)(M^2\eps)
 \sum_{m<n}(M^C\eps\delta_n\varsigma_m\langle t\rangle^{-2})
             (M\eps\delta_m)
 \lesssim\mathscr R^{(n)},\\
 |T_{742}|&\lesssim\sqrt{\E_f^{(n)}}\|\A_{\Res}^{(n)}h^{(n)}\|_2
 \||(\ell,\xi)|e^{c\lambda_n|(\ell,\xi)|^s}
       \widehat\Phi_{\Ext,\ell}^{[n-1]}(\xi)\|_{\ell^1_{\ell\ne0}L^1_\xi}
 \sum_{m<n}
 \|| (p,\zeta)|e^{c\lambda_n|(p,\zeta)|^s}
       \widehat f_p^{(m)}(\zeta)\|_{\ell^1_pL^1_\zeta}\\
 &\qquad\lesssim
 (M\eps\delta_n)(M^2\eps\delta_n)
 (M^C\eps\langle t\rangle^{-2})\sum_{m<n}M\eps\delta_m
 \lesssim\mathscr R^{(n)} .
\end{aligned}
 \label{eq:f-newton-coefficient-high-bounds}
\end{equation}

For the $T_{72}$ term in \eqref{eq:f-newton-stream-g-children}, $|\xi|\geq4$ and
$|(k,\eta-\xi)|\leq3|\xi|/16$ on its support.
The estimate \eqref{eq:reaction-multiplier-zero-mode}
therefore gives directly
\begin{equation}
\begin{aligned}
 |T_{72}|\lesssim{}&\sum_{m<n}\sum_k
 \iint_{\substack{(\eta,\xi)\in\mathbb R^2\\|\xi|\geq4}}
 |\eta-\xi|e^{c\lambda_n|(k,\eta-\xi)|^s}
 |\widehat f_k^{(m)}(\eta-\xi)|\\
 &\quad\times\Biggl[
 \sqrt{\frac{\partial_tw_k(\eta)}{w_k(\eta)}}
 \left|\sqrt{\A_k^{(n)}(\eta)\widetilde\A_k^{(n)}(\eta)}\widehat f_k^{(n)}(\eta)\right|
 \langle t\rangle^{1+s}\sqrt{\frac{\partial_tw_0(\xi)}{w_0(\xi)}}
 |\widetilde\A_0^{(n)}(\xi)\widehat g^{(n)}(\xi)|\\
 &\qquad+
 \sqrt{-\dot\lambda_n}\langle k,\eta\rangle^{s/2}
 |\A_k^{(n)}(\eta)\widehat f_k^{(n)}(\eta)|
 \frac{\langle\xi\rangle^{-s/2}}{\sqrt{-\dot\lambda_n}}
 |\A_0^{(n)}(\xi)\widehat g^{(n)}(\xi)|
 \Biggr]\,\dd\xi\dd\eta.
\end{aligned}
\label{eq:f-newton-coordinate-majorants}
\end{equation}
Since $\widetilde\A_0^{(n)}=\A_{\NR}^{(n)}$ and
$\A_0^{(n)}\leq2\A_{\NR}^{(n)}$, the definition \eqref{eq:CK-dg}
and \eqref{eq:lambda-derivative} give
\begin{equation}
\begin{aligned}
 \left\|\langle t\rangle^{1+s}\sqrt{\frac{\partial_tw_0}{w_0}}
 \widetilde\A_0^{(n)}\widehat g^{(n)}\right\|_{L^2(|\xi|\geq4)}
 &\lesssim\sqrt{\CK_{\partial g}^{(n)}},\\
 \left\|\frac{\langle\xi\rangle^{-s/2}}{\sqrt{-\dot\lambda_n}}
 \A_0^{(n)}\widehat g^{(n)}\right\|_{L^2(|\xi|\geq4)}
 &\lesssim\langle t\rangle^{a-s}\sqrt{\CK_{\partial g}^{(n)}}.
\end{aligned}
\label{eq:f-newton-coordinate-densities}
\end{equation}
Indeed, the ratios to the corresponding square roots of the $\CK$ 
integrands are bounded by $\langle\xi\rangle^s/|\xi|$ and
$2\langle t\rangle^{-1-s}(-\dot\lambda_n)^{-1}
\langle\xi\rangle^{s/2}/|\xi|^{1+s/2}\lesssim\langle t\rangle^{a-s}$.

For the $T_{72}$ term in \eqref{eq:f-newton-stream-g-children}, apply
Young's inequality to \eqref{eq:f-newton-coordinate-majorants},
using \eqref{eq:augmented-radius-coercivity} and
\eqref{eq:f-newton-coordinate-densities}. This gives \begin{equation}
 |T_{72}|\lesssim(1+\langle t\rangle^{a-s})
 \sqrt{\CK_{f,+}^{(n)}\CK_{\partial g}^{(n)}}
 \sum_{m<n}
 \||\zeta|e^{c\lambda_n|(p,\zeta)|^s}
                \widehat f_p^{(m)}(\zeta)\|_{\ell^1_pL^1_\zeta}\lesssim
 M\eps\sqrt{\CK_{f,+}^{(n)}\CK_{\partial g}^{(n)}}.
 \label{eq:f-newton-T92-bound}
\end{equation} 
The terms $T_{71},T_{72},T_{731},T_{732},T_{741},T_{742}$ complete the estimate for
$T_7$.

A standard application of Young's inequality for products yields 
\[
 CM^C\eps^2\delta_n^2
 \sum_{b<n}\delta_b^{-1}\varsigma_b\mathscr Q^{(b)}
 +CM^C\eps^4\delta_n^2\langle t\rangle^{-1-\delta}
 \leq C\mathscr N^{(n)}+C\mathscr R^{(n)}.
\]

Combining \eqref{eq:f-newton-T8-bound}, \eqref{eq:f-newton-T10-bound},
\eqref{eq:f-newton-T91-bound}, \eqref{eq:f-newton-T931-bound},
\eqref{eq:f-newton-T932-bound}, \eqref{eq:f-newton-coefficient-high-bounds},
and \eqref{eq:f-newton-T92-bound}, and using
$\sqrt{\CK_{f,+}^{(n)}\CK_{f,\partial g,h}^{(n)}}
\leq\CK_{f,+}^{(n)}
 +(\CK_{\partial g}^{(n)}+M^{-2}\CK_h^{(n)})/2$
gives
\begin{equation*}
\begin{aligned} 
 |\mathcal N_f^{(n)}|
 \leq\frac1{32}\CK_{f,+}^{(n)}
 +CM^C\eps\left(\CK_{\partial g}^{(n)}+M^{-2}\CK_h^{(n)}\right)
 +C\mathscr N^{(n)}+C\mathscr R^{(n)} .
\end{aligned}
\end{equation*} 
This concludes the proof of \eqref{eq:f-Newton-bound}.

\end{proof}

Combining the energy relation \eqref{eq:f-energy-identity} with
\eqref{eq:f-transport-bound}, \eqref{eq:f-reaction-bound}, and
\eqref{eq:f-Newton-bound} from Lemma~\ref{lem:f-energy-errors} yields
\begin{equation}
  \frac12\frac{\dd}{\dd t}\E_f^{(n)}
  +\frac34\CK_{f,+}^{(n)}
  \leq (-\dot\lambda_n)\E_f^{(n)}+C\mathscr R^{(n)}+C\mathscr N^{(n)}
   +CM^C\eps\left(\CK_{\partial g}^{(n)}
      +M^{-2}\CK_h^{(n)}+\CK_{\mathrm{aux}}^{(n)}\right).
  \label{eq:f-energy-closed}
\end{equation}

\subsection{Estimates for the coordinate quantities}
We collect estimates for the coordinate quantities $h^{(n)},g^{(n)}$.

\begin{lemma}[\(\E_h\) estimate]
\label{lem:h-energy-estimate}
Under the bootstrap assumptions,
\begin{equation}
  \frac12\frac{\dd}{\dd t}\E_h^{(n)}
  +\frac34\CK_h^{(n)}
  \lesssim
  \CK_{\partial g}^{(n)}+\CK_{\mathrm{aux}}^{(n)}
  +\mathscr R^{(n)}+\mathscr N^{(n)}.
  \label{eq:h-energy-bound}
\end{equation}
Here, the $\CK_{h}, \CK_{\pa g}, \CK_{\rm aux}$-terms are  defined in \eqref{eq:CK-dg}, \eqref{eq:CK-h} and \eqref{eq:CK-aux}. 
\end{lemma}

\begin{proof}The proof is an adaptation of the arguments in \cite{BM13}; we omit the details.
\end{proof}

\begin{lemma}[\(\E_{\pa g}\)-estimate]
\label{lem:dg-energy-estimate}
Under the bootstrap assumptions,
{\begin{align}
  \frac12\frac{\dd}{\dd t}\E_{\partial g}^{(n)}
  +\frac34\CK_{\partial g}^{(n)}
  +\frac{1-s}{2t}\E_{\partial g}^{(n)}
  &\lesssim
  M^C\eps\CK_{f,+}^{(n)}+M^C\eps M^{-2}\CK_h^{(n)}
  \notag\\
  &\quad{}
  +M^C\eps\avg t^{-1-a/2}\E_{g,\mathrm{low}}^{(n)}
  +\mathscr R^{(n)}+\mathscr N^{(n)}.
  \label{eq:dg-energy-bound}
\end{align}}
\end{lemma}

\begin{proof}The proof is an adaptation of the arguments in \cite{BM13}; we omit the details.
\end{proof}

\begin{lemma}[\(\E_{ g,\rm low}\)-estimate]
\label{lem:g-low-energy-estimate}
Under the bootstrap assumptions, after $\eps$ is chosen sufficiently small,
\begin{equation}
  \frac12\frac{\dd}{\dd t}\E_{g,\mathrm{low}}^{(n)}
  +\frac12\CK_{g,\mathrm{low};\lambda}^{(n)}
  +\frac{4-\beta_g}{4t}\E_{g,\mathrm{low}}^{(n)}
  \lesssim\mathscr R^{(n)}.
  \label{eq:g-low-energy-bound}
\end{equation}

\end{lemma}

\begin{proof} The proof is an adaptation of the arguments in \cite{BM13}; we omit the details.
\end{proof}

\section{Proof of the main theorem}
\label{sec:convergence}

We first close the energy and support bootstraps for finite auxiliary
systems.  Their terminal sums yield global Euler solutions for data
supported away from the walls, with bounds independent of the terminal
layer.  We then pass to the compact-support approximation and prove
Theorem~\ref{thm:main}.

For the convergence and coordinate reconstruction arguments, fix a common radius factor $c\in(0,1)$ for
Corollary~\ref{cor:elliptic-profile-output} and
Lemma~\ref{lem:newton-velocity-output}.
Keep this choice fixed when choosing $\lambda_\infty$ and throughout
the reconstruction. We now choose
the radius $\lambda_\infty$ in \eqref{eq:main-parameter-choices} so that
\begin{equation}
  0<2\lambda_\infty
  <c\inf_{n,t\geq1}\lambda_n(t).
  \label{eq:final-gevrey-radius}
\end{equation}

\subsection{Closing the energy bootstrap}
Recall $\E^{(n)}$, $\CK_{f,\partial g,h}^{(n)}$, and $\mathscr Q^{(n)}$
from \eqref{eq:combined-energy}, \eqref{eq:combined-CK}, and
\eqref{eq:newton-CK-coupling}, respectively.

\begin{lemma}[Auxiliary coercivity comparison]
\label{lem:auxiliary-coercivity-comparison}
For every $1\leq n\leq {\mathcal{N}}$ and $t\geq1$,
\begin{equation}
 \CK_{\mathrm{aux}}^{(n)}
 \le C\left\{\CK_{\partial g}^{(n)}
       +\avg t^{-1-a/2}\E_{g,\mathrm{low}}^{(n)}\right\}.
 \label{eq:auxiliary-coercivity-closure}
\end{equation}
\end{lemma}

\begin{proof}
On $|\eta|\geq1$, the factor $|\eta|$ in
$\partial_vg^{(n)}$ controls both integrands in
\eqref{eq:CK-aux}, and
$\widetilde {\A}_0^{(n)}={\A}_{\NR}^{(n)}$.  This gives
$C\CK_{\partial g}^{(n)}$.  On $|\eta|<1$, $w_{\NR}(t,\eta)=1$, so
$\partial_t w_{\NR}(t,\eta)/w_{\NR}(t,\eta)=0$.  For the remaining $\lambda$ term,
\eqref{eq:lambda-derivative} and $\beta_g=4-a/2$ give
\begin{equation*}
 (-\dot\lambda_n)\avg t^{2+2s-\beta_g}
 \lesssim\avg t^{-1-a/2};
\end{equation*}
the bounded Fourier multipliers at $|\eta|<1$ are controlled by
$\E_{g,\mathrm{low}}^{(n)}$.  This proves the claim.
\end{proof}

\begin{proposition}[Combined energy inequality]
\label{prop:combined-energy-inequality}
Choose \(M\) sufficiently large and then \(\eps\) sufficiently small.  Under
the bootstrap assumptions, every shell $1\leq n\leq {\mathcal{N}}$ satisfies
\begin{equation}
\begin{aligned}
 &\frac{\dd}{\dd t}\E^{(n)}(t)
  +\tfrac18\CK_{f,\partial g,h}^{(n)}(t)
  +\frac{1-s}{2t}\E_{\partial g}^{(n)}(t)
     +\frac a{8t}\E_{g,\mathrm{low}}^{(n)}(t)\\
 &\quad\leq2(-\dot\lambda_n(t))\E^{(n)}(t)
       +C\mathscr R^{(n)}(t)+C\mathscr N^{(n)}(t).
\end{aligned}
  \label{eq:combined-energy-inequality}
\end{equation}
for a fixed constant $C>0$.
\end{proposition}

\begin{proof} \space Recall the bounds \eqref{eq:f-energy-closed}, \eqref{eq:dg-energy-bound}, 
\eqref{eq:h-energy-bound} and
\eqref{eq:g-low-energy-bound}. Applying
Lemma~\ref{lem:auxiliary-coercivity-comparison}, $\avg t^{-1-a/2}\leq t^{-1}$ and
$(4-\beta_g)/4=a/8$ gives 
\begin{align*}
 &\tfrac12\frac{d}{dt}\E^{(n)}
  +\tfrac18\CK_{f,\partial g,h}^{(n)}
  +\frac{1-s}{4t}\E_{\partial g}^{(n)}
  +\frac a{8t}\E_{g,\mathrm{low}}^{(n)}
  +\tfrac12\CK_{g,\mathrm{low};\lambda}^{(n)}\\
 &\quad\leq(-\dot\lambda_n)\E_f^{(n)}
  +(2CM^C\eps+CM^{-2})
       \bigl(\CK_{f,\partial g,h}^{(n)}
                      +t^{-1}\E_{g,\mathrm{low}}^{(n)}\bigr)
  +C(\mathscr R^{(n)}+\mathscr N^{(n)}).
\end{align*}
This estimate and Lemma~\ref{lem:g-low-energy-estimate} yield the result. 
\end{proof}

Multiply \eqref{eq:combined-energy-inequality} by
$
 \exp\left(-2\int_1^t(-\dot\lambda_n(\tau))\dd\tau\right),$ 
which lies between $e^{-2\lambda_0}$ and $1$.
Since $a<1$ and $a\leq1-s$, the dissipation satisfies
\[
 \frac18\CK_{f,\partial g,h}^{(n)}
 +\frac{1-s}{2t}\E_{\partial g}^{(n)}
 +\frac a{8t}\E_{g,\mathrm{low}}^{(n)}
 \geq\frac a8\mathscr Q^{(n)},
\]
 Integration therefore gives
\begin{equation}
 \max\left\{
   \sup_{1\leq t\leq T}\E^{(n)}(t),
   \frac a8\int_1^T\mathscr Q^{(n)}(t)\dd t
 \right\}
 \leq e^{2\lambda_0}\left[
   \E^{(n)}(1)
   +C\int_1^T(\mathscr R^{(n)}+\mathscr N^{(n)})(t)\dd t
 \right].
 \label{eq:integrated-energy-inequality}
\end{equation}
The integrable remainder, the bootstrap assumption, and
\eqref{eq:energy-common-shell-sum} yield
\begin{equation}
 \int_1^T\mathscr R^{(n)}(t)\dd t
 \lesssim M^C\eps^3\delta_n^2,\qquad
 \int_1^T\mathscr N^{(n)}(t)\dd t
 \leq4M^{C+2}\eps^3\delta_n^2
       \sum_{m\leq n}\varsigma_m\delta_m
 \lesssim M^{C+2}\eps^3\delta_n^2.
 \label{eq:integrated-Newton-coupling-chain}
\end{equation}

The auxiliary comparison \eqref{eq:auxiliary-coercivity-closure}
also gives
\begin{equation}
 M^{-2}\int_1^T\CK_{\mathrm{aux}}^{(n)}(t)\dd t
 \leq CM^{-2}\int_1^T
   \left(\CK_{\partial g}^{(n)}
       +t^{-1}\E_{g,\mathrm{low}}^{(n)}\right)(t)\dd t
 \leq\int_1^T\mathscr Q^{(n)}(t)\dd t,
 \label{eq:integrated-auxiliary-recovery}
\end{equation}
where the fixed lower bound for $M$ includes $CM^{-2}\leq1$
for the constant in this comparison.
Consequently, the definition \eqref{eq:bootstrap-envelope}
and $\E^{(n)}(1)=\mathcal B_n(1)$ imply
\[
 \begin{aligned}
 \mathcal B_n(T)
 &\leq2\sup_{1\leq t\leq T}\E^{(n)}(t)
       +2\int_1^T\mathscr Q^{(n)}(t)\dd t\\
 &\leq\left(2+\frac{16}{a}\right)e^{2\lambda_0}
 \left[\mathcal B_n(1)
       +C\int_1^T(\mathscr R^{(n)}+\mathscr N^{(n)})(t)\dd t
 \right].
 \end{aligned}
\]
Since $2+16/a\leq32/a$, the previously fixed constant
$C_{\rm B}$ from \eqref{eq:bootstrap-constant} gives
\begin{equation}
 \mathcal B_n(T)
 \leq C_{\rm B}\mathcal B_n(1)
       +CM^{C+2}\eps^3\delta_n^2.
 \label{eq:energy-bootstrap-conclusion}
\end{equation}
Using \eqref{eq:initial-bootstrap-margin}, and decreasing
$\eps$ with $M$ fixed so that $CM^C\eps\leq1$, we obtain
\[
 \mathcal B_n(T)
 \leq\left(\frac1{100}+1\right)M^2\eps^2\delta_n^2
 \leq2M^2\eps^2\delta_n^2.
\]
This proves \eqref{eq:bootstrap-energy-improvement},
uniformly in $n,T$, and $\mathcal N$.
 We next prove the support improvement on the same finite interval.

\subsection{Finite-profile bounds and support improvement}

\begin{lemma}
\label{lem:finite-terminal-bounds}
Let a finite system solve \eqref{eq:n-system} on $[1,T]$
and satisfy \eqref{eq:bootstrap-energy-assumption}. Write
$(f,h,g,\phi)=(f^{[{\mathcal{N}}]},h^{[{\mathcal{N}}]},g^{[{\mathcal{N}}]},\phi^{[{\mathcal{N}}]})$.
The terminal fields solve \eqref{eq:profile-system} with the sum
of their actual auxiliary time-one data, and
\begin{align}
 &f\in C([1,T];\Gs{s}{\lambda_\infty}{\sigma}),\qquad
 h,g\in C([1,T];\Gs{s}{\lambda_\infty}{\sigma-6}),
 \label{eq:terminal-regularity}\\
 &\sum_{n=1}^{\mathcal{N}}\sup_{1\leq t\leq T}
 \left(\|f^{(n)}(t)\|_{\Gs{s}{\lambda_\infty}{\sigma}}
 +\|h^{(n)}(t)\|_{\Gs{s}{\lambda_\infty}{\sigma-6}}
 +\|g^{(n)}(t)\|_{\Gs{s}{\lambda_\infty}{\sigma-6}}\right)
 \lesssim M^C\eps.
 \label{eq:finite-sum-Gevrey-bound}
\end{align}
The constants are independent of ${\mathcal{N}}$ and $T$.
\end{lemma}

\begin{proof}
The lower bounds in \eqref{eq:main-multiplier}, the assumed energy bootstrap,
and the strict radius gap in \eqref{eq:final-gevrey-radius} \space bound the $n$th summand in \eqref{eq:finite-sum-Gevrey-bound} by
$CM^C\eps\delta_n$.
For $g^{(n)}$, use its low energy on $|\eta|\leq1$ and its
differentiated energy on $|\eta|>1$. Summability of $\delta_n$ proves
\eqref{eq:finite-sum-Gevrey-bound}, uniformly in $\mathcal N$.
Sobolev continuity and the Fourier tail at the strict radius gap give
\eqref{eq:terminal-regularity}. Lemma~\ref{prop:newton-equivalence}
identifies the terminal auxiliary solution.
For the original finite data in \eqref{eq:newton-data}, \eqref{eq:terminal-cutoff-identity} also
identifies its time-one datum with the prescribed physical datum.
This last identification is not asserted for the stronger-Gevrey
regularizations.
\end{proof}

Recall $\mathbf d(v)=\dist(v,\left\{L,R\right\})$ and $d_n=2^{-n}$
from \eqref{eq:wall-distance-scales}. We now prove the support estimate used in
the continuation argument below. The key point is to bound the vertical drift relative to the distance from the wall. The vertical component of the profile transport field is
\begin{equation}
  b^{[{\mathcal{N}}]}=(1+h^{[{\mathcal{N}}]})\partial_z\phi_{\neq}^{[{\mathcal{N}}]}+g^{[{\mathcal{N}}]}.
  \label{eq:normal-profile-drift}
\end{equation}
It vanishes at $L,R$ by the Dirichlet trace and the support of $g^{(n)}$.

\begin{lemma}[Normalized drift and support improvement]
\label{lem:support-improvement}
Let the finite system satisfy
\eqref{eq:bootstrap-energy-assumption}--\eqref{eq:bootstrap-support-assumption}
on $[1,T]$, with the initial support in \eqref{eq:support-scale}.
Then, for $1\leq t\leq T$,
\begin{equation}
  \norm{\frac{b^{[{\mathcal{N}}]}(t)}{\mathbf d}}_{L^\infty(\T\times(L,R))}
  \lesssim M^C\eps
    \left(\frac1{\avg t^2}+\frac1{\avg t^{1+s}}\right).
  \label{eq:normalized-drift}
\end{equation}
If $\eps$ is sufficiently small, then every increment satisfies
\begin{equation}
  \supp X^{(n)}(t)\subset
  \left\{v:\mathbf d(v)\geq\frac{5d_n}{88}\right\},
  \qquad X\in\left\{f,g,h\right\},\quad 1\leq t\leq T.
  \label{eq:global-support-bound}
\end{equation}
The drift bounds also hold for every prefix $1\leq m\leq\mathcal N$
in place of the terminal index $\mathcal N$.
\end{lemma}

\begin{proof}
The functions $\partial_z\phi_{\neq}^{[{\mathcal{N}}]}$ and $g^{[{\mathcal{N}}]}$ vanish at
both endpoints.  Integrating their vertical derivatives from the closer
endpoint gives
\[
 \|b^{[{\mathcal{N}}]}/\mathbf d\|_\infty
 \leq(1+\|h^{[{\mathcal{N}}]}\|_\infty)
       \|\partial_v\partial_z\phi_{\neq}^{[{\mathcal{N}}]}\|_\infty
       +\|\partial_vg^{[{\mathcal{N}}]}\|_\infty.
\]
Corollary~\ref{cor:elliptic-profile-output}, with one additional derivative,
and \eqref{eq:bootstrap-consequences-partial-g} prove
\eqref{eq:normalized-drift}; the same argument bounds $g^{[{\mathcal{N}}]}/\mathbf d$.
In particular
\begin{equation}
 \int_1^T(\|b^{[{\mathcal{N}}]}/\mathbf d\|_\infty+\|g^{[{\mathcal{N}}]}/\mathbf d\|_\infty)\dd t
 \leq CM^C\eps.\label{eq:integrated-normalized-drift}
\end{equation}
Choose the common majorant and moving radius
\[
 q(t)=CM^C\eps(\avg t^{-2}+\avg t^{-1-s}),\qquad
 R_n^{\mathrm{supp}}(t)=\frac{d_n}{16}\exp\left(-\int_1^tq\right).
\]
Choose the constant in $q$ to dominate both drift bounds, and
keep that choice fixed in $R_n^{\mathrm{supp}}$. Along characteristics of either vertical drift,
\begin{equation}
 \left|\frac{\dd}{\dd t}\log \mathbf d(V(t))\right|\leq q(t)
 \label{eq:relative-distance-characteristic}
\end{equation}
almost everywhere, including crossings of the midpoint by the
Lipschitz chain rule.  Thus data supported at distance at least $\frac{d_n}{16}$ from the walls, and sources
supported at distance at least $R_n^{\mathrm{supp}}(t)$, remain outside this moving barrier.
Induct on $n$, applying the barrier at each level in the order
$f^{(n)},g^{(n)},h^{(n)}$.  In \eqref{eq:n-system-f} the source contains a
vorticity factor from an earlier level; in \eqref{eq:Lg-expanded} it has
either $f^{(n)}$ or a vorticity factor from an earlier level.  The remaining sources
in \eqref{eq:n-system-h} contain already localized coordinate fields
or increments from earlier levels.  Since $\frac{d_m}{16}\geq \frac{d_n}{16}$ for $m<n$, each source is
outside the same level-$n$ barrier.  Transport and stretching do not
enlarge support.  Choosing $e^{-CM^C\eps}\geq10/11$ proves
\eqref{eq:global-support-bound}.
\end{proof}

\subsection{Construction and global closure of the finite systems}
\label{sec:finite-global-closure}

The intermediate partial sums solve the profile equations, but need not
be profiles of Euler solutions on the prescribed physical channel.
We first construct these auxiliary solutions. Euler compatibility will
be used only for the terminal datum and for reconstruction of the final
solution.

\begin{proposition}[Local construction of auxiliary profiles]
\label{prop:finite-profile-lwp}
Let $r\in\mathbb N$, $r\geq5$, and $T_0\geq1$. Here $T_0$ is the initial or restart time:
$T_0=1$ for the initial construction, and $T_0>1$ is allowed for continuation. Let
$(f_{\mathrm{in}},h_{\mathrm{in}},g_{\mathrm{in}})$ be smooth and vanish
in a neighborhood of the two endpoints of $I=[L,R]$, with
$1+h_{\mathrm{in}}\geq\kappa>0$. Define
\begin{equation}
 \mathfrak r_{\mathrm{in}}
 =T_0(1+h_{\mathrm{in}})\partial_vg_{\mathrm{in}}
   +(f_{\mathrm{in}})_0+h_{\mathrm{in}}.
 \label{eq:auxiliary-defect-definition}
\end{equation}
There is a unique smooth, boundary-flat solution of
\eqref{eq:profile-system} on $[T_0,T_1]$, $T_1>T_0$, whose length depends only on
$T_0,r,\kappa,I$ and
$\|(f_{\mathrm{in}},h_{\mathrm{in}},\mathfrak r_{\mathrm{in}})\|_{H^r}$.
It continues through a finite interval whenever $1+h$ remains bounded away from zero
and $(f,h,\mathfrak r)$ remains bounded in $H^r$. No coordinate-length
or wall-velocity compatibility condition is imposed.
\end{proposition}

\begin{proof}
The zero-mode profile equation is
\[
 (\partial_t+g\partial_v)f_0
 =-(1+h)\partial_v\Pzero(f_{\neq}\partial_z\phi_{\neq}).
\]
To recover the derivative of $g$ without losing a derivative of $h$,
we use the transported defect
\begin{equation}
 \begin{aligned}
 \mathfrak r&=t(1+h)\partial_vg+f_0+h,(\partial_t+g\partial_v)\mathfrak r=0,\\
 (\partial_t+g\partial_v)h&=\frac{\mathfrak r-f_0-h}{t},
 \quad \partial_vg=\frac{\mathfrak r-f_0-h}{t(1+h)},\quad g(t,L)=0.
 \end{aligned}
 \label{eq:auxiliary-defect-system}
\end{equation}
Indeed,
\[
 (\partial_t+g\partial_v)\bigl(t(1+h)\partial_vg\bigr)
 =-(1+h)\partial_vg+(1+h)\partial_v\Pzero(f_{\neq}\partial_z\phi_{\neq}),
\]
which cancels $(\partial_t+g\partial_v)(f_0+h)$. Thus the primitive for $g$ gains one derivative.
The reconstructed $g$ satisfies its profile equation because its
residual has zero $v$ derivative and vanishes at $L$. At the other wall,
$\partial_tg(t,R)=-(\partial_vg(t,R)+2/t)g(t,R)$, so its zero trace
persists.

For existence, define the fixed interval and prescribed shear by
\begin{equation}
 \begin{gathered}
 Y_{\mathrm{in}}(v)=-1+\int_L^v\frac{\dd w}{1+h_{\mathrm{in}}(w)},
 \quad \ell_{\mathrm{aux}}=\int_L^R\frac{\dd v}{1+h_{\mathrm{in}}(v)},
 \quad I_{\mathrm{aux}}=[-1,-1+\ell_{\mathrm{aux}}],\\
 V_{\mathrm{in}}=Y_{\mathrm{in}}^{-1},\qquad
 \\
 S_{\mathrm{aux}}(Y_{\mathrm{in}}(v))
 =v+T_0g_{\mathrm{in}}(v)+(1+h_{\mathrm{in}}(v))\partial_v(\phi_{\mathrm{in}})_0(v).
 \end{gathered}
 \label{eq:auxiliary-physical-shear}
\end{equation}
Then $S_{\mathrm{aux}}'(Y_{\mathrm{in}})=1+\mathfrak r_{\mathrm{in}}$.
The $H^r$ transport iteration, using the Dirichlet inverse to recover
$\psi$, gives a local
solution on $\mathbb T_{2\pi}\times I_{\mathrm{aux}}$ of
\begin{equation}
 \begin{gathered}
 \partial_t\omega+
 (S_{\mathrm{aux}}(y)e_x+\nabla^\perp\psi)\cdot\nabla\omega=0,
 \quad \Delta\psi=\omega,\quad\psi|_{\partial I_{\mathrm{aux}}}=0,\\
 \omega(T_0,x,y)
 =f_{\mathrm{in}}(x-T_0V_{\mathrm{in}}(y),V_{\mathrm{in}}(y)),
 \end{gathered}
 \label{eq:auxiliary-physical-equation}
\end{equation}
where $e_x=(1,0)$ and $\nabla^\perp\psi=(-\partial_y\psi,\partial_x\psi)$.
We now solve
\[
 \partial_t(tV)=S_{\mathrm{aux}}-\partial_y\psi_0,
 \qquad V(T_0)=V_{\mathrm{in}}.
\]
Averaged transport gives
$\partial_t(-\partial_y\psi_0)=\Pzero(\omega_{\neq}\partial_x\psi_{\neq})$.
The integration constant is zero because
$\int_{I_{\mathrm{aux}}}\partial_y\psi_0\,\dd y=0$ and
$\int\Pzero(\Delta\psi\,\partial_x\psi)\,\dd y=0$.
Consequently $V$ preserves the endpoints. After decreasing $T_1$, $\partial_yV>0$ on $[T_0,T_1]$,
and pullback by $z=x-tV(t,y)$, $v=V(t,y)$ gives
\eqref{eq:profile-system}, with $h=\partial_yV-1$, $g=\partial_tV$.

The primitive in \eqref{eq:auxiliary-defect-system} and the Sobolev
transport estimate give, on $[T_0,T_1]$, for some fixed integer $p_r$,
\begin{equation}
 \begin{aligned}
 \|g\|_{H^{r+1}}&\leq C_{r,T_0,\kappa}
 (1+\|h\|_{H^r})^{p_r}\|(f,h,\mathfrak r)\|_{H^r},\\
 \frac{\dd}{\dd t}\|(f,h,\mathfrak r)\|_{H^r}
 &\leq C_{r,T_1,\kappa}
 (1+\|(f,h,\mathfrak r)\|_{H^r})^{p_r}\|(f,h,\mathfrak r)\|_{H^r}.
 \end{aligned}
 \label{eq:auxiliary-sobolev-estimate}
\end{equation}
The $H^{r-1}$ difference estimate is taken on the common profile
interval $I$; the primitive controls $g_1-g_2$ in $H^r$, including
the term $(g_1-g_2)\partial_v\mathfrak r_2$.
Uniqueness, smooth persistence, and continuation follow.
\end{proof}

\begin{lemma}[Gevrey persistence for the auxiliary equation]
\label{lem:stronger-gevrey-channel-persistence}
Let $0<\rho\leq1$. A smooth boundary-flat solution of
\eqref{eq:auxiliary-physical-equation}, with Gevrey-$\rho$ initial
vorticity and fixed Gevrey-$\rho$ shear $S_{\mathrm{aux}}$, retains a
positive Gevrey-$\rho$ radius on each compact Sobolev-lifespan interval.
If the initial auxiliary coordinate $V_{\mathrm{in}}$ is
Gevrey-$\rho$ on $I_{\mathrm{aux}}$ and $\partial_yV_{\mathrm{in}}>0$,
then, writing $Y(t,\cdot)=V(t,\cdot)^{-1}$, the profiles
\[
 h=(\partial_yV-1)\circ Y,\qquad g=\partial_tV\circ Y,\qquad
 f(t,z,v)=\omega(t,z+tv,Y(t,v))
\]
retain a positive Gevrey-$\rho$ radius on each compact
Sobolev-lifespan interval on which $1+h$ stays bounded away from zero.
\end{lemma}

\begin{proof}
For $r\geq5$, we use 
\eqref{eq:factorial-Gevrey-definition} and define
\begin{equation}
 \|\omega\|_{\mathcal F_{\rho,\tau}^{r}}
 =\sum_{m\geq0}\frac{\tau^m}{(m!)^{1/\rho}}
             \|\nabla^m\omega\|_{H^r,\mathrm{sum}}.
 \label{eq:auxiliary-factorial-norm}
\end{equation}
Across a strict radius gap this norm has the form 
\eqref{eq:factorial-Fourier-equivalence}, with $\tau=(\rho\lambda)^{1/\rho}$.
Odd reflection of the flat vorticity gives the derivativewise Dirichlet
estimate, uniformly for interval lengths in a compact subset of $(0,\infty)$.
Choose
$\tau_{\mathrm{sh}}>0$ such that
\[
 K_{\mathrm{sh}}:=\|S_{\mathrm{aux}}\|_{H^{r+1}}
 +\sum_{j\geq1}\frac{\tau_{\mathrm{sh}}^{j-1}}{(j!)^{1/\rho}}
                    \|S_{\mathrm{aux}}^{(j)}\|_{H^r}<\infty.
\]
The shear commutator is
\[
 [\partial_x^{\alpha_1}\partial_y^{\alpha_2},S_{\mathrm{aux}}\partial_x]\omega
 =\sum_{j=1}^{\alpha_2}\binom{\alpha_2}{j}S_{\mathrm{aux}}^{(j)}
       \partial_x^{\alpha_1+1}\partial_y^{\alpha_2-j}\omega.
\] \space For $\tau\leq\tau_{\mathrm{sh}}$, factorial convolution bounds this
commutator by
$C_rK_{\mathrm{sh}}\tau\partial_\tau\|\omega\|_{\mathcal F_{\rho,\tau}^{r}}$. \space The nonlinear estimate uses the derivative gained by the Dirichlet
inverse and the vanishing transport wall flux. We omit the factorial Leibniz-rule calculation for the nonlinear transport term; together these estimates give
\begin{equation}
 \begin{aligned}
 D^+\|\omega\|_{\mathcal F_{\rho,\tau(t)}^{r}}
 &\leq C_r(\|\omega\|_{H^r}+K_{\mathrm{sh}})
              \|\omega\|_{\mathcal F_{\rho,\tau(t)}^{r}}\\
 &\quad+\left[\dot\tau+C_r\tau
       (K_{\mathrm{sh}}+\|\omega\|_{\mathcal F_{\rho,\tau(t)}^{r}})\right]
              \partial_\tau\|\omega\|_{\mathcal F_{\rho,\tau(t)}^{r}}.
 \end{aligned}
 \label{eq:auxiliary-factorial-energy}
\end{equation}
Here $D^+$ is the upper right derivative, and the estimate is first
applied to finite derivative sums. On a compact Sobolev interval $[T_0,T]$,
choose $0<\tau(T_0)<\tau_{\mathrm{sh}}$ with finite initial norm. \space Choose $\dot\tau/\tau$ to be $-2C_r$ times the sum of
$K_{\mathrm{sh}}$ and the Gronwall majorant from the first line of
\eqref{eq:auxiliary-factorial-energy}. This majorant is finite on
$[T_0,T]$, so comparison bounds all truncated norms and gives
$\inf_{[T_0,T]}\tau>0$, independently of the derivative cutoff.
Monotone convergence removes the truncation. For the profiles, the
coordinate equation gives
\begin{equation}
 V(t,y)=\frac{T_0}{t}V_{\mathrm{in}}(y)
 +\frac1t\int_{T_0}^t
   \bigl(S_{\mathrm{aux}}(y)-\partial_y\psi_0(q,y)\bigr)\,\dd q
 =:T_1+T_2.
 \label{eq:auxiliary-coordinate-persistence}
\end{equation}
The hypothesis on $V_{\mathrm{in}}$ controls $T_1$; the Gevrey bound for $\omega$ obtained from \eqref{eq:auxiliary-factorial-energy}
and the derivativewise Dirichlet estimate control $T_2$.
Since $\partial_yV$ is bounded away from zero, intrinsic Gevrey inversion
and composition give a positive radius for $Y$ and the profiles in the
statement, using
$\partial_tV=(S_{\mathrm{aux}}-\partial_y\psi_0-V)/t$.
\end{proof}

\begin{lemma}[Uniform local Gevrey bounds]
\label{lem:finite-profile-gevrey-entrance}
There are radii
\begin{equation}
 \lambda_0<\lambda_b<\lambda_a<3\lambda_0/2
 \label{eq:finite-entrance-radii}
\end{equation}
and $\tau\in(0,1]$, depending only on the fixed parameters, such that
every finite system of partial data obtained by summing \eqref{eq:newton-data} has
auxiliary solutions satisfying
\begin{align}
 \sup_{1\leq t\leq1+\tau}
 \|(f^{[m]},h^{[m]},g^{[m]})(t)\|_{\Gs{s}{\lambda_a}{\sigma+2}}
 &\leq C\eps,\qquad 1\leq m\leq {\mathcal{N}},
 \label{eq:finite-entrance-partial-bound}\\
 \sup_{1\leq t\leq1+\tau}
 \|(f^{(n)},h^{(n)},g^{(n)})(t)\|_{\Gs{s}{\lambda_a}{\sigma+2}}
 &\leq C\eps\delta_n,\qquad 1\leq n\leq {\mathcal{N}}.
 \label{eq:finite-entrance-increment-bound}
\end{align}
The increments are continuous with values in
$\mathcal G^{s,\lambda_b;\sigma+2}$, and
\begin{equation}
 \supp X^{(n)}(t)\subset
 \left\{v:\mathbf d(v)\geq\frac{5d_n}{96}\right\},
 \qquad X\in\{f,h,g\},\quad 1\leq t\leq1+\tau.
 \label{eq:finite-entrance-support}
\end{equation}
After fixing the lower bound for $M$ sufficiently large,
\begin{equation}
 \mathcal B_n(1+\tau)\leq C\eps^2\delta_n^2
              \leq4M^2\eps^2\delta_n^2,\qquad n\leq {\mathcal{N}}.
 \label{eq:finite-entrance-bootstrap}
\end{equation}
The same conclusions hold for smooth auxiliary partial data with
the same increment supports and the uniform initial bounds
\[
 \|(f_{\mathrm{in}}^{(n)},h_{\mathrm{in}}^{(n)},g_{\mathrm{in}}^{(n)})
       \|_{\Gs{s}{3\lambda_0/2}{\sigma+2}}
 \leq C\eps\delta_n,
\]
where $C$ is fixed independently of the family and the terminal layer.
\end{lemma}

\begin{proof}
Set $r=\lceil\sigma+2\rceil$, $\lambda_b=9\lambda_0/8$ and
$\lambda_a=5\lambda_0/4$, with fixed intermediate radii below
$3\lambda_0/2$. The initial radius reserve controls
$\partial_vg_{\mathrm{in}}$ in \eqref{eq:auxiliary-defect-definition}.
The product estimate therefore bounds the partial profiles $(f,h,\mathfrak r)$
by $C\eps$ and their consecutive differences by $C\eps\delta_n$.

In \eqref{eq:auxiliary-physical-shear},
$\ell_{\mathrm{aux}}=2+O(\eps)$ and $S_{\mathrm{aux}}-y=O(\eps)$
at the chosen intrinsic Gevrey radii. With $x=z+ty$, the scalar
$W(t,z,y)=\omega(t,z+ty,y)$ solves
\[
 \partial_tW+(S_{\mathrm{aux}}-y-\partial_y\Phi,
               \partial_z\Phi)\cdot\nabla_{z,y}W=0,
 \quad(\partial_z^2+(\partial_y-t\partial_z)^2)\Phi=W,
 \quad\Phi|_{\partial I_{\mathrm{aux}}}=0.
\]
The proof of Lemma~\ref{lem:short-time-sheared-elliptic} is unchanged for
$1\leq t\leq2$ and these interval lengths: the same trace cancellation
and uniform ellipticity apply. Using
\eqref{eq:auxiliary-factorial-energy} with exponent $s$ and drift
$S_{\mathrm{aux}}-y$ bounds the norm by $C\eps$ and keeps the radius at least
$e^{-C\eps\tau}$ times its initial value. The fixed composition gaps in
\eqref{eq:short-time-composition-estimate} then give the partial bound
at a radius strictly larger than $\lambda_a$ on a common interval.

Take differences on the common profile interval $I$, writing
$\delta f=f_1-f_2$, $\delta h=h_1-h_2$,
$\delta\mathfrak r=\mathfrak r_1-\mathfrak r_2$,
$\delta g=g_1-g_2$, and $\delta\phi=\phi_1-\phi_2$.
Let $\mathbf V_i$ be the full transport fields in
\eqref{eq:partial-transport-field} for these two profiles.
Subtract the two elliptic equations and apply
Lemma~\ref{lem:short-time-sheared-elliptic} at Sobolev order $r-1$.
The Gevrey product estimate and the primitive for $\delta g$ give,
below the uniform background radius,
\[
 \begin{aligned}
 \|\nabla\delta\phi\|_{\mathcal F_{s,\tau_*}^{r}}
 &\leq C\|\delta f\|_{\mathcal F_{s,\tau_*}^{r-1}}
 +C\eps\|\nabla\delta\phi\|_{\mathcal F_{s,\tau_*}^{r}}
 +C\eps\|\delta h\|_{\mathcal F_{s,\tau_*}^{r}},\\
 \|\delta g\|_{\mathcal F_{s,\tau_*}^{r}}
 +\|\mathbf V_1-\mathbf V_2\|_{\mathcal F_{s,\tau_*}^{r}}
 &\leq C\|(\delta f,\delta h,\delta\mathfrak r)
                                      \|_{\mathcal F_{s,\tau_*}^{r}}.
 \end{aligned}
\]
Absorb the first line's middle term with $C\eps\leq1/2$.
The coefficient derivative uses only
$\|\partial_v\delta h\|_{\mathcal F^{r-1}_{s,\tau_*}}
\leq\|\delta h\|_{\mathcal F^r_{s,\tau_*}}$.
Every derivative of a background field, including the one in
$\delta g\,\partial_v\mathfrak r_2$, is controlled using the fixed radius gap.
The transport difference estimate is consequently
\begin{equation}
 \begin{aligned}
 D^+\|(\delta f,\delta h,\delta\mathfrak r)
                         \|_{\mathcal F_{s,\tau_*(t)}^{r}}
 &\leq C\|(\delta f,\delta h,\delta\mathfrak r)
                         \|_{\mathcal F_{s,\tau_*(t)}^{r}}\\
 &\quad+(\dot\tau_*+C\tau_*)\partial_{\tau_*}
 \|(\delta f,\delta h,\delta\mathfrak r)
                         \|_{\mathcal F_{s,\tau_*(t)}^{r}}.
 \end{aligned}
 \label{eq:auxiliary-gevrey-difference}
\end{equation}
Choosing $\dot\tau_*=-2C\tau_*$ and shortening the common interval
gives \eqref{eq:finite-entrance-increment-bound} and continuity by
Gronwall.

The normal velocities vanish at the walls and are bounded by
$C\eps\mathbf d$. The common barrier
$d_ne^{-C\eps(t-1)}/16$, with $e^{-C\eps\tau}\geq5/6$, proves
\eqref{eq:finite-entrance-support}. Finally $\lambda_b>\lambda_0$
and $s>1/2$ absorb the Orr and polynomial factors; initial smallness
and continuity give \eqref{eq:finite-entrance-bootstrap}.
\end{proof}

\begin{lemma}[Stronger Gevrey approximation of the finite data]
\label{lem:stronger-gevrey-auxiliary-approximation}
Fix the Euler datum supported away from the walls and its finite terminal index ${\mathcal{N}}$.
There is a sequence of smooth auxiliary partial data, indexed by
$j\in\mathbb N$, at levels $1\leq m\leq {\mathcal{N}}$, whose components are
Gevrey-$\rho$, with the same $\rho>s$ as in
\eqref{eq:alpha-dagger}. Their consecutive increments have the
supports in \eqref{eq:support-scale} and satisfy the time-one bootstrap
bound \eqref{eq:bootstrap-energy-assumption} with a strict inequality uniformly in $j,{\mathcal{N}}$, and
converge to the original increments in
$\mathcal G^{s,3\lambda_0/2;\sigma+2}$. Neither the intermediate
regularizations nor their terminal regularization are required to be
Euler-compatible.
\end{lemma}

\begin{proof}
Smooth the zero extensions of the three time-one profiles
$(f_{\mathrm{in}},h_{\mathrm{in}},g_{\mathrm{in}})$ by
$e^{-\nu(k^2+\eta^2)}$ (with $k=0$ for $h,g$), and multiply by the
same cumulative cutoffs $\chi_m$, $1\leq m\leq {\mathcal{N}}$. Take consecutive
differences. The cutoffs preserve the support margins, while their
Gevrey-$\rho$ regularity and Gaussian smoothing give the stronger
Gevrey regularity. For fixed ${\mathcal{N}}$, multiplication by these finitely many
cutoffs is continuous across the available initial radius reserve,
so the increments converge in the stated norm. Choose
$\nu_j\downarrow0$ so that each of the finitely many increment
errors is at most $2^{-j}\eps\delta_m$ in that norm. The strict
time-one bounds are then uniform; stronger-Gevrey radii may depend
on $j,{\mathcal{N}}$.
\end{proof}

\subsubsection{Global bootstrap closure}
\begin{proof}[Proof of Proposition~\ref{prop:bootstrap-improvement}]
Fix ${\mathcal{N}}$ and approximation $j$ from
Lemma~\ref{lem:stronger-gevrey-auxiliary-approximation}, and let
$\tau_{{\mathcal{N}},j}$ be the minimum of the maximal Sobolev lifespans of the corresponding exact
profiles $(f_j^{[m]},h_j^{[m]},g_j^{[m]})$, $m\leq {\mathcal{N}}$.
Lemma~\ref{prop:newton-equivalence} gives the increment equations, and
Lemma~\ref{lem:finite-profile-gevrey-entrance} starts the bootstrap
uniformly. We first close the bootstrap for these approximate data, then take $j\to\infty$ with ${\mathcal{N}}$ fixed.

We first verify that the norms in the absorption estimates are finite. Fix $T<\tau_{{\mathcal{N}},j}$.  The regularized data satisfy
$Y_{\mathrm{in}}'=(1+h_{\mathrm{in}})^{-1}>0$; Gevrey inversion therefore
gives an intrinsically Gevrey-$\rho$ initial coordinate $V_{\mathrm{in}}$.
The physical initial vorticity and prescribed shear are then Gevrey-$\rho$
by \eqref{eq:auxiliary-physical-equation} and
\eqref{eq:auxiliary-physical-shear}.
By Lemma~\ref{lem:stronger-gevrey-channel-persistence}, the finitely many
profile fields have a common positive Gevrey-$\rho$ radius
$\mu_{{\mathcal{N}},j,T}$ on $[1,T]$.  Since $\rho>s>1/2$, for every fixed
$b,C,p>0$,
\[
 \sup_\xi\langle\xi\rangle^p
 e^{b|\xi|^s+C|\xi|^{1/2}-\mu_{{\mathcal{N}},j,T}|\xi|^\rho}<\infty.
\]
The Orr corrections are bounded by $Ce^{C|(k,\eta)|^{1/2}}$ and their
logarithmic weight derivatives lie in $[0,8]$. Dominated convergence
therefore gives continuous energies and locally integrable $\CK$ terms; their suprema and integrals are continuous
in the upper endpoint.
Decrease $\mu_{{\mathcal{N}},j,T}$ to include the extensions and their complete sources in
\eqref{eq:exterior-extension-definition}--\eqref{eq:interior-extension-definition}. Recall $\mathscr G_m^*$,
$*\in\left\{\Ext,\Int\right\}$, from
\eqref{eq:exterior-green-operator}--\eqref{eq:interior-green-operator},
and $\lambda_m^\Delta$ from \eqref{eq:kernel-localization-radius}.
For $0<\nu<\min(\mu_{{\mathcal{N}},j,T},\allowbreak\lambda_m^\Delta)$,
concavity and Young's inequality applied to
\eqref{eq:localized-green-kernel}--\eqref{eq:localized-green-kernel-delta}
give, for $p\geq0$,
\begin{align*}
 &\|\mathscr G_m^*q\|_{\Gs{\rho}{\nu}{p}}
  +\|\Delta_L\mathscr G_m^*q\|_{\Gs{\rho}{\nu}{p}}\\
 &\quad\leq C_p2^{2m}
  \big\|(1+\xi^2)\langle\xi\rangle^p
       e^{-(\lambda_m^\Delta-\nu)|\xi|^\rho}\big\|_{L^1_\xi}
       \|q\|_{\Gs{\rho}{\nu}{p}}<\infty.
\end{align*}
Here $k\neq0$ gives $k^2+|\xi-kt|^2\geq1$, and the zero-mode component vanishes.
Applying this bound to \eqref{eq:exterior-extension-definition}--\eqref{eq:interior-extension-definition},
using the endpoint identity \eqref{eq:exact-extension-endpoint}, proves finiteness
of the absorption norms; a further radius loss gives continuity.
The radii may depend on ${\mathcal{N}},j,T$. No zero extension of the Dirichlet stream
is used, and channel agreement is invoked only on the regions where the cutoffs equal one.
Let $[1,T_{{\mathcal{N}},j}^B)$ be the maximal bootstrap interval.  On every
$[1,T]$ with $T<T_{{\mathcal{N}},j}^B$, the component estimates and their sum have
constants independent of ${\mathcal{N}},j,T$.  For
$b_j^{[m]}=(1+h_j^{[m]})\partial_z\phi_{j,\neq}^{[m]}+g_j^{[m]}$,
the mean-value and Sobolev argument of
Lemma~\ref{lem:support-improvement} gives
\begin{equation}
 \left\|b_j^{[m]}(t)/\mathbf d\right\|_\infty
 +\left\|g_j^{[m]}(t)/\mathbf d\right\|_\infty
 \lesssim M^C\eps(\langle t\rangle^{-2}
                         +\langle t\rangle^{-1-s}),\qquad m\leq {\mathcal{N}}.
 \label{eq:finite-normalized-drift}
\end{equation}
The moving-barrier argument in Lemma~\ref{lem:support-improvement},
applied on $[1,T]$ to these same drifts, gives, after decreasing $\eps$,
\begin{equation}
 \operatorname{supp}X_j^{(m)}(t)\subset\left\{\mathbf d\geq5d_m/88\right\},
 \qquad X\in\left\{f,g,h\right\},\quad m\leq {\mathcal{N}}.
 \label{eq:finite-global-support-improvement}
\end{equation}
Together with \eqref{eq:energy-bootstrap-conclusion}, this strictly improves
both bootstraps, so continuity gives $T_{{\mathcal{N}},j}^B=\tau_{{\mathcal{N}},j}$.
The positive uniform $s$-radius bounds $f_j^{[m]},h_j^{[m]}$
in every fixed Sobolev norm and $g_j^{[m]}$ one order higher. On a finite
interval the defect therefore satisfies, for each fixed $r$,
\[
 \|\mathfrak r_j^{[m]}(t)\|_{H^r}
 \lesssim t(1+\|h_j^{[m]}(t)\|_{H^r})\|g_j^{[m]}(t)\|_{H^{r+1}}
       +\|f_j^{[m]}(t)\|_{H^r}+\|h_j^{[m]}(t)\|_{H^r}<\infty.
\]
Together with $1+h_j^{[m]}\geq1/2$, these are exactly the continuation
hypotheses of Proposition~\ref{prop:finite-profile-lwp}. They exclude
a finite $\tau_{\mathcal N,j}$. Hence
\[
 \mathcal B_{m,j}(\infty)\leq2M^2\eps^2\delta_m^2,
 \qquad m\leq {\mathcal{N}},
\]
uniformly in $j,{\mathcal{N}}$.

At fixed ${\mathcal{N}}$, the uniform positive $s$-radius bounds every Sobolev norm,
and the equations bound time derivatives one order lower. Compactness therefore
gives strong convergence in $C([1,T];H^r)$ for every fixed $r,T$ as
$j\to\infty$, after extraction. Continuity of the Dirichlet inverse
preserves the auxiliary equations; uniqueness identifies the
limits with the original finite systems, consistently in $T$ at fixed ${\mathcal{N}}$.
Pass first to bounded sets of Fourier frequencies in each nonnegative energy
and CK integral, then use monotone convergence and lower semicontinuity
of the time supremum. Strong $L^2$ convergence preserves support. Thus
\begin{equation}
 T_m^{\max}=\infty\quad(1\leq m\leq {\mathcal{N}}),\qquad
 \mathcal B_m(\infty)\leq2M^2\eps^2\delta_m^2.
 \label{eq:finite-global-bound}
\end{equation}
These bounds are independent of the finite terminal layer ${\mathcal{N}}$.

To recover the precision bounds for the original finite systems, use
the radii $\theta\lambda_n$ and derivatives $\theta\dot\lambda_n$,
$1/2\leq\theta<1$, keeping the Orr weights unchanged. The lossy bounds
make the reduced precision norms finite, and the proof of
Proposition~\ref{prop:precision-elliptic} applies uniformly in $\theta$
with the original CK densities on the right, since every reduced CK
integrand is bounded by its original counterpart. These densities are
integrable by \eqref{eq:finite-global-bound}; letting $\theta\uparrow1$
and using monotone convergence gives the precision bounds almost
everywhere in time and their integrated versions.
\end{proof}

\subsection{Coordinate reconstruction}
\begin{lemma}[Coordinate reconstruction and gauge compatibility]
\label{lem:coordinate-reconstruction}
Let $(f,h,g,\phi)$ be the terminal profile from
Lemma~\ref{lem:finite-terminal-bounds}, whose time-one datum is exactly
the physical datum supplied by Proposition~\ref{prop:short-time-weighted-propagation}.
Then there is a unique inverse coordinate $Y(t,v)$ such that, with
$v=V(t,y)=Y(t,\cdot)^{-1}(y)$,
\begin{equation}
 V(t,y)=y+\frac1t\int_0^t\avg{U^x}(\tau,y)\dd\tau,
 \label{eq:reconstructed-coordinate}
\end{equation}
and the pullbacks
\begin{equation*}
 \omega(t,x,y)=f(t,x-tV(t,y),V(t,y)),
 \qquad
 \psi(t,x,y)=\phi(t,x-tV(t,y),V(t,y))
\end{equation*}
solve the Euler system with the fixed zero-flux Dirichlet normalization.

\end{lemma}

\begin{proof}
We reconstruct the inverse coordinate by solving the transport equation
\begin{equation}
 \partial_tY+g\partial_vY=0,
 \qquad Y(1,v(1,y))=y,
 \label{eq:inverse-coordinate-transport}
\end{equation}
where $v(1,y)$ is the time-one coordinate from
\eqref{eq:coordinates}.  Differentiating the transport equation gives
$(\partial_t+g\partial_v)((\partial_vY)^{-1})
=(\partial_vY)^{-1}\partial_vg$, the same equation and initial
value as for $1+h$.  Hence
\begin{equation}
 (\partial_vY)^{-1}=1+h.\label{eq:reconstructed-jacobian}
\end{equation}

For the terminal profile fields in \eqref{eq:profile-system}, put
\begin{equation}
 U_0^x=-(1+h)\partial_v\phi_0.
 \label{eq:reconstructed-qQ}
\end{equation}

The zero-mode equations give $f_0=-(1+h)\partial_v U_0^x$ and
$(\partial_t+g\partial_v)f_0=-(1+h)\partial_v\Pzero(f_{\neq}\partial_z\phi_{\neq})$; the latter follows by periodic integration
by parts in the $f$ equation.  Using
$[(\partial_t+g\partial_v),\partial_v]=-\partial_vg\partial_v$ yields

\begin{equation}
 \partial_v((\partial_t+g\partial_v)U_0^x-\Pzero(f_{\neq}\partial_z\phi_{\neq}))=0.
 \label{eq:mean-momentum-up-to-constant}
\end{equation}
Since $\partial_t\partial_vY+\partial_v(g\partial_vY)=0$ and $g=0$ at both walls,
the constant in \eqref{eq:mean-momentum-up-to-constant} has
$(\partial_vY)$-weighted integral
$\partial_t\int(\partial_vY)U_0^x-\int(\partial_vY)\Pzero(f_{\neq}\partial_z\phi_{\neq})=0$:
the first integral vanishes by Dirichlet normalization, and the second
by integration by parts: after the change $(x,y)=(z+tv,Y(t,v))$,
\[
 \int_L^R(\partial_vY)\Pzero(f_{\neq}\partial_z\phi_{\neq})\,\dd v
 =\frac1{2\pi}\int_\Omega(\Delta\psi)\partial_x\psi\,\dd x\dd y=0.
\]
The horizontal integral vanishes by periodicity, and $\partial_x\psi=0$
on the walls makes the boundary term vanish. As $\partial_vY>0$, this proves
\begin{equation}
 (\partial_t+g\partial_v)U_0^x=\Pzero(f_{\neq}\partial_z\phi_{\neq}).\label{eq:mean-momentum-identity}
\end{equation}
The $g$ equation also gives $(\partial_t+g\partial_v)(v-Y+tg)=\Pzero(f_{\neq}\partial_z\phi_{\neq})$.
The two quantities agree at $t=1$ by \eqref{eq:short-time-g-data};
transport uniqueness proves
\begin{equation}
 U_0^x=v-Y+tg.\label{eq:coordinate-gauge-identity}
\end{equation}

Inverting $Y$ gives $\partial_tV=g(t,V)$, $\partial_yV=1+h(t,V)$; hence
$\partial_t[t(V-y)]=\avg{U^x}$, and the time-one normalization yields
\eqref{eq:reconstructed-coordinate}. The chain rule and
\eqref{eq:reconstructed-jacobian} recover Euler and $\Delta\psi=\omega$,
with unchanged Dirichlet trace and zero-flux normalization.

\end{proof}

For each finite terminal sum, \eqref{eq:global-support-bound} separates
the increments from the walls.  Reducing $\eps$,
\eqref{eq:bootstrap-consequences-partial-fh} gives $\|h^{[{\mathcal{N}}]}\|_\infty\leq1/2$,
so the terminal coordinate remains a diffeomorphism.  The next estimates
refer to $(f,h,g,\phi)=(f^{[{\mathcal{N}}]},h^{[{\mathcal{N}}]},g^{[{\mathcal{N}}]},\phi^{[{\mathcal{N}}]})$ and its
reconstructed Euler solution; all constants are independent of ${\mathcal{N}}$.
\subsection{Boundary weight, scattering, and damping}

\begin{lemma}[Recovery of the weighted boundary estimate]
\label{lem:recover-boundary-weight}
There is $c_\infty>0$, depending only on $c_*$ and the fixed parameters, such
that
\begin{equation}
  \sup_{t\geq1}\norm{e^{c_\infty \mathbf d(v)^{-\alpha_*}}f(t)}_{L^2(\T\times(L,R))}
  \lesssim\eps.
  \label{eq:profile-boundary-weight}
\end{equation}
Consequently, after decreasing $c_\infty$, the same bound holds for
$e^{c_\infty d(y)^{-\alpha_*}}\omega(t,x,y)$.
\end{lemma}

\begin{proof}
We combine the support separation with the decay of the increments. On the support of the $n$th increment,
\eqref{eq:global-support-bound} and \eqref{eq:support-scale} imply
$\mathbf d^{-\alpha_*}\leq C_12^{n\alpha_*}$, with a fixed $C_1$.  Hence
\begin{align*}
  \norm{e^{c_\infty \mathbf d^{-\alpha_*}}f(t)}_2
  &\leq\sum_{n=1}^{\mathcal{N}}
    e^{c_\infty C_12^{n\alpha_*}}\norm{f^{(n)}(t)}_2 \lesssim M\eps\sum_{n\geq1}
    \exp\bigl(-(c_*-c_\infty C_1)2^{n\alpha_*}\bigr)
  \lesssim\eps
\end{align*}
provided $c_\infty C_1<c_*/2$. Changing variables
proves the assertion in physical variables.
\end{proof}

\begin{lemma}[Scattering and velocity damping]
\label{lem:scattering-damping}
There is $f_\infty\in\Gs{s}{\lambda_\infty}{\sigma-6}$ such that
\begin{equation}
  \norm{f(t)-f_\infty}_{\Gs{s}{\lambda_\infty}{\sigma-6}}
  \longrightarrow0.
  \label{eq:profile-scattering-conclusion}
\end{equation}
Furthermore,
\begin{equation}
  \norm{\Pneq U^x(t)}_{L^2(\Omega)}
  \lesssim\frac{\eps}{\avg t},
  \qquad
  \norm{U^y(t)}_{L^2(\Omega)}
  \lesssim\frac{\eps}{\avg t^2}.
  \label{eq:physical-damping-conclusion}
\end{equation}
\end{lemma}

\begin{proof}
The profile equation, the Gevrey product estimate,
\eqref{eq:summed-profile-transport}, and
\eqref{eq:bootstrap-consequences-partial-g} imply
\begin{equation}
  \norm{\partial_tf(t)}_{\Gs{s}{\lambda_\infty}{\sigma-6}}
  \lesssim M^C\eps^2
    \left(\frac1{\avg t^2}+\frac1{\avg t^{\beta_g/2}}\right).
  \label{eq:integrable-profile-derivative}
\end{equation}
The six available Sobolev orders control the derivative and product in the profile equation.
Since $\beta_g/2=2-a/4>1$, integration gives convergence to $f_\infty$.

In physical variables,
\begin{equation}
  U^x=-(1+h)\pvt\phi,
  \qquad U^y=\partial_z\phi.
  \label{eq:physical-profile-velocity-relation}
\end{equation}
The Jacobian of the transformation from $(x,y)$ to $(z,v)$ is uniformly
comparable to one.  Therefore \eqref{eq:summed-physical-velocity}, the
product estimate, and \eqref{eq:physical-profile-velocity-relation} yield
\eqref{eq:physical-damping-conclusion}.
\end{proof}

\subsection{Approximation by data supported away from the walls}
\label{sec:outer-approximation}

We now remove the assumption that the initial support lies a positive distance from the
walls.  The approximation in this subsection is made in the physical
channel $\Omega=\mathbb T_{2\pi}\times[-1,1]$.  Its index $N$ labels
different physical initial data; the corresponding total number of layers is
$\mathcal{N}=N+J_s$.
The terminal index $\mathcal N$ is chosen so that $\chi_{\mathcal N}=1$
on the supports of all three time-one profiles. The approximation index $N$ and the fixed shift $J_s$
are specified in \eqref{eq:outer-cutoff-data} and
\eqref{eq:outer-terminal-index}, respectively; $N$ tends to infinity,
whereas $J_s$ is chosen uniformly for the family.

\begin{lemma}[Passage from finite support to boundary-flat data]
\label{lem:outer-compact-approximation}
The uniform finite-terminal estimates extend to every datum satisfying
\eqref{eq:data-assumption}.  More precisely, after fixing a smaller
positive output radius, the compactly supported approximations below
converge on every finite time interval to the Euler solution with the
original datum.  The uniform boundary weight, Gevrey bound, damping,
and scattering estimates pass to this solution, with constants and a
smallness threshold depending only on the parameters in
Theorem~\ref{thm:main}.
\end{lemma}

\begin{proof}

We reserve a fixed part of the initial Gevrey radius for the
physical cutoff and work at the Gevrey radius
\begin{equation}
 \bar\lambda_0=\frac34\lambda_0.
 \label{eq:outer-working-radius}
\end{equation}
With the Gevrey-$\rho$ cutoff $\gamma$ used in
\eqref{eq:exact-shell-partition},
where $0\leq\gamma\leq1$, $\gamma=0$ on $(-\infty,1]$ and
$\gamma=1$ on $[2,\infty)$, define, for $N\geq2$,
\begin{equation}
 \chi_N^{\rm phys}(y)=\gamma(2^Nd(y)),\qquad
 \omega_{{\rm in},N}=\chi_N^{\rm phys}\omega_{\rm in}.
 \label{eq:outer-cutoff-data}
\end{equation}

Thus we choose initial data $\omega_{{\rm in},N}$ supported in
$\{d\geq2^{-N}\}$ and equal to $\omega_{\rm in}$ on
$\{d\geq2^{1-N}\}$. By \eqref{eq:data-assumption}, the discarded collar satisfies
$\|(1-\chi_N^{\rm phys})\omega_{\rm in}\|_2
\leq\eps e^{-c_0 2^{(N-1)\alpha_*}}$.
Lemma~\ref{lem:localization} bounds the Gevrey cost of each cutoff
shell. Interpolating these shell bounds across the fixed radius gap
and then summing the tails gives the quantitative approximation estimates
\begin{equation}
 \|\omega_{{\rm in},N}-\omega_{\rm in}\|_{
                       \Gs{s}{3\bar\lambda_0}{2\sigma+8}}
 \leq C\eps e^{-c2^{N\alpha_*}},\qquad
 \|\omega_{{\rm in},N}\|_{\Gs{s}{3\bar\lambda_0}{2\sigma+8}}
 +\|e^{c_0d^{-\alpha_*}}\omega_{{\rm in},N}\|_2\leq C\eps.
 \label{eq:outer-data-convergence}
\end{equation}
Thus the finite-support argument applies at base radius $\bar\lambda_0$
with smallness parameter $C\eps$ and a threshold independent of $N$.

The bound $\|U_N^y/(1-y^2)\|_\infty\lesssim\eps$ on $0\leq t\leq1$,
from the proof of Proposition~\ref{prop:short-time-weighted-propagation}, preserves 
$\{d<e^{-C\eps}2^{-N}\}$ through $t=1$. There
$\langle U_N^x\rangle=c_{\pm,N}$ and $v_N(1,y)=y+c_{\pm,N}$;
\eqref{eq:short-time-h-data}--\eqref{eq:short-time-g-data} show that
all three time-one profiles vanish for $\dist(v,\partial I_N)<e^{-C\eps}2^{-N}$.
Fix $J_s\geq0$ with

\begin{equation}
 2^{-J_s-3}<e^{-C\eps_0},\qquad \mathcal N=N+J_s.
 \label{eq:outer-terminal-index}
\end{equation}
Here $\eps_0$ is the data-smallness threshold in
Theorem~\ref{thm:main} (see \eqref{eq:main-parameter-choices}), decreased once so that the finite-support
argument applies at base radius $\bar\lambda_0$ to data of size
$C\eps_0$ under Proposition~\ref{prop:bootstrap-improvement}.
Then $\chi_{\mathcal N}$ equals one on
their supports, so the terminal auxiliary datum is exactly the
physical time-one datum; lower auxiliary partial sums solve the
profile equations. The profile intervals
$I_N=[-1+c_{-,N},1+c_{+,N}]$ satisfy
$|c_{\pm,N}|\leq C\eps$ and $1\leq|I_N|\leq3$.
Translation and uniformly bounded dilation give common cutoff and
Dirichlet constants.

On the fixed physical channel, set
\begin{equation}
 \Phi_N(t,y)=\int_0^t\langle U_N^x\rangle(\tau,y)\,d\tau,
 \qquad W_N(t,x,y)=\omega_N(t,x+ty+\Phi_N(t,y),y).
 \label{eq:outer-physical-scattering-profile}
\end{equation}
Choose $\lambda_\infty>0$, independently of $N$, small enough that the
finite-data estimates and composition with the terminal coordinate
hold at $2\lambda_\infty$, with the reserve in
\eqref{eq:final-gevrey-radius}. For the terminal profiles,
$W_N(t,x,y)=f_N(t,x,V_N(t,y))$ and
$\partial_tW_N=(\partial_tf_N+g_N\partial_vf_N)(t,x,V_N(t,y))$.
Lemma~\ref{lem:recover-boundary-weight},
\eqref{eq:summed-profile-transport}--\eqref{eq:summed-physical-velocity},
\eqref{eq:bootstrap-consequences-partial-g}, and the fixed-radius
composition estimate give common $c_\infty>0$ and
\begin{align}
 \sup_{t\geq0}\|e^{c_\infty d^{-\alpha_*}}\omega_N(t)\|_2
 +\sup_{t\geq1}\|W_N(t)\|_{\Gs{s}{2\lambda_\infty}{\sigma}}
 &\leq C\eps,
 \label{eq:outer-uniform-physical-bounds}\\
 \|\Pneq U_N^x(t)\|_2+\langle t\rangle\|U_N^y(t)\|_2
 &\leq\frac{C\eps}{\langle t\rangle},\qquad t\geq0,
 \label{eq:outer-uniform-damping}\\
 \|\partial_tW_N(t)\|_{\Gs{s}{2\lambda_\infty}{\sigma-6}}
 &\leq C\eps^2
       \bigl(\langle t\rangle^{-2}+\langle t\rangle^{-\beta_g/2}\bigr),
       \qquad t\geq1.
 \label{eq:outer-uniform-time-tail}
\end{align}
The short-time bounds cover $0\leq t\leq1$; all constants are
independent of $N$.

Sobolev compactness and uniqueness give
$\omega_N\to\omega$ in $C([0,T];H^r(\Omega))$ for every fixed $r,T$.
Writing $\Phi=\int_0^t\langle U^x\rangle\,d\tau$ and
$W(t,x,y)=\omega(t,x+ty+\Phi(t,y),y)$, the strict radius gap gives
$W_N\to W$ in $C([1,T];\Gs{s}{\lambda_\infty}{\sigma})$.
Truncation of the exponential boundary weight, followed by monotone
convergence, passes its uniform bound to the limit; the Gevrey and damping bounds
pass by the same strong convergence.
Integrate \eqref{eq:outer-uniform-time-tail} before taking $N\to\infty$.
Since $\beta_g/2>1$, the resulting Cauchy estimate gives a limit
$\omega_\infty$ with
\begin{equation}
 \|W(t)-\omega_\infty\|_{\Gs{s}{\lambda_\infty}{\sigma-6}}
 \leq C\eps^2
       \bigl(\langle t\rangle^{-1}
                    +\langle t\rangle^{1-\beta_g/2}\bigr).
 \label{eq:outer-scattering-tail}
\end{equation}
Finally, $\Pneq\omega=\partial_xU^y-\partial_y\Pneq U^x$ and the
damping bound give weak convergence to zero against smooth test functions supported
in the interior; density and the uniform $L^2$ bound give weak convergence in $L^2(\Omega)$.
\end{proof}

\begin{proof}[Proof of Theorem~\ref{thm:main}] \space First suppose the datum is supported away from the walls.
Lemmas~\ref{lem:finite-terminal-bounds} and
\ref{lem:coordinate-reconstruction}, with Sobolev uniqueness, give the
global smooth Euler solution.  Proposition~\ref{prop:short-time-weighted-propagation}
and Lemma~\ref{lem:recover-boundary-weight} prove \eqref{eq:main-boundary}
after decreasing $c_\infty$ so that $0<c_\infty<c_0$
as in \eqref{eq:main-parameter-choices} and $c_\infty\leq c_1$.
Lemma~\ref{lem:scattering-damping} gives \eqref{eq:main-damping} for $t\geq1$;
on $[0,1]$, \eqref{eq:short-time-gevrey-conclusion} and Dirichlet
Biot--Savart give $\sup_t\|U(t)\|_2\lesssim\eps$.

Lemma~\ref{lem:outer-compact-approximation} extends these conclusions
to the prescribed boundary-flat datum. Its uniform physical bound
\eqref{eq:outer-uniform-physical-bounds} and time-tail estimate
\eqref{eq:outer-scattering-tail} give \eqref{eq:main-scattering};
the curl-identity argument at the end of its proof gives
$\Pneq\omega(t)\rightharpoonup0$ in $L^2(\Omega)$.
\end{proof}

\appendix

\section{Short-time propagation of the boundary weight}
\label{sec:short-time}

We prove Proposition~\ref{prop:short-time-weighted-propagation}: the boundary weight
persists until $t=1$, and the profile data satisfy the bounds needed
in Section~\ref{sec:decomposition}.

The linear Couette flow shears the Fourier support of $\omega$.
We therefore use
\begin{equation}
  F(t,z,y)=\omega(t,z+ty,y),
  \qquad
  \Phi(t,z,y)=\psi(t,z+ty,y),
  \qquad 0\leq t\leq1.
  \label{eq:short-time-sheared-unknowns}
\end{equation}
Under $x=z+ty$, the physical derivatives become
$\partial_x=\partial_z$ and $\partial_y|_x=\partial_y|_z-t\partial_z$.
Thus the transformed transport and Dirichlet Poisson equations are
\begin{equation}
 \partial_tF+\mathcal U\cdot\nabla_{z,y}F=0,
 \qquad
 \mathcal U=(-\partial_y\Phi,\partial_z\Phi),
 \qquad
 (\partial_z^2+(\partial_y-t\partial_z)^2)\Phi=F,
 \qquad \Phi|_{y=\pm1}=0.
 \label{eq:short-time-sheared-system}
\end{equation}
Thus $\nabla_{z,y}\cdot\mathcal U=0$ and
$\mathcal U^y|_{y=\pm1}=0$.

We first establish the boundary elliptic and Gevrey estimates used below.
Here $\Omega=\mathbb T_{2\pi}\times[-1,1]$ as in \eqref{eq:euler},
$\alpha=(\alpha_1,\alpha_2)$ is a derivative multi-index, and
$|\alpha|=\alpha_1+\alpha_2$.  For a smooth scalar function $Q$ on the closed strip, put
\[
 \norm{\nabla^m Q}_{H^r,\mathrm{sum}}
 =\sum_{|\alpha|=m}\binom m\alpha
   \norm{\partial^\alpha Q}_{H^r(\Omega)},
 \qquad \binom m\alpha=\frac{m!}{\alpha_1!\alpha_2!}.
\]
For a vector or derivative tensor, we also sum these norms over
its components; equivalently, derivative indices are ordered and summed
in $\ell^1$. In one dimension, the degreewise norm consists of a single term.
For $J\in\mathbb N_0$ define
\begin{equation}
 \norm Q_{\mathfrak G_J^{s,\lambda;r}}^2
 =\sum_{m=0}^J\frac{(s\lambda)^{2m/s}}{(m!)^{2/s}}
   \norm{\nabla^m Q}_{H^r,\mathrm{sum}}^2,
 \qquad
 \norm Q_{\mathfrak G^{s,\lambda;r}}
 =\sup_J\norm Q_{\mathfrak G_J^{s,\lambda;r}}.
 \label{eq:factorial-Gevrey-definition}
\end{equation}
This choice agrees, across any strict radius gap, with the Fourier
frequency convention $|k,\eta|=|k|+|\eta|$. To see this without changing
the radius, let $Q$ be boundary-flat and write $\widetilde Q$ for its zero
extension. With $\xi=(k,\eta)$ and $|\xi|=|k|+|\eta|$, the binomial
identity, the triangle inequality, and Cauchy--Schwarz give
\[
 \bigl\||\xi|^m\widehat{\widetilde Q}\bigr\|_2
 \leq\norm{\nabla^mQ}_{L^2,\mathrm{sum}}
 \leq\sqrt{m+1}\,
       \bigl\||\xi|^m\widehat{\widetilde Q}\bigr\|_2.
\]
Indeed, $\sum_{|\alpha|=m}\binom m\alpha|\xi^\alpha|=|\xi|^m$,
and the sum of the squares of these nonnegative terms is at most the
square of their sum. Stirling's formula, with maximizing index
$m=s\lambda u^s+O(1)$, yields, for $u\geq0$,
\begin{equation}
 \begin{split}
 \sum_{m\geq0}\frac{(s\lambda)^{2m/s}u^{2m}}{(m!)^{2/s}}
 &\approx_{s,\lambda}e^{2\lambda u^s}(1+u)^{s/2-1},\\
 \sum_{m\geq0}(m+1)
       \frac{(s\lambda)^{2m/s}u^{2m}}{(m!)^{2/s}}
 &\lesssim_{s,\lambda}
       e^{2\lambda u^s}(1+u)^{3s/2-1}.
 \end{split}
 \label{eq:factorial-Fourier-calibration}
\end{equation}
Thus, for $0<\lambda_-<\lambda<\lambda_+$ and fixed real $r\geq0$,
\begin{equation}
 \norm Q_{\Gs{s}{\lambda_-}{r}}
 \leq C\norm Q_{\mathfrak G^{s,\lambda;r}}
 \leq C\norm Q_{\Gs{s}{\lambda_+}{r}}.
 \label{eq:factorial-Fourier-equivalence}
\end{equation}
The strict gaps absorb the displayed polynomial factors. We make no claim about the Fourier norm
of the zero extension of an arbitrary strip function;
coordinate maps and Dirichlet gradients are measured in intrinsic strip norms.

We also recover fixed Sobolev orders across a strict radius gap.
For $0<\lambda_-<\lambda$ and each fixed integer $\ell\geq0$,
\[
 \sup_{m\geq0}
 \left(\frac{\lambda_-}{\lambda}\right)^{m/s}
 \left(\frac{(m+\ell)!}{m!}\right)^{1/s}<\infty.
\]
Hence
\[
 \frac{(s\lambda_-)^{m/s}}{(m!)^{1/s}}
     \norm{\nabla^{m+\ell}Q}_{L^2,\mathrm{sum}}
 \leq C_{\ell,s,\lambda_-,\lambda}
 \frac{(s\lambda)^{(m+\ell)/s}}{((m+\ell)!)^{1/s}}
     \norm{\nabla^{m+\ell}Q}_{L^2,\mathrm{sum}}.
\]
Since
\[
 \norm{\nabla^mQ}_{H^r,\mathrm{sum}}
 \leq C_r\sum_{\ell=0}^{\lceil r\rceil}
             \norm{\nabla^{m+\ell}Q}_{L^2,\mathrm{sum}},
\]
squaring and summing gives the Sobolev-order shift.
At cutoff $J$ this uses input degrees through $J+\lceil r\rceil$;
we include these degrees or apply the shift after monotone passage to the full norm.
For flat scalars, zero extension commutes with integer derivatives, so
Plancherel and this shift prove the first inequality in
\eqref{eq:factorial-Fourier-equivalence}. Restriction of each full-space
$H^r$ derivative norm proves the reverse inequality. Constants depend
on the fixed gaps and orders, not on the cutoff.

\begin{lemma}[Short-time elliptic estimate on the strip]
\label{lem:short-time-sheared-elliptic}
Let $0\leq t\leq1$, $r\in\mathbb N_0$, $J\in\mathbb N_0$, and $\lambda>0$.
If $Q$ is boundary-flat and $\Phi$ solves
\begin{equation}
  (\partial_z^2+(\partial_y-t\partial_z)^2)\Phi=Q,
  \qquad \Phi|_{y=\pm1}=0,
  \label{eq:short-time-flat-elliptic-problem}
\end{equation}
then, uniformly in $t$ and $J$,
\begin{equation}
 \norm{\nabla_{z,y}\Phi}_{\mathfrak G_J^{s,\lambda;r+1}}
 +\norm{\nabla_{z,y}^2\Phi}_{\mathfrak G_J^{s,\lambda;r}}
 \leq C_r\norm{Q}_{\mathfrak G_J^{s,\lambda;r}}.
 \label{eq:short-time-factorial-elliptic-estimate}
\end{equation}
Consequently, for $0<\lambda'<\lambda$ and real $r\geq0$,
\begin{equation}
 \norm{\nabla_{z,y}\Phi}_{\mathfrak G^{s,\lambda';r+1}}
 +\norm{\nabla_{z,y}^2\Phi}_{\mathfrak G^{s,\lambda';r}}
 \leq C_{r,\lambda-\lambda'}\norm{Q}_{\Gs{s}{\lambda}{r}}.
 \label{eq:short-time-sheared-elliptic-estimate}
\end{equation}
\end{lemma}

\begin{proof}[Proof sketch]

Boundary flatness and the Dirichlet condition give, at $y=\pm1$,
\[
 (\partial_y-t\partial_z)^{2\ell}\Phi=0,\qquad
 (\partial_y-t\partial_z)^{2\ell+1}\Phi
   =(-1)^\ell\partial_z^{2\ell}\partial_y\Phi.
\]
Expand $\partial_y=(\partial_y-t\partial_z)+t\partial_z$ in the traces of
$\partial_z^p\partial_y^q\Phi$. Testing the commuted equation against
$\partial_z^2\partial_z^p\partial_y^q\Phi$ produces a real multiple of
\[
 \operatorname{Re}\int_{\mathbb T_{2\pi}}
 \partial_z^{p+q}\partial_y\Phi\,
 \overline{\partial_z^{p+q+1}\partial_y\Phi}\,\dd z=0
 \qquad (y=\pm1).
\]
For $q=0$, the Dirichlet trace vanishes directly.
Horizontal Poincar\'e controls the nonzero modes. For the zero mode,
$p>0$ gives zero; $p=q=0$ uses the Dirichlet estimate; and $p=0,q\geq1$
gives $\partial_y^{q+1}\Phi=0$ at the walls, permitting vertical Poincar\'e.
Thus, uniformly in $p,q$ and $0\leq t\leq1$,
\[
 \|\nabla\partial_z^p\partial_y^q\Phi\|_2
 +\|\nabla^2\partial_z^p\partial_y^q\Phi\|_2
 \leq C\|\partial_z^p\partial_y^qQ\|_2.
\]
Sum the ordered derivatives with their factorial weights, including the
fixed integer Sobolev orders. Monotone convergence and the Sobolev-order shift across the
strict radius gap give the real-order estimate. Only the flat source is zero-extended.
\end{proof}

\begin{lemma}[Short-time Gevrey propagation]
\label{lem:short-time-gevrey-propagation}
Let $r>4$ and let $F$ solve
\eqref{eq:short-time-sheared-system} with boundary-flat initial datum.
There is $\eps_1=\eps_1(s,\lceil r\rceil,\lambda_0)>0$ such that
\begin{equation}
  \norm{F(0)}_{\Gs{s}{3\lambda_0}{r}}\leq\eps\leq\eps_1
  \quad\Longrightarrow\quad
  \sup_{0\leq t\leq1}
  \norm{F(t)}_{\Gs{s}{(5/2)\lambda_0}{r}}
  \leq C\eps.
  \label{eq:short-time-sheared-gevrey-bound}
\end{equation}
All normal and tangential traces of $F$ continue to vanish at the walls.
\end{lemma}

\begin{proof}[Proof sketch]

A Sobolev Friedrichs construction gives the solution and preserves flatness.
Fix $J\in\mathbb N_0$ and a positive differentiable radius $\lambda_J(t)$.
Transport cancellation, Lemma~\ref{lem:short-time-sheared-elliptic}, and
the normalized Leibniz convolution give
\begin{equation}
\begin{aligned}
 \tfrac12\partial_t\norm F_{\mathfrak G_J^{s,\lambda_J;\lceil r\rceil}}^2
 &\leq\left(\frac{\dot\lambda_J}{s\lambda_J}
      +C\norm F_{\mathfrak G_J^{s,\lambda_J;\lceil r\rceil}}\right)
 \sum_{m=0}^Jm\frac{(s\lambda_J)^{2m/s}}{(m!)^{2/s}}
       \norm{\nabla^mF}_{H^{\lceil r\rceil},\mathrm{sum}}^2\\
 &\quad+C\norm F_{\mathfrak G_J^{s,\lambda_J;\lceil r\rceil}}^3.
\end{aligned}
 \label{eq:short-time-factorial-energy}
\end{equation}
Choose
$\dot\lambda_J=-2Cs\lambda_J
\norm F_{\mathfrak G_J^{s,\lambda_J;\lceil r\rceil}}$
and $\lambda_J(0)=23\lambda_0/8$.
The input radius and Sobolev shift give an initial norm at most $C_1\eps$.
For small $\eps$, the differential inequality yields
\begin{gather*}
 \norm F_{\mathfrak G_J^{s,\lambda_J;\lceil r\rceil}}\leq2C_1\eps,\\
 \lambda_J(t)=\frac{23\lambda_0}{8}
 \exp\!\left(-2Cs\int_0^t
    \norm{F(\tau)}_{\mathfrak G_J^{s,\lambda_J(\tau);\lceil r\rceil}}
       \dd\tau\right)>\frac{21\lambda_0}{8}.
\end{gather*}
Apply monotone convergence at radius $21\lambda_0/8$, then apply
\eqref{eq:factorial-Fourier-equivalence} to obtain $5\lambda_0/2$.
\end{proof}

\begin{lemma}[Near-identity Gevrey change of variables]
\label{lem:short-time-gevrey-composition}
Use the factorial Gevrey norm \eqref{eq:factorial-Gevrey-definition}
and the endpoint interval \eqref{eq:fixed-interval}.
Let $a,b$ be one-dimensional Gevrey functions with
$\norm{a}_{\mathfrak G^{s,(19/8)\lambda_0;r}}+
 \norm{b'-1}_{\mathfrak G^{s,(19/8)\lambda_0;r-1}}\leq\delta_{\mathrm{comp}}$,
where $r>4$ and $\delta_{\mathrm{comp}}$ measures the size of the coordinate
perturbation and is assumed sufficiently small.  Assume
\(1/2<b'<3/2\) and that \(b\) maps the endpoints to \(L,R\).  Define
\begin{equation*}
 \Theta:\mathbb T_{2\pi}\times[-1,1]
 \longrightarrow\mathbb T_{2\pi}\times[L,R],
 \qquad \Theta(z,y)=(z+a(y),b(y)).
\end{equation*}
If $Q_{\rm tar}$ is boundary-flat on the target strip and $Q_{\rm src}$ is
boundary-flat on the source strip, then
\begin{equation}
 \norm{Q_{\rm tar}\circ\Theta}_{\Gs{s}{2\lambda_0}{r-4}}
 +\norm{Q_{\rm src}\circ\Theta^{-1}}_{\Gs{s}{2\lambda_0}{r-4}}
 \leq C\left(
   \norm{Q_{\rm tar}}_{\mathfrak G^{s,(19/8)\lambda_0;r}}
  +\norm{Q_{\rm src}}_{\mathfrak G^{s,(19/8)\lambda_0;r}}
 \right).
 \label{eq:short-time-composition-estimate}
\end{equation}
Moreover \(\dist(b(y),\left\{L,R\right\})\approx d(y)\), and the same change of
variables preserves an exponential boundary weight after decreasing its
constant.
\end{lemma}

\begin{proof}[Proof sketch]

Both compositions preserve flatness. We estimate them at radius
$35\lambda_0/16$, then pass to $2\lambda_0$.
Sum derivative indices by output row, so the identity has
first-derivative size one. For ordered positive partitions
$j_1+\cdots+j_p=m$, factorial normalization leaves
\[
 \left(\frac{p!\prod_i j_i!}{m!}\right)^{1/s-1}\leq1.
\]
A positive power-series majorant bounds composition with the forward map;
differentiating $b\circ b^{-1}=\mathrm{Id}$ gives the inverse majorant.
Strict gaps recover the fixed Sobolev orders and Fourier radius $2\lambda_0$,
with constants and a smallness threshold depending on these gaps and orders,
not on the cutoff. Constant vertical translations are unrestricted.
Endpoint distance and Jacobian comparison preserve the weakened boundary weight. 
\end{proof}

\begin{lemma}[Exponential weighted Hardy inequality]
\label{lem:exponential-weighted-Hardy}
Let \(\alpha>0\), \(0<c'<c\), and \(G\in H^1(0,1)\) with \(G(0)=0\).
Then
\begin{equation}
 \int_0^1e^{2c'r^{-\alpha}}|G(r)|^2\dd r
 \leq C_{c,c',\alpha}
 \int_0^1e^{2cr^{-\alpha}}|G'(r)|^2\dd r.
 \label{eq:exponential-weighted-Hardy}
\end{equation}
The analogous estimate holds at the right endpoint.
\end{lemma}

\begin{proof}
For \(0<r<1\), Cauchy--Schwarz gives
\begin{equation}
 |G(r)|^2\leq
 \left(\int_0^r e^{-2cs^{-\alpha}}\dd s\right)
 \left(\int_0^r e^{2cs^{-\alpha}}|G'(s)|^2\dd s\right).
 \label{eq:Hardy-pointwise-gap}
\end{equation}
The first factor is bounded by
\(Cr^{1+\alpha}e^{-2cr^{-\alpha}}\), by one integration by parts after the
change of variables \(u=s^{-\alpha}\).  Multiply by
\(e^{2c'r^{-\alpha}}\), integrate in \(r\), and use Tonelli's theorem.
The resulting inner integral is bounded uniformly because
\(c-c'>0\).  This proves \eqref{eq:exponential-weighted-Hardy}; replacing
\(r\) by \(1-r\) proves the estimate at the other endpoint.
\end{proof}

\begin{proposition}[Short-time weighted propagation]
\label{prop:short-time-weighted-propagation}
There are constants $c_1=c_1(c_0,\alpha_*)>0$ and $C>0$ such that, if
\eqref{eq:data-assumption} holds and $\eps$ is sufficiently small, then
\begin{align}
  \sup_{0\leq t\leq1}
  \norm{e^{c_1d(y)^{-\alpha_*}}\omega(t)}_{L^2(\Omega)}
  &\leq C\eps,
  \label{eq:short-time-conclusion}\\
  \sup_{0\leq t\leq1}\norm{F(t)}_{
       \Gs{s}{(5/2)\lambda_0}{2\sigma+8}}
  +\sup_{0\leq t\leq1}\norm{\omega(t)}_{H^3}
  &\leq C\eps.
  \label{eq:short-time-gevrey-conclusion}
\end{align}
At $t=1$, the functions $(f_{\mathrm{in}},h_{\mathrm{in}},g_{\mathrm{in}})$ defined by
\eqref{eq:profile-f-phi} obey the following
target-coordinate bounds:
\begin{equation}
  \norm{f_{\mathrm{in}}}_{\Gs{s}{2\lambda_0}{2\sigma+4}}
  +\norm{h_{\mathrm{in}}}_{\Gs{s}{2\lambda_0}{2\sigma+4}}
  +\norm{g_{\mathrm{in}}}_{\Gs{s}{2\lambda_0}{2\sigma+4}}
  +\norm{e^{c_2\mathbf d^{-\alpha_*}}(f_{\mathrm{in}},h_{\mathrm{in}},g_{\mathrm{in}})}_{L^2}
  \leq C\eps
  \label{eq:short-time-profile-data}
\end{equation}
for some $c_2>0$.  These bounds verify the hypotheses used in
Section~\ref{sec:decomposition}: Lemma~\ref{lem:localization}
applies directly to all three profiles; the estimate using the radius reserve
in Appendix~\ref{sec:localized-data-proof} yields \eqref{eq:shell-data}.
\end{proposition}

\begin{proof}[Proof sketch]

Gevrey propagation at $r=2\sigma+8$, tangency and Dirichlet
elliptic regularity give $\|U^y/(1-y^2)\|_\infty\lesssim\eps$.
Choose bounded increasing truncations with
$0\leq r\Theta_M'(r)\leq\Theta_M(r)$ and $\Theta_M(r)\uparrow r$.
For sufficiently small $\eps$,
\begin{align*}
 &(\partial_t+U^y\partial_y)
   \big[c_0(1-t/2)\Theta_M((1-y^2)^{-\alpha_*})\big]\\
 &\qquad\leq c_0
 \left(-\tfrac12+2\alpha_*
          \left\|\frac{U^y}{1-y^2}\right\|_\infty\right)
 \Theta_M((1-y^2)^{-\alpha_*})\leq0.
\end{align*}
Integrate the resulting weighted transport inequality before applying monotone
convergence. Since $1-y^2\leq2d(y)$, this proves
\eqref{eq:short-time-conclusion} with $c_1=c_02^{-1-\alpha_*}$.
Differentiating
\eqref{eq:coordinates} gives
\begin{align}
 f_{\mathrm{in}}(z,v(1,y))&=F(1,z+v(1,y)-y,y),
 \label{eq:short-time-f-data}\\
 h_{\mathrm{in}}(v(1,y))&=-\int_0^1\avg\omega(\tau,y)\dd\tau=\partial_yv(1,y)-1,
 \label{eq:short-time-h-data}\\
 g_{\mathrm{in}}(v(1,y))
 &=\avg{U^x}(1,y)-\int_0^1\avg{U^x}(\tau,y)\dd\tau,
 \label{eq:short-time-g-data}\\
 \partial_y[g_{\mathrm{in}}(v(1,y))]
 &=-F_0(1,y)+\int_0^1F_0(\tau,y)\dd\tau.
 \label{eq:short-time-dg-data}
\end{align}
The gauge fixes $v(1,-1)+1=c_-$, $|c_-|\lesssim\eps$, and
$g_{\mathrm{in}}(L)=g_{\mathrm{in}}(R)=0$.
Apply Lemma~\ref{lem:short-time-gevrey-composition} with
$a(y)=y-v(1,y)$ and $b(y)=v(1,y)$, using
$5\lambda_0/2>19\lambda_0/8>2\lambda_0$.
Averaging, the weighted Hardy inequality, and endpoint distance comparison give the
weighted bounds; the identities also preserve flatness.
\end{proof}

\section{The localized initial-data estimates}
\label{sec:localized-data-proof}

We prove \eqref{eq:shell-data} with a fixed radius reserve, directly for
the localizations of the three time-one profiles in \eqref{fingindef}.
The constants depend on the fixed parameters, not on the terminal layer.

\begin{proof}
Let $q$ denote any of $f_{\mathrm{in}},h_{\mathrm{in}},g_{\mathrm{in}}$.
Choose a common $c>0$ below the fixed decay coefficients in the
three weighted-profile estimates and in the cutoff Fourier estimate,
and keep it fixed throughout this proof.
The Fourier decay of $\Psi_n$ and the weighted time-one bound
\eqref{eq:short-time-profile-data} give, for $n\geq1$,
\[
 \|\Psi_nq\|_2\leq C\eps e^{-c2^{n\alpha_*}},\qquad
 \|\Psi_nq\|_{\Gs{s}{2\lambda_0}{2\sigma+4}}
       \leq C\eps e^{C2^{n\alpha_\dagger}}.
\]
For this interpolation argument retain the value of $c$ in the
first bound, common to the three profiles, and choose
$0<c_{\rm in}\leq c/8$ once and for all.
In the second bound one first uses the Gevrey product estimate and then
optimizes $2\lambda_0|\eta|^s-c(d_n|\eta|)^\rho$; the resulting loss is
of order $2^{n\rho s/(\rho-s)}=2^{n\alpha_\dagger}$.
Polynomial factors are absorbed by increasing the fixed constant.
Fourier interpolation with exponents $1/4$ and $3/4$ now yields
\begin{align}
 \|\Psi_nq\|_{\Gs{s}{3\lambda_0/2}{\sigma+3}}
 &\leq C\|\Psi_nq\|_2^{1/4}
       \|\Psi_nq\|_{\Gs{s}{2\lambda_0}{2\sigma+4}}^{3/4}
   \leq C\eps e^{-\frac c4 2^{n\alpha_*}+C2^{n\alpha_\dagger}}
   \leq C\eps e^{-c_{\rm in}2^{n\alpha_*}}.
 \label{eq:strict-shell-data}
\end{align}
Here $\sigma+3\leq\tfrac34(2\sigma+4)$, and
$\alpha_*>\alpha_\dagger$ absorbs the localization loss.
The reserve $c/4-c_{\rm in}\geq c/8>0$ makes the final bound uniform
in $n$; the finitely many initial shells are included in $C$.  The first
increment $\chi_1q=(\Psi_0+\Psi_1)q$ obeys the same estimate after increasing
the constant.  Lowering the radius to $\lambda_0$ proves
\eqref{eq:shell-data}; it also controls any fixed finite number of derivatives
needed in the local construction.  The support assertion
\eqref{eq:support-scale} follows from the zero region of $\chi_n$.
Finally, telescoping the partition \eqref{eq:exact-shell-partition} and using
\eqref{eq:terminal-cutoff-identity} recovers the time-one datum.
\end{proof}

\section{\texorpdfstring{The weight and structural multiplier estimates}{The weight and structural multiplier estimates}}
\label{sec:orr-multiplier-appendix}
We collect the weight construction and structural comparisons used with the multipliers defined in Subsection~\ref{sec:weights}.

\subsection{Critical intervals and weights}
\label{sec:orr-weights}
We use the construction of \cite[(3.5)--(3.10)]{BM15}, with
$C_\kappa=1$ and the correction exponent $\mu=54$ fixed in
\eqref{eq:weight-parameters}. The scalar weights are even in $\eta$;
$w_{-k}(t,-\eta)=w_k(t,\eta)$, but $w_k$ need not be even at fixed $k\ne0$.
Set $N(\eta)=\lfloor\sqrt{\eta}\rfloor$ for $\eta\geq0$.
For every \(k\in\mathbb Z\) with
\(1\leq |k|\leq N(|\eta|)\) and
\(\operatorname{sgn}k=\operatorname{sgn}\eta\), define
\begin{align}
 t_{k,\eta}&=\frac{|\eta|}{|k|}-\frac{|\eta|}{2|k|(|k|+1)},
 \qquad t_{0,\eta}=2|\eta|,\label{eq:critical-endpoints}\\
 I_{k,\eta}&=[t_{|k|,|\eta|},t_{|k|-1,|\eta|}].
 \label{eq:critical-intervals}
\end{align}
For all other \((k,\eta)\), set \(I_{k,\eta}=\varnothing\).  
The full cutoff $N(\eta)$ is used in the recursion below.
For the echo interactions, we retain only the critical intervals with
left endpoint at least $2\sqrt{|\eta|}$:
\begin{equation}
 \mathbf I_{k,\eta}=
 \begin{cases}
 I_{k,\eta},&I_{k,\eta}\ne\varnothing,\quad
             t_{|k|,|\eta|}\geq2\sqrt{|\eta|},\\
 \varnothing,&\text{otherwise}.
 \end{cases}
 \label{eq:restricted-critical-intervals}
\end{equation}
The largest retained positive index is $\sqrt{\eta}/2+O(1)$.
If \(|\eta|<1\), all the weights below are set equal to one.  For
\(\eta\geq1\), let
\begin{equation}
 \begin{aligned}
  b_{1,\eta}=1-\frac1\eta,
  \qquad
  &b_{k,\eta}=\frac{2(k-1)}k\left(1-\frac{k^2}{\eta}\right)
       \quad(k\geq2),\\
  &a_{k,\eta}=\frac{2(k+1)}k\left(1-\frac{k^2}{\eta}\right)
       \quad(k\geq1).
 \end{aligned}
  \label{eq:weight-shape-parameters}
\end{equation}
Starting from $w_{\NR}(2\eta,\eta)=1$, define, for $1\leq k\leq N(\eta)$,
\begin{equation}
 w_{\NR}(t,\eta)=
 \begin{cases}
 \left(\dfrac{k^2}{\eta}[1+b_{k,\eta}|t-\eta/k|]\right)^{C_\kappa}
 w_{\NR}(t_{k-1,\eta},\eta),&\eta/k\leq t\leq t_{k-1,\eta},\\[3pt]
 (1+a_{k,\eta}|t-\eta/k|)^{-1-C_\kappa}
 w_{\NR}(\eta/k,\eta),&t_{k,\eta}\leq t\leq\eta/k,
 \end{cases}
 \label{eq:wNR-right}
\end{equation}
\begingroup
\expandafter\edef\csname @currentlabel\endcsname{\theequation}
\label{eq:wNR-left}
\endgroup
\begin{equation}
 \frac{w_{\Res}(t,\eta)}{w_{\NR}(t,\eta)}=
 \frac{k^2}{\eta}
 \begin{cases}
 1+b_{k,\eta}|t-\eta/k|,&\eta/k\leq t\leq t_{k-1,\eta},\\
 1+a_{k,\eta}|t-\eta/k|,&t_{k,\eta}\leq t\leq\eta/k.
 \end{cases}
 \label{eq:wR-right}
\end{equation}
\begingroup
\expandafter\edef\csname @currentlabel\endcsname{\theequation}
\label{eq:wR-left}
\endgroup
Extend both \(w_{\NR}\) and \(w_{\Res}\) constantly to
\(t\leq t_{N(\eta),\eta}\), set \(w_{\NR}=w_{\Res}=1\) for
\(t\geq2\eta\), and extend both weights evenly to negative \(\eta\).
The choices in \eqref{eq:weight-shape-parameters} make the recursive pieces
continuous at their endpoints.  Finally, we define the mode-dependent weight by
\begin{equation}
 w_k(t,\eta)=
 \begin{cases}
 w_{\Res}(t,|\eta|),&t\in I_{k,\eta},\\
 w_{\NR}(t,|\eta|),&t\notin I_{k,\eta}.
 \end{cases}
 \label{eq:mode-dependent-weight}
\end{equation}

\subsection{Structural estimates for the multipliers}
\label{sec:multiplier-toolkit}

We record the projections and weight comparisons used in the elliptic
and energy estimates, with the radii and multipliers of
Subsection~\ref{sec:weights}.

\paragraph{Frequency projections.}
We use the Littlewood--Paley convention of \cite[Appendix~A.1]{BM15}.
Choose an even smooth cutoff $\psi:\R\to[0,1]$, nonincreasing in $|\nu|$,
with $\psi(\nu)=1$ for $|\nu|\leq1/2$ and $\psi(\nu)=0$ for
$|\nu|\geq3/4$.  Put $\mathbf D=\{2^j:j=-1,0,1,\ldots\}$ and
\begin{equation}
 \pi_{1/2}(\nu)=\psi(\nu),\qquad
 \pi_N(\nu)=\psi(\nu/(2N))-\psi(\nu/N),\quad N\in\mathbf D,\ N\geq1.
 \label{eq:explicit-LP-symbols}
\end{equation}
These nonnegative symbols sum to one.  Their supports satisfy
\begin{equation}
 \operatorname{supp}\pi_{1/2}\subset\left\{|\nu|\leq3/4\right\},\qquad
 \operatorname{supp}\pi_N\subset\left\{N/2\leq|\nu|\leq3N/2\right\},\quad N\geq1.
 \label{eq:explicit-LP-annuli}
\end{equation}
All unrestricted dyadic sums run over $\mathbf D$.  The strict truncation
$F_{<N/8}=\sum_{M\in\mathbf D,\,M<N/8}F_M
=\sum_{M\in\mathbf D,\,M\leq N/16}F_M$ is zero if the sum is empty.
For $N\geq8$ its symbol is $\psi(8\nu/N)$, and hence
\begin{equation}
 \operatorname{supp}\widehat F_{<N/8}\subset\left\{|\nu|\leq3N/32\right\},\qquad
 \frac{|\nu_{\mathrm{low}}|}{|\nu_{\mathrm{high}}|}\leq\frac3{16}.
 \label{eq:explicit-LP-low-high-support}
\end{equation}
The first nonempty low--high pair has high label $N=8$.
The complementary comparable pairs are $N/8\leq N'\leq8N$.
We use $\nu=\eta$ for the vertical projections in the precision elliptic
estimates and $\nu=(k,\eta)$, with $|\nu|=|k|+|\eta|$, for the full
energy projections; in the latter formulas $\psi(\nu)$ means $\psi(|\nu|)$.
Their restriction to horizontal mode zero gives the scalar energy projections.
Thus the precision and energy decompositions both use the same strict
$<N/8$ convention.  The full symbols may have corners on the axes.
We use only nonnegativity, the partition of unity, the support bounds,
and the $L^2$ and finite-overlap estimates from Plancherel; no derivative
of a projection symbol or kernel moment is used.

\paragraph{Weight comparisons and product estimates.}

Constants involving Gevrey radii or Sobolev factors may depend on the
fixed parameters $s,a,\mu,\lambda_0,\sigma$ and cutoffs, but not on time or
Newton and operator indices.  The raw-weight comparisons retain the stated dependence on the
parameters.  We use comparisons (3.16), (3.18), (3.31), (3.34), and
(3.39b), together with the early-time finite-difference argument in
\cite[Lemma~3.7]{BM15}, expressed in the notation of the Newton scheme.  The radius normalization used below relies on
$a\leq2s-1$, which gives $\avg t^{-s}\lesssim\sqrt{-\dot\lambda_n}$.
The raw weights have $C_\kappa=1$, hence growth parameter
$4(1+2C_\kappa)=12$ in \cite[Lemma~3.1]{BM15}.  They are unchanged by the
larger correction exponent $\mu=54$ in \eqref{eq:weight-parameters} and
\eqref{eq:J-mode}.  Combining their comparisons with the square-root
exponentials in ${\J}_k$ gives the
$e^{C\mu|(k-\ell,\eta-\xi)|^{1/2}}$ bounds below.

\begin{lemma}[Growth and logarithmic derivatives of the weights]
\label{lem:weight-growth}
For all \(t\geq1\), \(k\in\mathbb Z\), and
\(\eta\in\mathbb R\),
\begin{equation}
  0<w_k(t,\eta),\,w_{\NR}(t,\eta),\,w_{\Res}(t,\eta)\leq1,
  \qquad
  w_k^{-1}+w_{\NR}^{-1}+w_{\Res}^{-1}
  \lesssim e^{6|\eta|^{1/2}}=e^{(\mu/9)|\eta|^{1/2}}.
  \label{eq:weight-total-growth}
\end{equation}
If \(t\in\mathbf I_{k,\eta}\) and \(k\neq0\), then
\begin{equation}
  \frac{\partial_tw_k(t,\eta)}{w_k(t,\eta)}
  \approx
  \frac{\partial_tw_{\NR}(t,\eta)}{w_{\NR}(t,\eta)}
  \approx
  \frac{\partial_tw_{\Res}(t,\eta)}{w_{\Res}(t,\eta)}
  \approx\frac1{1+|t-\eta/k|}.
  \label{eq:weight-log-derivative}
\end{equation}
On every unrestricted critical interval, the corresponding upper bounds in
\eqref{eq:weight-log-derivative} remain valid.
\end{lemma}

\begin{proof}[Proof sketch] \space Multiply and differentiate the interval formulas as in
\cite[Lemmas~3.1 and 3.3, Remark~8]{BM15}, with $C_\kappa=1$.
For $x\geq1$ the initial inverse weight is
$[x^{\lfloor\sqrt x\rfloor}/(\lfloor\sqrt x\rfloor!)^2]^3$.
It is continuous at the square endpoints and its logarithmic derivative
is at most $3/\sqrt x$; hence every inverse weight is at most
$e^{6(\sqrt{|\eta|}-1)}$ for $|\eta|\geq1$.
The correction exponent $\mu=54$ does not change these raw weights.
Only the restricted intervals admit the two-sided derivative comparison;
on unrestricted terminal intervals we use the upper bound.
\end{proof}

\begin{lemma}[Separation of critical intervals]
\label{lem:critical-separation}
Let $\alpha\geq1$, $\alpha^{-1}|\xi|\leq|\eta|\leq\alpha|\xi|$, and
$t\in I_{k,\eta}\cap I_{\ell,\xi}$. Then $|k|\approx_\alpha|\ell|$.
If $k\ne\ell$, either
\begin{equation}
 \left|t-\frac\eta k\right|\geq\frac{|\eta|}{10\alpha k^2},
 \qquad\left|t-\frac\xi\ell\right|\geq\frac{|\xi|}{10\alpha\ell^2},
 \label{eq:far-from-both-resonances}
\end{equation}
or
\begin{equation}
 |\eta-\xi|\gtrsim_\alpha\frac{|\eta|}{|\ell|}.
 \label{eq:vertical-frequency-separation}
\end{equation}
\end{lemma}

\begin{proof}
This is the critical-interval trichotomy of \cite[Lemma~3.2]{BM15},
expressed in our notation for the unrestricted intervals.
\end{proof}

\begin{lemma}
\label{lem:weight-comparison}
There is \(c_w>0\) such that, for all \(t\geq1\) and \(\eta,\xi\in\mathbb R\),
\begin{equation}
  \frac{w_{\NR}(t,\xi)}{w_{\NR}(t,\eta)}
  +\frac{w_{\NR}(t,\eta)}{w_{\NR}(t,\xi)}
  \lesssim e^{c_w|\eta-\xi|^{1/2}}.
  \label{eq:wNR-comparison}
\end{equation}
For all $k,\ell\in\mathbb Z$,
\begin{equation}
  \frac{{\J}_k(t,\eta)}{{\J}_\ell(t,\xi)}
  \lesssim
  \left(1+\bm1_{\left\{t\in\mathbf I_{k,\eta}\right\}}
    \frac{|\eta|}{k^2(1+|t-\eta/k|)}\right)
  e^{C\mu|(k-\ell,\eta-\xi)|^{1/2}}.
  \label{eq:J-ratio}
\end{equation}
\begingroup
\expandafter\edef\csname @currentlabel\endcsname{\theequation}
\label{eq:raw-J-general-ledger}
\endgroup
The parenthetical factor is taken to be one if $k=0$ or $k=\ell$, if
$t\notin\mathbf I_{k,\eta}$, or, for comparable $|\eta|,|\xi|$, if
\[
 t\in(\mathbf I_{k,\eta}\setminus\mathbf I_{k,\xi})
       \cup(I_{k,\eta}\setminus I_{k,\xi}).
\]
\end{lemma}

\begin{proof}[Proof sketch] \space We follow the proof of \cite[Lemmas~3.5--3.6 and Remark~10]{BM15},
comparing first the scalar weights and then the two terms in ${\J}_k$.
In the improved cases, we retain the distinction between restricted
and unrestricted intervals.
\end{proof}

\begin{lemma}[Same-mode reduced-weight comparison]
\label{lem:reduced-J}
For every fixed $k\ne0$ and all $t\geq1$,
\begin{equation}
 \frac{\widetilde {\J}_k(t,\eta)}{\widetilde {\J}_k(t,\xi)}
 +\frac{\widetilde {\J}_k(t,\xi)}{\widetilde {\J}_k(t,\eta)}
 \leq C e^{C\mu|\eta-\xi|^{1/2}}.
 \label{eq:reduced-J-ratio}
\end{equation}
We do not deduce this conclusion from the corresponding estimate for
${\J}_k$, which contains the additive term $e^{\mu|k|^{1/2}}$.
\end{lemma}

\begin{proof}
We use the fixed-mode comparison of the raw weights in the proof of
\cite[(3.31)]{BM15}, before adding the horizontal term in ${\J}_k$.
Indeed,
\[
 \frac{\widetilde {\J}_k(t,\eta)}{\widetilde {\J}_k(t,\xi)}
 =\frac{w_k(t,\xi)}{w_k(t,\eta)}
       e^{\mu(\sqrt{|\eta|}-\sqrt{|\xi|})}.
\]
The scalar comparison and explicit same-mode resonant factor give
both directions.
\end{proof}

\begin{lemma}[Raw-weight and coupled CK estimates]
\label{lem:raw-BM-ledger}
For $|\eta|\geq1$, the logarithmic derivative
$\frac{\partial_tw_k(t,\eta)}{w_k(t,\eta)}$ vanishes outside the scalar
active range $[t_{N(|\eta|),|\eta|},2|\eta|]$; it need not vanish
outside the mode's own critical interval.  For $|\eta|<1$ it is zero.
The following statements hold almost everywhere in $t$.
\begin{enumerate}[(i)]
 \item For every fixed $\alpha\geq1$, if $\alpha^{-1}|\xi|\leq|\eta|\leq\alpha|\xi|$, then
 \begin{equation}
  \sqrt{\frac{\partial_tw_\ell(t,\xi)}{w_\ell(t,\xi)}}
  \lesssim_\alpha \avg{\eta-\xi}
   \left(\sqrt{\frac{\partial_tw_k(t,\eta)}{w_k(t,\eta)}}
         +\frac{|\eta|^{s/2}}{\avg t^s}\right).
  \label{eq:raw-w-endpoint-transfer}
 \end{equation}
 \item 
In the reverse imbalanced case
 $t\in I_{\ell,\xi}\setminus I_{k,\eta}$,
 \begin{equation}
  \frac{{\J}_k(t,\eta)}{{\J}_\ell(t,\xi)}
  \lesssim
  \frac{\ell^2(1+|t-\xi/\ell|)}{|\xi|}
  e^{C\mu|(k-\ell,\eta-\xi)|^{1/2}}.
  \label{eq:raw-J-reverse-ledger}
 \end{equation}
 We do not require $|\eta|$ and $|\xi|$ to be comparable here.  The extension
  beyond the comparable-frequency case of \cite[(3.32)]{BM15} is proved below.
 \item If the frequencies are comparable, \(k\neq\ell\), and
 \begin{equation}
  t\in\mathbf I_{k,\eta}\cap\mathbf I_{k,\xi},
  \label{eq:BM-common-index-set}
 \end{equation}
 then
 \begin{equation}
  \frac{{\J}_k(t,\eta)}{{\J}_\ell(t,\xi)}
  \lesssim
   \frac{|\eta|}{k^2}
   \sqrt{\frac{\partial_tw_k(t,\eta)}{w_k(t,\eta)}}\sqrt{\frac{\partial_tw_\ell(t,\xi)}{w_\ell(t,\xi)}}
   e^{C\mu|(k-\ell,\eta-\xi)|^{1/2}}.
  \label{eq:raw-J-two-CK-ledger}
 \end{equation}
\end{enumerate}
The same pointwise assertions hold after a low--high paraproduct switch.
After a high--low switch, the alternative
\eqref{eq:raw-J-two-CK-ledger} places the $\CK$-term on the high
coefficient.  These are the only two allocations used in the reaction,
transport, and precision elliptic estimates.
\end{lemma}

\begin{proof}
The interval comparisons are \cite[Lemma~3.4(ii), (3.16),
(3.30)--(3.34), and Remark~10]{BM15}.  In
\eqref{eq:BM-common-index-set}, both intervals have the numerator index
$k$, exactly as in the published (3.34).  The inactive and terminal cases
use (3.16), including its full frequency-gap factor.

To extend the published estimate (3.32) beyond comparable frequencies,
we estimate the two terms of ${\J}_k$.  On $I_{\ell,\xi}\setminus I_{k,\eta}$,
\begin{align*}
 \frac{{\J}_k(t,\eta)}{{\J}_\ell(t,\xi)}
 &\leq\frac{w_{\Res}(t,|\xi|)}{w_{\NR}(t,|\xi|)}
 \left[\frac{w_{\NR}(t,|\xi|)}{w_{\NR}(t,|\eta|)}
 e^{\mu(\sqrt{|\eta|}-\sqrt{|\xi|})}
 +w_{\NR}(t,|\xi|)e^{\mu(\sqrt{|k|}-\sqrt{|\xi|})}\right]\\*
 &\lesssim\frac{\ell^2(1+|t-\xi/\ell|)}{|\xi|}
 e^{C\mu|(k-\ell,\eta-\xi)|^{1/2}}.
\end{align*}
Here \eqref{eq:wNR-comparison}, $w_{\NR}\leq1$, and
$\sqrt{|k|}\leq\sqrt{|k-\ell|}+\sqrt{|\xi|}$ give the second line,
since $1\leq|\ell|\leq\sqrt{|\xi|}$.  This also covers opposite signs
and $k=0$.  Finally, \eqref{eq:lambda-derivative} and $a\leq2s-1$ give
the radius conversion used in the elliptic estimates:
\begin{equation}
 \frac{|(k,\xi)|^{s/2}}{\avg t^s}
 \lesssim(-\dot\lambda_n)^{1/2}|(k,\xi)|^{s/2}.
 \label{eq:w-endpoint-paid-by-lambda}
\end{equation}
\end{proof}

\begin{lemma}
\label{lem:zero-mode-endpoint-dichotomy}
Let $k\neq0$, $|\eta|\approx|\xi|$, and
$|(k,\eta-\xi)|\leq c|\xi|$ for a fixed $c\in(0,1)$.  Recall that
$w_0=w_{\NR}$, and index both restricted intervals from
\eqref{eq:restricted-critical-intervals} by $k$:
\[
 \chi_D=\bm1_{\left\{t\in\mathbf I_{k,\eta}\cap\mathbf I_{k,\xi}\right\}},
 \qquad \chi_*=1-\chi_D.
\]
Then, 
\begin{align}
 \chi_D\frac{{\J}_k(t,\eta)}{{\J}_0(t,\xi)}
 &\lesssim \chi_D\frac{|\eta|}{k^2}
 \sqrt{\frac{\partial_tw_k(t,\eta)}{w_k(t,\eta)}\frac{\partial_tw_0(t,\xi)}{w_0(t,\xi)}}
 e^{C\mu|(k,\eta-\xi)|^{1/2}},
 \label{eq:zero-mode-J-CK-dichotomy}\\
 \chi_*\frac{{\J}_k(t,\eta)}{{\J}_0(t,\xi)}
 &\lesssim \chi_*e^{C\mu|(k,\eta-\xi)|^{1/2}}.
 \label{eq:zero-mode-J-improved-dichotomy}
\end{align}
\end{lemma}

\begin{proof}
We apply the comparisons \cite[(3.31), (3.34)]{BM15}, using $w_0=w_{\NR}$.
\end{proof}

\begin{lemma}[Orr--elliptic balance]
\label{lem:orr-elliptic-balance}
Recall the full and restricted critical intervals from
Appendix~\ref{sec:orr-weights} and the flat elliptic symbol
$\ell^2+|\xi-\ell t|^2$.
Let $\ell\ne0$ and $t\geq1$.
Away from the critical interval,
\begin{equation}
  \bm 1_{\left\{t\notin I_{\ell,\xi}\right\}}\frac{\avg{\ell,\xi}}{\ell^2+|\xi-\ell t|^2}
  \lesssim
  \frac{\avg{\ell,\xi}^{s}}{\avg t^{2s}}.
  \label{eq:nonresonant-elliptic-symbol}
\end{equation}

\end{lemma}

\begin{proof}[Proof sketch] \space This is the nonresonant symbol comparison used in
\cite[Sections~4.2 and 6.1]{BM15}.  Away from $I_{\ell,\xi}$,
$\ell^2+|\xi-\ell t|^2\gtrsim
\avg{\ell,\xi}^{1-s}\avg t^{2s}$; the opposite-sign case follows
directly from $|\xi-\ell t|=|\xi|+|\ell|t$.
\end{proof}

We use the Gevrey product estimate of \cite[Lemma~A.3]{BM15}.
The scalar resonant weights satisfy the following comparisons.

\begin{lemma}[Scalar resonant envelope and finite-index comparison]
\label{lem:scalar-resonant-envelope}
For $m\in\mathbb Z$ and almost every $t\geq1$,
\begin{equation}
 w_{\Res}(t,\eta)=\min_{m\in\mathbb Z}w_m(t,\eta),
 \qquad
 \frac{\partial_tw_m(t,\eta)}{w_m(t,\eta)}\lesssim \frac{\partial_tw_{\Res}(t,\eta)}{w_{\Res}(t,\eta)}.
 \label{eq:scalar-resonant-envelope}
\end{equation}
Moreover, for every fixed $\alpha\geq1$, if
$\alpha^{-1}|\xi|\leq|\eta|\leq\alpha|\xi|$, then
\begin{align}
 \frac{{\J}_{\Res}(t,\eta)}{{\J}_{\Res}(t,\xi)}
 +\frac{{\J}_{\Res}(t,\xi)}{{\J}_{\Res}(t,\eta)}
 &\lesssim_\alpha
 e^{C_\alpha\mu|\eta-\xi|^{1/2}},
 \label{eq:scalar-resonant-J-comparison}\\
 \sqrt{\frac{\partial_tw_{\Res}(t,\eta)}{w_{\Res}(t,\eta)}}
 &\lesssim_\alpha\avg{\eta-\xi}^{C_\alpha}
 \left[
   \sqrt{\frac{\partial_tw_{\Res}(t,\xi)}{w_{\Res}(t,\xi)}}
   +\frac{|\xi|^{s/2}}{\avg t^s}
 \right].
 \label{eq:scalar-resonant-derivative-comparison}
\end{align}
\end{lemma}

\begin{proof}[Proof sketch] \space On each active interval $I_{j,|\eta|}$, the defining formulas give
\[
 \frac{j^2}{|\eta|}
 \leq\frac{w_{\Res}(t,\eta)}{w_{\NR}(t,\eta)}\leq1,\qquad
 \frac{\partial_tw_m(t,\eta)}{w_m(t,\eta)}
 \leq2\frac{\partial_tw_{\Res}(t,\eta)}{w_{\Res}(t,\eta)}.
\]
At the interval endpoints the ratio is one; outside the active range
the two scalar weights agree and their derivatives vanish.
Thus $w_{\Res}$ is the minimum, attained by a nonzero mode.
Lemma~\ref{lem:reduced-J} gives
\[
 \J_{\Res}(t,\eta)
 =\max_{m\ne0}\widetilde\J_m(t,\eta)
 \leq Ce^{C\mu|\eta-\xi|^{1/2}}
       \max_{m\ne0}\widetilde\J_m(t,\xi)
 =Ce^{C\mu|\eta-\xi|^{1/2}}\J_{\Res}(t,\xi).
\]
Interchanging $\eta,\xi$ proves \eqref{eq:scalar-resonant-J-comparison}.
For \eqref{eq:scalar-resonant-derivative-comparison}, the case of a zero output
derivative is immediate. Otherwise take the active signed mode $m$
at $\eta$, apply \eqref{eq:raw-w-endpoint-transfer} at $(m,\eta),(m,\xi)$,
then use \eqref{eq:scalar-resonant-envelope} at $\xi$.
The mode comparison includes unequal active indices and terminal endpoints.
\end{proof}

\begin{lemma}[Transfer from a mode weight to the resonant coefficient weight]
\label{lem:elliptic-J-transfer}
For $k\in\mathbb Z$ and $\xi,\xi'\in\mathbb R$, suppose
\begin{equation}
  |(k,\xi')|\leq\frac14|\xi-\xi'|.
  \label{eq:J-transfer-region}
\end{equation}
Then, for almost every $t\geq1$,
\begin{align}
  {\J}_k(t,\xi)
  &\lesssim {\J}_{\Res}(t,\xi-\xi')
    e^{C\mu|(k,\xi')|^{1/2}},
  \label{eq:J-transfer-to-resonant}\\
  \sqrt{\frac{\partial_tw_k(t,\xi)}{w_k(t,\xi)}}
  &\lesssim e^{C\mu|(k,\xi')|^{1/2}}
  \left[
    \sqrt{\frac{\partial_tw_{\Res}(t,\xi-\xi')}{w_{\Res}(t,\xi-\xi')}}
    +\frac{|\xi-\xi'|^{s/2}}{\avg t^s}
  \right].
  \label{eq:w-transfer-to-resonant}
\end{align}
\end{lemma}

\begin{proof}
We apply the scalar envelope to the coefficient-weight transfer of
\cite[Section~4.2.2, (4.8)]{BM15}. The hypothesis gives
\[
 \frac34|\xi-\xi'|\leq|\xi|\leq\frac54|\xi-\xi'|,
 \qquad |k|\leq\frac14|\xi-\xi'|.
\]
By \eqref{eq:scalar-resonant-envelope},
\eqref{eq:scalar-resonant-J-comparison}, and $w_{\Res}\leq1$,
\begin{equation*}
 {\J}_k(t,\xi) \leq {\J}_{\Res}(t,\xi)+e^{\mu|k|^{1/2}} \lesssim {\J}_{\Res}(t,\xi-\xi') e^{C\mu|(k,\xi')|^{1/2}}.
\end{equation*}
For the logarithmic derivative, first use the envelope bound
\eqref{eq:scalar-resonant-envelope}, then apply
\eqref{eq:scalar-resonant-derivative-comparison} with $\alpha=4/3$.
Its polynomial factor $\avg{\xi'}^C$ is bounded by
$Ce^{C\mu|(k,\xi')|^{1/2}}$, which proves
\eqref{eq:w-transfer-to-resonant}.
\end{proof}

\vspace{2 mm}

\noindent \textbf{Acknowledgments:} SI is grateful for support from NSF CAREER DMS2442781 and NSF DMS2306528. SH is supported by NSF DMS2406293.

\vspace{2 mm}

\noindent \textbf{AI Statement:} Publicly available versions of ChatGPT (Sol and Astra, within Codex) assisted with language editing, organization, bibliographic checks, and checking calculations, based on a 63 page draft prepared entirely by the authors containing a complete proof. The authors take full responsibility for the manuscript in its present form.

\footnotesize\par\bigskip
\Needspace{5\baselineskip}
\phantomsection

\end{document}